\documentclass[12pt,twoside,bezier]{article}

\usepackage{amssymb,amsthm}

\usepackage{imakeidx}
\usepackage{enumerate}

\newcommand{\iv}{^{-1}}
\newcommand{\vk}{ van Kampen }
\newcommand{\ct}{ contiguity }
\newcommand{\topp}{{\bf top}}

\newcommand{\bott}{{\bf bot}}

\newcommand{\oo}{{\bf 0}}

\newcommand{\la}{\langle}
\newcommand{\ra}{\rangle}

\newcommand{\tool}{\stackrel{\ell}{\too} }

\newcommand{\ttt}{{\cal T}}

\newcommand{\too}{\to}

\newcommand{\bb}{{\cal B}}

\newcommand{\mmm}{\mathbf{M}}
\newcounter{ppp}
\newcounter{pdten}
\newcounter{pdeleven}

\usepackage{imakeidx}
\makeindex[name=g,title=Subject index]
\usepackage[columns=2]{idxlayout}

\begin{document}
\theoremstyle{plain}
\newtheorem{theo}{Theorem}[section]
\newtheorem{lm}[theo]{Lemma}
\newtheorem{cy}[theo]{Corollary}
\newtheorem{df}[theo]{Definition}
\newtheorem{remark}[theo]{Remark}
\newtheorem{prop}[theo]{Proposition}
\newtheorem{prob}[theo]{Problem}

%\usepackage{imakeidx}
%\makeindex[name=g,title=Subject index]
%\usepackage[columns=2]{idxlayout}

\title{Subgroups of groups finitely presented in Burnside varieties}

 \author{A.Yu. Olshanskii, \thanks{The author was supported in part by the NSF grants DMS-1500180 and DMS-1901976.}}

\date{}
\maketitle

\tableofcontents

\begin{abstract} For all sufficiently large odd integers $n$, the following version of Higman's embedding theorem is proved in the variety ${\cal B}_n$ of all groups satisfying the identity $x^n=1$. A finitely generated group $G$ from ${\cal B}_n$
has a presentation $G=\langle A\mid R\rangle$ with a finite
set of generators $A$ and a recursively enumerable set $R$ of defining relations if and only if it is a subgroup of a group $H$ finitely
presented in the variety ${\cal B}_n$. It follows that there is a 'universal' $2$-generated finitely
presented in ${\cal B}_n$ group containing isomorphic copies of all finitely
presented in ${\cal B}_n$ groups as subgroups.
\end{abstract}

\medskip

{\bf Key words:} generators and relations in groups, finitely and recursively presented groups, Burnside
variety, van Kampen diagram, Turing machine, S-machine 

\medskip 

{\bf AMS Mathematical Subject Classification:} 20F05, 20F06, 20F50, 20F65, 20E06, 20E10, 03D10, 20F10, 03D25

\section{Introduction}

The celebrated theorem of G. Higman \cite{H} asserts that a finitely generated group is recursively presented if and only if it is a subgroup of a finitely presented group. This theorem detects a deep connection between the
logic concept of recursiveness and properties of
finitely presented groups.

The Higman embedding $G \hookrightarrow H$ of a recursively presented group $G$ makes the finitely presented group $H$ quite large, saturating it with free subgroups. Therefore
there were no nontrivial analogies of Higman's theorem in proper subvarieties of the variety of all groups. (The statement is trivial if every finitely generated group of a variety $\cal V$ is finitely presented in $\cal V$.)
O. Kharlampovich finished the paper \cite{Kh}
on the following question:

{\it Is every finitely generated group with the identity $x^n=1$ and recursively enumerable set of defining relations embeddable into a group given by this identity and a finite set of relations ($n$ is a large number) ? In other words, is the analogy of Higman's embedding theorem true in the variety ${\cal B}_n$ ?}

Even earlier the same problem was formulated in \cite{KT} (problem 12.32) by S.V. Ivanov :

{\it Prove an analogue of Higman’s theorem for the Burnside variety  ${\cal B}_n$ of groups
of odd exponent $n >> 1$, that is, prove that every recursively presented group of
exponent $n$ can be embedded in a finitely presented (in ${\cal B}_n$) group of exponent $n$.}

Recall that the \index[g]{Burnside variety} {\it Burnside variety} of exponent $n\ge 1$, denoted by ${\cal B}_n$, is the class
of all groups satisfying the identical relation $x^n=1$. So ${\cal B}_n$ is the class of groups
having finite exponent dividing $n$.  Every $m$-generated group of this variety is isomorphic to a
factor group of the $m$-generated free group $B(m,n)$ in ${\cal B}_n$, which in turn is isomorphic to the factor group of
a free group $F(m)$ of rank $m$ modulo the normal subgroup
$F_m^n$ generated by all powers $g^n$, where $g\in F_m$. (If two words from $F(m)$ have equal images in the \index[g]{free Burnside group} {\it free Burnside group} $B(m,n)$, we say that these
words are equal {\it modulo} \index[g]{Burnside relations} {\it  the Burnside relations.})

The group
$G$ is called \index[g]{finitely presented group in ${\cal B}_n$} {\it finitely presented} in the variety ${\cal B}_n$
if $G\cong B(m,n)/N$ for some $m\ge 1$, where $N$ is a normal closure of a finite subset $R\subset B(m,n)$, in other words, $G$ satisfies the identity $x^n=1$ and have a presentation $\langle a_1,\dots, a_m\mid R\cup S\rangle$ in the class of all groups, where $R$ is a finite set of relators and every relator from
$S$ is an $n$-th power in the free group $F(a,\dots, a_m)$.

\begin{theo} \label{HE} There is a constant $C$ such that for every
odd integer $n\ge C$, the following is true.

A finitely generated group $G$ satisfying the identity $x^n=1$
has a presentation $G=\langle A\mid R\rangle$ with a finite
set of generators $A$ and a recursively enumerable set $R$ of defining relations if and only if it is a subgroup of a group $H$ finitely
presented in the variety ${\cal B}_n$.
\end{theo}

The embedding property of Theorem \ref{HE} holds for groups $G$ with countable sets of generators as well. Furthermore, the group $H$ can be chosen with two generators. More precisely, we have

\begin{cy}\label{twog} There is a constant $C$ such that for every
odd integer $n\ge C$, the following is true.

Let a group $G$ from the variety ${\cal B}_n$
be given by a countable set of generators
$x_1, x_2,\dots$ and a recursively enumerable
set of defining relations $R$ in these generators. Then $G$ is a subgroup of a 2-generated group $E$ finitely presented in ${\cal B}_n$.
\end{cy}

\proof The group $G$ is isomorphic to the factor group $B(\infty,n)/N$, where $B(\infty,n)$ is the free Burnside group of exponent $n$
with countable set of free generators $x_1, x_2,\dots$ and $N$ is the normal closure of the
set $R$ in $B(\infty,n)$.

We will use a theorem of D. Sonkin
(Theorem 2.2 \cite{S}). According to this theorem the 2-generated free Burnside group $B(2,n)$ with
free basis $(a,b)$ has a subgroup $K$ with the
following properties.

The subgroup $K$ is generated by a countable set of words
$w_1, w_2,...$ in the generators $a,b$ such that this set is recursive, $K$ is a free Burnside group of exponent $n$ with free generators $w_1, w_2,\dots$, and $K$ has \index[g]{Congruence Extension Property} 
{\it Congruence Extension Property} (CEP).  CEP
means that every normal subgroup of $K$ is
the intersection of a normal subgroup of $B(2,n)$ with $K$. Equivalently, a normal closure in $K$ of arbitrary subset $S\subset K$
is the intersection of the normal closure of $S$ in $B(\infty,n)$ with $K$.

Now the group $G$ can be presented as the factor group $K/M$, where $M$ is the normal closure in $K$ of a set $S$, where $S=\{r(w_1,w_2,\dots) \mid r(x_1,x_2,\dots)\in R\}$.
So the set $S$ is a recursively enumerable
set of words in the generators $a,b$.

By CEP, we have $M=K\cap P$ for the normal closure $P$ of $S$ in $B(2,n)$. Hence we obtain that $G\cong K/M=K/K\cap P\cong KP/P$, where $KP/P$ is a subgroup of the two generated group
$B(2,n)/P$ and $P$ is the normal closure
in $B(2,n)$ of the recursively enumerable set of words $S$. By Theorem \ref{HE}, the group
$B(2,n)/P$ is, in turn, a subgroup of a finitely presented in ${\cal B}_n$ group $H$.

It remains to embed the group $H$ in a 2-generated finitely presented in ${\cal B}_n$
group. If $H\cong B(m,n)/L$, then we repeat the above trick using that $B(m,n)$ is embeddable in $B(2,n)$ as a CEP-subgroup. Indeed, the group $B(m,n)$ is a retract of $B(\infty,n)$, and so a CEP-subgroup. Since the CEP property is transitive, $B(m,n)$ is isomorphic to a
CEP subgroup $K'$ of $B(2,n)$.

The normal subgroup $L$ becomes a normal closure of a finite set $R'$ in $K'$. Therefore $L$ is the
intersection of $K'$ with the normal closure
$Q$  of $R'$ in $B(2,n)$.
Hence the group $H$ is embeddable in the $2$-generated finitely presented in ${\cal B}_n$ group $E=B(2,n)/Q$, as required.
\endproof

The next statement is similar to another Higman's theorem \cite{H}.

\begin{cy} \label{all} For every sufficiently large odd integer $n$,
there exists a 2-generated finitely presented in ${\cal B}_n$ group $E$ containing, as subgroups with pairwise
trivial intersections, isomorphic copies of all recursively presented groups $\{G_i\}_{i=1}^{\infty}$ of exponents dividing $n$.
\end{cy}

\proof The set of finite presentations $G_i=\langle A_i\mid  R_i\rangle$ in ${\cal B}_n$ is recursively enumerable. One can
assume that the finite sets $A_i$-s are disjoint and obtain a recursive presentation
$\langle X \mid \cup_i R_i\rangle$, where $X=\cup_i A_i=\{x_1,x_2,\dots\}$. By Dyck's
lemma the group $G$ defined by this presentation in ${\cal B}_n$ admits a retraction $\rho_i$
onto every $G_i$, where $\rho_i$ is identical on $A_i$
and maps all other generators from $X$ to $1$.
Therefore $G$ contains an isomorphic copy
of every group $G_i$. So does the $2$-generated finitely presented group $E\ge G$ given by Corollary \ref{twog}. By Theorem \ref{HE}, every recursively presented group from ${\cal B}_n$ is also embeddable in $E$.
\endproof

The next corollary shows that there are finitely presented in ${\cal B}_n$ groups with undecidable word problem. Recall that the first examples of a finitely presented groups with undecidable word problem were constructed by
W.W. Boone \cite{B} and P.S. Novikov \cite{N}.
M.Sapir \cite{Sa} obtained the first example of a finitely presented in a variety $\cal V$ group of finite exponent with undecidable word
problem. Then O. Kharlampovich \cite{Kh} found finitely presented in varieties ${\cal B}_n$ groups with undecidable word problem, provided $n=pr$, $p$ was an odd prime, and either $r$ had an odd divisor $\ge 665$ or $r\ge 2^{48}$. Sapir and Kharlampovich referred to the property that non-cyclic free Burnside groups of large exponent $n$ are infinite (\cite{NA, A, I94}), and both of them simulated the work of Minskiy's machine in a normal subgroup of smaller exponent. Therefore the exponent $n$ was a composite number in their examples. We obtain new examples now, in particular, for all large prime exponents:

\begin{cy}  For every sufficiently large
odd integer $n\ge C$, there exists a 2-generated finitely presented in ${\cal B}_n$ group with undecidable word problem.
\end{cy}

\proof We will use the properties of the subgroup $K\le B(2,n)$ defined in the proof
of Corollary \ref{twog}, and chose a recursively enumerable but non-recursive subset $T$ in the set of free generators  $W=\{w_1, w_2,\dots\}$. Let $N$ be the normal closure of $T$ in $B(2,n)$, and so the group $G=B(2,n)/N$ is recursively presented. A word $w_i\in W$ in the generators $a,b$ of $B(2,n)$
is trivial in $G$ iff it belongs to $N$.
%where $M$ is the normal closure of $T$ in $K$.
By CEP, this condition is equivalent to the property $w_i\in T$, because
$W$ is the free basis of $K$. Since the set $T$
is not recursive, the word problem in $G$ is undecidable. By Corollary \ref{twog}, the finitely generated group $G$ is a subgroup of a finitely presented in ${\cal B}_n$ group $H$ with two generators. Hence the word problem in
$H$ is undecidable too.
\endproof

The following proposition will be considered at the end of this paper since the proof of it is based on a few lemmas needed for the main theorem. Recall that a subgroup $A\le B$ is called \index[g]{Frattini embedded subgroup} {\it Frattini embedded} in the group $B$ (see \cite{Th}) provided
two elements of $A$ are conjugate in $B$ if and only if they are conjugate in $A$.

\begin{prop} \label{prop12} The embeddings $G\hookrightarrow H$, $G\hookrightarrow E$ and $G_i\hookrightarrow E$ given by Theorem \ref{HE}, Corollary \ref{twog} and Corollary \ref{all} enjoy the
following additional properties.

(1) The group $G$ (every group $G_i$) is a CEP subgroup in $H$ (resp., in $E$).

(2) The group $G$ (every group $G_i$) is Frattini embedded in $H$ (resp., in $E$).
\end{prop}

W.W. Boone and G. Higman  obtained a pure algebraic characterization of groups with decidable word problem. They proved in \cite{BH} that the word problem is decidable in a finitely generated group $G$ iff $G$ is embeddable in a simple subgroup of a finitely presented group. Similar characterization can be obtained in the varieties ${\cal B}_n$ ($n$ is a large odd). We do not include it in the present paper since the proof is based on a modification of the more elaborate techniques from \cite{book} (in particular, from Chapter 11).

\medskip

There are two main ingredients in the proof of Theorem \ref{HE}: diagram analysis 
of the consequences of defining relations and programming of machines producing 
the set of relations of the group $G$. The work of a Turing machine or of an S-machine 
has finite program; a proper interpretation of their commands in terms of finite sets of 
group relations leads to the proof of the classical Higman's theorem. In the present 
paper, we should make the "machine relations" compatible with Burnside relations so that
their union does not kill nontrivial elements of $G$. 

Recall that the first finitely generated infinite groups satisfying the Burnside identity
$x^n=1$ for a large odd exponent $n$ were constructed by P.S. Novikov and S.I. Adian
(see \cite{N59}, \cite{NA}, \cite{A}). Much shorter proof, although with a worse estimate
for the exponent $n$, was given in \cite{O82} (see \cite{book} for other results). 
Among the extensions and the applications of our approach, we emphasize the development 
presented in the joint paper of Olshanskii and Sapir \cite{OS03}, which is the important
tool of the current exposition too. 

As in \cite{book} and \cite{OS03}, the concept of A-map is permanently used in the
diagram part of this paper. To define it, we use so called contiguity submaps $\Gamma$ between 
two cells (or faces, regions) $\Pi_1$ and $\Pi_2$ of a planar map,
i.e. very thin submaps in comparison with the perimeters $|\partial\Pi_1|$ and $|\partial\Pi_2|$. 
Loosely speaking, the A-property says that $|\partial\Pi_1| << |\partial\Pi_2|$ if  $\Gamma$
is not too short in comparison with $|\partial\Pi_1|$. (See accurate definitions and basic
properties of $A$-maps in Section \ref{amaps}). Note that there are no groups or van 
Kampen diagrams in Chapter 5 of \cite{book}; the properties of A-maps obtained there
used global estimates for surface maps on regarded as metric spaces. Lemmas \ref{four}
and \ref{0cont} borrowed from \cite{book} make redundant the local analysis of cancellations
and inductive classification of periodic words typical for \cite{NA}, \cite{A}. In particular,
Lemma \ref{0cont} helps to prove that the free Burnside group $B(m,n)$ is infinite. This follows
from well-known Prochet - Thue - Morse examples \cite{P,T,M} of infinite sets of aperiodic 
words in 2-letter alphabet. 

We start with a version of Higman's embedding of a recursively presented group $G$ in a finitely presented group $\tilde G$ .
Then we impose Burnside relations $w^n=1$ on  $\tilde G$: at first we consider the words in so called
$a$-generators of $\tilde G$, then inductively introduce the periods $w$ depending on $a$- and $\theta$-generators, finally the periods depend on all the generators of $\tilde G$. At each of these steps,
the van Kampen diagram with standard metric do not satisfy the hyperbolic properties of $A$-map since the
inductively defined groups $G(i)$ have many abelian subgroups. However they becomes $A$-maps with respect
of special metrics. In particular, as in \cite{OS03}, we prescribe length  $0$ to every $a$-letter
when introducing periods $w$ depending on $a$- and $\theta$-letters only, and even $\theta$-letters
get length $0$ when all the generators are involved.

It was proved in \cite{OS03} that the inductive transitions ${\bf G}(0)\to\dots\to {\bf G}(i)\to\dots\to {\bf G}(\infty)$
with surjective homomorphism are successful, i.e. the diagrams over ${\bf G}(\infty)$ are $A$-maps, if the original presentations of ${\bf G}(0)$ and of $\tilde G = {\bf G}(1/2)$ have Properties (Z1), (Z2), (Z3)
(see Subsection \ref{axioms}).  Property (Z1) imposes some restrictions on the relations of ${\bf G}(0)$ implying
that the corresponding cells can be further regarded having rank $0$ in diagrams. Property (Z2) is imposed
on hub relations used for the Higman's embedding. (This property is empty, if the set of hubs empty.)
Property (Z3) is related to the known property of periodic words (see \cite{FW}) saying that simple periods 
of the same long periodic word are equal up to cyclic permutation. (Various inductive generalizations of such properties of periodic words were obtained in \cite{NA, A}; for more direct generalizations formulated in different terms, see \cite{book, OS03}.)
Property (Z3) reduces the information we need on periodic words to the following condition.

(*) {\it Let $W$ be a word of positive length which is not a conjugate of a shorter word in ${\bf G}(1/2)$,
and $X$ be a $0$-word (i.e. every letter of $X$ has length $0$). Assume that the word $W^{-4}XW^4$ is also equal in ${\bf G}(0)$ to a $0$-word $Y$. Then $(WX)^n=W^n$ in ${\bf G}(0)$.} 

(The exponent $4$ can be replaced with $3$, but here we use the exponent from \cite{OS03}.)

To obtain Condition (*), we simplify the diagram $\Delta$ corresponding to the equality $W^{-4}XW^4=Y$.
Changing the pair $(W,X)$ by conjugate pairs in ${\bf G}(0)$, we step-by-step remove unwanted cells
until $\Delta$ becomes a diagram over the group $M$ whose relations simulate the work of the machine
$\bf M$ recognizing the relations of $G$. After further transformations one can find trapezia in $\Delta$
corresponding to some computations of ${\bf M}$. The labels of "horizontal" lines of a trapezium are the words $W_i$-s in a computation $W_0\to W_1\to\dots\to W_k$, where $W_{i}$ is obtained from $W_{i-1}$ by an application of
some rule $\theta_i$, and the product $\theta_1\dots\theta_k$ is called the history $h$ of this computation. It turns out that Condition (*) is a consequence of the properties of $M$-computations. For example:

(**) {\it If two computations with histories $h$ and $ghg^{-1}$ both start and end with the same word $W$
(under some restrictions), then $g$ is equal, modulo the Burnside relations, to a word $g'$ such that there 
is a computation with history $g'$  starting and ending with the word $W$.}

The designing of a machine $\bf M$ recognizing the recursively enumerable set of defining relations
of given group $G$ and satisfying conditions like (**), is the central point of this paper. The crucial issue
and the important difference in comparison with \cite{OS03} is the following. In the beginning, we know nothing about the Turing machine enumerating the defining relations of $G$, except for the {\it existence} of such a machine,
whereas the S-machine $\cal M$ from \cite{OS03} was explicitly defined to produce the words of the form $w^n$
in some alphabet. (Note that S-machines invented by M. Sapir in \cite{SBR} are more suitable for many group-theoretic applications than classical Turing machines working with positive words only, see \cite{OS01}.)

Now we start with a Turing machine $M_0$ having a few additional properties (e.g., no configuration can appear
in a computation of $M_0$ two times). Then each command of $M_0$ should be transformed in a set of rules of
an $S$-machine. Unfortunately, known constructions become insecure when one interprets  the computations of
 S-machine by trapezia, because the cells corresponding to the relations of $G$ and to the Burnside
relations can be inserted between the horizontal bands, which distorts the computation. In \cite{OS03},
we had very concrete S-machine which produced arbitrary word $w$, and a after number of properly defined rewritings
one obtained the required relator $w^n$. Those  rewritings were organized and controlled so that 
arbitrary interference of unwanted cells  replaced $w^n$ with some word $(w')^n$, which was also a relator.
But now $M_0$ looks like a black box device without information about its internal workings.

Therefore we suggest a new way of replacement of a Turing machine $M_0$ with an equivalent S-machine $M_1$. For every command $\theta$ of $M_0$, we construct an S-machine $M(\theta)$, and for every
configuration $\cal W$ of $M_0$, we define a configuration $W=F({\cal W})$ of $M_1$, where the inductive definition of the function $F$ is based on aperiodic endomorphisms $\psi_j$ of 2-generated free semigroup.
For example, if $\theta$ inserts a letter $a_j$ at the end of some tape, then the canonical word
of $M(\theta)$ replaces the whole tape subword $w$ of $W$ with a word $\phi_j(w)={\bf w_j}\psi_j(w)$,
where ${\bf w}_j$ is an aperiodic word with a small cancellation property.
The S-machine $M_1$ is  a union of all $M(\theta)$-s. 

The machine $M_1$ has a number of helpful properties.
For instance, if some tape has a word $w$ in the beginning of a computation and the word $w'$ at the end,
then either $w'$ and $w$ are equal in the free group or they have different canonical images in the free Burnside group. This and other features make $M_1$ interference-free, i.e. arbitrary replacements
of the tape words with equal words modulo the Burnside relations, during a computation changes the output with
an equal word modulo the Burnside relations too.

Many copies of $M_1$ gives us $M_2$, since one needs this modification to obtain Property (Z2) later. Then we should add a  machine
$T$ which translates the language recognized by $M_0$ into the language recognized by $M_1$.
Such a modification adds one more tape which does not belong to $M_2$, but $T$ just cleans 
up this tape, and so it is easy to control possible unwanted ``noise'' during a computation.

Note that $M_1$ and the main machine $\bf M$ do not satisfy the usual definition of $S$-machine
from \cite{SBR}, \cite{OS01}, \cite{OS19}. So we have to give  a more general definition, where the
tape alphabets of $\theta$-admissible words are replaced with finitely generated subgroups (depending 
on the rule $\theta$) of free groups. As a consequence, some rules cannot be written in "`inserting"' form
now, and a rule $\theta$ may insert a word in sectors locked by $\theta$.

In sections \ref{gd} and \ref{more}, we get back to groups and diagrams and obtain the required embedding 
$G\hookrightarrow H$ starting with the Higman-type embedding $G\hookrightarrow \tilde G$, where the finite set of defining relations of $\tilde G$ is based on the rules of the machine $\bf M$.
It turns out that the canonical mapping of $G$ to the Burnside quotient $H$  of $\tilde G$ is injective, which completes the proof of Theorem \ref{HE}. Section \ref{more} is closer to Sections 8 - 10 of \cite{OS03}, although some proofs are different and shorter now.

\section{Maps and diagrams}
\label{amaps}

Planar maps and \vk diagrams are standard tools to study the presentations of
complicated groups (see, for example, \cite{LS} and \cite{book}).
Here we present the main concepts related to maps and diagrams. A
\vk diagram is a labeled map, so we start with discussing maps.

\subsection{Graphs and maps}
\label{maps}

We are using the standard definition of a \index[g]{graph}{\em
graph}. In particular, every edge has a direction, and every edge $\bf e$
has an inverse edge ${\bf e}^{-1}$ (having the opposite direction). For every
edge ${\bf e}$, \index[g]{${\bf e}_-$ and ${\bf e}_+$}${\bf e}_-$ and ${\bf e}_+$ are the beginning and the
end vertices of the edge.

A \index[g]{map on a surface} {\em map} on a surface $X$ (in this paper, on a disk or on an annulus) is simply a finite connected graph drawn on $X$ which subdivides this surface into polygonal
2-cells (= cells). Edges do not have labels, so no group presentations are
involved in studying maps. On the other hand the properties of maps
help finding the structure of \vk diagrams because every diagram
becomes a map after we remove all the labels.

Let us recall the necessary definitions from \cite{book}.

A map is called \index[g]{graded map}{\em graded} if every cell $\Pi$ is
assigned a non-negative number, $r(\Pi)$, its \index[g]{rank of a cell}{\em
rank}. A map $\Delta$ is called a \index[g]{rank of map}{\em map of rank}
$i$ if all its cells have ranks $\le i$. Every graded map $\Delta$
has a \index[g]{type of a map}{\it type} $\tau (\Delta)$ which is the vector
whose first coordinate is the maximal rank of a cell in the map,
say, $r$, the second coordinate is the number of cells of rank $r$
in the map, the third coordinate is the number of cells of maximal rank
$<r$ and so on. We compare types lexicographically.

Let $\Delta$ be a graded map on a surface $X$. The cells of rank
$0$ in $\Delta$ are called \index[g]{$0$-cell} $0$-{\em cells}. Some of the
edges in $\Delta$ are called \index[g]{$0$-edge} $0$-{\em edges}. Other
edges and cells will be called \index[g]{positive edge}{\em positive}. For every path ${\bf p}$ in $\Delta$, there are the beginning and the end vertices ${\bf p}_-$ and $ {\bf p}_+$. We
define the \index[g]{length of a path} {\em length of a path} ${\bf p}$ in a map as the number of
positive edges in it. The length of a path ${\bf p}$ is denoted by
$|{\bf p}|$, in particular, $|\partial(\Pi)|$ is the \index[g]{perimeter} {\it perimeter}
perimeter of a cell $\Pi$ and $|\partial(\Delta)|$ is the perimeter of a disk map $\Delta$.

We assume that \index[g]{Properties (M1), (M2), (M3)}

\begin{itemize}
\item[(M1)] the inverse edge of a $0$-edge is also a $0$-edge,
\item[(M2)] if $\Pi$ is a $0$-cell, then we have either $|\partial\Pi|=0$ or $|\partial\Pi|=2$,
\item[(M3)] for every cell $\Pi$ of rank $r(\Pi)>0$ the length
$|\partial(\Pi)|$ of its contour is positive.
\end{itemize}

\subsection{Bands}\label{bds}

\label{prelim}

A \index[g]{disk diagram} Van Kampen {\it diagram} over a presentation $\la {\cal X}\ |\
{\cal R}\ra$ is a planar map with edges labeled by elements from
${\cal X}^{\pm 1}$ (and therefore paths labeled by words), such that the label of the contour of each
cell belongs to ${\cal R}$ up to taking cyclic shifts and inverses
\cite{book}. In $\S 11$ of \cite{book}, 0-edges of \vk diagrams
have label 1, but in \cite{OS03} and in this paper, 0-edges will be labeled also by
the so called \index[g]{0-letter} {\em 0-letters}. All diagrams which we
shall consider in this paper, when considered as planar maps, will
obviously satisfy conditions (M1), (M2), (M3) from Section
\ref{maps}. By \index[g]{van Kampen Lemma}{van Kampen Lemma}, a word $w$ over ${\cal X}$ is
equal to 1 modulo ${\cal R}$ if and only if there exists a \vk
disk diagram over ${\cal R}$ with boundary label $w$. Similarly
(see \cite{LS}, \cite{book}) two words $w_1$ and $w_2$ are conjugate modulo $R$
if and only if there exists an annular diagram with the internal
counterclockwise contour labeled by $w_1$ and external
counterclockwise contour labeled by $w_2$.

There may be many
diagrams with the same boundary label. Sometimes we can reduce the
number of cells in a diagram or ranks of cells by replacing a
2-cell subdiagram by another subdiagram with the same boundary
label. The simplest such situation is when there are two positive cells in
a diagram which have a common edge and are mirror images of each
other, i.e. the labels of the boundary paths of these cells starting with this edge are equal. In this case one can {\em cancel} these two cells (see
\cite{LS}). In order strictly define this cancellation (to
preserve the topological type of the diagram), one needs to use
the so called $0$-refinement. The $0$-refinement consists of
adding $0$-cells with (freely trivial) boundary labels of the form
$1aa^{-1}$ or $111$ (see \cite{book}, $\S 11.5$). A diagram having no
such mirror pairs of positive cells is called \index[g]{reduced diagram} {\it reduced}.

As in the book \cite{R} and our the papers \cite{SBR}, \cite{BORS} and
\cite{O97}, one of the main tools to study \vk diagrams are
bands and annuli.

Let $S$ be a subset of ${\cal X}$. An \index[g]{band} $S$-{\it band} $\bb$
is a sequence of cells $\pi_1,...,\pi_n$, for some $n\ge 0$ in a
\vk diagram such that

\begin{itemize}
\item Each two consecutive cells in this sequence have a common edge
labeled by a letter from $S$.
\item Each cell $\pi_i$, $i=1,...,n$ has exactly two $S$-edges
(i.e. edges labeled by a letter from $S$).
\item If $n=0$, then the boundary of $S$ has form ${\bf ee}^{-1}$ for
an $S$-edge ${\bf e}$.
\end{itemize}

\setcounter{pdten}{\value{ppp}} Figure \thepdten\ illustrates this
concept. In this Figure, edges ${\bf e}, {\bf e}_1,...,{\bf e}_{n-1},{\bf f}$ are
$S$-edges, the lines $l(\pi_i,{\bf e}_i), l(\pi_i,{\bf e}_{i-1})$ connect
fixed points in the cells with fixed points of the corresponding
edges.

\bigskip
\unitlength=0.90mm \special{em:linewidth 0.4pt}
\linethickness{0.4pt}
\begin{picture}(149.67,30.11)
\put(19.33,30.11){\line(1,0){67.00}}
\put(106.33,30.11){\line(1,0){36.00}}
\put(142.33,13.11){\line(-1,0){35.67}}
\put(86.33,13.11){\line(-1,0){67.00}}
\put(33.00,21.11){\line(1,0){50.00}}
\put(110.00,20.78){\line(1,0){19.33}}
\put(30.00,8.78){\vector(1,1){10.33}}
\put(52.33,8.44){\vector(0,1){10.00}}
\put(76.33,8.11){\vector(-1,1){10.33}}
\put(105.66,7.78){\vector(1,2){5.33}}
\put(132.66,8.11){\vector(-1,1){10.00}}
\put(16.00,21.11){\makebox(0,0)[cc]{${\bf e}$}}
\put(29.66,25.44){\makebox(0,0)[cc]{$\pi_1$}}
\put(60.33,25.44){\makebox(0,0)[cc]{$\pi_2$}}
\put(133.00,25.78){\makebox(0,0)[cc]{$\pi_n$}}
\put(145.33,21.44){\makebox(0,0)[cc]{${\bf f}$}}
\put(77.00,25.11){\makebox(0,0)[cc]{${\bf e}_2$}}
\put(122.00,25.44){\makebox(0,0)[cc]{${\bf e}_{n-1}$}}
\put(100.33,25.44){\makebox(0,0)[cc]{$S$}}
\put(96.33,30.11){\makebox(0,0)[cc]{$\dots$}}
\put(96.33,13.11){\makebox(0,0)[cc]{$\dots$}}
\put(26.66,4.78){\makebox(0,0)[cc]{$l(\pi_1,{\bf e}_1)$}}
\put(52.66,4.11){\makebox(0,0)[cc]{$l(\pi_2,{\bf e}_1)$}}
\put(78.33,4.11){\makebox(0,0)[cc]{$l(\pi_2,{\bf e}_2)$}}
\put(104.66,3.78){\makebox(0,0)[cc]{$l(\pi_{n-1},{\bf e}_{n-1})$}}
\put(134.00,3.11){\makebox(0,0)[cc]{$l(\pi_n,{\bf e}_{n-1})$}}
\put(88.66,30.11){\makebox(0,0)[cc]{${\bf q}_2$}}
\put(88.66,13.11){\makebox(0,0)[cc]{${\bf q}_1$}}
\put(49.33,25.11){\makebox(0,0)[cc]{${\bf e}_1$}}
\put(33.33,21.11){\circle*{0.94}}
\put(46.66,21.11){\circle*{1.33}}
\put(60.33,21.11){\circle*{0.94}}
\put(73.66,21.11){\circle*{1.33}}
\put(95.00,20.56){\makebox(0,0)[cc]{...}}
\put(19.33,21.00){\line(1,0){14.00}}
\put(129.67,20.67){\line(1,0){12.67}}
\put(19.00,13.00){\vector(0,1){17.00}}
\put(46.67,13.00){\vector(0,1){17.00}}
\put(73.67,13.00){\vector(0,1){17.00}}
\put(117.00,13.00){\vector(0,1){17.00}}
\put(117.00,20.67){\circle*{1.33}}
\put(129.67,20.67){\circle*{0.67}}
\put(142.33,13.00){\vector(0,1){17.00}}
\put(142.33,20.67){\circle*{1.33}}
\put(19.00,21.00){\circle*{1.33}}
\put(14.00,8.33){\vector(1,1){11.00}}
\put(10.33,5.00){\makebox(0,0)[cc]{$l(\pi_1,{\bf e})$}}
\put(148.00,8.00){\vector(-1,1){10.33}}
\put(154.67,3.11){\makebox(0,0)[cc]{$l(\pi_n,{\bf f})$}}
\end{picture}

\begin{center}
\nopagebreak[4] Fig. \theppp.

\end{center}
\addtocounter{ppp}{1}

The broken line formed by the lines $l(\pi_i,{\bf e}_i)$,
$l(\pi_i,{\bf e}_{i-1})$ connecting points inside neighboring cells is
called the \index[g]{connecting line} {\em connecting line}
$\bb$. The $S$-edges ${\bf e}$ and ${\bf f}$ are called the \index[g]{start/end edges of a band} {\em
start} and {\em end} edges of the band.

The counterclockwise boundary of the subdiagram formed by the
cells $\pi_1,...,\pi_n$ of $\bb$ has the form ${\bf e}^{-1} {\bf q}_1{\bf f} {\bf q}_2^{-1}$.
We call ${\bf q}_1$ the \index[g]{bottom/top of a band} {\em bottom} of $\bb$ and ${\bf q}_2$
the {\em top} of $\bb$, denoted $\bott(\bb)$ and
$\topp(\bb)$. (Both of them are the \index[g]{side of a band} {\it sides}
of $\bb$.)

We say that two bands $\cal S$ and $\cal T$ \index[g]{crossing bands}{\em cross} if their connecting
lines cross (see
\setcounter{pdeleven}{\value{ppp}} Figure \thepdeleven b). We say that a band is an \index[g]{annulus}{\em annulus}
if its connecting line is a closed curve. In this case the start
and the end edges of the band coincide (see
\setcounter{pdeleven}{\value{ppp}} Figure \thepdeleven a).

\bigskip
\begin{center}
\unitlength=1.5mm
%%\special{em:linewidth 0.4pt}
\linethickness{0.4pt}
\begin{picture}(101.44,22.89)
\put(30.78,13.78){\oval(25.33,8.44)[]}
\put(30.78,13.89){\oval(34.67,18.00)[]}
\put(39.67,9.56){\line(0,-1){4.67}}
\put(25.67,9.56){\line(0,-1){4.67}}
\put(18.56,14.89){\line(-1,0){5.11}}
\put(32.56,7.11){\makebox(0,0)[cc]{$\pi_1=\pi_n$}}
\put(21.00,7.33){\makebox(0,0)[cc]{$\pi_2$}}
\put(15.89,11.78){\makebox(0,0)[cc]{$\pi_3$}}
\put(43.44,12.89){\line(1,0){4.67}}
\put(43.00,8.44){\makebox(0,0)[cc]{$\pi_{n-1}$}}
\put(19.22,10.89){\line(-1,-1){4.00}}
\put(25.89,20.44){\circle*{0.00}}
\put(30.78,20.44){\circle*{0.00}}
\put(30.33,1.56){\makebox(0,0)[cc]{a}}
\put(35.44,20.44){\circle*{0.00}}
\put(62.78,4.89){\line(0,1){4.67}}
\put(62.78,9.56){\line(1,0){38.67}}
\put(101.44,9.56){\line(0,-1){4.67}}
\put(101.44,4.89){\line(-1,0){38.67}}
\put(82.11,16.22){\oval(38.67,13.33)[t]}
\put(62.78,15.78){\line(0,-1){7.33}}
\put(101.44,16.44){\line(0,-1){8.22}}
\put(82.00,12.22){\oval(27.78,14.67)[t]}
\put(67.89,12.89){\line(0,-1){8.00}}
\put(95.89,13.11){\line(0,-1){8.22}}
\put(67.89,13.78){\line(-1,1){4.67}}
\put(95.67,14.00){\line(1,1){5.11}}
\put(73.67,9.56){\line(0,-1){4.67}}
\put(89.67,9.56){\line(0,-1){4.67}}
\put(76.56,21.33){\circle*{0.00}}
\put(81.44,21.33){\circle*{0.00}}
\put(86.11,21.33){\circle*{0.00}}
\put(65.22,12.22){\makebox(0,0)[cc]{$\pi_1$}}
\put(65.22,6.89){\makebox(0,0)[cc]{$\pi$}}
\put(71.00,6.89){\makebox(0,0)[cc]{$\gamma_1$}}
\put(92.78,7.11){\makebox(0,0)[cc]{$\gamma_m$}}
\put(98.56,7.11){\makebox(0,0)[cc]{$\pi'$}}
\put(98.33,12.67){\makebox(0,0)[cc]{$\pi_n$}}
\put(94.11,20.22){\makebox(0,0)[cc]{$\cal S$}}
\put(85.89,6.89){\makebox(0,0)[cc]{$\cal T$}}
\put(76.56,6.89){\circle*{0.00}} \put(79.44,6.89){\circle*{0.00}}
\put(82.33,6.89){\circle*{0.00}}
\put(83.44,1.56){\makebox(0,0)[cc]{b}}
\end{picture}
\end{center}
\begin{center}
\nopagebreak[4] Fig. \theppp.

\end{center}
\addtocounter{ppp}{1}

\subsection{Bonds and contiguity submaps}

 \label{bcs}

Now let us introduce the main concepts from Chapter 5 of
\cite{book}.  Consider a map $\Delta$. A path in $\Delta$ is
called \index[g]{reduced path} {\em reduced} if it does not contain
subpaths of length 2 which are homotopic to paths of length 0 by
homotopies involving only 0-cells. We shall usually suppose that
the contour of $\Delta$ is subdivided into several subpaths called
\index[g]{section of a boundary}{\em sections}. We usually assume that sections are
reduced paths.

By property (M2), if a 0-cell has a positive edge, it has exactly two
positive edges. Thus we can consider bands of 0-cells having
positive edges. The start and end edges of such bands are called
\index[g]{adjacent edges} {\em adjacent}. (We obtain a simple generalization of the
cancellation of a pair of positive cells in a diagram replacing the common edge
of their boundaries with two adjacent edges.)

Let ${\bf e}$ and ${\bf f}$ be adjacent edges of $\Delta$ belonging to two
positive cells $\Pi_1$ and $\Pi_2$ or to sections of the contour
of $\Delta$. Then there exists a band of $0$-cells connecting
these edges. The union of the cells of this band is called a
\index[g]{$0$-bond}{\em $0$-bond} between the cells $\Pi_1$ and $\Pi_2$
(or between a cell and a section of the contour of $\Delta$, or
between two sections of the contour). 
The contour of the $0$-bond has the form ${\bf p}\iv {\bf e}{\bf s} {\bf f}\iv$, where
${\bf p}$ and ${\bf s}$ are paths of length 0 because every $0$-cell has at
most two positive edges.

Now suppose we have chosen two pairs $({\bf e}_1, {\bf f}_1)$ and $({\bf e}_2, {\bf f}_2)$
of adjacent edges such that ${\bf e}_1, {\bf e}_2$ belong to the
contour of a cell $\Pi_1$ and ${\bf f}_1^{-1}, {\bf f}_2^{-1}$ to the contour of a
cell $\Pi_2$. Then we have two bonds $E_1$ and $E_2$. If $E_1=E_2$
(that is ${\bf e}_1={\bf e}_2, {\bf f}_1={\bf f}_2$) then let $\Gamma=E_1$. Now let
$E_1\ne E_2$ and ${\bf z}_1{\bf e}_1{\bf w}_1{\bf f}_1^{-1}$ and ${\bf z}_2{\bf e}_2{\bf w}_2{\bf f}_2^{-1}$  be the
contours of these bonds. Further let ${\bf y}_1$ and ${\bf y}_2$ be subpaths
of the contours of $\Pi_1$ and $\Pi_2$ (of $\Pi_1$ and ${\bf q}$, or
${\bf q}_1$ and ${\bf q}_2$, where ${\bf q}$, ${\bf q}_1$ and ${\bf q}_2$ are some segments of
the boundary contours) where ${\bf y}_1$ (or ${\bf y}_2$) has the form
${\bf e}_1{\bf p}{\bf e}_2$ or ${\bf e}_2{\bf p}{\bf e}_1$ (or $({\bf f}_1{\bf u}{\bf f}_2)^{-1}$ or $({\bf f}_2{\bf u}{\bf f}_1)^{-1}$).
If ${\bf z}_1{\bf y}_1{\bf w}_2{\bf y}_2$ (or ${\bf z}_2{\bf y}_1{\bf w}_1{\bf y}_2$) is a contour of a disk
submap $\Gamma$ which does not contain $\Pi_1$ or $\Pi_2$, then
$\Gamma$ is called a $0$-{\em contiguity submap} of $\Pi_1$ to
$\Pi_2$ (or $\Pi_1$ to ${\bf q}$ or of ${\bf q}_1$ to ${\bf q}_2$). The contour of
$\Gamma$ is naturally subdivided into four parts. The paths ${\bf y}_1$
and ${\bf y}_2$ are called the {\em contiguity arcs}. We write
${\bf y}_1=\Gamma \bigwedge \Pi_1$, ${\bf y}_2=\Gamma \bigwedge \Pi_2$ or
${\bf y}_2=\Gamma \bigwedge {\bf q}$. The other two paths are called the {\em
side arcs} of the \ct submap.

The ratio $|{\bf y}_1|/|\partial(\Pi_1)|$ (or $|{\bf y}_2|/|\partial(\Pi_2)|$
is called the {\em degree of contiguity} of the cell $\Pi_1$ to
the cell $\Pi_2$ or to ${\bf q}$ (or of $\Pi_2$ to $\Pi_1$). We denote
the degree of contiguity of $\Pi_1$ to $\Pi_2$ (or $\Pi_1$ to ${\bf q}$)
by $(\Pi_1,\Gamma,\Pi_2)$ (or $(\Pi_1,\Gamma,{\bf q})$). Notice that
this definition is not symmetric and when $\Pi_1=\Pi_2=\Pi$, for
example, then $(\Pi,\Gamma,\Pi)$ is a pair of numbers.

We say that two contiguity submaps are \index[g]{disjoint contiguity submaps}{\em disjoint} if none of
them has a common point with the interior of the other one, and
their contiguity arcs do not have common edges.

Now we are going to define $k$-bonds and $k$-contiguity submaps
for $k>0$. In these definitions we need a fixed number
$\varepsilon$, $0<\varepsilon< 1$.

Let $k>0$ and suppose that we have defined the concepts of
$j$-bond and $j$-contiguity submap for all $j<k$. Consider three
cells $\pi, \Pi_1, \Pi_2$ (possibly with $\Pi_1=\Pi_2$) satisfying
the following conditions:

\begin{enumerate}
\item $r(\pi)=k, r(\Pi_1)>k, r(\Pi_2)>k$,
\item there are disjoint submaps $\Gamma_1$ and $\Gamma_2$ of
$j_1$-contiguity of $\pi$ to $\Pi_1$ and of $j_2$-contiguity of
$\pi$ to $\Pi_2$, respectively, with $j_1<k, j_2<k$, such that
$\Pi_1$ is not contained in $\Gamma_2$ and $\Pi_2$ is not
contained in $\Gamma_1$,
\item $(\pi,\Gamma_1,\Pi_1)\ge\varepsilon$,
$(\pi,\Gamma_2,\Pi_2)\ge\varepsilon$.
\end{enumerate}

Then there is a minimal submap $E$ in $\Delta$ containing $\pi,
\Gamma_1, \Gamma_2$. This submap is called a \index[g]{$k$-bond}{\em
$k$-bond between $\Pi_1$ and $\Pi_2$ defined by the contiguity
submaps $\Gamma_1$ and $\Gamma_2$ with principal cell $\pi$} (see
\setcounter{pdeleven}{\value{ppp}} Figure \thepdeleven).

The {\em contiguity arc} $q_1$ of the bond $E$ to $\Pi_1$ is
$\Gamma_1\bigwedge \Pi_1$. It will be denoted by $E\bigwedge
\Pi_1$. Similarly $E\bigwedge\Pi_2$ is by definition the arc
${\bf q}_2=\Gamma_2\bigwedge\Pi_2$. The contour $\partial(E)$ can be
written in the form ${\bf p}_1{\bf q}_1{\bf p}_2{\bf q}_2$ where ${\bf p}_1$, ${\bf p}_2$ are called
the {\em side arcs} of the bond $E$.

\unitlength=1mm \special{em:linewidth 0.4pt} \linethickness{0.4pt}
\begin{picture}(90.33,55.67)
\put(67.33,40.33){\oval(46.00,9.33)[b]}
\put(44.33,42.00){\line(0,1){3.00}}
\put(44.33,46.67){\line(0,1){3.33}}
\put(90.33,42.00){\line(0,1){3.00}}
\put(90.33,46.67){\line(0,1){3.33}}
\put(47.17,52.33){\oval(5.00,5.33)[lt]}
\put(86.17,51.67){\oval(8.33,8.00)[rt]}
\put(67.50,23.17){\oval(29.67,11.00)[]}
\put(58.67,28.67){\line(0,1){7.00}}
\put(74.33,28.67){\line(0,1){7.00}}
\put(67.33,5.33){\oval(46.00,10.67)[t]}
\put(60.67,17.67){\line(0,-1){7.00}}
\put(71.67,10.67){\line(0,1){7.00}}
\put(58.00,47.00){\makebox(0,0)[cc]{$\Pi_1$}}
\put(68.67,38.33){\makebox(0,0)[cc]{${\bf q}_1$}}
\put(62.33,7.33){\makebox(0,0)[cc]{${\bf q}_2$}}
\put(80.00,4.00){\makebox(0,0)[cc]{$\Pi_2$}}
\put(71.83,10.67){\line(0,1){7.00}}
\put(71.67,17.50){\line(1,0){5.50}}
\put(75.42,23.08){\oval(14.17,11.17)[r]}
\put(74.50,28.83){\line(1,0){2.67}}
\put(77.17,28.83){\line(0,0){0.00}}
\put(74.50,28.83){\line(0,1){6.67}}
\put(58.50,35.67){\line(0,-1){7.00}}
\put(58.50,28.67){\line(-1,0){0.83}}
\put(58.17,23.17){\oval(11.33,11.33)[l]}
\put(58.33,17.50){\line(1,0){2.17}}
\put(60.50,17.50){\line(0,-1){6.83}}
\put(60.50,15.33){\vector(0,-1){2.00}}
\put(74.50,30.33){\vector(0,1){2.67}}
\put(66.00,10.67){\vector(-1,0){2.67}}
\put(65.00,35.67){\vector(1,0){2.83}}
\put(70.00,27.17){\makebox(0,0)[cc]{${\bf v}_1$}}
\put(68.83,19.83){\makebox(0,0)[cc]{${\bf v}_2$}}
%\bezier{108}(63.67,35.67)(60.33,20.50)(67.00,10.67)
\put(61.00,23.00){\makebox(0,0)[cc]{$t$}}
\put(76.00,23.00){\makebox(0,0)[cc]{$\pi$}}
\end{picture}

\begin{center}
\nopagebreak[4] Fig. \theppp.

\end{center}
\addtocounter{ppp}{1}

Bonds between a cell and a section of the contour of $\Delta$ or
between two sections of the contour are defined in a similar way.

Now let $E_1$ be a $k$-bond and $E_2$ be a $j$-bond between two
cells $\Pi_1$ and $\Pi_2$, $j\le k$ and either $E_1=E_2$ or these
bonds are disjoint. If $E_1=E_2$ then $\Gamma=E_1=E_2$ is called
the $k$-contiguity submap between $\Pi_1$ and $\Pi_2$ determined
by the bond $E_1=E_2$. If $E_1$ and $E_2$ are disjoint then the
corresponding \index[g]{contiguity submap/subdiagram}{\em contiguity submap} $\Gamma$ is
defined as the smallest submap containing $E_1$ and $E_2$, bounded by these submaps and  segments of $\partial\Pi_1$ and $\partial\Pi_2$, and not
containing $\Pi_1$ and $\Pi_2$ (see
\setcounter{pdeleven}{\value{ppp}} Figure \thepdeleven).

\unitlength=1mm \special{em:linewidth 0.4pt} \linethickness{0.4pt}
\begin{picture}(135.00,76.00)
\put(73.50,68.00){\oval(123.00,16.00)[b]}
\put(73.50,13.00){\oval(123.00,20.00)[t]}
\put(35.00,41.00){\oval(26.00,16.00)[]}
\put(107.00,41.00){\oval(28.00,16.00)[]}
\put(30.00,54.00){\oval(6.00,4.00)[]}
\put(43.00,54.00){\oval(6.00,4.00)[]}
\put(102.50,54.00){\oval(7.00,4.00)[]}
\put(114.50,54.00){\oval(7.00,4.00)[]}
\put(29.00,27.00){\oval(6.00,4.00)[]}
\put(41.00,27.00){\oval(6.00,4.00)[]}
\put(102.00,27.00){\oval(6.00,4.00)[]}
\put(114.00,27.00){\oval(6.00,4.00)[]}
\put(28.00,29.00){\line(0,1){4.00}}
\put(31.00,33.00){\line(0,-1){4.00}}
\put(39.00,29.00){\line(0,1){4.00}}
\put(42.00,33.00){\line(0,-1){4.00}}
\put(28.00,25.00){\line(0,-1){2.00}}
\put(31.00,23.00){\line(0,1){2.00}}
\put(39.00,25.00){\line(0,-1){2.00}}
\put(43.00,23.00){\line(0,1){2.00}}
\put(28.00,49.00){\line(0,1){3.00}}
\put(32.00,52.00){\line(0,-1){3.00}}
\put(44.00,52.00){\line(1,-2){3.00}}
\put(41.00,56.00){\line(-1,4){1.00}}
\put(45.00,60.00){\line(-1,-4){1.00}}
\put(29.00,56.00){\line(-1,4){1.00}}
\put(33.00,60.00){\line(-1,-2){2.00}}
\put(101.00,56.00){\line(-1,4){1.00}}
\put(107.00,60.00){\line(-3,-4){3.00}}
\put(113.00,56.00){\line(0,1){4.00}}
\put(116.00,60.00){\line(0,-1){4.00}}
\put(99.00,33.00){\line(1,-2){2.00}}
\put(104.00,29.00){\line(1,4){1.00}}
\put(100.00,25.00){\line(0,-1){2.00}}
\put(100.00,23.00){\line(0,0){0.00}}
\put(104.00,23.00){\line(0,1){2.00}}
\put(113.00,23.00){\line(0,1){2.00}}
\put(115.00,25.00){\line(0,-1){2.00}}
\put(113.00,29.00){\line(-1,2){2.00}}
\put(114.00,33.00){\line(1,-4){1.00}}
\put(101.00,52.00){\line(0,-1){3.00}}
\put(104.00,49.00){\line(0,1){3.00}}
\put(41.00,52.00){\line(-5,-3){5.00}}
\put(113.00,52.00){\line(-4,-3){4.00}}
\put(116.00,52.00){\line(2,-3){4.00}}
\end{picture}
\begin{center}
\nopagebreak[4] Fig. \theppp.

\end{center}
\addtocounter{ppp}{1}
 The \index[g]{contiguity arcs}{\em contiguity arcs} ${\bf q}_1$ and ${\bf q}_2$ of $\Gamma$ are
the intersections of $\partial(\Gamma)$ with $\partial(\Pi_1)$ and
$\partial(\Pi_2)$. The contour of $\Gamma$ has the form
${\bf p}_1{\bf q}_1{\bf p}_2{\bf q}_2$ where ${\bf p}_1$ and ${\bf p}_2$ are called the
{\em side arcs} of $\Gamma$. 
The ratio $|{\bf q}_1|/|\partial(\Pi_1)|$ is called the
{contiguity degree}{\em contiguity degree} of $\Pi_1$ to $\Pi_2$ with
respect to $\Gamma$ and is denoted by
\index[g]{$(\Pi_1,\Gamma,\Pi_2)$} $(\Pi_1,\Gamma,\Pi_2)$. If $\Pi_1=\Pi_2=\Pi$ then
$(\Pi,\Gamma,\Pi)$ is a pair of numbers.

Contiguity submaps of a cell to a section of the contour and
between sections of the contour are defined in a similar way. (In
this paper we do not use notion  "degree of contiguity of a
section of a contour to" anything.)

We are going to write ``bonds" and ``\ct submaps" in place of
``$k$-bonds" and ``$k$-contiguity submaps". Instead of writing
``the \ct submap $\Gamma$ of a cell $\Pi$ to ...", we sometimes
write ``the $\Gamma$-\ct of $\Pi$ to ...".

The above definition involved the standard decomposition of the
contour of a \ct submap $\Gamma$ into four sections ${\bf p}_1{\bf q}_1{\bf p}_2{\bf q}_2$
where \index[g]{wedge}${\bf q}_1=\Gamma\bigwedge \Pi_1$,
${\bf q}_2=\Gamma\bigwedge \Pi_2$ (or ${\bf q}_2=\Gamma\bigwedge {\bf q}$ if ${\bf q}$ is
a section of $\partial(\Delta)$), we shall write
\index[g]{$\partial(\Pi_1,\Gamma,\Pi_2)$} ${\bf p}_1{\bf q}_1{\bf p}_2{\bf q}_2=\partial(\Pi_1,\Gamma,\Pi_2)$, and so on.

As in \cite{book} (see $\S 15.1$) we fix certain real numbers
\index[g]{parameters} {\it parametes} $$\iota << \zeta << \varepsilon <<\delta << \gamma <<
\beta << \alpha$$ between 0 and 1 where "$<<$" means ``much
smaller". Here ``much" means enough to satisfy all the
inequalities in Chapters 5, 6 of \cite{book}. We also set
$$\bar\alpha=\frac{1}{2}+\alpha, \;\;\bar\beta=1-\beta\;\;,
\bar\gamma=1-\gamma,\;\; h=\delta\iv,\;\; n=\iota\iv$$.

\subsection{Condition A}

\label{conditionA}

The set of positive cells of a map $\Delta$ is denoted by
$\Delta(2)$ and as before the length of a path in $\Delta$ is the
number of positive edges in the path. A path ${\bf p}$ in $\Delta$ is
called \index[g]{geodesic path}{\em geodesic} if $|{\bf p}|\le |{\bf p}'|$ for any path ${\bf p}'$
combinatorially homotopic to ${\bf p}$.

The condition A has three parts:\index[g]{condition A}

\begin{itemize}
\item[A1] If $\Pi$ is a cell of rank $j>0,$ then
$|\partial(\Pi)|\ge nj$.

\item[A2] Any subpath of length $\le \max(j,2)$ of the contour of
an arbitrary cell of rank $j$ in $\Delta(2)$ is geodesic in
$\Delta$.

\item[A3] If $\pi, \Pi\in \Delta(2)$ and $\Gamma$ is a \ct submap
of $\pi$ to $\Pi$ with $(\pi,\Gamma,\Pi)\ge \varepsilon$, then
$|\Gamma\bigwedge\Pi|<(1+\gamma)k$, where $k=r(\Pi)$.
\end{itemize}

A map satisfying conditions A1, A2, A3 will be called an
\index[g]{$A$-map} $A$-{\em map}. As in \cite{book}, Section 15.2, a
(cyclic) section ${\bf q}$ of a contour of a map $\Delta$ is called a
\index[g]{smooth section}{\em smooth section of rank} $k > 0$ if:

1) every subpath of ${\bf q}$ of length $\le \max(k,2)$  is
geodesic in $\Delta$;

2) for each \ct submap $\Gamma$ of a cell $\pi$ to ${\bf q}$ satisfying
$(\pi, \Gamma,{\bf q})\ge \varepsilon$, we have $|\Gamma\bigwedge
{\bf q}|<(1+\gamma)k$.

\begin{remark}\label{geo} This definition implies that every geodesic section ${\bf q}$  of positive length in
$\partial\Delta$ is smooth of rank $k=|{\bf q}|$.
\end{remark}

It is shown in \cite{book}, $\S\S 16,  17$, that $A$-maps have
several ``hyperbolic" properties. (They show that $A$-maps of rank $i$ are hyperbolic spaces
with hyperbolic constant depending on $i$ only.) The next three lemmas are Lemma 15.8,
Theorem 16.1 and Theorem 16.2 from \cite{book}.
(The proof of Lemma \ref{0cont} is contained in the proof of Theorem 2 \cite{book}.)

\begin{lm} \label{bara} In an arbitrary A-map $\Delta$, the degree of contiguity of an arbitrary cell $\pi$ to an arbitrary cell $\Pi$ or to an arbitrary smooth section $\bf q$ of the contour $\partial\Delta$ via arbitrary contiguity submap is less than $\bar\alpha$. $\Box$
\end{lm}

\begin{lm}\label{four} (a) Let $\Delta$ be a disk $A$-map whose contour is subdivided into at most $4$ sections ${\bf q^1}$, ${\bf q^2}$,
${\bf q^3}$, ${\bf q^4}$ and $r(\Delta)>0$. Then there exists a positive cell $\pi$ and disjoint contiguity submaps $\Gamma_1$, $\Gamma_2$, $\Gamma_3$, $\Gamma_4$ of $\pi$ to
${\bf q^1},\dots, {\bf q^4}$ (some of these submaps may be absent) such that the sum of corresponding contiguity degrees $\sum_{i=1}^4(\pi, \Gamma_i, {\bf q}_i)$ is greater than $\bar\gamma=1-\gamma$.

(b)Let $\Delta$ be an annular $A$-map whose contours are subdivided into at most $4$ sections ${\bf q^1}$, ${\bf q^2}$,
${\bf q^3}$, ${\bf q^4}$ regarded as cyclic or ordinary paths, and $r(\Delta)>0$. Then there exists a positive cell $\pi$ and disjoint contiguity submaps $\Gamma_1$, $\Gamma_2$, $\Gamma_3$, $\Gamma_4$ of $\pi$ to
${\bf q^1},\dots, {\bf q^4}$ (up to three of these submaps may be absent) such that the sum of corresponding contiguity degrees $\sum_{i=1}^4(\pi, \Gamma_i, {\bf q}_i)$ is
greater than $\bar\gamma$. $\Box$
\end{lm}

\begin{lm}\label{0cont} Let $\Delta$ be a disk or annular $A$-map. Assume that there is a cell $\Pi$ in $\Delta$ and a contiguity submap $\Gamma$ of $\Pi$ to a section $q$ of the boundary with $(\Pi,\Gamma,{\bf q})\ge \varepsilon$. Then there is a cell $\pi$ of $\Delta$ and a contiguity submap $\Gamma'$ of $\pi$ to $\bf q$ such that $r(\Gamma')=0$ and
$(\pi,\Gamma',{\bf q})\ge\varepsilon$. $\Box$
\end{lm}

\section{Factor groups of finite exponent}

\label{segregation}

\subsection{Axioms}
\label{axioms}

Recall that in Chapter 6 of \cite{book}, it is proved that for any
sufficiently large odd $n$, there exists a presentation of the
free Burnside group $B(m,n)$ with $m\ge 2$ generators, such that
every reduced diagram over this presentation is an $A$-map.
The method from \cite{book} was generalized in \cite{OS03} to
Burnside factors of non-free groups in the following way.

In \cite{book}, 0-edges are labeled by 1 only, and (trivial)
0-relators have the form $1\cdot\dots\cdot 1$ or $a\cdot 1\cdot
... \cdot 1\cdot a^{-1} \cdot 1\cdot...\cdot 1$. In \cite{OS03}, we enlarged
the collection of labels of 0-edges and the collection of
0-relations. Namely, we divide the generating set of our
group into the set of 0-letters and the set of positive letters.
Then in the diagrams, 0-letters label 0-edges. We shall also add
defining relators of our group which have the form $uavb^{-1}$
where $u,v$ are words consisting of 0-letters only (we call such
words 0-{\it words}), and $a,b$ are a positive letters, as well as
defining relations without positive letters, to the list of zero
relations. The corresponding cells in diagrams will be called
$0$-cells.

Let ${\bf G}$ be a group given by a presentation $\langle {\cal X}\ |
\ {\cal R}\rangle$. Let ${\cal X}$ be a symmetric set (closed
under taking inverses). Let ${\cal Y}\subseteq {\cal X}$ be a
symmetric subset which we shall call the {\it set} of {\em
non-zero letters} or $Y$-{\em letters}. All other letters in
${\cal X}$ will be called \index[g]{$0$-letter} $0$-{\em letters}.

By the \index[g]{$Y$-length}$Y$-{\em length} of a word $A$, written
\index[g]{$|A|_Y$} $|A|_Y$ or just $|A|$, we mean the number of occurrences
of \index[g]{$Y$-letter}$Y$-letters in $A$. The $Y$-length $|g|_Y$ of an
element $g\in \bf G$ is the length of the $Y$-shortest word $w$
representing element $g$. A word $w$ in generators is called
\index[g]{minimal word} {\it minimal} (\index[g]{cyclically minimal 
word}{\it cyclically minimal}) in a group if
it is not equal (resp., not conjugate) in this group to a shorter word.

Let \index[g]{${\bf G}(\infty)$} ${\bf G}(\infty)$ be the group given by the presentation $\langle
{\cal X}\ |\ {\cal R}\rangle$ inside the variety of Burnside groups of
exponent $n$. We will choose a presentation of the group ${\bf
G}(\infty)$ in the class of all groups, that is we add enough
relations of the form $A^n=1$ to ensure that ${\bf G}(\infty)$ has
exponent $n$. This infinite presentation will have the form
$\langle {\cal X} \ | \ {\cal R}(\infty)= \cup {\cal S}_i\rangle$
where sets ${\cal S}_i$ will be disjoint, and relations from
${\cal S}_i$ will be called relations of {\it rank } $i$. Unlike
\cite[Chapter 6]{book}, $i$ runs over $0, 1/2, 1, 2,3,\dots.$
($i=0,1,2,\dots$ in \cite[Chapter 6] {book}.) The set ${\cal R}$
will be equal to ${\cal S}_0\cup {\cal S}_{1/2}$. The detailed
description of sets ${\cal S}_i$ will be given below. We shall
call this presentation of ${\bf G}(\infty)$ a \index[g]{graded presentation}{\it
graded presentation}.

As in \cite{book}, we proved in \cite{OS03}
that (reduced) diagrams over ${\cal R}(\infty)$ are A-maps.
This goal was achieved if
${\cal R}$ satisfies conditions \index[g]{Conditions (Z1), (Z2), (Z3)}(Z1), (Z2), (Z3)
presented below. While listing these conditions, we also fix some
notation and definitions.

\begin{itemize}
\item[(Z1)] The set ${\cal R}$ is the union of two disjoint subsets
${\cal S}_0={\cal R}_0$ and ${\cal S}_{1/2}$. The group $\langle
{\cal X}\ |\ {\cal S}_0\rangle$ is denoted by \index[g]{${\bf G}(0)$}${\bf G}(0)$ and is
called the group of rank $0$. We call relations from ${\cal S}_0$
\index[g]{relations of rank $0$} relations of rank 0. The relations from ${\cal S}_{1/2}$ have rank
$1/2.$ The group $\langle {\cal X}\ |{\cal R}\rangle$ is denoted
by \index[g]{${\bf G}(1/2)$} ${\bf G}(1/2).$
\begin{itemize}
\item[(Z1.1)]The set ${\cal S}_0$ consists of all relations from ${\cal R}$ which have
$Y$-length $0$ and all relations of ${\cal R}$ which have the form
(up to a cyclic shift) $ay_1by_2^{-1}=1$ where $y_1,y_2\in {\cal Y}^{+}$,
$a$ and $b$ are of $Y$-length $0$ (i.e. \index[g]{$0$-word}$0$-words).

\item[(Z1.2)] The set ${\cal S}_0$ implies all Burnside relations $u^n=1$ of $Y$-length 0.
\end{itemize}
\end{itemize}
The subgroup of ${\bf G}(0)$ generated by all $0$-letters is
called the \index[g]{$0$-subgroup }{\em $0$-subgroup} of ${\bf G}(0)$.
Elements from this subgroup are called \index[g]{$0$-elements}$0$-{\em
elements}. Elements which are not conjugates of elements from the
$0$-subgroup are called \index[g]{essential element}{\em essential}. An
essential element $g$ from ${\bf G}(0)$ is called \index[g]{cyclically $Y$-reduced element}{\em
cyclically $Y$-reduced} if a shortest word $A$ representing $g$ is cyclically $Y$-reduced, i.e. no cyclic permutation of $A$ has a subword of lengths $2$ equal to a 0-word. It is shown in \cite{OS03} (Lemma 3.5) that an element
$g$ is cyclically $Y$-reduced if and only if it is cyclically
minimal (in rank 0), i.e. it is not a conjugate of a $Y$-shorter
element in ${\bf G}(0)$.

Notice that condition (Z1.1) allows us to consider $Y$-bands in
\vk diagrams over ${\cal S}_0$. Maximal $Y$-bands do not
intersect.

\begin{itemize}
\item[(Z2)] The relators of the set ${\cal S}_{1/2},$ will be called
\index[g]{hub}{\it hubs}. The corresponding cells in \vk diagrams
are also called hubs. They satisfy the following properties
\begin{itemize}
\item[(Z2.1)] The $Y$-length of every hub is at
least $n$.
\item[(Z2.2)] Every hub is linear in $\cal Y,$ i.e. contains
at most one occurrence $y^{\pm 1}$ for every letter $y\in \cal Y.$
\item[(Z2.3)]
Assume that each of words $v_1w_1$ and $v_2w_2$ is a cyclic
permutation of a hub or of its inverse, and
$|v_1|\ge\varepsilon|v_1w_1|.$ Then an equality $u_1v_1=v_2u_2$
for some 0-words $u_1,$ $u_2,$ implies in ${\bf G}(0)$ equality
$u_2w_1=w_2u_1$ (see \setcounter{pdeleven}{\value{ppp}} Figure
\thepdeleven).

\unitlength=.85mm \special{em:linewidth 0.4pt}
\linethickness{0.4pt}
\begin{picture}(113.00,74.00)
\put(40.00,25.00){\framebox(60.00,20.00)[cc]{}}
\put(71.50,57.00){\oval(83.00,24.00)[]}
\put(71.50,14.00){\oval(83.00,22.00)[]}
\put(40.00,38.00){\vector(0,-1){7.00}}
\put(100.00,40.00){\vector(0,-1){8.00}}
\put(70.00,25.00){\vector(1,0){10.00}}
\put(62.00,45.00){\vector(1,0){12.00}}
\put(85.00,69.00){\vector(-1,0){10.00}}
\put(86.00,3.00){\vector(-1,0){12.00}}
\put(37.00,34.00){\makebox(0,0)[cc]{$u_1$}}
\put(71.00,22.00){\makebox(0,0)[cc]{$v_1$}}
\put(69.00,42.00){\makebox(0,0)[cc]{$v_2$}}
\put(103.00,36.00){\makebox(0,0)[cc]{$u_2$}}
\put(70.00,65.00){\makebox(0,0)[cc]{$w_2$}}
\put(102.00,6.00){\makebox(0,0)[cc]{$w_1$}}
\end{picture}

\begin{center}
\nopagebreak[4] Fig. \theppp.

\end{center}
\addtocounter{ppp}{1}

\end{itemize}
\end{itemize}

A word $B$ over ${\cal X}$ is said to be {\it
cyclically minimal in rank} $1/2$ if it is not a conjugate in
${\bf G}(1/2)$  of a word of smaller $Y$-length.  An element $g$
of ${\bf G}(0)$ is called {\em cyclically minimal in rank} $1/2$
if it is represented by a word which is cyclically minimal in rank
$1/2$.

\begin{itemize}
\item[(Z3)]  With every essential element $g\in {\bf G}(0)$ we
associate a subgroup $\oo(g)\le {\bf G}(0)$, normalized by $g$. If
$g$ is cyclically $Y$-reduced, let \index[g]{{\bf 0}(g)}$\oo(g)$ be the
maximal subgroup consisting of $0$-elements which contains $g$ in
its normalizer. By Lemmas  3.3 (2) and 3.4
(3) from \cite{OS03}, arbitrary essential element $g$ is equal to a product
$vuv^{-1}$ where $u$ is cyclically $Y$-reduced. In this case we
define $\oo(g)=v\oo(u)v^{-1}.$ It is proved in Lemma
3.7 \cite{OS03} that $\oo(g)$ is well defined.

\begin{itemize}
\item[(Z3.1)]
Assume that $g$ is an essential cyclically minimal in rank $1/2$
element of ${\bf G}(0)$ and for some $x\in {\bf G}(0)$, both $x$
and $g^{-4}xg^4$ are 0-elements. Then $x\in \oo(g).$

\item[(Z3.2)] For every essential cyclically minimal in rank $1/2$ element
$g\in {\bf G}(0)$, there exists a $0$-element $r$ such that $gr^{-1}$ commutes with
every element of ${\bf 0}(g)$ \footnote{We noticed in Section 3.1 \cite{OS03} that (Z3.2) can 
be replaced with the weaker condition (Z3.2'): For every essential cyclically minimal in rank $1/2$ element
$g\in {\bf G}(0)$, the extension of $\oo(g)$ by the automorphism
induced by $g$ (acting by conjugation) satisfies the identity
$x^n=1$. Condition (Z2.3') goes back to S.V. Ivanov's paper \cite{I92}.}

\end{itemize}
\end{itemize}

The following statement is given by Lemma 3.12 and Lemma 3.13 in \cite{OS03}.

\begin{lm} \label{star1}
Suppose that ${\cal R}$ satisfies conditions (Z1.1), (Z1.2), and (Z3.2). Then for every essential element $g\in {\bf G}(0)$ which
is cyclically minimal in rank $1/2$, and every $x\in \oo(g)$ we
have $(gx)^n=g^n$ and $g^nx=xg^n$. $\Box$
\end{lm}

In \cite{OS03}, the approach from \cite[Chapter
6]{book} was  adapted to diagrams over a presentation satisfying
(Z1), (Z2), (Z3). The construction is defined below .

\subsection{The construction of a graded presentation}

\label{construction}

Here we recall several concepts analogous to the concepts in
\cite{book}, chapter 6. They were introduced in \cite{OS03} by induction on the rank $i$.

We say that a certain statement {\it holds in rank $i$} if it
holds in the group \index[g]{${\bf G}(i)$}${\bf G}(i)$ given by ${\cal
R}_i=\cup_{j=0}^i
 {\cal S}_j$.

Recall that  a word $w$  is \index[g]{periodic word} $A$-{\em periodic} if it
is a subword of some power of the word $A.$

{\em Up to the end of this section, the $Y$-length of a word $W$
will be called simply \index[g]{length of word}{\em length} and will be denoted by $|W|$.}

  Let us define simple words of rank $i=0,1/2,1,2,...$ and periods of rank
$i=1,2,3,\dots$. Let ${\cal S}_0={\cal R}_0$ and ${\cal S}_{1/2}$
be the set of hubs. Thus the groups ${\bf G}(0)$ and ${\bf
G}(1/2)={\bf G}$ are defined by relations of the sets ${\cal R}_0$
and ${\cal R}_{1/2}={\cal S}_0\cup{\cal S}_{1/2},$ respectively. A
word $A$ of positive length is said to be \index[g]{simple in rank $0$ or $1/2$ word}{\it simple} in rank 0
(resp. 1/2) if it is not a conjugate in rank 0 (resp. 1/2) either
of a word of length 0 or of a product of the form $B^l P$ for an essential $B$, where
$|B|<|A|$ and $P\in \oo(B)$ in rank 0.

\begin{remark} \label{0simple} It follows from the definition
that a simple word in rank $0$ is not a conjugate in rank 0 of any
shorter word. Also if two words of the same
length are conjugate in rank 0 and one of them
is simple in rank 0, then the other one is simple in rank 0 as well.
\end{remark}

Suppose that $i\ge 1$, and we have defined the sets of relators
\index[g]{${\cal R}_j$}${\cal R}_j$, $j<i$, and the corresponding groups ${\bf
G}(j)$ of ranks $j<i$. Suppose also that we have defined simple
words of ranks $j$, $j<i$, and periods of rank $j$, $1\le j< i$.

For every $i = 1/2,1,2,3,..$ let $i_-$ denote $i-1$ if $i>1$,
$1/2$ if $i=1$ or $0$ if $i=1/2$. Similarly $i_+$ is $1/2$ if
$i=0$ or $1$ if $i=1/2$ or $i+1$ if $i\ge 1$.

For every $i\ge 1$ consider the set ${\cal X}'_i$ of all words of
length $i$ which are simple in rank $i_-$, and the
\index[g]{$\sim_{i_-}$} equivalence $\sim_{i_-}$ given in Lemma \ref{Equival}
(we apply it for the smaller rank here).  Now choose a set of
representatives \index[g]{${\cal X}_i$}${\cal X}_i$ of the equivalence classes
of ${\cal X}'_i$. The words of ${\cal X}_i$ are said to be
\index[g]{period of rank $i$}{\it periods} of rank $i.$ The set of words ${\cal
R}_i$ (defining the group ${\bf G}(i)$ or rank $i$) is the union
of ${\cal R}_{i_-}$ and \index[g]{${\cal S}_i$}${\cal S}_i= \{A^n, \ A\in {\cal
X}_i\}$.

Let $A$ be a word of positive length. We say that $A$ is
\index[g]{simple in rank $i$ word}{\it simple} in rank $i\ge 1$ if it is not a
conjugate in rank $i$ either of a word of length 0 or of a word of
the form $B^lP,$ where $B$ is a period of rank $j\le i$ or an
essential word with $|B|<|A|$, and $P$ represents a word of the
subgroup $\oo(B).$

\begin{lm}\label{Equival} (\cite{OS03}, Lemma 3.24) The following relation $\sim_i$ is an
equivalence on the set of all simple in rank $i$ words:  $A\sim_i
B$ by definition, if there are words $X, P, R$, where $P\in
\oo(A)$, $R\in\oo(B)$ (in rank 0), such that $AP=XB^{\pm
1}RX^{-1}$ in rank $i$. $\Box$
\end{lm}

Let \index[g]{${\cal R}(\infty)$}${\cal R}(\infty)=\cup_{i=0}^\infty {\cal R}_i$. The set
${\cal R}(\infty)$ defines the group ${\bf G}(\infty)$, which is the
inductive (= direct) limit of the sequence of epimorphisms  ${\bf G}(0)\to {\bf G}(1/2)\to {\bf G}(1)\to\dots\to  {\bf G}(i)\to\dots$

If $\Pi$ is a cell of a diagram over ${\cal R}(\infty)$ with the
boundary labeled by a word of the set ${\cal S}_{j},$ then, by
definition, $r(\Pi)=j.$ The boundary label of every cell of a rank
$\ge 1$  has a period $A$ defined up to
cyclic shifts.

For $j\ge 1$, a pair of distinct cells $\Pi_1$ and $\Pi_2$ of rank
$j$ of a diagram $\Delta$ is said to be a \index[g]{$j$-pair of cells} $j$-{\it
pair}, if their counterclockwise contours ${\bf p}_1$ and ${\bf p}_2$ are
labeled by $A^n$ and $A^{-n}$ for a period $A$ of rank $j$ and
there is a path ${\bf t}=({\bf p}_1)_- -({\bf p}_2)_-$ without self-intersections
such that the label $Lab({\bf t})$ is equal, in rank $j_-$, to an element of
$\oo(A).$ Then $Lab({\bf t})$ and $A^n$ commute in rank $j_-$ by Lemma
\ref{star1}, and so the subdiagram with contour ${\bf p}_1{\bf t}{\bf p}_2{\bf t}^{-1}$
can be replaced in $\Delta$ by a diagram of rank $j_-$. As a
result, we obtain a diagram with the same boundary label as
$\Delta$ but of a smaller type.

Similarly, a pair of hubs $\Pi_1,$ $\Pi_2$ with counterclockwise
contours ${\bf p}_1$, ${\bf p}_2$ forms a \index[g]{$\frac 12$-pair}$\frac 12$-{\it pair}, if
the vertices $({\bf p}_1)_-$ and $({\bf p}_2)_-$ are connected by a path ${\bf t}$
without self-intersections such that the label of the path
${\bf p}_1{\bf t}{\bf p}_2{\bf t}^{-1}$ is equal to 1 in ${\bf G}(0).$ Consequently, any
diagram can be replaced by a reduced diagram having no $j$-pairs, $j=1/2,
1, 2,...$, i.e. by a  \index[g]{$g$-reduced diagram}{\it g-reduced} diagram with the
same boundary label(s)\footnote{Recall that in \cite{book}, $\S
13.2,$ such diagrams are called {\em reduced}. But in this paper
we shall consider different kinds of reduced diagrams, so we call
diagrams without $j$-pairs g-reduced; ``g" stands for ``graded".}.

The following statement is a part of Proposition 3.19  \cite{OS03}.

\begin{lm}\label{mainprop} If the presentation of a group ${\bf G}$ satisfies
properties (Z1), (Z2), (Z3) then
all g-reduced diagrams over the presentation
$\langle {\cal X\mid R}(\infty)\rangle$
are $A$-maps. Every relator from ${\cal S}_i$, $i\ge
1$, is of the form $A^n$ for a cyclically reduced word $A$ of
$Y$-length $i$. The group ${\bf G}(\infty)$ defined by this presentation is the factor group ${\bf
G}/\la g^n, g\in {\bf G}\ra$. $\Box$
\end{lm}

\begin{remark}\label{fb} If ${\cal Y}={\cal X}$, ${\cal S}_0={\cal S}_{1/2}=\emptyset$, then ${\bf G}$ is the free group with basis $\cal X$ and trivial subgroups ${\bf 0}(g)$, and so Conditions
(Z1), (Z2), (Z3) obviously hold. In this case, the factor group ${\bf G}(\infty)={\bf
G}/\la g^n, g\in {\bf G}\ra$ is the free Burnside group $B({\cal X})$ of exponent $n$ with free set of generators $\cal X$. Therefore $g$-reduced diagrams over ${\bf G}(\infty)=B({\cal X})$
are $A$-maps.
\end{remark}

The following two lemmas will be  explicitly used in Section \ref{more}.

\begin{lm}\label{3.25} (\cite{OS03}, Lemma 3.25).
Every word $X$ is conjugate in rank $i\ge 0$ either of a $0$-word or of a word $A^lP$, where
$|A|\le |X|$, $A$ is either a period of rank $j\le i$ or a simple in rank $i$ word and $P$ represents in rank $0$ an element of the subgroup ${\bf 0}(A)$.

If $X$ is not simple in rank $0$ and is not conjugate to a $0$-word in rank $0$, then
$|A|<|X|$.  $\Box$
\end{lm} 

Note that the second sentence of Lemma \ref{3.25} is not formulated in Lemma 3.25 \cite{OS03}, but %immediately
it follows from the
first phrase of the proof.

\begin{lm} \label{3.32} (\cite{OS03}, Lemma 3.32).
Let $\Delta$ be a $g$-reduced diagram of rank $i$ with the contour ${\bf p}_1{\bf q}_1{\bf p}_2{\bf q}_2$, where $Lab({\bf q}_1)$ and
$Lab({\bf q}_2)$ are periodic words with simple in rank $i$ periods $A$ and $B$, resp., and
$|A|\ge |B|$. Assume that $|{\bf p}_1|, {\bf p}_2|<\alpha|B|$, $|{\bf q}_1|>\frac{3}{4}\delta^{-1}|A|$
and $|{\bf q}_2|>\delta^{-1}|B|$. Then the word $A$ is a conjugate in rank $i$ of a product $B^{\pm 1}R$, where $R$ represents an element of the subgroup ${\bf 0}(B)$, and $|A|=|B|$.
Moreover, if the words $Lab({\bf q}_1)$ and
$Lab({\bf q}_2)^{-1}$ start with $A$ and $B^{-1}$, resp., then $A$ is equal in rank $i$ to $Lab({\bf p}_1)^{-1}B^{\pm 1}R\; Lab({\bf p}_1)$.
 $\Box$ \end{lm} 

\section{Aperiodic endomorphisms}\label{apend}

The rules of the S-machines introduced in this paper will use words from an infinite  set  \index[g]{set of words ${\bf W}$} {${\bf W}$
of positive words in the 2-letter alphabet $\{a,b\}$. This set has  the following properties.

\medskip
\index[g]{Properties (A), (B), (C) of the set ${\bf W}$}
(A) For arbitrary two distinct words $u$ and $v$ from ${\bf W}$, their maximal common prefix (common suffix) has length at most $\frac 15 \min (|u|, |v|)$ (where $|w|$ is the length of a word $w$). In particular, the subgroup generated by $\bf W$ in the free group
$F=F(a,b)$ is a free group with basis $\bf W$.

(B) Suppose in some product of words from ${\bf W}^{\pm 1}$ after all cancellations, we have a non-empty
$B$-periodic subword $V$ of length $\ge 2|B|$. Then
either $B$ is freely conjugate of a some product of
words from ${\bf W}^{\pm 1}$ or $|V|\le 11|B|$.

(C) The subgroup $\cal B$ generated by $\bf W$ in the free Burnside
quotient $B(a,b)=F/F^n$ of $F$ is the free
Burnside group of exponent $n$ with basis $\bf W$.

\medskip

For every large enough odd $n$, such infinite set ${\bf W}$ was constructed by D. Sonkin in
\cite{S}.

Choose now $({\bf w}_1, {\bf w}_1(a) ,{\bf w}_1(b)), ({\bf w}_2, {\bf w}_2(a), {\bf w}_2(b)), \dots$, disjoint triples of distinct words from  ${\bf W}$.
With every triple $({\bf w}_j, {\bf w}_j(a), {\bf w}_j(b))$ we associate
the subgroup \index[g]{$H(j)$} $H(j)=gp\{{\bf w}_j(a), {\bf w}_j(b)\}$, the left coset \index[g]{$K(j)$} $K(j)= {\bf w}_j H(j)$ and the endomorphism \index[g]{$\psi_j$} $\psi_j$ of the free group $F(a,b)$ to $F(a,b)$ (and of the free Burnsine group of exponent $n$ $B(a,b)$ to $B(a,b)$) such that $$\psi_j(a)={\bf w}_j(a)\;\; and \;\;\psi_j(b)={\bf w}_j(b)$$

\begin{remark} \label{mono} It follows from Property (A) of the set $\bf W$ (from (C))
that every $\psi_j$ is a monomorphism of $F(a,b)$
(of $B(a,b)$) and $\psi_j (F(a,b))\cap \psi_{j'} (F(a,b))=\{1\}$ (resp., $\psi_j (B(a,b))\cap \psi_{j'} (B(a,b))=\{1\}$) if $j\ne j'$. In other words, $H(j)\cap H(j')= \{1\}$, and the images of these two subgroups in
$B(a,b)$ have trivial intersection too.
\end{remark}

Let \index[g]{$\phi_j$} $\phi_j$ be the mapping  $F(a,b)\to F(a,b)$ and  $B(a,b)\to B(a,b)$ given by the formula $\phi_j(u)={\bf w}_j\psi_j(u)$. Property (A) of the set $\bf W$ implies
$|u|<|\phi_j(u)|$ for every $u\in F(a,b)$.

\begin{remark}\label{disj} By Remark \ref{mono},
every mapping $\phi_j$ is injective both on $F(a,b)$
and on $B(a,b)$. Besides we have $\phi_j (F(a,b))\cap \phi_{j'} (F(a,b))=\emptyset$ (resp., $\phi_j (B(a,b))\cap \phi_{j'} (B(a,b))=\emptyset$) if $j\ne j'$ by Property
(A) (resp., by (C)) of the set $\bf W$. In other words, $K(j)\cap K(j')=\emptyset$, and the images of these two cosets in
$B(a,b)$ have empty intersection too.
\end{remark}

\begin{lm} \label{produ} Assume that $w_0\in F(a,b)$
and for $s=1,\dots, t$, we have $w_i =\phi_{j_i}^{\epsilon_i}(w_{i-1})$, for some mappings $\phi_{j_i}$-s, where $\epsilon_i=\pm 1$ and we write $w_i =\phi_{j_i}^{-1}(w_{i-1})$ if $w_{i-1} =\phi_{j_i}(w_{i})$. If $w_t=w_0$, then the reduced form of
the word $\phi_{j_1}^{\epsilon_1}\dots\phi_{j_t}^{\epsilon_t}$
is empty. \end{lm}

\proof
Assuming that the word $\phi_{j_1}^{\epsilon_1}\dots\phi_{j_t}^{\epsilon_t}$ is reduced we will induct on $t$ with obvious statement for $t=0$.
By Remark \ref{disj}, we cannot have $\epsilon_{i-1}=1$ and $\epsilon_i=-1$
for any $i$.
Thus, we may assume that $\epsilon_1=\dots=\epsilon_l=-1$
for some $l\le t$ and $\epsilon_{l+1}=\dots=\epsilon_t=1$. If $1\le l<t$,
then by Remark \ref{disj}, the equality $w_0=w_t$
implies $\phi_{i_1}=\phi_{j_t}$. It follows that
$u_1=u_{t-1}$. So the inductive conjecture can be
applied to the series $u_1,\dots, u_{t-1}$.

It remains to assume that $\epsilon_1=\dots=\epsilon_t=1$. In this case $t=0$,
since for $t>0$, we have $|u_0|<|u_1|<\dots< |u_t|$, a contradiction.
\endproof

\begin{lm}\label{posit} Let $w_0,\dots,w_t$ be words representing
elements of the free group $F(a,b)$, where
for every $i=1,\dots, t$, we have either
$w_i=\phi_j(w_{i-1})$ or $w_{i-1}=\phi_j(w_i)$ for
some $j=j(i)$.
Then there is a sequence of positive words $u_0,\dots,u_t$ such that

(1) $u_i=\phi_j(u_{i-1})$ or, respectively, $u_{i-1}=\phi_j(u_i)$ for the same $\phi_j$ as above,
$i=1,\dots,t$;

(2) $u_l=1$ (empty word) for some $l\le t$, and every $u_i$ is
an image of $1$ under a product of some $\psi_j$-s;

(3) if $w_0=w_t$, then $u_0=u_t$.
\end{lm}

\proof We will induct on the length $t$ of the sequence
with obvious base $t=0$.

Assume $t>0$ and for some $i$ we have $\phi_l(w_{i+1})=w_i=\phi_j(w_{i-1})$. Then $l=j$ by Remark \ref{disj}, and
$w_{i+1}=w_{i-1}$ in $F(a,b)$.
Hence one can apply the inductive conjecture to
the sequence $w_0,\dots,w_{i-1}, w_{i+2}\dots w_t$ of length $t-2$, obtain a sequence $u_0,\dots,u_{i-1}, u_{i+2}\dots u_t$ with conditions (1)-(3), and
define $u_i=\phi_j(u_{i-1})$, $u_{i+1}=u_{i-1}$.
obtaining the required sequence of length $t$.

Therefore one may assume that for some $l\in [0;t]$, the mappings act as follows: $w_l\mapsto\dots\mapsto w_0$
and $w_l\mapsto\dots\mapsto w_t$. Then we define $u_l=1$, and by induction $u_{i+1}=\phi_{j_{i+1}}(u_i)$ for $i\ge l$
and $u_{i-1}=\phi_{j_i}(u_i)$ for $i\ge l$. Clearly,
we obtain Properties (1) and (2). It remains to prove
Property (3).

Suppose $1\le l<t$. If
$w_0=\phi_j(w_1)$ and $w_t=\phi_{j'}(w_{t-1})$, then  $j=j'$ by Remark \ref{disj}. Then the equality  $w_0=w_t$ implies $w_1=w_{t-1}$ in $F(a,b)$. Since we may
assume by induction that $u_1=u_{t-1}$, we obtain
$u_0=\phi_j(u_1)=\phi_j(u_{t-1})=u_1$, as required.

It remains to  assume that all the mappings are arranged as $w_0\mapsto\dots\mapsto w_t$. Then $|w_0|<|w_1|<\dots <|w_t|$, contrary to
the assumption $w_0=w_t$.
Thus, Property (3) is proved.
\endproof

\begin{lm} \label{eo} Let $w_0,\dots,w_t$ be words representing
elements of the free group $F(a,b)$, where
for every $i=1,\dots, t$, we have either
$w_i=\phi_j(w_{i-1})$ or $w_{i-1}=\phi_j(w_i)$ for some $j=j(i)$.
Then either $w_0$ and $w_t$ are
equal in $F(a,b)$ or these two words represent different
elements of  $B(a,b)$.
\end{lm}

\proof Assume that for some $i$, we have $\phi_l(w_{i+1})=w_i=\phi_j(w_{i-1})$. Then $l=j$ by Remark \ref{disj} and
$w_{i+1}=w_{i-1}$ in $F(a,b)$.
Hence inducting on $t$, one can replace $t$ with $t-2$.

Thus, we may assume that for some $l\in [0;t]$, the mappings are directed as follows: $w_l\mapsto\dots\mapsto w_0$
and $w_l\mapsto\dots\mapsto w_t$. Suppose $1\le l<t$. Then
$w_0=\phi_j(w_1)$ and $w_t=\phi_j(w_{t-1})$ for the same $\phi_j$ by Remark \ref{disj}. Then the equality $w_1=w_{t-1}$ implies $w_0=w_t$ in $F(a,b)$
and $w_1\ne w_{t-1}$ implies $w_0\ne w_t$ in $B(a,b)$
since the mapping $\phi_j$ is injective on $B(a,b)$.
Hence to complete the proof, one can consider a shorter sequence $w_1,\dots, w_{t-1}$ again.

Ir remains to assume that all the mappings are directed as $w_0\mapsto\dots\mapsto w_t$. Under this
assumption, we will prove that
$w_0\ne w_t$ in $B(a,b)$ if $t>0$. Proving by contradiction, we may assume that we have a counter-example with the shortest first word $w_0$. It follows that the word $w_0$ is minimal in $B(a,b)$, i.e. it is not equal to a shorter word.
Indeed if $w_0=v_0v_0'$, where $|v_0|<|w_0|$
and $v'_0$ is trivial in $B(a,b)$, then applying the first mapping $\phi_{j_1}$
we would obtain $v_0v_0'\to {\bf w}_{j_1}\psi(v_0)_{j_1}\psi_{j_1}(v_0')$,
which is equal to ${\bf w}_{j_1}\psi_{j_1}(v_0)$  in $B(a,b)$
since $\psi_{j_1}$ is a homomorphism of $B(a,b)$.
Therefore $\phi_{j_1}(w_0)= \phi_{j_1} (v_0)$, and by induction,
$\phi_{j_t}(\dots(\phi_{j_1}(w_0)))=\phi_{j_t}(\dots(\phi_{j_1}(v_0)))$ in $B(a,b)$,
and $v_0$ would begin a shorter counter-example.

Let $w_t=\phi_j(w_{t-1})$.
If $w_0$ is also a $\phi_j$-image in $F(a,b)$, i.e. $w_0=\phi_j(w)$ for some word $w$, then
we have a sequence $w\mapsto w_0\mapsto\dots\mapsto w_{t-1}$ of the same length but with $|w|<|w_0|$ and $w=w_{t-1}$ in $B(a,b)$ by Remark \ref{disj}, a contradiction.
Hence we may assume further that $w_0$ is not a $\phi_j$-image in $F(a,b)$.

The equality $w_0=w_t$ in $B(a,b)$ provides us with a $g$-reduced
van Kampen diagram $\Delta$ over the presentation
of $B(a,b)$ defined in Section \ref{construction} (see Remark \ref{fb}) with boundary path ${\bf pq}$, where
$Lab({\bf p})\equiv w_0$, $Lab({\bf q})\equiv w_t^{-1}$. (We use symbol $\equiv$ for letter-by-letter equality of words.) Since
$w_t$ is a $\phi_j$-image, there is a subdiagram $\Gamma$ of minimal type with contour ${\bf pr}$, where $Lab({\bf r})^{-1}$
is a $\phi_j$-image in $F(a,b)$.
$\Gamma$ has at least one cell of positive rank,
 because $w_0$ is not a $\phi_j$-image.

 By Lemma \ref{four}, there is a cell $\Pi$ in $\Gamma$ such that
 the sum of degrees of contiguity of $\Pi$ to ${\bf p}$ and to
 ${\bf r}$ is greater than $\bar\gamma$. However the degree of contiguity of $\Pi$ to the geodesic path ${\bf p}$ cannot
 exceed $\bar\alpha$ by Remark \ref{geo} and Lemma \ref{bara}. Hence the degree
 of contiguity of $\Pi$ to ${\bf r}$ is $>\bar\gamma-\bar\alpha>1/3$.
 By Lemma \ref{0cont}, we have a cell $\pi$ of positive rank in $\Gamma$ with immediate contiguity to  ${\bf r}$ of degree
 $\ge\varepsilon$ (see \setcounter{pdeleven}{\value{ppp}} Figure
\thepdeleven).

% This is a LaTeX picture output by TeXCAD.
% File name: [pr.pic].
% Version of TeXCAD: 4.3
% Reference / build: 30-Jun-2012 (rev. 105)
% For new versions, check: http://texcad.sf.net/
% Options on the following lines.
%\grade{\on}
%\emlines{\off}
%\epic{\off}
%\beziermacro{\on}
%\reduce{\on}
%\snapping{\off}
%\pvinsert{% Your \input, \def, etc. here}
%\quality{8.000}
%\graddiff{0.005}
%\snapasp{1}
%\zoom{4.0000}
\unitlength 1mm % = 2.845pt
\linethickness{0.4pt}
\ifx\plotpoint\undefined\newsavebox{\plotpoint}\fi % GNUPLOT compatibility


\begin{center}
\nopagebreak[4] Figure \theppp
\end{center}
\addtocounter{ppp}{1}

 Since (1) the boundary label of $\pi$ is a power $A^n$
 for some cyclically reduced word $A$, (2) the label of
 ${\bf r}^{-1}$ is a reduced form of a word from the coset $K(j)={\bf w}_jH(j)$ of $H(j)$, and (3) $\varepsilon n>30$, Properties (A),(B) of the set $\bf W$ imply that $A^{15}$
 is a subword of a word from the subgroup $H(j)$, and 
 $A$ is freely conjugate to a word from the subgroup $H(j)$
  of $F(a,b)$, i.e. a cyclic permutation $A'$ of $A$ belongs to $H(j)$. Property $(A)$ of the set $\bf W$  
	guarantees that a subpath ${\bf x}$ of ${\bf r}$  labeled by
 $(A')^{\pm 14}$ can be replaced with the complement of
 ${\bf x}$ in $\partial \pi$ so that for the obtained path ${\bf r}'$,
 we again have $Lab({\bf r'}^{-1})\in K(j)$. However the type
  of the subdiagram $\Gamma'$ with contour ${\bf pr}'$ is strictly less
 than the type of $\Gamma$.
  Thus, the lemma is proved by contradiction.
 \endproof
Recall that a subgroup $L$ of a group $K$ is called \index[g]{malnormal subgroup}{\it malnormal} if
$xLx^{-1}\cap L =\{1\}$ for every $x\in K\backslash L$. This property is transitive, i.e., if $M$ is a malnormal subgroup of $L$ and $L$ is malnormal in $K$, then $M$ is malnormal in $K$.
\begin{lm} \label{maln} The canonical image  $\tilde H(j)$ of the subgroup $H(j)\le F(a,b)$ in the group
$B(a,b)$ is malnormal in $B(a,b)$.
\end{lm}
\proof Assume we have an equality $z^{-1}uz=v$ in $B(a,b)$, where the words $u$ and $v$ represent elements
of $H(j)$ with non-trivial images in $B(a,b)$. We claim that $z$ represents an element from $\tilde H(j)$ too.
 
According to Remark \ref{fb}, we have a diagram $\Delta$ over the
 presentation of $B(a,b)$ with four boundary sections labeled
 by $z, u, z^{-1}$ and $v$, resp. Identifying two of them, we
 obtain an annular diagram $\Delta_0$ with two boundary sections
 ${\bf p}_1$ and ${\bf p}_2$ labeled by $u$ and $v$, where the vertex $o_1=({\bf p}_1)_-=({\bf p}_1)_+$ is connected with the vertex $o_2=({\bf p}_2)_-=({\bf p}_2)_+$ by a simple path $\bf q$ labeled by the word $z$.
 If $\Delta_0$ is not a $g$-reduced diagram, then one can make reductions
 preserving the labels of ${\bf p}_1$ and ${\bf p}_2$ and changing $\bf q$ by ${\bf q}'$, where the labels of $\bf q$ and ${\bf q}'$ are equal modulo the
 relations of $B(a,b)$ (see Section 13 in \cite{book}). So we may
 regard $\Delta_0$ as a $g$-reduced diagram.

 Now we will induct on the number $k$ of positive cells in $\Delta_0$. If $k=0$, then we have the conjugation of $u$ and
 $v$ in the free group $F(a,b)$. So
 the statement
 follows from Property (A) of the set $\bf W$.

 If $k\ge 1$, then by Lemmas \ref{four} (b) and \ref{0cont}, we get
 a positive cell $\pi$ in $\Delta_0$ with immediate conjugacy
 to ${\bf p}_1$ or to ${\bf p}_2$ of degree $\ge\varepsilon$, and as in the proof of Lemma \ref{eo}, one can
 change ${\bf p}_1$ (or ${\bf p}_2$) by a homotopic path ${\bf p}'_1$ such that
 the word $Lab ({\bf p}'_1)$ belongs to $H(j)$, it is nontrivial in $B(a,b)$, the annular diagram
 bounded by ${\bf p}'_1$ and ${\bf p}_2$ has $k-1$ positive cells, and the
 path ${\bf q'}$ connecting $({\bf p}'_1)_-$ and ${\bf p}_2$ is homotopic to ${\bf q}$.
 So by the inductive conjecture,
 $z$ belongs to $\tilde H(j)$.
 \endproof

 \section{Machines} \label{mach}

\subsection{Turing machines}

We will use a model of {\it recognizing} Turing machine ($TM$) close to the model from \cite{SBR}.
A  $TM$ with $k$ tapes and $k$
heads is  a tuple $$ M= \langle A, Y, Q,
\Theta, \vec s_1, \vec s_0 \rangle$$ where $A$ is the input alphabet,
$Y=\sqcup_{i=1}^k Y_i$ is the tape alphabet, $Y_1 \supset A,$
$Q=\sqcup_{i=1}^k Q_i$ is the set of states of the machine, $\Theta$ is a set of commands, $\vec s_1$ is
the $k$-vector of start states, $\vec s_0$ is the $k$-vector of
accept states. ($\sqcup$ denotes the disjoint union.) The sets $Y, Q, \Theta$
are finite.

We
assume that the machine normally starts working with
states of the heads forming the vector $\vec s_1$,  with the head
placed at the right end of each tape, and accepts if it reaches the
state vector $\vec s_0$.  In general, the machine can be turned on
in any configuration and turned off at any time.

A {\em configuration} of tape number $i$ of a $M$ is a word
$u q v$ where $q\in Q_i$ is the current state of the head,
$u$ is the word to the left of the head, and $v$ is the word to the
right of the head, $u,v$ are words in the alphabet $Y_i$.
A tape is {\em empty} if $u$, $v$ are empty words.

A \index[g]{configuration of Turing machine}{\em configuration} $U$ of the machine $M$ is a word
$$\alpha_1U_1\omega_1\alpha_2U_2\omega_2... \alpha_kU_k\omega_k$$
where $U_i$ is the configuration of tape $i$, and the end-markers $\alpha_i, \omega_i$
of the $i$-th tape are special separating symbols.

An \index[g]{input configuration}{\it input configuration} $w(u)$ is a configuration, where all tapes,
except for the first one, are empty, the configuration of the first
tape (let us call it the {\it input tape}) is of the form $uq$,
$q\in Q_1$, $u$ is a word in the alphabet $A$,
and the states form the start vector $\vec s_1$. The
\index[g]{accept configuration} {\em accept configuration} is the configuration where the state
vector is $\vec s_0$, the accept vector of the machine, and all
tapes are empty.

To every $\theta\in \Theta,$ there corresponds a \index[g]{command of Turing machine} {\it command} (marked by the
same letter $\theta$), i.e., a pair of sequences of words
 $[V_1,...,V_k]$ and $[V'_1,...,V'_k]$
 such that for each $j\le k,$ either both $V_j=uqv$ and $V'_j=u'q'v'$ are configurations of the tape number $j,$ or
 $V_j =\alpha_j qv$ and $V'_j = \alpha_j q'v',$
 or $V_j = uq\omega_j$ and $V'_j = u'q' \omega_j ,$ or $V_j = \alpha_j q\omega_j$ and $V'_j = \alpha_j q' \omega_j $ ($q,q'\in Q_j$ ).

In order to execute this command, the machine checks if every $V_i$
is a subword of the current configuration of the machine,
and if this condition holds the machine replaces $V_i$ by $V'_i$
for all $i=1,\dots,k.$ Therefore we also use the notation:
$\theta: [V_1\to V'_1,\dots, V_k\to V'_k],$  where $V_j \to V'_j$ is
called the $j$-th {\it part} of the command $\theta.$

 Suppose we have a sequence of configurations $w_0,...,w_t$ and a word
 $H= \theta_1\dots\theta_{t}$ in the alphabet $\Theta,$
such that for every $i=1,..., t$, the machine passes from $w_{i-1}$ to
$w_i$ by applying the command $\theta_i$. Then the sequence
$w_0\to w_1\to\dots\to w_t$ is said to be a \index[g]{computation of Turing machine} {\it computation} with
\index[g]{history of computation} {\it history} $H.$ In this case we shall write $w_0\cdot H=w_t.$
 The number $t$ will be called the {\em time} or the \index[g]{length of computation} {\em length} of the computation.

A configuration $w$ is called \index[g]{accepted configuration} {\em accepted} by a machine $M$ if
there exists at least one computation which starts with $w$ and
ends with the accept configuration.

A word $u$ in the input alphabet $A$ is said to be \index[g]{accepted word} {\em accepted} by the machine if the
corresponding input configuration is accepted. (A configuration with
the vector of states $\overrightarrow s_1$
is never accepted if it is not an input configuration.)  The set of all
accepted words over the alphabet $A$ is called the \index[g]{language recognized by Turing machine}{\em language  ${\cal L}_M$
recognized by the machine $M$}.

We do not only consider {\em deterministic} Turing machines, for example,
we allow several transitions with the same left side.  A Turing machine $M$ is called \index[g]{symmetric machine}{\it symmetric} if with every command $\theta: [V_1\to V'_1,\dots, V_k\to V'_k],$  it contains the
\index[g]{inverse command}{\em inverse} command $\theta^{-1}: [V'_1\to V_1,\dots, V'_k\to V_k].$ Given a deterministic machine $M$, one can extend the set of commands by adding $\theta^{-1}$
for every $\theta\in \Theta$. The obtained (non-deterministic)
Turing machine $Sym(M)$ is called the symmetrization of $M$.

A computation of a machine is called \index[g]{reduced computation} {\it reduced} if the history $H$ of it admits
no cancellations, i.e. $H$ has no 2-letter subwords of the form $\theta\theta^{-1}$. Clearly, every computation can be made reduced (without changing the start or end configurations of the computation) by removing consecutive mutually inverse rules.

We will use the following additional properties of Turing machines; most of them were formulated
in \cite{birget, BORS, O12}

\begin{lm} \label{bos} For every deterministic Turing machine $M$ recognizing a language $\cal L$
there exists a Turing machine $M_0$ with the following properties.
\begin{enumerate}
\item The language recognized by $M_0$ is $\cal L$.
\item $M_0$ is symmetric.
\item If ${\cal C}: w_0\to\dots\to w_t$ is a reduced
computation of $M_0$ and $w_0\equiv w_t$,
then $t=0$ .

\item If a computation ${\cal C}: w_0\to\dots\to w_t$ of $M_0$ with length $t>0$ starts and ends with input configurations of $M_0$ then both
    $w_0$ and $w_t$ are accepted by $M_0$.

\item The machine $M_0$ accepts only when all tapes are empty.
\item Every command of $M_0$ or its inverse
has one of the following
forms for some $i$
\begin{equation}\label{eq57}
[q_1\omega_1\to q_1\omega_1, ...,q_{i-1}\omega_{i-1}\to q_{i-1}\omega_{i-1}, q_i\omega_i\to aq_i\omega_i, q_{i+1}\omega\to q_{i+1}\omega_{i+1},...]
\end{equation}
\begin{equation}\label{eq56}
[q_1\omega_1\to q_1\omega_1, ...,q_{i-1}\omega_{i-1}\to q_{i-1}\omega_{i-1}, q_i\omega_i\to q'_i\omega_i, q_{i+1}\omega\to q_{i+1}\omega_{i+1},...]
\end{equation}
\begin{equation}\label{eq58}
[q_1\omega\to q_1\omega_1,...,q_{i-1}\omega_{i-1}\to q_{i-1}\omega_{i-1}, \alpha_i q_i\omega_i\to \alpha_i q_i'\omega_i, q_{i+1}\omega_{i+1}\to q_{i+1}\omega_{i+1},...]
\end{equation}
where $a$ belongs to the tape alphabet of tape $i$,
$q_j, q_j'$ are state
letters of tape $j$.
\item The letters used on different tapes are from disjoint alphabets.
This includes the state letters.
\end{enumerate}
\end{lm}

\begin{remark} \label{odna} (1) Only one part (number $i$) of the commands (\ref{eq57}, \ref{eq56}, \ref{eq58}) is changing. 

(2) If the head observes the right markers
at the beginning of a computation, it will observe the right markers during
the whole computation. If the head does not observe the right markers at
the beginning, then no
command is applicable, and so the computation is trivial.
\end{remark}

\proof Properties (1), (2), and (5) are provided by Lemma  2.3 \cite{O12} for
the machine $M_0=Sym(M)$.
%, where $M$ can be  assumed deterministic.

The same lemma 2.3 \cite{O12} (Properties (b),(c)) says that the history of every reduced computation $\cal C$ of $M_0$
has the form $H=H_1H_2^{-1}$, where $H_1$ and $H_2$ are histories of $M$,
and $M$ can be chosen so that the same configuration cannot occur twice in a computation of $M$. Hence to prove Property (3) by contradiction, one may assume that both
$H_1$ and $H_2$ are non-empty. Since $H_1\ne H_2$, we have
$H_1=H_2H'$ (or $H_2=H_1H'$) for some non-empty $H'$ since the machine $M$ is deterministic and the computations $\cal C$ and ${\cal C}^{-1}$ start with the same configuration. But then $w_0\cdot H_2\equiv (w_0\cdot H_2)H',$ a contradiction, since $H'$ is a history of $M$-computation.

To prove Property (4), we first modify the deterministic machine $M$ and obtain a deterministic Turing machine $M'$ working as follows.

After the initial command (or commands) changing the states of the heads only (but preserving tape words), the work of $M'$ is subdivided in three Steps.

\medskip

Step 1. $M'$ copies the input word from the first tape to two auxiliary tapes with numbers $k+1$ and $k+2$ by deleting letter-by letter from the first tape and inserting the copies of erased letters to the last two tapes. Connecting rule is  applied when the first
tape is empty. It changes the states of all the heads.

Step 2. $M'$ works as $M$ but uses the $k+1$-st tape
as the input one. It keeps the $k+2$-d tape unchanged. The copy  of the accept command of $M$ connects this Step with the next one. So the tapes
$1, 2,\dots, k+1$ are empty when this command applies.

Step 3. $M'$ letter-by-letter erases the content
of the tape number $k+2$ and finally accepts.

\medskip

Let $M''=Sym(M')$. By Property (3), we have $w_0\ne w_t$ and as there, the history is $H=H_1H_2^{-1}$, where $H_1$ and $H_2$ are histories of $M'$-computations. Here both $H_1$ and $H_2$ must contain the commands of
both Steps 1 and 2 since $w_0\ne w_t$ and the definition of Step 1 imply that $w_0 \cdot H'\ne w_t\cdot H''$ for arbitrary histories without commands of Step 2. Furthermore,
both $H_1$ and $H_2$ have to involve the commands of Step 3, because Step 2 does not
change the (different) copies of the input words obtained on the tape number $k+2$. So by the
definition  of $M'$, both $w_0$
and $w_t$ are accepted words for $M$, and
therefore for $M'$ and $M''$.

Property (6) is obtained in Lemma 3.1 \cite{SBR}
by means of division of every tape in two parts, which does not affect
(1) - (5), and so we get a modification $M'''$. However the analogs of (\ref{eq57}, \ref{eq58}) are a bit weaker
there, e.g. the analog of (\ref{eq57}) in \cite{SBR} allows to change the states of
many heads:
\begin{equation}\label{eq59}
[q_1\omega\to q_1'\omega, ...,q_{i-1}\omega\to q_{i-1}'\omega, aq_i\omega\to q_i'\omega, q_{i+1}\omega\to q_{i+1}'\omega,...]
\end{equation}

But this can be easily improved. For example, if one has $q_1\to q_1'$ and
$q_2\to q_2'$ for a command $\theta$, then it is possible to introduce
a new state letter $q$ and replace $\theta$ with three commands $\theta_1$:
$q_1\to q$, $q_2\to q_2$, $\theta_2$: $q\to q,  q_2\to q'_2$, and
$\theta_3$: $q\to q'_1, q'_2\to q'_2$  changing each only one state letter.
Since the new letter $q$ is involved in these three commands only,
it  is easy to see that the modified machine
keeps the properties (1) - (5). The same trick works if one head inserts/deletes a letter and the same head (or another one ) changes the state.
Therefore Property (6) follows by induction.

To obtain (7), one can  just prescribe the number of the tape to the tape and state letters as the additional index and change the commands of the machine accordingly. The machine $M_0$ is built.
\endproof

\begin{remark}\label{pair} (1) The Steps 1 - 3 of the Turing machines $M''=Sym (M')$ constructed in the proof of Lemma \ref{bos} (4) (which are also the Steps of $M_0$) will be applied in the proof of Lemma \ref{inin}.

(2) The rewriting commands of Step 1
have the form $$\theta_a: aq_1\omega_1\to q_1\omega_1,\dots, \alpha_{k+1}q_{k+1}\to \alpha_{k+1}q_{k+1}a, \alpha_{k+2}q_{k+2}\to \alpha_{k+2}q_{k+2}a$$ for every $a\in A$. When we
replace such a command with several commands $\theta_{a,1},\dots, \theta_{a, m}$ using
the trick from the proof of item (6) of Lemma \ref{bos},
we introduce a number of auxiliary states of the heads, where every state is used for the transition from
some $\theta_{a,i}$ to $\theta_{a,i+1}$ only for some $i$ and does not used in any other command.

Therefore if an $M_0$-computation ${\cal C}: w_0\to\dots\dots\to w_t$ has no those auxiliary state letters in $w_0$ and in $w_t$, then every
command $\theta_{a,i}$ can occur in the history
$H$ of $\cal C$ only in the block-subword $b=(\theta_{a,1}\dots \theta_{a, m})^{\pm 1}$. In particular, such a block has a unique command $\theta'(b)^{-1}$ (command $\theta''(b)$) of the form (\ref{eq57}) that deletes the letter $a$ from the first sector (resp., inserts the copies of $a$ in sectors
$k+1$ and $k+2$.
\end{remark}

In cases (\ref{eq57}), (\ref{eq57}), and (\ref{eq58}), we say that
the head $Q_i$ is \index[g]{working head} {\it working} for $\theta$, and so there is exactly one working head for every command of $M_0$.

\subsection{S-machines (modified) as rewriting systems}\label{SM}

There are several equivalent definitions of  S-machines (see \cite{S}). However we need a more general definition in comparison with that was used in \cite{SBR}, \cite{OS06} or \cite{OS19}, and below the finite sets of letters
$Y(\theta)$ are replaced with finitely generated subgroups of free groups. Besides, some parts of the rules $\theta$ having
form $aqb\to cq'd$ cannot be written as $q\to (a^{-1}c)q'(db^{-1})$, because now the words $a,b,c,d$ in the tape alphabet are not necessarily belong to $Y(\theta)$. In particular, some rules may now insert/deletes words in locked
sectors.

A ``hardware'' of an $S$-machine $\bf S$ is a pair $(Y,Q),$ where $Q=\sqcup_{i=1}^s Q_i$ and $Y= \sqcup_{i=0}^s Y_i$ for some $s\ge 2$. 
We always set $Y_{s}=Y_0=\emptyset$ and if $Q_{s}=Q_0$ (i.e., the indices of $Q_i$ are counted modulo $s$), then we say that $\bf S$ is a \index[g]{circular machine} {\em circular} $S$-machine.

The elements from $Q$ are called \index[g]{state letters}{\it state letters} or \index[g]{$q$-letters} $q$-{\it letters} of an S-machine, the elements from $Y$ are \index[g]{tape letters} {\it tape letters} or \index[g]{$a$-letters} $a$-{letters} of an S-machine. The sets $Q_i$ (resp.
$Y_i$) are called \index[g]{part of state/tape letters}{\it parts} of state and tape letters of $Q$ (resp. $Y$). The number of $a$-letters (resp. $q$-letters, $\theta$-letters)  in a word $W$ is called the \index[g]{$a$-length, $q$-length, $\theta$-length} $a$-length
(resp., $q$-length, $\theta$-length) of $W$ denoted by $|W|_a$ ($|W|_q$, $|W|_{\theta}$, resp.).

The language of \index[g]{admissible word}{\it admissible words}
 of an S-machine $\bf S$ (the language of $\bf S$-admissible words)
consists of reduced words $W$ of the form
\begin{equation}\label{admiss}
q_1u_1q_2\dots u_{r-1} q_{r},
\end{equation}
where $r\ge 1$, every $q_i$  is a state letter from some part $Q_{j(i)}^{\pm 1}$, $u_i$ are reduced group words in the alphabet of tape letters of the part $Y_{k(i)}$ and for every $i=1,...,r-1$, one of the following holds:

\begin{itemize}
\item If $q_i$ is from $Q_{j(i)}$ then $q_{i+1}$ is either from $Q_{j(i)+1}$ or is equal to $q_i\iv$ and $k(i)=j(i)+1$.

\item If $q_i\in Q_{j(i)}\iv$ then $q_{i+1}$ is either from $Q_{j(i)-1}^{-1}$ or is equal to $q_i\iv$ and $k(i)=j(i)$.
\end{itemize}
Every subword $q_iu_iq_{i+1}$ of an admissible word  will be called the  {\em $Q_{j(i)}^{\pm 1}Q_{j(i+1)}^{\pm 1}$-sector} of that word. An admissible word may contain many $Q_{j(i)}^{\pm 1}Q_{j(i+1)}^{\pm 1}$-sectors.

Usually parts of the set $Q$ of state letters are denoted by capital letters. For example, a part $P$ would consist of
letters $p$ with various indices.

If an admissible word $W$
has the form (\ref{admiss}), $W=q_1u_1q_2u_2...q_r,$
and $q_i\in Q_{j(i)}^{\pm 1},$
$i=1,...,r$, $u_i$  are
group words in tape letters, then we shall say that the \index[g]{base} {\it base} of an admissible word  $W$ is the word
$Q_{j(1)}^{\pm 1}Q_{j(2)}^{\pm 1}...Q_{j(r)}^{\pm 1}$. Here $Q_i$ are just symbols which denote the corresponding parts of the set of state letters. Note that, by the definition of admissible words, the base is not necessarily a reduced word.

Instead of saying that the parts of the set of state letters of $\bf S$ are $ Q_1, Q_2,... , Q_s$ we will write
that \index[g]{standard base} {\em the standard base}
of the $S$-machine is $Q_1...Q_s$.

The software of an $S$-machine with the standard base $Q_1...Q_s$ is a set of \index[g]{rules of S-machine} {\em rules} $\Theta$.  Every $\theta\in \Theta$
is given by

(1) a  family \index[g]{$Y_i(\theta)$} \{$Y_i(\theta)\}_{i=0}^s$ 
of finitely generated subgroups $Y_i(\theta)\le F(Y_i)$
of free groups and

(2) a sequence $$[a_1q_1b_1\to a'_1q_1'b'_1,...,a_sq_sb_s\to a'_sq_s'b'_s],$$ where $q_i\in Q_i$,
$a_i, a'_i$ are reduced words from the group $F(Y_{i-1})$, $b_i, b'_i$ are reduced words from $F(Y_i)$
(recall that $Y_{0}=Y_s=\emptyset$)
for every $i$.

Each component $a_iq_ib_i\to a'_iq_i'b_i$ is called a \index[g]{part of rule} {\em part} of the rule.
In most cases it will be clear what the sets $Y_i(\theta)$ are.
By default we assume that $Y_i(\theta)$ is the whole free group  $F(Y_i)$.

Every rule $\theta=[a_1q_1b_1\to a'_1q'_1b'_1,...,a_sq_sb_s\to a'_sq'_sb'_s]$ has the \index[g]{inverse rule} inverse one
$$\theta\iv=[a'_1q_1'b'_1\to a_1\iv q_1b_1\iv,..., a'_sq_s'b'_s\to a_s\iv q_sb_s]$$ It is always the case that
$Y_i(\theta\iv)=Y_i(\theta)$.

 Thus the set of rules
$\Theta$ of an S-machine is divided into two disjoint parts, \index[g]{$\Theta^+$}
$\Theta^-$} $\Theta^+$ and
$\Theta^-$ such that for every $\theta\in \Theta^+$, we have $\theta\iv\in
\Theta^-$ and for every $\theta\in\Theta^-$, we have $\theta\iv\in\Theta^+$ (in particular $\Theta\iv=\Theta$, that is any $S$-machine is symmetric).

The rules from $\Theta^+$ (resp. $\Theta^-$) are called \index[g]{positive, negative rules} {\em
positive} (resp. {\em negative}).

 To \index[g]{application of a rule} apply a rule  $\theta=[a_1q_1b_1\to a'_1q'_1b'_1,...,a_sq_sb_s\to a'_sq'_sb'_s]$ to an admissible word $p_1u_1p_2u_2...p_r$ (\ref{admiss})
where each $p_i\in Q_{j(i)}^{\pm 1}$, means

\begin{itemize}
\item to check for every subword $p_iu_ip_{i+1}$
    if the following word belongs to the subgroup $Y_{j(i)}(\theta)$:

    (a) $b_{j(i)}^{-1}u_ia_{j(i)+1}^{-1}$,
    when $p_i\in Q_{j(i)}$ and $p_{i+1}\in Q_{j(i)+1}$,

     (b)  $a_{j(i)+1}u_ia_{j(i)+1}^{-1}$, when $p_{i+1}\in Q_{j(i)+1}$ and $p_{i}=p_{i+1}^{-1}$,

     (c)  $b_{j(i)}^{-1}u_ib_{j(i)}$, when $p_{i}\in Q_{j(i)}$ and $p_{i+1}=p_{i}^{-1}$,
   
and if this property holds,

\item replace each $p_i=q_{j(i)}^{\pm 1}$ by $(a_{j(i)}^{-1}a'_{j(i)}q'_{j(i)}b'_{j(i)}b_{j(i)}^{-1})^{\pm 1}$
\item if the resulting word is not reduced or starts (ends) with $a$-letters, then reduce the word and trim the first and last $a$-letters to obtain an admissible word again.
\end{itemize}

It follows from the definitions that if $W'$ is obtained by applying a rule $\theta$ (we write $W'=W\cdot \theta$), then the inverse rule is applicable to $W'$ and $W'\cdot \theta^{-1}=W$.

\begin{remark} The application of every rule $\theta$ of $\bf S$ is effective since the membership problem is decidable for any finitely generated subgroup of a free group.
\end{remark}

\begin{remark}\label{qtoq'a} If $a,b,a',b'\in Y_{i}(\theta)$, then the component
$aq_ib\to a'q_i'b'$ of a rule $\theta$ can be written in the equivalent form $q_i\to cq'_id$, where
$c$ and $d$ are the reduced forms of the products $a'a^{-1}$ and $b'b^{-1}$, respectively.
\end{remark}

For example, applying the rule $[q_1\to a\iv q_1'b, q_2\to cq_2'd]$ to the admissible word $q_1b\iv q_2dq_2\iv q_1\iv$ (where $a,b,c,d$ are just letters) we first obtain the word
$$a\iv q_1'bb\iv cq_2'ddd\iv (q_2')\iv c\iv b\iv(q_1')\iv a,$$ then after trimming and reducing we obtain $$q_1'cq_2'd (q_2')\iv c\iv b\iv(q_1')\iv$$

 If a rule $\theta$ is applicable to an admissible word $W$, we say that $W$ belongs to the \index[g]{domain of a rule}{\em domain} of $\theta$. In other words, $W$ is \index[g]{$\theta$-admissible word} $\theta$-{\it admissible}.

We call an admissible word with the standard base a   \index[g]{configuration of S-machine} \emph{configuration} of an S-machine.

We usually (but not always) assume that every part $Q_i$ of the set of state letters contains a  \index[g]{start/end state letters} {\em start state letter} and an
{\em end state letter}. Then a configuration is called a   \index[g]{start/end configurations} \emph{start} (\emph{end}) configuration if all state letters in it are start (end) letters. As Turing machines, some S-machines are  {\em recognizing a language}. In that case we choose an  \index[g]{input sector} {\em input} sector, say, the $Q_1Q_2$-sector, of every configuration. The $a$-projection of that sector (i.e. the word obtaining by deleting of $q$-letters) is called the {\em input} of the configuration. In that case, the end configuration with empty $a$-projection is called the  \index[g]{accept configuration} {\it accept configuration}.

A  \index[g]{computation of S-machine}{\em computation} of  {\it length} $t\ge 0$ is a sequence of admissible words $W_0\to \dots\to W_t$ such that for
every $0=1,..., t-1$ the S-machine passes from $W_i$ to $W_{i+1}$ by applying one
of the rules $\theta_i$ from $\Theta$.  The word $h=\theta_1\dots\theta_t$ is called the  \index[g]{history of computation}{\it history}
of the computation. Since $W_t$ is determined by $W_0$ and the history $h$, we use notation $W_t=W_0\cdot h$ and say that $W_0$ belongs to the domain of $h$.

If by means a computation, the $S$-machine can take an input configuration with input $u$ to the accept configuration, we say that the word $u$ is  \index[g]{accepted word} {\em accepted} by the $S$-machine. We define  \index[g]{accepted configuration} {\em accepted} configurations (not necessarily start configurations) similarly. All accepted words $u$ form the \index[g]{language recognized by S-machine} {\it language recognized by the S-machine}.

A computation is called  \index[g]{reduced computation} {\em reduced} if its history is a reduced word. Clearly, every computation can be made reduced (without changing the start or end configurations of the computation) by removing consecutive mutually inverse rules. The domain of a reduced form of a history $h$ can be larger than the domain of $h$.

\begin{lm} \label{gen2}
Suppose the base of an admissible word $W$ is $Q_{i}Q_{i}\iv$
(resp., $Q_i\iv Q_i$). Assume that each rule $\theta$ of a
computation $W\equiv quq\iv\to \dots\to W'\equiv q'u'(q')\iv$
(resp., $W\equiv q^{-1}uq\to \dots\to W'\equiv (q')^{-1}u'q'$) with history $h$
has a component $a_{\theta,i}q_ib_{\theta,i}^{-1}\to a'_{\theta,i}
q_i'(b'_{\theta,i})^{-1},$
Assume that
the mapping $\theta\mapsto b_{\theta,i}^{-1} b'_{\theta,i}$
($\theta\mapsto a_{\theta,i}^{-1}a'_{\theta,i}$) extends to a monomorphism $\lambda$ of the free Burnside group $\cal B$ on $\theta$-letters of the history
to the free Burnside group $B=B(Y_{i+1})$ with basis $Y_{i+1}$ (to  $B=B(Y_i)$). Then $u'$ is equal to $u$ modulo the Burnside relations if and only if
the $\lambda$-image of $h$ belong to the centralizer of $u$ in $B$.
\end{lm}

\proof We consider the base $Q_{i}Q_{i}\iv$ only. By induction on $t$, we see that
$u'$ is equal to the word $v^{-1}uv$ in $B$, where $v=\lambda(h)$, which proves the lemma. %and so  $v$ is 
\endproof

If for some rule $\theta=[a_1q_1b_1\to a'_1q_1'b'_1,...,a_sq_sb_s\to a'_sq'_sb'_s]$ of an S-machine $\bf S$, the subgroup $Y_{i}(\theta)$ is trivial
then we say that $\theta^{\pm 1}$  \index[g]{locked sector}{\it locks} the $Q_iQ_{i+1}$-sector. In that case
we denote the $i$-th part of the rule as follows:  \index[g]{$\tool$} $a_iq_i\tool a_iq_i'$.
Since the admissible words are reduced, the definition of rule application implies

\begin{lm}\label{qqiv}
If the $i$-th component of the rule $\theta$ has the form $v_iq_i\tool v'_iq_i',$
then the
base of any $\theta$-admissible word cannot have subwords $Q_iQ_i\iv$ or $Q_{i+1}^{-1}Q_{i+1}.$
\end{lm} $\Box$

\begin{remark} The above definition of S-machines can be compared with the definition of multi-tape Turing machines.
The main differences are that the heads  of Turing machines are near-sighted, but the heads  of  S-machines are farsighted : they do not "see" the nearest tape letters  but know if the whole tape words belong to distinguished cosets associated with the subgroups
$Y_i(\theta)$.

Since S-machines are symmetric, for every computation ${\cal C}: W_0\to W_1\to \dots\to W_t$, there is an inverse one
${\cal C}^{-1}: W_t\to W_{t-1}\to \dots\to W_0$. An S-machine  can work with words containing negative letters, and words with "non-standard" order of state letters. In particular,
the {\it mirror copy} of the computation ${\cal C}$ is the computation $W_0^{-1}\to W_1^{-1}\to \dots\to W_t^{-1}$ with the same history. Many properties of computations obviously hold for their mirror copies.
\end{remark}

\subsection{Primitive machines}

For a set of letters
$A,$ let $A'$ and $A''$ be  disjoint  copies of $A$, the maps $a\mapsto a'$ and $a\mapsto a''$ identify $A$ with $A'$ and $A'',$ resp. Let \index[g]{$\overrightarrow Z$ and $\overleftarrow Z$} $\overrightarrow Z=\overrightarrow Z(A)$ and $\overleftarrow Z=\overleftarrow Z(A)$ be the $S$-machines
with tape alphabet $A'\sqcup A''$, state alphabet $\{L\}\cup P\cup \{R\},$ where
$P=\{p(1),p(2),p(3)\}$ and the following positive $S$-rules. For  \index[g]{$\overleftarrow Z$} $\overleftarrow Z$ we have the rules

$$\xi_1(a)=[L\to L, p(1)\to (a')^{-1}p(1)a'', R\to R],\;\; a\in A$$
{\em Comment:} The head can move from right to left, replacing the word in alphabet $A'$ by its copy in the alphabet $A''$. (However, the inverse rules can move the head from right to left. Also an application of the rule can insert/delete two letters, one letter from the left and one letter from the right of the $p$-letter.)

$$\xi_2:=[L\tool L,  p(1)\to p(2), R\to R]$$
{\em Comment:} When the head meets $L$, it turns into $p(2)$.

For \index[g]{$\overrightarrow Z$}$\overrightarrow Z$, we define the rules

$$\xi_3(a)=[L\to L, p(2)\to a' p(2)(a'')^{-1}, R\to R]$$
{\em Comment:}
%Since $\xi_3(a)^{-1}=[L\to L, p(2)\to a' %p(2)(a'')^{-1}, R\to R]$,
The head $p(2)$ can move from left to right, replacing the word in $A''$ by its copy in $A'$.
$$\xi_4=[L\to L, p(2)\tool p(3), R\to R]$$
{\em Comment:} When the head reaches the right end of the tape, it turns into $p(3)$.

\begin{remark} \label{aiv} For every $a\in A$, $i=1,3$, it will be convenient to denote $\xi_i(a)\iv$ by $\xi_i(a\iv)$. It is clear from the definition $\xi_i(a)$ that this does not lead to a confusion.\end{remark}

\begin{lm}\label{pR} Assume that ${\cal C}:\;W_0\to ...\to W_t$ is a reduced computation of $\overrightarrow Z$  (of $\overleftarrow Z$) with history $h$. Suppose the base $B=B({\cal C})$ of $\cal C$ contains $LP$ or $PR$.
Then $W_i\ne W_j$ if $0\le i<j\le t$.

Furthermore, if $W_i$ contains the subword $Lu_ip$ (subword $pu_iR$) for some tape word $u_i$ and the corresponding subword of $W_j$ is
$Lu_jp$ (resp., $pu_jR$) for the same $p$-letter $p$, then the history
of the computation $W_i\to\dots\to W_j$ is trivial in the free Burnside group of exponent $n$  provided $u_j$ is equal to $u_i$ modulo the Burnside relations.
 \end{lm}

\proof It suffices to assume that $B\equiv LP$ or $B\equiv PR$.

Let the machine $\overrightarrow Z$ start with a word $W_i\equiv Lu_ip(2)$, ends with $W_j\equiv  Lu_jp(2)$, where $u_i, u_j$ are words   in $A'$-letters, and does not apply the rule $\zeta_4$. Then $u_j$ is freely equal to $u_iu$, where $u$ is the copy $a'_{i+1}\dots a'_j$ of the reduced history $h_{ij}\equiv \xi_3(a_{i+1})...\xi_3(a_j)$ of this computation, and so $W_j\ne W_i$. If $W_i\equiv Lu_ip(2)$ and $W_j\equiv  Lu_jp(3)$, then $W_i$ and $W_j$ have just different state letters.

It is not possible that $W_k\equiv Lu_kp(3)$ for $0<k<t$ since in this case both transitions $W_{k-1}\to W_{k}$ and $W_{k+1}\to W_{k}$ were the applications of the rule $\zeta_4$, and the history $h$ were not reduced. If $t=2$ and both $W_0$ and $W_2$ contain $p(3)$ we have the same contradiction, since the history
should be equal to $\zeta_4^{-1}\zeta_4$. If $t>2$ and  $W_0\equiv Lu'_ip(3)$ and $W_t\equiv  Lu_jp(3)$,
then $u_j\ne u_i$ as in the previous paragraph since $W_1\equiv Lu_ip(2)$ and $W_{t-1}\equiv  Lu_jp(2)$
in this case.

Thus, the statement of the lemma is proved for the machine   $\overrightarrow Z$ and the subbase $LP$.
The proof of the Lemma is finished since other properties are similar. \endproof

We also need two simple machines depending on the mappings
 $\psi_j$ and $\phi_j$ defined in Section \ref{apend}. Now $A=\{a,b\}$.  The machine \index[g]{$\overleftarrow Z(\phi_j)$} $\overleftarrow Z(\phi_j)$ replaces a word with the $\phi_j$-image of it. As $\overleftarrow Z=\overleftarrow Z(A)$, the $S$-machines $\overleftarrow Z(\phi_j)$
has tape alphabet $A'\sqcup A''$, state alphabet $\{L\}\cup P\cup \{R\},$ where
$P=\{p(1),p(2), p(3)\}$ and the following positive $S$-rules depending on $\psi_j$.

$$\xi'_1(c)=[L\to L, p(1)\to (c')^{-1}p(1){\bf w}''_j(c), R\to R], \;\;c\in A=\{a,b\},$$
where ${\bf w}''_j(c)$ is the copy of ${\bf w}_j(c)$ in the alphabet $A''$, and the finitely generated subgroup $A''(\xi'_1(c))\le F(A'')$ corresponding
to this rule is the the copy $H''(j)$ of the subgroup $H(j)\le F(A)$ in $F(A'')$,  i.e. this rule cannot
be applied unless the words in the $PR$-sectors (and in $PP^{-1}$-, $R^{-1}R$-sectors if any occur)
belong to $H''(j)$.

$$\xi'_2:=[L\tool L,  p(1)\to {\bf w}'_jp(2), R\to R],$$
where ${\bf w}_j'$ is the copy of the word ${\bf w}_j$ in the alphabet
$A'$, $A'(\xi'_2)=\{1\}$ and $A''(\xi'_2)= H''(j).$
So the rule $\xi'_2$ (the rule $(\xi'_2)^{-1}$) is applicable only if the sector $LP$ is empty
(resp., if the tape word of this sector is ${\bf w}_j'$).

{\em Comment:} The difference in comparison with $\overleftarrow Z=\overleftarrow Z(A)$ is that the replacements $a'\mapsto \psi_j(a')$
given by the rules $\xi'_1(a)$ ($a\in A$) use the monomorphism $\psi_j$, and the rule $\xi'_2$ inserts
the word ${\bf w}_j'$ into the $LP$-sector .

The machine \index[g]{$\overrightarrow Z(\phi_j)$} $\overrightarrow Z(\phi_j)$ has the positive rules

$$\xi'_3(c)=[L{\bf w}'_j\to L{\bf w}'_j, p(2){\bf w}''_j(c)\to
{\bf w}'_j(c)p(2), R\to R],$$
where $A'(\xi'_3(c))=H'(j)$ and $ A''(\xi'_3)=H''(j)$.

The next rules is

$$\xi'_4 = [L{\bf w}'_j\to L{\bf w}'_j,  p(2)\tool p(3), R\to R]$$

Here $A'(\xi'_4)=H'(j)$, and $A''(\xi'_4) =\{1\}$.

{\it Comment:} The work of the machine $\overleftarrow Z(\phi_j)$ followed by the work of $\overrightarrow Z(\phi_j)$ replaces the configuration $Lup(1)R$ with $L\phi_j(u)p(3)R$.

We have the following analog of Lemma \ref{pR}.

\begin{lm}\label{Ephi} Suppose ${\cal C}:\;W_0\to ...\to W_t$ is a reduced computation of $\overleftarrow Z(\phi_j)$
or of $\overrightarrow Z(\phi_j)$ for some $j$. Assume that the base $B=B({\cal C})$ of $\cal C$ contains $LP$ or $PR$.
Then $W_i\ne W_j$ if $0\le i<j\le t$.

Furthermore, if $W_i$ contains the subword $Lu_ip$ (subword $pu_iR$) for a tape word $u_i$ and the corresponding subword of $W_j$ is
$Lu_jp$ (resp., $pu_jR$) for the same $p$-letter $p$, then the history
of the computation $W_i\to\dots\to W_j$ is trivial in the free Burnside group of exponent $n$ provided $u_j$ is equal to $u_i$ modulo the Burnside relations.
\end{lm}

\proof It suffices to repeat the argument of
the proof of Lemma \ref{pR} but instead of
a nontrivial copy of the word $a'_{i+1}\dots a'_j$, we now have
a (copy of the) nontrivial $\psi_j$-image of such word for the machine $\overleftarrow Z(\phi_j)$ since the homomorphism $\psi_j$ is injective by Remark \ref{mono}.\endproof

\begin{lm}\label{unl} Let $QQ'$ be a sector of the standard base of one of the machines $\overleftarrow Z$, $\overrightarrow Z$ $\overleftarrow Z(\phi_j)$, $\overrightarrow Z(\phi_j)$ and $Y'$ or $Y''$ be the tape subalphabet of this sector. Assume that a rule
$\eta$ does not lock the sector $QQ'$ and have
a word $xuy$ in the domain, where $x\in Q$,
$y\in Q'$. Then for
every rule $\xi$ of the same machine having at least one word $xvy$ in the domain,
we have
$Y'(\xi)\le Y'(\eta)$ (resp., $Y''(\xi)\le Y''(\eta)$).
\end{lm}

\proof This property can be checked by the inspection of every rule. For instance, if $QQ'=LP$ and
$\eta=\xi'_1(c)$, then $\xi$ is either
$\xi'_1(d)$ for a letter $d\in \{a^{\pm 1}, b^{\pm 1}\}$, or $\xi=\xi'_2$. In both cases the required inequality follows from the definition
of these rules.
\endproof

\begin{remark} \label{proj} Let ${\cal C}: W_0\to W_1\to\dots\to W_s$ be a computation
of $\overleftarrow Z$ or of $\overrightarrow Z$, or
of $\overrightarrow Z(\phi_j)$ with standard base, and $u'_l$, $v''_l$ be the tape words of the $LP$- and $PR$-sectors of the word $W_l$, $l=0,\dots,s$. Let $\pi: F(A')\to F(A'')$ be an isomorphism, such that
$\pi(a')=a''$ for every $a\in A$. Then it follows from the definition of $Z$-machines that the reduced forms of the words $\pi(u'_j)v''_j$ are equal for all $l=0,\dots,s$.

We will refer to this property as \index[g]{projection argument} {\it `projection argument'}. For the computations of the
machine $\overleftarrow Z(\phi_j)$ having no rules $(\xi'_2)^{\pm 1}$ in the history, it states that  the reduced form of $\psi_j(u_l)v_l$
is unchanged in a computation with standard base,
where $u_l$ and $v_l$ copy the words $u'_l$ and $v''_l$
in the alphabet $A=\{a,b\}$.

\end{remark}

\begin{remark} \label{wrong} Assume that $W_0\equiv Lup(1)vR$ and a reduced computation
${\cal C}: W_0\to W_1\to\dots\to W_s$ of $\overleftarrow Z$ starts with $W_0\to W_1=W_0\cdot\xi_1(a)$, $a\in A^{\pm 1}$, where the last letter of $u$ is not $a'$
and the word $v$ does not start with $a''$. Then we have $W_1\equiv Lu(a')^{-1}p(1)a''vR$ with reduced subwords $u_1\equiv u(a')^{-1}$ and
$v_1\equiv a''v$. Since $u_1$ and $v_1$ are non-empty, the word $W_1$ is in
the domains of the rules $\xi_1(b)$ ($b\in A^{\pm 1}$) only, and $b\ne a^{-1}$ in the transition
$W_1\to W_2=W_1\cdot \xi_1(b)$ since $\cal C$ is a reduced computation. Hence
the products $u_2\equiv u_1(b')^{-1}$ and $v_2\equiv b''v_1$ are reduced too,
and so on. Therefore all the rules in the history of $\cal C$ are of the type
$\xi_1(c)$ for $c\in A^{\pm 1}$ and both $LP$- and $PR$-sectors of $W_s$
are non-empty. Under above assumption, we say that the application of
the first rule $\xi_1(a)$ to $W_0$ is \index[g]{wrong application of a rule} {\it wrong}.

Similarly we define wrong applications of $\xi_3$-rules of $\overrightarrow Z$.

The definition of a wrong application of a rule $\xi'_1(c)$ of the machine
$\overleftarrow Z(\phi_j)$ is  similar, but now we assume that $c'$ is not the last letter of $u$ and
the word $v$ is a copy of the $\psi_j$-image of a reduced word which does not start with $c^{-1}$.
The definition of a wrong application of a rule $\xi'_3(c)$ of the machine
$\overrightarrow Z(\phi_j)$ is also similar:
we assume that $u$ (resp., $v$) is a a copy of $\phi_j$-image (a copy of $\psi_j$-image) of a reduced word which does not end with $c^{-1}$
(does not start with $c$).
\end{remark}

\begin{lm} \label{zto} Assume that for some rule $\xi$
of $\overleftarrow Z$ or of $\overrightarrow Z$, or of
$\overleftarrow Z(\phi_j)$, or of $\overleftarrow Z(\phi_j)$,
the canonical image $\tilde A$ of the subgroup $A'(\xi)$ (or of the subgroup $A''(\xi)$) in $B(a',b')$ (resp., in $B(a'',b'')$ contains
nontrivial in $\tilde A$ words $u$ and $z^{-1}uz$ for some $z\in B(a',b')$
(resp. $z\in B(a'',b'')$. Then the word $z$ belongs in $\tilde A$ too.
\end{lm}
\proof We consider (non-trivial) $A'(\xi)$ only. The statement is obvious if $\tilde A$ is the whole
$B(a',b')$. Otherwise $A'(\xi)$ is equal to a subgroup $H'(j)$,
and therefore the statement follows from Lemma \ref{maln}.
\endproof

\section{Encoding of the rules of Turing machines by primitive machines}\label{simu}

Let $M_0= \langle A, Y(0), Q,
\Theta, \vec s_1, \vec s_0 \rangle$ be a Turing machine satisfying the conditions
(1) - (7) of Lemma \ref{bos} with the tape alphabet $ Y(0)=\cup_{i=1}^{k} Y_i(0)$; it recognizes a language $\cal L$ in the input alphabet $A$. We will construct an \index[g]{$M_1$} $S$-machine $M_1$ that simulates the work of $M_0$ by applying
the mappings $\phi_j$ instead of the command of $M_0$.

Recall that $M_0$ is symmetric, and every rule $\theta$ has at most one active
 head inserting/deleting a letter or changing the state. Separating the commands of $M_0$ into two parts, we say that inserting commands are positive. If a command does not change tape words, then we may randomly choose the sign for this command and
 the opposite sign for the inverse command. We also may assume that there
 are no duplicate commands, i.e. every command of $M_0$ is uniquely determined
 by the set of parts of it, and there is no command $\theta=\theta^{-1}$.

If the tape subalphabet $ Y_i(0)$ of $M_0$  has letters $a_{i1},..., a_{i,k_i}$, then we chose disjoint  triples of different words  $({\bf w}_{i1}, {\bf w}_{i1}(a), {\bf w}_{i1}(b))..., ({\bf w}_{i,k_i}, {\bf w}_{i,k_i}(a), {\bf w}_{i,k_i}(b))$
in the alphabet
 $\{a(i),b(i)\}$ in accordance with Section \ref{apend},
 where the words  are taken from the
  (infinite) set $\bf W$ with Properties (A)-(C), and for every $a_{ij}$, we define the endomorphisms  \index[g]{$\psi_{ij}$} $$\psi_{ij}: a(i)\mapsto  {\bf w}_{ij}(a),\;\;
b(i)\mapsto {\bf w}_{ij}(b))$$ of the groups $F(a(i),b(i))$
and of $B(a(i),b(i))$, and define the mappings \index[g]{$\phi_{ij}$}$\phi_{ij}$ on these groups by the rule $\phi_{ij}(u)={\bf w}_{ij}\psi_{ij}(u)$.

 The tape alphabet of $M_1$ is $Y=\sqcup_{i=1}^{k} Y_i$, where $Y_i=Y_i'\sqcup Y''_i$ with $Y_i'=\{a'(i), b'(i)\}$ and  $Y_i''=\{a''(i), b''(i)\}$. The standard base of $M_1$ is $$Q_0P_1Q_1P_2\dots P_kQ_k,$$ where $\sqcup_{i=1}^k Q_i$ is the set of states of $M_0$, $Q_0=\{q_0\}$, and
$$P_i=\{p_i, p_{\theta,i}(1), p_{\theta,i}(2),  p_{\theta,i}(3)\mid \theta\in \Theta^+, i=1,\dots,k\}$$

 For every positive command $\theta$ of $M_1$
we introduce a set of positive rules $\Xi^+(\theta)$ of $M_1$. The S-mashine
\index[g]{$M(\theta)$} $M(\theta)$ (part of $M_1$) with the set of rules \index[g]{$\Xi^+(\theta)$} $\Xi^+(\theta)$ works as follows.

If $\theta$ is given by formula (\ref{eq57}), i.e. it has a unique working head $Q_i$ with part $q_i\omega\to a_{ij}q_i\omega$, where $a_{ij}\in  Y_i(0)=\{a_{i1},\dots, a_{i,k_i}\}$, then
$M(\theta)$ has the following set of rules $\Xi^+(\theta)$.

\medskip

(1) $p_i\tool p_{\theta,i}(1)$, $p_j\tool p_j$ for $j\ne i$, $q_s\to q_s$
for every $q_s$ from the left-hand side of (\ref{eq57}), i.e., the  $p$-letter of the working sector gets
a $\theta$-index.

{\bf Defining the rules of S-machines below we often omit their idling parts of the form $q\to q$.}

\medskip

(2) This gives start to the work of the
copy
of $\overleftarrow Z (\phi_{ij})$ working with the alphabet $Y_i$ (instead of $A$) and having the standard base $Q_{i-1}P_iQ_i$ (instead of $LPR$). When the
state letter $p_{\theta,i}(1)$ meets the left state letters $q_{i-1}$,
the copy of the rule $\xi'(2)$ is applicable, i.e.
\begin{equation}\label{12}
p_{\theta,i}(1)\to {\bf w}'_{ij}p_{\theta,i}(2),\;\; and \;\; also\;\;  p_{i-1}\tool p_{\theta,i-1}(1),
\end{equation}
which activates $P_{i-1}$-head (the head $P_k$ if $i=1$, i.e. the $i$-index is taken  modulo $k$).

\medskip

(3) The machine $\overleftarrow Z(\theta,i-1)$ (a copy of $\overleftarrow Z$
working on $Y_{i-1}$) is switched on now. When $p_{\theta,i-1}(1)$ meets $q_{i-2}$,
a copy of the rule $\xi(2)$ applicable, and by definition, this rule
should also switch on the machine $\overleftarrow Z(\theta,i-2)$ on the base $Q_{i-3}P_{i-2}Q_{i-2}$, i.e. we
have similar work subsequently on the bases between $Q_{i-1}$ and $Q_i$,
$Q_{i-2}$ and $Q_{i-1}$, $\dots$ $Q_0$ and $Q_1$, $Q_{k-1}$ and $Q_{k}$, $\dots$ $Q_i$ and $Q_{i+1}$.
So  machine $\overleftarrow Z (\phi_{ij})$ (for some $j$) works only between $Q_{i-1}$ and $Q_i$, in other intervals, p-letters just subsequently run from right to left.

\medskip

(4) The last rule changing $p_{\theta,i+1}(1) $ by $p_{\theta,i+1}(2)$ switches on the machine $\overrightarrow Z(\theta,i+1)$, the letter
$p_{\theta,i+1}(2)$ runs until it reaches $q_{i+1}$
and becomes $p_{\theta,i+1}(3)$, simultaneously
giving start to $\overrightarrow Z(\theta,i+2)$,
whose work follows by the subsequent work of
$\overrightarrow Z(\theta,i+3),\dots, \overrightarrow Z(\theta,i-1)$ ($i$-index is taken modulo $k$), $\overrightarrow Z(\phi_{ij})$, i.e. again, the copy of the
machine $\overrightarrow Z(\phi_{ij})$ works
between $Q_{i-1}$ and $Q_i$ only.

\medskip

(5) Finally, the rule with parts $p_{\theta,i}(3)\to p_i$ ($i=1,\dots, k$) erases the $\theta$-index in every $p$-letter and makes possible the work of  a machine
$M(\theta')$ (but $M(\theta')$ cannot start working unless  $\theta'$ can be applied for the states of heads $(q_1,\dots, q_k)$).

\medskip

{\it Comment}. The described work of $M(\theta)$ with the standard base replaces
the tape word $u$ in the $Q_{i-1}P_i$-sector with  $\phi'_{ij}(u)$ (where $\phi'_{ij}$ copies $\phi_{ij}$ in the alphabet $Y_i'$) if the command $\theta$
of $M_0$ inserts the letter $a_{ij}$ from the left
of $q_i$. Indeed, at Step (2) we can obtain the copy of $\phi'_{ij}(u)$ in the alphabet $Y''_i$, and the machine $\overrightarrow Z(\phi_{ij})$ rewrites it in the alphabet $Y_i'$ at Step (4). The ''useless'' work of $Z$-machines in other sectors will be exploited for the control of computations, especially if the base is not standard.

\medskip

If the command $\theta$ of $M_0$ has type (\ref{eq56}), then the work  of the machine $M(\theta)$ differs 
in comparison with the case (\ref{eq56}) as follows. At Steps (2) and (4),
the work of machines  $\overleftarrow Z(\phi_{ij})$ and  $\overrightarrow Z(\phi_{ij})$ is replaced
with the work of $\overleftarrow Z(\theta,i)$ and $\overrightarrow Z(\theta,i)$, respectively; the
rule changing $p_{\theta,i+1}(1) $ by $p_{\theta,i+1}(2)$ also changes $q_i$ by $q'_i$.

\medskip

{\it Comment} The described work of $M(\theta)$ just replaces $q_i$ with $q'_i$.

\medskip

If the command $\theta$ of $M_0$ has type (\ref{eq58}), then the machine $M(\theta)$ works as follows.

(1')  $q_{i-1}\tool q_{i-1},\;\;p_i\tool p_{\theta,i}(2), \;\; p_{i-1}\tool p_{\theta,i-1}(1)$,

\medskip

{\it Comment} This rule checks that there are no letters between $q_{i-1}$ and $p_i$ and immediately
 switches on the machine $\overleftarrow Z_{\theta,i-1}$.

\medskip

(2') Then as above the machines $\overleftarrow Z_{\theta,i-1}$,..., $\overleftarrow Z_{\theta,1}$,
$\overleftarrow Z_{\theta,k}$,...,$\overleftarrow Z_{\theta,i+1}$ work subsequently until the rule
$q_i\tool q'_i$, $p_{\theta,i+1}(1)\to p_{\theta,i+1}(2)$
completes moving of the $p$-letters to the left and
changes the $Q_i$-letter.

\medskip

(3')
The obtained set of state
letters makes possible the work of the machine $\overrightarrow Z_{\theta,i+1}$. In turn, this
machine will switches on the machine $\overrightarrow Z_{\theta,i+2}$ after the rule $ p_{\theta,i}(2)\tool p_{\theta,i}(3)$, and so on. Finally, the
rule $ p_{\theta,i}(2)\tool p_{\theta,i}(3)$
follows by the rule $p_{\theta,j}(3)\tool p_j$, where the latter one is
applicable simultaneously for all $j=1,\dots,k$.

\medskip

\begin{df} For every positive word $u$ in the alphabet $Y_i(0)$ of $M_0$, we define
the word \index[g]{$f_i$} $f_i(u)$ in the alphabet $Y'_i=\{a'(i), b'(i)\}$. If $u=1$ (i.e. $u$ is empty), then we define $f_i(u)=1$,
and by induction, if $u\equiv va_{ij}$ for some $a_{ij}\in Y_i(0)$, then
$$f_i(u)=\phi'_{ij}(f_i(v))={\bf w}'_{ij}\psi'_{ij}(f_i(v)),$$ where ${\bf w}'_{ij}, \psi'_{ij}, \phi'_{ij}$ copy ${\bf w}_{ij}, \psi_{ij}, \phi_{ij}$, resp., in the alphabet $Y'_i$.
\end{df}

\begin{lm} \label{smallcan} If $f_i(w)=\phi_{ij}(u)$ for some words $w$ and $u$, then $u$ is a positive
word and the last letter of $w$ is $a_{ij}$.
\end{lm}
\proof It follows from  Property (A) of the set of positive words $\bf W$ that  $\phi_{ij}$-images of  words $u$ are positive for positive words $u$ only. The second statement
is true by the definition of $f_i(v)$ and Remark \ref{disj}.
\endproof

For every configuration ${\cal W}=\alpha_1u_1q_1\omega_1\dots \alpha_ku_kq_k\omega_k$ of the Turing machine $M_0$, we define the configuration
\index[g]{$F({\cal W})$} $W=F({\cal W})\equiv q_0 f_1(u_1)p_1q_1\dots q_{k-1} f_k(u_k)p_kq_k$
of the S-machine $M_1$. The start (stop)
configuration of $M_1$ is the $F$-image of the start (stop) configuration of $M_0$.

\begin{lm} \label{tuda} If ${\cal W}_0\to\dots\to {\cal W}_s$ is a computation
of $M_0$, then there is a computation
$F({\cal W}_0)\to\dots\to F({\cal W}_1)\to\dots\to F({\cal W}_{s-1})\to\dots \to F({\cal W}_s)$  of $M_1$.
\end{lm}

\proof Arguing by induction, we may assume that $s=1$. Let ${\cal W}_1={\cal W}_0\cdot\theta$. If $\theta=\theta_1$ is a positive
 command of type (\ref{eq57}), and so it changes the $i$-th sector inserting a letter $a_{ij}$:
$\alpha_iu_iq_i\omega_i\to \alpha_iu_ia_{ij}q_i\omega_i$, then the work of the S-machine $M(\theta)$ defined above
can transform the configuration $F({\cal W}_0)=\dots q_{i-1}f(u_i)p_iq_i\dots $
to $$\dots q_{i-1}{\bf w}'_{ij}\psi'_{ij}(f_i(u_i))p_iq_i\dots= \dots q_{i-1}\phi'_{ij}(f_i(u_i))p_iq_i\dots=
\dots q_{i-1}f_i(u_ia_{ij})p_iq_i\dots = F({\cal W}_1),$$
as required. If $\theta^{-1}$ is a positive command of type (\ref{eq57}),
then the application of $\theta$ is  possible only if
the last letter of $u_i$ is $a_{ij}$, that is $u_i\equiv v_ia_{ij}$ for some word $v_i$.
Therefore there is a work of $M(\theta^{-1})$ transforming $F({\cal W}_1)=\dots q_{i-1}f_i(v_i)p_iq_i\dots $ to
$F({\cal W}_0)$. The inverse computation
takes $F({\cal W}_0)$ to $F({\cal W}_1)$, as required.

If $\theta$ has type (\ref{eq56}), then ${\cal W}_1$ is obtained from ${\cal W}_0$ by replacement of $q_i$ with $q'_i$
(or vice versa).
Therefore the same is true for the words $F({\cal W}_0)$ and $F({\cal W}_1)$, and the machine $M(\theta)$ can execute 
the similar replacement as well.

If $\theta$ is of type (\ref{eq58}), the
statement is also clear since the canonical work of $M(\theta)$ is possible because the word $u_i$ is empty and so $f_i(u_i)=1$. The machine $M(\theta)$ just replaces $q_i$ with $q'_i$ or vice versa. \endproof

\begin{lm}\label{canon} Let ${\cal C}: W_0\to\dots\to W_s$ ($s>0$) be a reduced
computation of $M(\theta)$ with standard base for some positive command $\theta$ of $M_0$. Suppose the state letters of the configurations $W_0$ and $W_s$ have
no $\theta$-indices while other configurations have $\theta$-indices. Then the computation $\cal C$ is unique, i.e. there is no other reduced computation of $M(\theta)$ starting with $W_0$ and ending with $W_s$.
If $Q_i$ is the working head for $\theta$, then
the only $Q_{i-1}P_i$- and $P_iQ_i$-sectors of $W_0$ and $W_s$ are different, namely,

(1) for a command $\theta$ given by (\ref{eq57}), these sectors of
$W_0$ are $q_{i-1}up_iq_i$ for some reduced word $u$ in the alphabet $\{a'(i),b'(i)\}$
(and some $q_{i-1}\in Q_{i-1}, q_i\in Q_i$), and $W_s$ has the sectors $q_{i-1}\phi'_{ij}(u)p_iq_i$,
or vice versa;

(2) for a command $\theta$ given by (\ref{eq56}), these sectors of
$W_0$ are $q_{i-1}up_iq_i$ for some reduced word $u$ in the alphabet $\{a'(i),b'(i)\}$, and $W_s$ has the sectors $q_{i-1}up_iq'_i$, or vice versa;

(3) for a command $\theta$ given by (\ref{eq58}), the sectors of
$W_0$ are $q_{i-1}p_iq_i$ and $W_s$ has the sectors $q_{i-1}p_iq'_i$, or vice versa.

In particular,
if $W_0=F({\cal W}_0)$ for some ${\cal W}_0$, then $W_s=F({\cal W}_0\cdot\theta^{\pm 1})$.
\end{lm}

\proof Assume that the command $\theta$ is given by (\ref{eq57}). The history of $\cal C$ can start from the rule $p_i\tool p_{\theta,i}(1)$ ((1) in the definition of $M(\theta)$) or with
the inverse rule with $p_{\theta,i}(3)\tool p_i$
($i=1,\dots,k$, (5) in the definition of $M(\theta)$) since the state letters of $W_0$ have
no $\theta$-indices. We consider the first case only.

Since $\cal C$ is reduced, the next rule has to be
the rule of the form $\xi'_1(c)$ of a machine $\overleftarrow Z(\phi_{ij})$ for some $j\in\{1,2\}$. Here $c'=c(ij)$ is the last letter $a(i)^{\pm 1}$ or $b(i)^{\pm 1}$)
of $u$ since otherwise the $p$-letter would
never loose its $\theta$-index by Remark \ref{wrong}. Hence this rule moves the
letter $p_{\theta,i}(1)$ left replacing $c'$ with the copy of the word ${\bf w}_{ij}(c)^{\pm 1}$ in the alphabet $Y''_i$. For the same
reason the subsequent rules  move it
until it meets $q_{i-1}$. Remark \ref{wrong}
now implies that no rules of type $\xi'_1(c)$ are applicable anymore. Hence
the next rule has to be of type (\ref{12}). It changes $p_{\theta,i}(1)$ by ${\bf w}'_{ij}p_{\theta,i}(2)$ and switches on the machine  $\overleftarrow Z(\theta,i-1)$. Then the rules of type
$\xi_1(c)$ of $\overleftarrow Z(\theta,i-1)$ move left the $P_{i-1}$-head until it meets $q_{i-2}$ and turns into
$p_{\theta,i-1}(2)$ switching on the machine
$\overleftarrow Z(\theta,i-2)$, and so on,
until $\overrightarrow Z(\theta,i-1)$ completes work, swithes on the machine  $\overrightarrow Z(\phi_{ij})$ and the final rule executes the transformation $p_{\theta,j}(3)\to p_j$ for all
$j=1,\dots,k$.

Note that the rules of $\overleftarrow Z(\phi_{ij})$ later-by-letter
replace  the word $u$ with the $\psi''_{ij}$-image of it and insert ${\bf w}'_{ij}$ from the left of the $P$-head, while
the rules of $\overrightarrow Z(\phi_{ij})$
just replace subwords $v''$ from the subgroup $H''(ij)$ by the copies $v'\in H'(ij)$, where the subgroup \index[g]{ $H(ij)$} $H(ij)$ is generated by ${\bf w}_{ij}(a)$ and ${\bf w}_{ij}(b)$. In other sectors, the words have been replaced with copies twice, and so Property (1),  of the lemma
follows.

If $\theta$ has type (\ref{eq56}), then the proof can be obviously modified.  If $\theta$ has type (\ref{eq58}), the above argument works as well, because the canonical work of $M(\theta)$ verifies that the $Q_{i-1}P_i$-sector has no tape letters.

The last statement follows from the previous one and the inductive definition of the functions $f_i$ and $F$.
\endproof

\begin{remark} \label{cano} By
Lemma \ref{canon}, the computation $\cal C$
starting with $W_0$, ending with $W_s$ and satisfying the assumptions of this lemma is unique, but it can start either from the
rule of one of the types (1) and (1') or from a rule of one of the
types  $(5)$ and (3') (as in the definition of $M(\theta)$ for a positive rule $\theta$). In the former case, we say that $\cal C$ is  \index[g]{canonical computation
of $M(\theta)$} the canonical computation
of $M(\theta)$, and in the later case, $\cal C$ is the canonical computation of $M(\theta^{-1})$. Thus, for any $\theta$ (positive or negative), the canonical computation of $M(\theta)$ is determined by
the command $\theta$ and the first admissible word $W_0$.
\end{remark}

\begin{lm} \label{Mt} %The described work of the
%machine $M(\theta)$ restricted to a base $Q_{i-1}P_iQ_i$ implies
%that  
A computation ${\cal C}: W_0\to\dots\to W_s$ of the machine $M(\theta)$ with a base $Q_{i-1}P_iQ_i$  has the following property. If  the state letters of
the words $W_0$ and $W_s$ have no $\theta$-indices but other words $W_k$ have $\theta$-indices, then the $Q_{i-1}P_i$-sector words $w_0$ and $w_s$  of $W_0$ and $W_s$ satisfy one of the
conditions: either $w_s\equiv w_0$
or $w_s=\phi'_{ij}(w_0)$, or $w_0=\phi'_{ij}(w_s)$ for some mapping $\phi'_{ij}$.
\end{lm}

\proof If the working head of the command $\theta$ of $M_0$ is $Q_{i'}$ with $i'\ne i$, then only
rules of $\overleftarrow Z(\theta,i)$ and $\overrightarrow Z(\theta,i)$ can change
the $Q_{i-1}P_i$- and $P_iQ_i$-sectors, and $w_s\equiv w_0$ by the projection argument (Remark \ref{proj})
since the sector $P_iQ_i$ is locked by the rules having no $\theta$-indices. 

Assume now that $Q_i$ is the working command of $M_0$ and, as in the proof of Lemma \ref{canon},
assume that we have the part $p_i\tool p_{\theta,i}$ in the first rule of $\cal C$. Then as there,
we see that the machine $\overleftarrow Z(\phi_{ij})$ has to start and complete its work switching
on the machine $\overleftarrow Z(\theta,i-1)$, and so on. However we have no $P_{i-1}$ (and many other heads)
in the base now. So after an application of a number of rules of $\overleftarrow Z(\theta,i-1)$, which change
nothing since the base is $Q_{i-1}P_iQ_i$, the machine $\overleftarrow Z(\phi_{ij})$ can be switched on again.
(Also such a return work can happen after an idling work of $Z(\theta,i-2)$ and return work of $Z(\theta,i-1)$, and so on.)

So the computation $\cal C$ can finish either with the work of $\overrightarrow Z(\phi_{ij})$ as in Lemma \ref{canon}, and we get $w_s=\phi'_{ij}(w_0)$, or the last admissible word $w_s$ looses $\theta$-indices after
a return work of $\overleftarrow Z(\phi_{ij})$ starting with the rule inverse to (\ref{12}). But in the 
latter case all the rules of this second switching on of $\overleftarrow Z(\phi_{ij})$ has to be
uniquely determined (for the similar reason as for the first switching on), and the second work of
$\overleftarrow Z(\phi_{ij})$ is just inverse to the first one. Hence we obtain $w_s=w_0$, and the lemma is proved. \endproof

A computation ${\cal C}: W_0\to\dots\to W_s\to \dots\to W_t$ is called the \index[g]{product of computations}{\it product} ${\cal C}'{\cal C}''$ of two
subcomputations  ${\cal C}': W_0\to\dots\to W_s$ and
${\cal C}'': W_s\to\dots\to W_t$.

\begin{lm} \label{suda} Let ${\cal C}:W_0\to\dots\to W_s$ be a reduced computation of the S-machine $M_1$, where $W_0=F({\cal W}_0)$ for some configurations ${\cal W}_0$ of $M_0$ and the state letters
of $W_s$ have no $\theta$-indices.
Then there is a reduced
computation ${\cal W}_0\to\dots\to {\cal W}_t$ of the Turing machine $M_0$
with a history $H=\theta_1\dots\theta_t$
such that
${\cal C}={\cal C}_1\dots{\cal C}_t$, where ${\cal C}_i$ is the canonical computation
of $M(\theta_i)$ starting and ending with configurations $F({\cal W}_{i-1})$ and $F({\cal W}_i)$,
and $W_s=F({\cal W}_t)$.
\end{lm}

\proof We can decompose ${\cal C}={\cal C}_1\dots{\cal C}_t$, where ${\cal C}_l:
W_{i_{l-1}}\to\dots \to W_{i_l}$ is a computation, without $\theta$-indices in
the state letters of the first configuration
$W_{i_{l-1}}$ and the last one $W_{i_l}$ only. Then by Lemma \ref{canon}, ${\cal C}_i$ is a canonical  computation
of a single S-machine $M(\theta_i)$
since the $\theta$-index can be changed
only after a rule erasing $\theta$-index.
The product
$\theta_1\dots\theta_t$ is reduced, because
by Remark \ref{cano} the equality $\theta_i = \theta_{i-1}^{-1}$ would imply that the
histories of the computations ${\cal C}_{i-1}^{-1}$ and ${\cal C}_i$ start with
the same rule, contrary to the assumption that $\cal C$ is a reduced computation.

By the assumption, $W_0=F({\cal W}_0)$. Assume by induction that for $l>0$ we have a computation ${\cal W}_0\to\dots\to {\cal W}_{l-1}$ of $M_0$ starting with ${\cal W}_{0}$
and ending with ${\cal W}_{l-1}$ such that $W_{i_{l-1}}=F({\cal W}_{l-1})$.

Assume first that $\theta_i$ has type (\ref{eq57}).
If $\theta_i$ is a positive rule, we can define ${\cal W}_l={\cal W}_{l-1}\cdot\theta_i$ and
$W_{i_l}=F({\cal W}_{l})$ by Lemma \ref{canon} and the inductive definition of the functions $f_i$ and $F$.

   If $\theta_i$ is negative, the word $W_{i_{l-1}}$ is obtained from $W_{i_l}$
   after a canonical computation of $M(\theta_i^{-1})$ by Lemma \ref{canon}.
   If $W_{i_l}$ has a subword $q_{i-1}up_iq_i$,
   then the configuration $W_{i_{l-1}}=F({\cal W}_{l-1})$ contains the subword $q_{i-1}\phi_{ij}(u)p_iq_i$
   by Lemma \ref{canon} (1).
   By the inductive conjecture,
   the subword $\phi_{ij}(u)$ of $W_{i_{l-1}}$ is equal to $f_i(v)$ for some $v$. Hence by 
	Lemma \ref{smallcan},
   the word $v$ ends with the letter $a_{ij}$ i.e. $v\equiv v'a_{ij}$, and so $\phi_{ij}(u)=f_i(v)=\phi_{ij}(f_i(v'))$, which implies $u\equiv f_i(v')$ since the mapping $\phi_{ij}$ is injective  by Remark \ref{disj}. It follows that ${\cal W}_{l-1}$ belongs to the domain of the rule $\theta_i$, and we define ${\cal W}_l={\cal W}_{l-1}\cdot\theta_i$. Then by the
   definition of $\theta_i$, the word
   ${\cal W}_l$ has the subword $p_i v' q_i\omega_i$, and the equality $u\equiv f_i(v')$ implies that $W_{i_l}=F({\cal W}_l)$, as desired.

   If $\theta_i$ is a command of type (\ref{eq56}), then both $\theta$ and $M(\theta)$ just replace $q_i$ with $q'_i$, and the statement becomes clear. If $\theta_i$ is a command of type (\ref{eq58}), then the word ${\cal W}_l$
   has the required form as well by Lemma
   \ref{canon} (2), since $1$ is the $f_i$-image of the empty word.
\endproof

   \begin{lm}\label{langu} If $\cal L$ is the language recognized by $M_0$, then the language recognized  by $M_1$ is $f_1({\cal L})=\{f_1(u)\mid u\in {\cal L}\}$.
   \end{lm}

   \proof Given an accepting computation
   $$\alpha_1uq_1\omega_1\alpha_2q_2\omega_2\alpha_kq_k\omega_k\to
   \dots\to\alpha_1q'_1\omega_1\alpha_2q'_2\omega_2\alpha_kq'_k\omega_k$$
of $M_0$, Lemma \ref{tuda} provides us an $M_1$-accepting computation
$$q_0f_1(u)p_1q_1p_2q_2\dots p_kq_k\to
   \dots\to q'_0p_1q'_1p_2q'_2\dots p_kq'_k$$
   since $f_i$-image of the empty word is $1$. Therefore the word $f_1(u)$ is
   accepted by $M_1$.

   Let now a word $U$ in the alphabet $\{a'(1), b'(1)\}$ be accepted by a computation ${\cal C}: W_0\to\dots\to W_t$ of $M_1$,
   where $W_0\equiv q_0Up_1q_1\dots$.
   Note that the accept word $W_t$ is of the form $F({\cal W})$ since $1$ is
   the $f_i$-image of the empty word.
   Applying Lemma \ref{suda} to the computation ${\cal C}^{-1}$, we see that
   that $U=f_1(u)$ for some positive word $u$.
   Now the same lemma applied to $\cal C$,
   provides us with an accepting computation
   ${\cal W}_0\to\dots\to{\cal W}_t$ of $M_0$, where ${\cal W}_0$ is of the form
   $\alpha_1uq_1\omega_1\alpha_2q_2\omega_2\alpha_kq_k\omega_k$,
   i.e. the word $u$ is accepted by $M_0$.
   \endproof

\begin{lm} \label{pos} There are no reduced computations  ${\cal C}:W_0\to\dots\to W_s$
of the S-machine $M_1$ with standard base such that $s>0$,  $W_0\equiv W_s$, and the state letters
of $W_0$  have no $\theta$-indices.
\end{lm}

\proof Proving by contradiction, we consider the factorization ${\cal C}={\cal C}_1\dots{\cal C}_m$, where ${\cal C}_l: W_{i_{l-1}}\to\dots \to W_{i_l}$ and each ${\cal C}_l$ is a canonical computation of some $M(\theta_l)$. We have $m\ge 1$, and it follows from Remark \ref{cano} that the product $\theta_1\dots\theta_m$ is reduced.

 If a head $Q_i$ is the working head for a rule
 $\theta_l$ of the form (\ref{eq57}), then by Lemma \ref{canon}, the tape word $w_{i_{l-1}}$ in
 the sector $Q_{i-1}P_i$ of $W_{i_{l-1}}$ is replaced
 with a word  $w_{i_l}=\phi'_{ij}(w_{i_{l-1}})$ in $W_{i_l}$ (or vice versa) by $M(\theta_l)$. If the head $Q_i$ is not working for the command $\theta_l$ or $\theta_l$ has type (\ref{eq58}), then the restriction of the canonical computation ${\cal C}_l$ to
sectors $Q_{i-1}P_iQ_i$ does not change the tape words of these
sectors by the projection argument (i.e. $W_{i_{l-1}}$ has the same sectors $Q_{i-1}P_iQ_i$ as $W_{i_l}$).

Thus by Lemma \ref{Mt},  we have a sequence of words $w_0, w_{i_1},\dots, w_{i_m}$ in the sectors $Q_{i-1}P_i$, where for every pair $w_{i_{k-1}}, w_{i_k}$
we have  either $w_{i_{k-1}}=w_{i_k}$ or $w_{i_k}$ is a $\phi'_{ij}$-image
of $w_{i_{k-1}}$ for some mapping $\phi'_{ij}$ of
$F(a'(i),b'(i))$, or vice versa. Then
  we can obtain a sequence of positive words $u_0, u_{i_1},\dots, u_{i_m}$ with properties (1) - (3) of
  Lemma \ref{posit}. In particular, $u_0=u_{i_m}$ if $w_0=w_{i_m}$.

  It follows that all the computations ${\cal C}_1,\dots, {\cal C}_m$ can be modified and replaced with canonical computations ${\cal C}'_1,\dots, {\cal C}'_m$, where the restriction of ${\cal C'}_l$ to the sector $Q_{i-1}P_i$ transforms the word $U_{i_{l-1}}$ to $U_{i_l}$, where every $U_{i_k}$ is obtained from $W_{i_k}$ by replacement of the subword $w_{i_k}$ with $u_{i_k}$.

  The modified computation ${\cal C}'_l$ differs from ${\cal C}_l$ not only in the sectors $Q_{i-1}P_iQ_i$; but in the other sectors we just
   vary the number of the rules changing nothing there, since the modifying rules work in
  the  sectors $Q_{i-1}P_iQ_i$ only.
  So one can subsequently make such changes for every $i$ and obtain an $M_1$-computation ${\cal C'}: V_0\to\dots\to V_t$ with $V_0\equiv V_t$. Here $V_0=F({\cal W})$
for some $M_0$-admissible word $\cal W$ by Lemma \ref{posit}, because every
tape subword $v_i$ of $V_0$ is the image of $1$ under a product of mappings of the form $\phi'_{ij}$.

Note that after the modification of $\cal C$, the product ${\cal C}'_1\dots {\cal C}'_m$ is reduced since ${\cal C}'_l$ and ${\cal C}_l$ start/end with the same rule inserting/erasing the $\theta$-index.
Since $V_0=F({\cal W})$,  we obtain by Lemma   \ref{suda}, a reduced computation
of $M_0$ with positive length and equal the first and the last configurations.  But this contradicts Property (3)
of Lemma \ref{bos}.
\endproof

 \begin{lm} \label{inin} Let ${\cal C}:W_0\to\dots \to W_s$ be a reduced $M_1$-computation with $s>0$, standard base and  history $h$. Assume that both words $W_0$ and $W_s$ are input configurations. Then
both of them are $M_1$-accepted and uniquely determined
by $h$.
\end{lm}
\proof By Lemma \ref{pos}, we have $W_s\ne W_0$, and so the corresponding input words $u_0$ and $u_s$ are different.

Let ${\cal C}={\cal C}_1\dots {\cal C}_t$ be the factorization, where every factor ${\cal C}_l$ is the canonical computation of some machine $M(\theta)$. The product $H=\theta_1\dots\theta_t$ is reduced by Remark \ref{cano}.

Since $u_0\ne u_s$, there are commands in $H$ sending a letter from the input tape to the tapes $k-1$ and $k$ (Step 1
according to Remark \ref{pair} (1), but we now subtract $2$ from the numbers of auxiliary tapes). Let ${\cal D}={\cal C}'_1\dots {\cal C}'_l$ be
the (first) maximal subcomputation of  $\cal C$ corresponding to the commands  of this kind,
where every ${\cal C}'_r$ ($1\le r\le l$) is
a product of the factors ${\cal C}_m$-s corresponding to a block of $\theta$-commands
of $M_0$, according to the definition of block
given in Remark \ref{pair}.

Since the sector $Q_{k-1}P_{k}$ is empty in the beginning of $\cal D$, the computation
${\cal C}'_1$ replaces it with $\phi'_{kj}(1)={\bf w}'_{ij}$ for some $j$ by Lemma \ref{canon} (1).
Similarly ${\cal C}'_{2}$ can replace the word ${\bf w}'_{ij}$ with a longer $\phi'_{kj'}$ image of it (not by
$\phi_{kj'}^{-1}$-image by Remark \ref{mono}, since $j'\ne j$, because ${\cal C}'_1$ and  ${\cal C}'_2$ do not correspond to mutually inverse subwords of $H$).
It follows  that
one never obtains $1$ again in the sector $Q_{k-1}P_{k}$ during $\cal D$, and so $\cal D$ is not a suffix of
$\cal C$, because $W_s$ has  this sector empty.

Since the transition  of $M_0$ to Step 2 is possible with empty
first tape only (after a command of type (\ref{eq58})), we have nothing in
the sector $Q_0P_1$ at the end of the computation $\cal D$, according to the
definitions of the machines $M(\theta)$-s
(see the rules (1') - (3')). Hence the input word $u_0$, and so the input configuration
$W_0$ is completely determined by the histories of $\cal D$ and $\cal C$. Furthermore by Lemma \ref{suda}, the input configuration $W_0$ has the form
$F({\cal W}_0)$ for an input word of $M_0$, because
${\cal D}^{-1}$ starts with the empty word $1=f_1(1)$ in $Q_0P_1$-sector. We have a similar form for $W_s$, which is also determined
by the history of $\cal C$.

Therefore we can apply Lemma \ref{suda} to the entire computation $\cal C$ now. The obtained reduced computation ${\cal W}_0\to\dots\to {\cal W}_t$ of $M_0$ has
accepted configurations ${\cal W}_0$ and ${\cal W}_t$ by Lemma \ref{bos} (4).
Therefore by Lemma \ref{langu},
the word $W_0$ (and similarly, the word $W_s$) is accepted by $M_1$, as required.
\endproof

\section{Main machine}\label{mainm}
\subsection{Translator machine}

An auxiliary input-output \index[g]{translator machine $T$} S-machine $T$ is designed for translation of arbitrary   language $\cal L$ of positive words in a finite alphabet $A=\{a_1,\dots,a_k\}$ into the language $f({\cal L})$.

The tape alphabet of $T$ is $A\sqcup Y'\sqcup Y''$, where $Y'=\{a',b'\}$
and $Y''=\{a'',b''\}$.  The set of state
letters of $T$ is
 $\{K\}\sqcup \{L\}\sqcup P\sqcup \{R\}$, where
$P=\{p, p_{a}(1), p_{a}(2), p_{a}(3)\mid a\in A\}$.

The input configuration of $T$ is $KuLpR$,
where the input word $u$ is a reduced word in the alphabet $A^{\pm 1}$, and the output
configuration is $KLvpR$, where the output word $v$ is a word in $Y'$.

The set o positive rules of $T$ consists
of the following rules of 'submachines'
\index[g]{$T(a_j)$} $T(a_j)$-s,  corresponding to  letters $a_j\in A$, $j=1,\dots, k$.

(1) The rule $$K\to K, \;L\to L,\; p\tool p_{a_j}(1),\; R\to R$$
prescribes an $a$-index $a_j$ to $p$ and gives start to

(2) the work of the machine $\overleftarrow Z(a_j)$, which is a copy
of $\overleftarrow Z(\phi_{1j})$ working with the alphabet $Y=Y'\sqcup Y''$, i.e. the $P$-head starts moving left replacing every letter in the $LP$-sector
with the $\psi''_{1j}$-image of the copy of it in the $PR$-sector.  (So the subgroup $Y''(\xi'_{a_j,1})$ associated with every
rule $\xi'_{a_j,1}$ is $H''(1j)$, the image of $F(a'',b'')$ under the homomorphism $\psi''_{1j}$.) When the
state letter $p_{a_j}(1)$ meets $L$,

(3) the next rule $\xi'(2)=\xi'_2(a_j)$ is applicable:

$$K\to K,\; a_jL\tool L,\; p_{a_j}(1)\to {\bf w}'_{1j}p_{a_j}(2),\; R\to R$$

(4) Then the machine $\overrightarrow Z(a_j)$ (a copy of $\overrightarrow Z(\phi_{1j})$
working on $Y$) is switched on.

When $p_{a_j}(2)$ meets $R$,
a `copy' of the rule  $\xi'(4)$ is applicable:
$p_{a_j}(2)\tool p_{a_j}(3)$, and the rule $p_{a_j}(3)\tool p$ removes the $a$-index..

\medskip

{\it Comment.}  The work of the machine $T(a_j)$ in sectors  $LPR$ is similar
to the work of $M_1$ (moreover, it is a copy of the work of $M(\theta)$ restricted to the base $Q_{i-1}P_iQ_i$ if the command $\theta$ of $M_0$ adds $a_j$ to the  $a$-word on tape number $i$) , i.e. it can replace the tape word $u$ of the  $LP$
sector with the copy of $\phi'_{1j}(u)$ in the alphabet $Y'$.
Such a computation also deletes the last tape letter in the sector $KL$ if
it ends with $a_j$.

\begin{remark} \label{can} The definition of \index[g]{canonical computation of  $T(a_j)$} canonical computations $C$ of
the machine $T(a_j)$ is the same as the
definition of the canonical work of $M(\theta)$ in Remark \ref{cano}. The role of $\theta$-index is played by
$a$-{\it index} $a_j$ of $p$-letters. We say that the canonical computation of $T(a_j)$
starts with the application of the rule
$p\tool p_{a_j}(1)$ and the canonical computation of $T(a_j^{-1})$
starts with the application of the rule
$p\tool p_{a_j}(3)$.
\end{remark}

\begin{lm}\label{T} There exists a computation $W_0\to\dots\to W_s$ of $T$ starting with an input configuration $W_0\equiv KuLpR$ and ending with
an output configuration $KLvpR$ if and only if $u$ is a positive
word in the alphabet $A$ and $v=f_1(u^*)$, where $u^*\equiv a_{i_t}\dots a_{i_1}$ is the mirror copy of $u\equiv a_{i_1}\dots a_{i_t}$.
\end{lm}

\proof If $u\equiv a_{i_1}\dots a_{i_t}$ a positive word, then starting with the input  configuration, after the canonical work of the machine $T(a_{i_t})$, we
obtain the configuration $Ka_{i_1}a_{i_2}\dots a_{i_{t-1}}L\phi_{1,i_t}(1)pR$. Then one can subsequently switch on
the machines $T(i_{t-1}), \dots, T(i_1)$ and obtain $$KL\phi_{1,i_1}( \phi_{1,i_2}(\dots \phi_{1,i_t}(1)\dots))pR=KLf_1(u^*)pR$$
by the inductive definition of $f_1$.

Conversely, given a reduced $T$-computation ${\cal C}: W_0\equiv KuLpR\to\dots\to KLvpR$,
we consider the sectors $LPR$, where the machine $T$ works as the machine
$M(\theta)$ for a command $\theta$ of the form (\ref{eq58}).

Therefore the argument of Lemma \ref{suda}
gives $v=f_1(w)$ for some (positive) $w$ and
the factorization ${\cal C= C}_1\dots{\cal C}_t$, where ${\cal C}_l$ is the canonical work of some $T(a_{j_l}^{\pm 1})$. Hence this factorization
corresponds to a reduced word $u_0$ in the
alphabet $A$, and since every $T(a_j^{\pm 1})$
multiplies the sector $KL$ by $(a_j)^{\mp 1}$ from the right,
we have $u_0\equiv u^*$, because the word $u$
is erased in the final configuration.

Let us prove now that there are no negative letters in $u\equiv a_{i_1}\dots a_{i_t}$. Assume $a_{i_t}=a_j^{-1}$ for some $j$. Then the computation should start
with the canonical work of $T(a_j^{-1})$ inserting $a_j$ in the sector $KL$, since otherwise it will increase the
length of the sector $KL$ forever by Remark \ref{wrong}. So $T(a_j^{-1})$ replaces $1$ in the sector $LP$ with the inverse image under $\phi'_{j,1}$, which does not exists since arbitrary $\phi'_{j,1}$-image has positive length. Thus, we have  $a_{i_t}=a_j$.  If $a_{i_{t-1}}=a_{j'}^{-1}$,
then $j'\ne j$, and similarly we obtain $\phi_{j,1}(1)=\phi_{j',1}(*)$, contrary to Remark \ref{disj}. Proceeding by induction, we see that
the whole word $u$ is positive and $w\equiv u_0\equiv u^*$.
\endproof

\subsection{The machine $M$}\label{M}

Let $L$ be a large integer. Consider now $L-1$ copies of the machine $M_1$, denote them by $M_{1}^{(i)}$, $i=1,\dots,L-1$,
 We denote the state and
tape letters of \index[g]{$M_{1}^{(i)}$} $M_{1}^{(i)}$ accordingly, by adding superscript $(i)$
to all tape and state letters
 and to all rules.
 Let \index[g]{$\Xi (M_1)$} $\Xi (M_1)$
be the set
of positive rules of $M_1$  \label{Bst} and $B$ be the standard base of $M_1$,
\index[g]{ $B^{(i)}$} $B^{(i)}$ be the copy of the word $B$ with new superscript $(i)$
added to all letters.
We now consider the $S$-machine \index[g]{$M_2$} $M_2$ with the rules $$\xi(M_2)=[\xi^{(1)},..,\xi^{(L-1)}] ,\;\;
\xi\in \Xi (M_1)$$
(we shall denote $\xi(M_2)$ by $\xi$ too) and the standard base
\begin{equation}\label{baza}
\{t^{(1)}\}B^{(1)}\{t^{(2)}\} B^{(2)} ...\{t^{(L-1)}\}B^{(L-1)}\{t^{(L)}\},
\end{equation}
where $t{(i)}$-s are just separating lettrs. The start/stop words are defined accordingly (every
letter in the standard base is replaced by the corresponding letter in the stop word of $M_1^{(i)}$).

For every admissible word $W$ of $M_1$ with the standard base we denote by \index[g]{$W(M_2)$} $W(M_2)$ the corresponding admissible word $t^{(1)}W^{(1)}t^{(2)}W^{(2)}\dots$   of
$M_2$ with the standard base (of $M_2$), where $W^{(i)}$ is the $i$-th copy of $W$. By definition, $W(M_2)$ is an input
of $M_2$ if $W$ is an input
word of $M_1,$ and so $M_2$ has $L-1$ input sectors.

 The letters in the copy \index[g]{$W^{(i)}$ } $W^{(i)}$ of the word $W$ are
equipped with the extra superscript $(i).$ Thus every $a$-letters
and every $q$-letter
of $M_2$ has
this extra index. We call it  the \index[g]{superscript of a letter} {\it superscript}
of a letter.

By definition the part of a rule $\xi$ of $M_2$ coincides with corresponding part of some rule $\xi^{(i)}$, except for $t$-letters $t^{(1)},\dots t^{(L)}$, which
are separating but never working: $t^{(i)}\to t^{(i)}$ for every rule $\xi$.

 Notice that for every rule $\xi$ of $M_2$ and every admissible word $W$ of $M_1$ with the
standard base of $M_1$, we have $W\cdot \xi=W'$ if and only if $W(M_2)\cdot \xi=W'(M_2).$ So there is no sense of
the obtained extension of $M_1$ if we are interested in
the properties of computations, but the inequality $L>>1$
will be helpful for the study of groups and diagrams
associated with machines.

Next, we want to connect $M_2$ to the translator machine, and at first we extend
the base of $M_2$. By definition,  the standard base of the S-machine \index[g]{$M_3$} $M_3$ is
obtained from the standard base of $M_2$ by  adding  one $t$-letter $t^{(0)}$ from the left: $\{t^{(0)}\}\{t^{(1)}\}B^{(1)}\dots $
But $M_3$ works
exactly as $M_2$ since every rule  $\xi$
of $M_2$ is extended by trivial part $t^{(0)}\tool t^{(0)}$. 
i.e. the new sector $t^{(0)}t^{(1)}$ is always locked by $M_3$.
However now we are able to connect $M_3$ with the translator $T$.

Let \index[g]{$M_4$} $M_4$ be the machine with the standard base
of the same form as $M_3$, that is $$\{t^{(0)}\}\{ t^{(1)}\} B^{(1)}\{t^{(2)}\}B^{(2)}\dots \{t^{(L)}\}$$
The analogy  $\bar\xi$ of every rule $\xi$ of $T$ works in the sector $t^{(0)} t^{(1)}$ as the rule $\xi$ in the sector $KL$,
and $\bar\xi$ copies in every two sectors $Q_0^{(i)}P_1^{(i)}Q_1^{(i)}$
($i=1,\dots, L-1$), the work of $\xi$ in the sectors $LPR$. Thus $M_4$ is built to {\it simultaneously} translate the content of the sector $t^{(0)} t^{(1)}$ to each of the input sectors $Q_0^{(i)}P_1^{(i)}$ of the machines $M_2$.

To connect the S-machine $M_4$ and $M_3$ we
introduce the \index[g]{rule $\xi(43)$} connection rule $\xi(43)$ which locks all the sectors of $M_3$ except for the input sectors $Q_0^{(i)}P_1^{(i)}$ and replaces all state letters of $M_4$ with the corresponding
letters of $M_3$. (However we will not introduce more indices, since it will be clear from the context if a rule
$\xi$ belongs to $M_3$ or to $M_4$.)

Denote by \index[g]{$M_5$} $M_5$ the obtained union of the machines $M_4$ and $M_3$. The input configuration of $M_5$ is the input configuration of $M_4$ (and so the input sector of $M_5$ is $t^{(0)}t^{(1)}$) and the accept configuration of $M_5$ is the accept configuration of $M_3$. The main \index[g]{machine $M$} machine $M$ is the circular form of $M_5$: we identify
$t^{(L)}=t^{(0)}$, and the superscripts are taken
modulo $L$ for $M$. So the standard base $B(M)$ of $M$ is the standard base $B(M_3)$ of $M_3$ 
%or any cyclic permutation of $B(M_3)$. 
However computations of
$M$ can have arbitrary long finite base, for example, $B^{(L-1)}t^{(L-1)}B^{(L)}t^{(0)}t^{(1)} B^{(1)}\dots$.

\begin{lm}\label{unll} Let $Q'Q''$ be a sector with a tape alphabet $\bar Y$ of the standard base of $M$.
Assume that there exist two words $q'uq''$ and $q'vq''$ with base $Q'Q''$ belonging to the domains of two rules $\eta$ and $\xi$, resp., and $\eta$ does not lock the sector $Q'Q''$.
Then

(1)we have $\bar Y(\xi)\le \bar Y(\eta)$;

(2) 
%if neither $Q'$ nor $Q''$ is a $t$-letter and 
if $\xi$ locks the sector $Q'Q''$,
then the application of $\xi$ changes either state $q'$ or $q''$.
\end{lm}

\proof The inclusion is obvious if $\bar Y(\xi)$ is the trivial subgroup or $\bar Y(\eta)=\bar Y$. Otherwise both rules
$\xi$ and $\eta$ belong to a copy of one of the
machines $\overleftarrow Z$, $\overrightarrow Z(\phi_j)$ $\overleftarrow Z(\phi_j)$, $\overrightarrow Z(\phi_j), T(a_j)$. In the later case $Q', Q''\ne K$ since $\bar Y(\xi)$ is nontrivial subgroup and $\bar Y(\eta)\ne \bar Y$, and so the machine $T(a_j)$ copies
 $\overrightarrow Z(\phi_j)$ and $\overleftarrow Z(\phi_j)$ in the sector $Q'Q''$. (See Comment after the definition of $T(a_j)$.) Therefore Statement (1)
follows from Lemma \ref{unl}. Statement (2) can be checked by inspection of the list of $M$-rules.
(For example,  for the pair $\xi_1(c), \xi_2$ and $QQ'=Q_{i-1}^{(j)}P_i^{(j)}$, the rule $\xi_2$ changes the $p$-letter.)
\endproof

Let $\cal L$ be a recursively enumerable set of positive
words in some alphabet $A=\{a_1,..,a_k\}$.
Then the
set ${\cal L}^*= \{w^*\mid w\in {\cal L}\}$ (where $w^*\equiv a_{i_t}\dots a_{i_1}$ if $w\equiv a_{i_1}\dots a_{i_t}$) is
recursively enumerable as well. Therefore there is a Turing machine $M_0$ recognizing the language ${\cal L}^*$, and one
may assume that $M_0$ satisfies the conditions (1) - (7) from
Lemma \ref{bos}. We now assume that the construction of $M_1$, of other auxiliary S-machines, and of the main machine $M$ are based of the Turing machine $M_0$ recognizing the language ${\cal L}^*$.

So by Lemma
\ref{langu}, the machine $M_1$ recognizes the
language $f_1({\cal L}^*)$. It follows from the definition
of the machines $M_2$ and $M_3$ that the
language of their accepted words is also $f_1({\cal L}^*)$.
The difference in comparison with $M_1$ is that the standard bases of $M_2$ and $M_3$ have $L-1$ input sectors, and an input configuration
is accepted iff all the input words are copies (with $L-1$ different
superscripts) of the same word from $f_1({\cal L}^*)$. This observation leads to

\begin{lm} \label{inH} The language recognized by the machine
 $M$ is $\cal L.$
 \end{lm}

 \proof If $u\in {\cal L}$, then by Lemma \ref{T}, there is a computation (of the translator machine $T$ and) of the machine $M_4$ transforming the
 input configuration $t_0ut_1\dots$ of $M$ with input $u$ in the input configuration of $ M_3$ with input words $f_1({\cal L}^*)$ (with different superscripts),
 which is accepted by $ M_3$ as we noticed above. Therefore the word $u$ is accepted by $M$.

 Conversely, assume that an input configuration with input  $u$
 is accepted by a reduced computation ${\cal C}: W_0\to\dots\to W_s$ of $\bf M$.
 The history $h$ of $\cal C$ has the form $h_1\xi(43)h_2$,
 where $h_1$ is the history of the translator machine $T$.
 Assume that $h_2$ contains the rule $\xi(43)^{-1}$, and so
 $h_2=h_3\xi(43)^{-1}h_4$, where $h_3$ is the history of the
 machine $M_3$. Then by Lemma \ref{inin}, the subcomputation
 with history $h_3$ starts with an accepted input configuration
 of the machine $M_3$.  Thus the last configuration of the
 subcomputation with history $h_1$
 has the form $t_0t_1vp_1q_1\dots$, where $v= f_1(w^*)$ for some $w\in {\cal L}$ by Lemma \ref{langu}.
 Moreover by Lemma \ref{T}, we have $w\equiv u$, since the
 mapping $f_1$ is injective. So $u\in \cal L$, as desired.

 Now assume that $h_2$ does not contain rule $\xi(43)^{-1}$, i.e.
 it is the history of an accepting computation of $M_3$. Then
 we obtain the subword $t_0t_1vp_1q_1$ in the beginning of
 the subcomputation with the history $h_2$ as in the previous
 paragraph, and so $u\in \cal L$ again.
 \endproof

 \begin{lm} \label{coeq} Let ${\cal C}: W_0\to W_1\to \dots\to W_s$ be a reduced computation of $M$ with a base $Q_{i-1}^{(j)}P_i^{(j)}Q_i^{(j)}$ and %cyclically simple
a history $h$.
Assume that the state p-letters of $W_0$ and $W_s$ have neither
$\theta$-indices no $a$-indices, and the tape words $u_0$
and $u_s$
of the $Q_{i-1}^{(j)}P_i^{(j)}$-sectors of the words $W_0$ and $W_s$ are
 equal modulo the Burnside relations.
Then $u_s\equiv u_0$.
\end{lm}

\proof The history $h$ of the computation $\cal C$ can be factorized as $(\chi_0) h_1\chi_1h_2\chi_2 \dots h_t(\chi_t)$, where $\chi_j =\xi(43)^{\pm 1}$ or $\chi_j$ is empty  and each $h_l$ corresponds to the canonical work of $M(\theta^{\pm 1})$ for some command $\theta$ of $M_0$ or to the canonical work of $T(a_j^{\pm 1})$. Let $W_{i_l}=W_0(\chi^{\pm 1})h_1(\chi^{\pm 1})h_2(\chi^{\pm 1}) \dots h_l(\chi^{\pm 1})$. If $u_0,\dots, u_t$ are tape words of the sector $Q_{i-1}^{(j)}P_i^{(j)}$ of the words $W_0=W_{i_0},\dots, W_{i_t}$, then by
Lemma \ref{Mt}, for every pair $(u_{l-1}, u_l)$, we
have either $u_l=\phi'_{ij}(u_{l-1})$ for some mapping $\phi_{ij}$ or $u_{l-1}=\phi'_{ij}(u_{l})$, or $u_l\equiv u_{l-1}$. Since $u_s$ and $u_0$ are equal modulo the Burnside relations, we have $u_s\equiv u_0$
by Lemma \ref{eo}, and the statement of the lemma follows.
\endproof

\begin{lm} \label{ztoo} Assume that for some rule $\xi$ of
of $M$ and the tape alphabet $\bar Y$ of some sector,
the canonical image $\tilde Y$ of the subgroup $\bar Y(\xi)$ in the free Burnside group $B(\bar Y)$ contains
both nontrivial in $\tilde Y$ words $u$ and $z^{-1}uz$ for some $z\in B(\bar Y)$.
 Then the word $z$ belongs to $\tilde Y$ too.
\end{lm}

\proof The statement is obvious if $Y(\xi)=\bar Y$. Otherwise $\xi$
is a copy of a rule of one of the machines $\overleftarrow Z$, $\overrightarrow Z$,
$\overleftarrow Z(\phi_j)$, $\overleftarrow Z(\phi_j)$. Therefore
the statement follows from Lemma \ref{zto}.
\endproof

 \section{Interference of Burnside relations}\label{interf}

In this section, we want to show that changing tape words by equal words in the free Burnside group, one
cannot essentially spoil computations of $M$. From now, the language ${\cal L}$ is the language of positive words in the generators of the group $G$, which represent the identity
of $G$. Since the \index[g]{group $G$} group $G$ from the formulation of Theorem \ref{HE} is recursively presented group of finite
exponent, this language is recursively enumerable and
arbitrary relation $w=1$ of $G$ is a consequence of the
relations with left-hand sides from $\cal L$, i.e.
$G=\langle A\mid {\cal L}\rangle$.

We say that two words in the tape subalphabet $\bar Y$ of some sector
are \index[g]{congruent tape words} {\it congruent} if they represent the same element of the group $G$ if $\bar Y$ is the input alphabet of $M$, and they
represent the same element of the free Burnside group of exponent $n$ with basis $\bar Y$ if $\bar Y$ is a  tape alphabet of a non-input sector.

An admissible word is said to be \index[g]{regular admissible word} {\it regular} if it has no subwords $qwq^{-1}$, where $q$ is a $q$-letter and $w$ is an $a$-word congruent to $1$. An $M$-computation $W_0\to\dots\to W_s$
is called \index[g]{regular computation} {\it regular} if $W_0$ (or, equivalently some $W_j$ for $0\le j\le s$) is regular.

\begin{df} We say that two regular
$M$-admissible words $W$ and $W'$ are \index[g]{congruent admissible words} congruent ($W'\cong W$) if they
have the same vector of state letters
$q(1)\dots q(l)$ and their tape words $w_j$, $w'_j$
between $q(j)$ and $q(j+1)$ are congruent
for every $j=1,\dots, l-1$.
\end{df}

\begin{remark} \label{Wcong} Note that if two words $W$ and $W'$ are congruent and $\eta$-admissible for
some rule $\eta$, then the words $W\cdot\eta$ and $W'\cdot\eta$ are also congruent, because
multiplication of the congruent tape
words by the same words from the left/right
preserves the congruence.
\end{remark}

\begin{df} \label{Grewr} We say that ${\cal C}: W_0\to\dots\to W_s$
is a \index[g]{quasi-computation} {\it quasi-computation} with history
$h=\eta_1\dots\eta_s$  if every word $W_r$ is regular and $\eta_{r+1}$-admissible ($r=0,\dots,s-1$),
and for every $r=1,\dots,s-1,$  we have $W_r\cong W_{r-1}\cdot\eta_l$ and $W_s=W_{s-1}\eta_s$. If the history $h$ is reduced, then the quasi-computation $\cal C$ is called reduced.
\end{df}

It follows from Remark \ref{Wcong} that for every quasi-computation $W_0\to\dots\to W_s$ with history $h$,
there is a reduced quasi-computation $W_0\to\dots\to W'_{s'}$ whose history is the reduced form of $h$ and  $W'_{s'}\cong W_s$. The \index[g]{inverse quasi-computation}  {\it inverse quasi-computation} ${\cal C}^{-1}$ with history $h^{-1}$ is $V_0\to\dots\to V_s$, where $V_0\equiv W_s$, $V_s\equiv V_0$, and $V_i\equiv W_{s-i-1}\cdot\eta_{s-i}$ for $i=1,\dots s-1$.

\begin{lm} \label{sw} Let ${\cal C}: W_0\to W_1\to W_2$ be a reduced quasi-computation with history $h=\eta_1\eta_2$ and base $Q'Q''$ for some base
letters $Q'$ and $Q''$. If the rule $\eta_2$
(the rule $\eta_1$) does not lock this sector, then there is a regular $M$-computation ${\cal C}': W'_0\to W'_1\to W'_2$ (resp.,  $W'_2\to W'_1\to W'_0$) with history $h$ (with history $h^{-1}$),
where $W_0'\equiv W_0$, $W'_1\equiv W'_0\cdot\eta_1$, and $W'_2\cong W_2$
(resp., $W'_2\equiv W_2$, $W'_1\equiv W'_2\cdot\eta_2^{-1}$, and $W'_0\cong W_0$).

\end{lm}
\proof Let $\bar Y$ be the tape alphabet
of the sector $Q'Q''$. Since the word $W_1$ is $\eta_2$-admissible and has the same state letter as the word $W_0\cdot\eta_1$ (which is $\eta_1^{-1}$-admissible), we have $\bar Y(\eta_1^{-1})\le \bar Y(\eta_2)$ by Lemma
\ref{unll}. Therefore the word $W'_1=W_0\cdot\eta_1$
is $\eta_2$-admissible. Since it is congruent to $W_1$, the word $W_2'=W'_1\cdot\eta_2$ is congruent to $W_1\cdot\eta_2=W_2$ by Remark \ref{Wcong}.

The second statement of the lemma 
follows from the first one
after one replaces $\cal C$ with ${\cal C}^{-1}$.
\endproof

\begin{lm} \label{exten} Let ${\cal C}:W_0\to\dots\to W_s$ be a reduced computation  of $M$ with base $t^{(i)}B^{(i)}t^{(i+1)}$,
(which is the subbase of the standard base of $M$),
where $t_i\ne t_0$. Assume that $|W_0|_a=|W_s|_a=0$ and each of $W_0$, $W_s$ has the vector of state letters, which is a part of either start
or end configuration of $M$.
Then $\cal C$ is a restriction of a quasi-computation
${\cal D}: V_0\to\dots\to V_s$ with standard base, where each of $V_0$, $V_s$
is either an accepted input  configuration of $M$ or the accept configuration,
%$|V_0|_a=|V_s|_a=0$
and the restriction of $\cal D$ to the subbase $t_1\dots t_L$ is a computation.
\end{lm}

\proof We obviously can extend $\cal C$ to the base
$t^{(1)}B^{(1)}t^{(2)}\dots t^{(L-1)}B^{(L-1)}t^{(L)}$ since
the rules of $M$ act uniformly on each subbase $t^{(k)}B^{(k)}t^{(k+1)}$, $k\ne 0$. We will use the same
notation $\cal C$ for such an extension. So it remains to extend
$\cal C$ to the sector $t_0t_1$.

The history $h$ can be written as $h=h_1\chi_1h_2\chi_2\dots h_l$, where
$\chi_k=\chi(43)^{\pm 1}$ ($k=1,\dots,l-1$),
and $h_1,\dots,h_l$ are the histories of either $M_3$ or $M_4$. We will induct on the number $l$, and so we make the assumption of the lemma weaker allowing the last admissible word
be in the domain of  rule $\chi(43)^{\pm 1}$ too (which unlock sectors $Q_0^{(j)}P_1^{(j)}$). So by the inductive conjecture, we can extend the computation with history
$h_1\chi_1h_2\chi_2\dots h_{l-1}$ and obtain either the empty
sector $t_0t_1$ at the end on the quasi-computation or a word $u$ from the language $\cal L$ in this sector if 
$h_{l-1}$ is a history of $M_4$.

If $h_l$ is a history of $M_3$, then before the aplication of $\chi_l$ , one can replace
$u$ with the congruent empty word since $u=1$ in $G$ the further extension
is obvious since every rule of $M_3$ locks the sector
$t_0t_1$.

Let now  $h_l$ be a history of $M_4$. We first assume that $\chi_{l-1}=
\chi(43)^{-1}$. Then either  $h_{l-1}^{-1}$ is a history of an accepting
computation ${\cal E}: U\to \dots\to U'$ of $M_3$ (if $l=2$)
or it is a computation
of $M_3$ starting and ending with input configurations.
In both cases $U$ is an accepted input configuration
of $M_3$ (see Lemma \ref{inin} for the second case).
Thus, the input tape word $v$ of $U$ (in the sector $Q_0^{(1)}P_1^{(1)}$) belongs to the language $f_1({\cal L}^*)$ (see Section \ref{M}).
So there is $u\in \cal L$ such that $v=f_1(u^*)$,
and by Lemma \ref{T}, one can extend the subcomputation
of $\cal C$ with history $h_l$ to the sector $t_0t_1$,
and obtain at the end the accepted input configuration $t_0ut_1\dots t_L$. 
%If $u$ is not empty, then one can replace it with an
%empty word and obtain the required quasi-computation. (We use that every word from $\cal L$ is trivial in %the group $G$.)

It remains to assume that  $l=1$.
Then $h_l$
is a product $g_1\dots g_k$, where
every $g_r$ is the history of the canonical work of a submachine
$T(a_{i_r}^{\pm 1})$ for some letter $a_{i_r}$.
Note that $g_{s-1}$ and $g_s$ cannot be associated with $a$ and $a^{-1}$ for some $a\in Y_0^{\pm 1}$ since the computation $\cal C$ is reduced. Hence the word $a_{i_1}^{\pm 1}\dots a_{i_k}^{\pm 1}$ is reduced to.
The canonical work of every machine $T(a_{i_r})$ corresponds to the application
of the mapping $\phi_{1,i_r}$ to the sector $Q_{i-1}^{(j)}P_i^{(j)}$. Hence if $k\ge 1$, then we have a contradiction
with Lemma \ref{produ}. If $k=0$, then the sector $t_0t_l$ can be left empty under the computation with history $h_l$, as required.
\endproof

We say that a history $h$ is \index[g]{stable history} {\it stable} if all the rules of $h$ do not change the states of the heads of the S-machine $M$. According to the definition of $M$ this means that all rules from $h$ are  the copies of  the rules $\xi_1(*)$ or all of them are copies of the rules $\xi_3(*)$ of $Z$-machines, or all the rules are copies of the rules $\xi_1'(*)$ of
the machine $\overleftarrow Z(\phi_{ij})$, or all the rules are copies of the rules $\xi_3'(*)$ of
the machine $\overrightarrow Z(\phi_{ij})$. (Only one working $p$-head from
the standard base of $M_1$ or of $T$ is moving, although
the standard base of $M$ has $L-1$ copies of such a moving head).

A reduced history is called \index[g]{simple history} {\it simple} if it has no non-empty maximal stable subhistories $h$, such that the word $h$ is trivial in the free Burnside group of exponent $n$.

\begin{lm} \label{si} For every quasi-computation $W_0\to\dots\to W_s$ with a history $h$, there is a quasi-computation
$W'_0\to\dots\to W'_{s'}$ with a simple history $h'$ such that $W'_0\cong W_0$, $W'_{s'}\cong W_s$, and
$h'$ is equal to $h$ modulo the Burnside relations.
\end{lm}
\proof Assume that $h$ is not simple, i.e. a nontrivial maximal stable sub-quasi-computation ${\cal D}: W_j\to\dots\to W_l$ has history $h_0$ trivial in the free Burnside group. Since $\cal D$ is stable,  the words $W_j$ and $W_l$ have the same vector of state letters. They also have congruent tape subword in each sector. Indeed, if
the working head $P_i$ of $\cal D$ works as the running head in one of the machines $\overleftarrow Z(\theta,i)$, $\overrightarrow Z(\theta,i)$, then it consequently
inserts/deletes the copies of letters from $h_0$ in the tape alphabets from the left and from the right, and so $\cal D$
does not change sector words modulo the Burnside relations. If $P_i$ works as in $\overrightarrow Z(\phi_{ij})$, then every rule from $h_0$ multiplies the sector tape words by the (copy of the) corresponding generator of the group $H(j)$.
So $\cal D$ changes tape words by equal words modulo the Burnside relations too. If $P_i$ works as in $\overleftarrow Z(\phi_{ij})$, then one may assume that $\cal D$ corresponds to the positive command $\theta$
and use both above arguments.

Since we have $W_l\cong W_j$, the sub-quasi-computation $\cal D$ can be removed from $\cal C$. The obtained quasi-computation can be made reduced. This procedure proves the lemma.
\endproof

Two reduced
histories $h\equiv h_1g_1h_2g_2\dots g_{l}h_{l}$  and $h'\equiv h'_1g_1h'_2g_2\dots g_{l}h'_{l}$, where $h_i$, $h'_i$ are all maximal stable subhistories (here $h_1$ or/and $h_l$ can be empty),  are called \index[g]{congruent histories} {\it congruent}
if $h'_i=h_i$ ($i=1,\dots, l$) in the free Burnside group of exponent $n$.

\begin{lm} \label{lock1} Let ${\cal C}:W_0\to W_1\to
W_2$ be a reduced quasi-computation with a history $h=\eta_1\eta_2$ and the base
$(Q_{i-1}^{(j)}P_i^{(j)}Q_i^{(j)})^{\pm 1}$ 
for some $i\in \{1,\dots,k\}$, $j\in \{1,\dots,L-1\}$.
Suppose that
the rule $\eta_2$ locks one of the sectors $Q_{i-1}^{(j)}P_i^{(j)}$, $P_i^{(j)}Q_i^{(j)}$, but not both.

Then there is a reduced regular $M$-computation $\cal C'$ starting with $W_0$ whose simple history $h'$
 is congruent to $h$  and the last  word of $\cal C'$ is congruent
 to $W_2$.
\end{lm}
\proof
(a) We will omit superscripts and assume that $\eta_2$ belongs to $\overleftarrow Z(\phi_{ij})$ and locks the sector $Q_{i-1}P_i$.
Then we assume first that $W_0\cdot\eta_1=
 q_{i-1}up_{\theta,i}(1)vq_i$ for some $a$-words $u$ and $v$.
 The congruent
 word is $W_1\equiv  q_{i-1}p_{\theta,i}(1)v'q_i$ (and so it
 is $\eta_2$-admissible). We consider
 an extension $\cal D$ of the computation $W_0\to W_0\cdot\eta_1$
 by a stable computation of $\overleftarrow Z(\phi_{ij})$ with the moving state letter
 $p_{\theta,i}(1)$, which erases the word
 $u$ from the left. Since $u$ is equal to the empty word in the free Burnside group,
 the same extension replaces $v$ by a
 word $v''v$, where $v''$ is trivial modulo the Burnside relations since $v''$ is (a copy of) the image of $u$ under a homomorphism $\psi'_{ij}$.

If $\cal E$ is the reduced form of $\cal D$, we have that the
 history of $\cal E$ is congruent to $h$ and the last admissible word is equal to
 (the reduced form of) $q_{i-1}p_{\theta,i}(1)v''vq_i$ The only difference in comparison with $W_1$ is the occurrence $v''v$ instead of $v'$.
 Since $v''v$ and $v'$ are equal modulo the Burnside relations, it remains to obtain
 ${\cal C}'$ applying Lemma \ref{sw} to the
 quasi-computation in sector $P_iQ_i$ whose history is product of the last rule of $\cal E$ and $\eta_2$.

 If $W_0\cdot\eta_1=
 q_{i-1}up_{\theta,i}(2)vq_i$ and therefore we have $W_1=q_{i-1}{\bf w}'_{ij}p_{\theta,i}(2)v'q_i$,
  where $u={\bf w}'_{ij}u'$, $u'\in H'(ij)$, and $u$ is equal to ${\bf w}'_{ij}$
  modulo the Burnside relations. Then instead of
  $\overleftarrow Z(\phi_{ij})$ in the previous paragraph, now we exploit the
  machine $\overrightarrow Z(\phi_{ij})$, which however
  moves left erasing the word $u'$ and copying this
  word from the right of the $P_i$-head in the $P_iQ_i$-sector.

 If the rule  $\eta_2$ locks the sector  $P_iQ_i$,
 one can also argue as above. For example, if
 $W_0\cdot\eta_1=
 q_{i-1}up_{\theta,i}(1)vq_i$ and $v$ is trivial in the
 free Burnside group, then by Property (C), it is a trivial
 word in the free Burnside subgroup $H''(ij)$ freely generated by the words ${\bf w}''_{i}(a)$ and ${\bf w}''_{ij}(b)$. Hence there is
 an extension $\cal D$ of the computation $W_0\to W_0\cdot\eta_1$
 by a stable computation with the moving state letter
 $p_{\theta,i}(1)$ of $\overleftarrow Z(\phi_{ij})$, which erases the word
 $v$. The history of this moving is trivial modulo the
 Burnside relation, and so the word $u$ will be replaced
 by a congruent word, and one can finish as above but using
 that the sector $Q_{i-1}P_i$ is unlocked by $\eta_2$ now.

 If $\eta_2$ belongs to a copy of $\overleftarrow Z$ and $\overrightarrow Z$, then in the sectors
 $Q_{i-1}P_iQ_i$, we have the work of 'more primitive'
 machines copying, and the proof is easier.
\endproof

\begin{remark} \label{W0H} In the first (the second) variant of Lemma \ref{sw} and in Lemma \ref{lock1}, the constructed computation $\cal C'$ depends on the first word $W_0$ (resp., the last word $W_s$) and the history of $\cal C$ only. In particular, the replacement
of $W_1$ with a congruent word can change the quasi-computation $\cal C$, but cannot change the computation ${\cal C}'$.
\end{remark}

\begin{lm} \label{QPQ} Let ${\cal C}:W_0\to\dots\to
W_s$ be a quasi-computation with
simple
history $h=\eta_1\dots\eta_s$ and the base $Q_{i-1}^{(j)}P_i^{(j)}Q_i^{(j)}$ for some $i\in \{1,\dots,k\}$, $j\in \{1,\dots,L-1\}$.
Suppose that

(1) exactly one rule of the history locks both $Q_{i-1}^{(j)}P_i^{(j)}$- and $P_i^{(j)}Q_i^{(j)}$-sectors, namely, the
 first rule $\eta_1$ or

(2) exactly two rules of the history lock both $Q_{i-1}^{(j)}P_i^{(j)}$- and $P_i^{(j)}Q_i^{(j)}$-sectors, namely, the
 first rule $\eta_1$ and the last rule $\eta_s$.

 Then there is a reduced regular $M$-computation $\cal C'$ starting with $W_0$ whose (simple) history $h'$
 is congruent to $h$  such that the last admissible word of $\cal C'$,

 (a) under the assumption (1), is congruent
 to $W_s$ and

 (b) the under assumption (2), is $W_s$.

 \end{lm}

 \proof We will omit the superscripts $(j)$ in the proof. Inducting on $t$, we want for every prefix $h(t)$ of $h$ of length $t$, $1\le t\le s$, to construct a regular reduced computation ${\cal C}'(t)$ of $M$ such that the history $h'(t)$ of ${\cal C}'(t)$ is congruent to $h(t)$ and the last admissible word of ${\cal C}'(t)$ is congruent to $W_t$.
 Since under assumption (2), the last rule of of $h=h(s)$ locks
 all the sectors, the last word of the computation with history $h'(s)$ must
 be equal to $W_s$. If $t=1$ this property
 is trivial with $h'(1)=h(1)$.

 Assume that $t\ge 2$ and the required ${\cal C}'(t-1): V_0\to\dots\to V_{i_{t-1}}$ is
 constructed, where $V_{i_{t-1}}\cong W_{t-1}$. If the rule $\eta_{t}$ unlocks both sectors
 $Q_{i-1}P_i$ and $P_iQ_i$, then
 we can apply Lemma \ref{sw} to the quasi-computation $V_{i_{t-1}-1}\to W_{t-1}\to W_t$ and obtain the required computation ${\cal C}'(t):
 V_0\to\dots\to V_{i_{t-1}}\to W'_t$, where
 $W'_t\cong W_t$ with history $h'(t)$ congruent to $h(t)$.
  Similarly, if $\eta_t$ locks only one of the sectors,
   then
  we can apply Lemma \ref{lock1} to the quasi-computation
  $V_{i_{t-1}-1}\to W_{t-1}\to W_t$ and obtain ${\cal C}'(t)$.

Assume now that $\eta_t$ locks both sectors.
 It follows from the assumption (2) that
 $V_{i_{s-1}}\equiv q_{i-1}upvq_i \cong W_{s-1}$,
 where both $u$ and $v$ are congruent to the empty word since the rule $\eta_s$ locks both
 sectors $Q_{i-1}P_i$ and $P_iQ_i$ and the word $W_{s-1}$ is $\eta_s$-admissible.
Therefore the words $u$ and $v$ are freely trivial by Lemma \ref{Mt} and Lemma \ref{eo}. We  obtain
the required computation ${\cal C}'$ extending ${\cal C}'(s-1)$ with $V_{i_{s-1}}\to V_{i_{s-1}}\cdot\eta_s\equiv W_s$.
\endproof

\begin{df}We call the base $Q(1)Q(2)\dots Q(k)$ of a computation or a quasi-computation ${\cal C}: W_0\to W_1\to \dots\to W_t$  \index[g]{revolving base} {\it revolving} if $k>1$ and $Q(1)=Q(k)$.
\end{df}

If the computation $\cal C$ with revolving base has history $h$, then for every $j=1,\dots k$ there is the computation (quasi-computation) ${\cal C'}: W'_0\to W'_1\to \dots\to W'_t$ with history $h$ and
the revolving base  $Q(j)Q(j+1)\dots Q(k)Q(2)\dots Q(j-1)Q(j)$, where every word $W'_i$ is obtained from $W_i$
by the corresponding  cyclic permutation of the sectors; $\cal C'$ is a \index[g]{cyclic permutation of a computation with revolving base} {\it cyclic permutation} of $\cal C$.
So $Q(1)Q(2)\dots Q(k)$ can be regarded as the circular base of $\cal C$.

 \begin{lm} \label{mnogo} Let ${\cal C}:W_0\to\dots\to
W_s$ be a quasi-computation with a revolving base,
with simple history $h=\eta_1\dots\eta_s$, and let the copies of the word $(Q_{i-1}^{(j)}P_i^{(j)}Q_i^{(j)})^{\pm 1}$ (with equal or different superscripts but with the same $i\ge 1$) occur in the cyclic base of ${\cal C}$  $m>0$ times. Suppose
that at least one rule $\eta_r$ locks both $Q_{i-1}^{(j)}P_i^{(j)}$- and $P_i^{(j)}Q_i^{(j)}$-sectors.
Then there is a quasi-computation ${\cal C'}: W'_0\to\dots\to
W'_{s'}$ such that

 (a) the simple history $h'$ of $\cal C'$ is congruent to $h$;

 (b) the restriction of $\cal C'$ to every subbase
 of the form $(Q_{i-1}^{(j)}P_i^{(j)}Q_i^{(j)})^{\pm 1}$
 is a reduced regular $M$-computation;

 (c) if the restriction $\cal D$ of $\cal C$ to a subbase $Q'Q''$, where neither $Q'$  nor $Q''$ is equal to $P_{i}^{(*)}$, is a computation of $M$, then the restriction $\cal D'$ of $\cal C'$ to this subbase is also a regular computation of $M$;

(d) the first and the last admissible words of $\cal C'$ are congruent to $W_0$ and $W_s$, respectively;

(e) if the state letters of an admissible word $W'_l$ have neither $\theta$-indices nor $a$-indices,
then either the $(Q_{i-1}^{(*)}P_i^{(*)})^{\pm 1}$-sectors of $W'_l$ are all empty
or they are all non-trivial in the free Burnside group of exponent $n$.
\end{lm}

\proof Since the rule $\eta_r$ locks the sectors $(Q_{i-1}^{(j)}P_i^{(j)})^{\pm 1}$ and   $(P_i^{(j)}Q_{i}^{(j)})^{\pm 1}$, the letters $P_i^{(j)}$ (with arbitrary
superscript $(j)$) can occur in the revolving base (taken up to cyclic permutations) solely between $(Q_{i-1}^{(j)})^{\pm 1}$ and $(Q_{i}^{(j)})^{\pm 1}$
by Lemma \ref{qqiv}.

Assume first that

(1) exactly one rule of the history locks both $Q_{i-1}^{(j)}P_i^{(j)}$- and $P_i^{(j)}Q_i^{(j)}$-sectors, namely, the
 first rule $\eta_1$, or

(2) exactly two rules of the history lock both $Q_{i-1}^{(j)}P_i^{(j)}$- and $P_i^{(j)}Q_i^{(j)}$-sectors, namely, the
 first rule $\eta_1$ and the last rule $\eta_s$.

Since in the beginning all $m$
pairs of the sectors are empty, the computation ${\cal C}(t)$,
we defined by induction in the proof of Lemma \ref{QPQ},
is the same (up to superscripts) for all pairs of such sectors. (See Remark \ref{W0H}.) Therefore
the statements (a), (b), and (d) follow from Lemma \ref{QPQ}.

Consider now a sector $Q'Q''$, where neither $Q'$  nor $Q''$ is equal to $P_{i}^{(*)}$ and the restriction of $\cal C$ to this sector is an $M$-computation. Then the modification
of the work of the $P_{i}^{(*)}$-heads needed for
the inductive construction of ${\cal C}'(t)$, does
not change the configurations in the $Q'Q''$-sector at all
by the definition of $M$; it can just extend the
 computation in this sector by a longer trivial computation.
 This proves the statement (c) under the above assumption.

 Assume now that the history $h$ of $\cal C$ has  $l>0$  rules $\chi_1,\dots,\chi_l$ simultaneously locking all
 $Q_{i-1}^{(*)}P_i^{(*)}$- and $P_i^{(*)}Q_i^{(*)}$-sectors. So, $h\equiv
 h_0\chi_1h_1\dots\chi_l h_l$, where the subhistories
 $h_0,\dots,h_l $ contain no such rules for given $i$.
 Then we consider the quasi-computations with subhistories
 $h_0\chi_1$, $\chi_1h_1\chi_2$,... and apply the
 statements of Lemmas \ref{QPQ} and \ref{mnogo} under assumptions (1) or (2)
 formulated above. (More precisely, the assumption (1)
 works for the inverse subcomputation with
 history $\chi_1^{-1}h_0^{-1}.$)

 To obtain Property (e), we apply Lemma \ref{Mt} to the subcomputation of ${\cal C'}^{\pm 1}$ (in the sectors $(Q_{i-1}^{(j)}P_i^{(j)}Q_i^{(j)})^{\pm 1}$) starting with the rule locking both sectors and ending with $W'_l$. Then the statement
 follows from Lemma \ref{eo}.
 Thus, the lemma is proved.\endproof

  \begin{lm} \label{open} Let ${\cal C}:W_0\to\dots\to
W_s$ be a quasi-computation with a revolving base,
a simple
history $h=\eta_1\dots\eta_s$, and let the copies of the word
(with equal or different superscripts but with the same subscript $i$) occur in the base of ${\cal C}$  $m>0$ times. Suppose the base of $\cal C$ has no subwords of the form $(P_i^{(j)})^{\pm 1}(P_i^{(j)})^{\mp 1}$
and at least one rule of $h$
locks one of the sectors $Q_{i-1}^{(j)}P_i^{(j)}$ and $P_i^{(j)}Q_i^{(j)}$ but no  rule
locks both these sectors.

Then there is a reduced $M$-computation $\cal C'$  with Properties (a) - (e) of Lemma \ref{mnogo}.
\end{lm}

\proof It follows from the assumption of the lemma that every occurrence of $(P_i^{(j)})^{\pm 1}$ in the revolving base
appears in a subword  $(Q_{i-1}P_iQ_i)^{\pm 1}$. We first assume that there is a rule locking $Q_{i-1}P_i$-sectors in $h$ and there is a rule
locking $P_iQ_i$-sectors. The word $W_0$ has  subwords of the form $(q_{i-1}upvq_i)^{\pm 1}$ (where we omit the superscripts and the indices at $p$). For different subbases $Q_{i-1}P_iQ_i$ of $W_0$, the pairs $(u,v)$ can be different
but
the first components of such pairs are congruent.
Indeed, all of them are empty in the word $W_l$ obtained
after the application of the rule locking the $Q_{i-1}P_i$-sectors, and the application of arbitrary
rule preserves congruence of words by Remark \ref{Wcong}.
 The second components are also congruent by the same reason.

For the beginning, we can therefore replace all such pairs in different subwords with the pair $(u,v)$, changing
the admissible word $W_0$ by a congruent word, and so the word $W_1$ of the quasi-computation can be left unchanged. Indeed, this can be done by Lemma \ref{sw} if the rule $\eta_1$ does not lock $Q_{i-1}P_i$-sector ($P_iQ_i$-sector); otherwise the words in
all $Q_{i-1}P_i$-sectors (resp., all $P_iQ_i$-sectors)
of $W_0$, are copies of the same word.

Now one can repeat the inductive construction of ${\cal C}'(t)$ and $\cal C'$  given in the 
proof of Lemma \ref{QPQ}.
Since the initial subwords of $W_0$ with the base of the form $Q_{i-1}P_iQ_i$ are just copies of each other,
the computations ${\cal C}'(t)$
can be built simultaneously for all subbases $(Q_{i-1}P_iQ_i)^{\pm 1}$
by Remark \ref{W0H}. So we obtain Properties (a) - (d) of Lemma \ref{mnogo} arguing as in Lemma \ref{QPQ}.

 If one of the two sectors $Q_{i-1}P_i$ and $P_iQ_i$ is unlocked by all the rules from $h$, then there is no need to modify $W_0$ making the $a$-words in such sector equal. 
%The cases (2), (3), and (4) are similar.

 To obtain Property (e) from Lemma \ref{mnogo}, we consider a word $W'_l$ without
 $\theta$-and $a$-indices and two computations ${\cal C}_1: W'_l\to\dots\to W'_s$ and ${\cal C}_2: W'_l\to\dots\to W'_0$.
 If a $Q_{i-1}^{(*)}P_i^{(*)}$-sector of $W'_l$ contains a non-empty
 word $w$ trivial in the free Burnside group of exponent $n$,
 then we can replace $W'_l$ by a word $W''_l$ with
 empty $Q_{i-1}^{(*)}P_i^{(*)}$-sectors and replace $W'_{l}\cdot \eta'_{l+1}$ with a congruent word $W''_l\cdot\eta'_{l+1}$.
 So, we obtain a quasi-computation starting with $W''_l$. Then again,
 this quasi-computation gives a computation ${\cal C'}_1$ starting with
 $W''_l$. Similarly  we obtain a computation ${\cal C'}_2$.
 The computation $({\cal C'}_2)^{-1}{\cal C'}_1$ has Property (e), because by Lemmas \ref{Mt} and \ref{eo}, any word $W_k$ without $\theta$-indices in this computation has empty $Q_{i-1}^{(*)}P_i^{(*)}$-sector word $w$
 provided $w$ is trivial modulo the Burnside relations.
 \endproof

\begin{lm} \label{beza} Let ${\cal C}:W_0\to\dots\to
W_s$ be a quasi-computation with a revolving base having no subwords of the form $(P_i^{(j)})^{\pm 1}(P_i^{(j)})^{\mp 1}$ and with a
simple
history $h=\eta_1\dots\eta_s$.  Then there is a quasi-computation $\cal C'$ starting with a word congruent to $W_0$, ending with a word congruent to $W_s$, whose simple history $h'$ is congruent to $h$, and the restriction of $\cal C'$ to
every sector $Q'Q''$ is a regular computation provided $Q'\notin\{t^{(0)}, (t^{(1)})^{-1}\}$.

If the base (or a cyclic permutation of it)
has a subword $Q_{i-1}^{(j)}P_i^{(j)}Q_i^{(j)}$ and $\cal C$ has an admissible word $W_r$ without $\theta$-indices
and $a$-indices,
then in every such admissible word $W'_r$ of $\cal C'$, every
$Q_{i-1}^{(*)}P_i^{(*)}$-sector is either empty or non-trivial
in the free Burnside group of exponent $n$.
\end{lm}

 \proof If we have a 2-letter subbase of the form $Q_{i}Q_{i}^{-1}$ or
 $Q_{i}^{-1}Q_{i}$, then such sectors cannot be locked by Lemma \ref{qqiv}, and the restrictions of $\cal C$ to such base can be replaced with computations having the same history, due
 to Lemma \ref{sw}.

 By the definition of admissible words and the assumption of the lemma, a base letter $(P_i^{(j)})^{\pm 1}$ can occur only in the subbases $(Q_{i-1}^{(j)}P_i^{(j)}Q_{i}^{(j)})^{\pm 1}$ of the revolving base.

If the history $h$ of $\cal C$ has rules locking the sectors
 $Q_{i-1}^{(j)}P_i^{(j)}$ and rules locking the
 sectors $P_i^{(j)}Q_{i}^{(j)}$, then 
%the subbases (3) and (4) are impossible, and one
 %has $(Q_{i-1}^{(j)}P_i^{(j)}Q_{i}^{(j)})^{\pm 1}$ in (2). Then 
the restrictions of ${\cal C}$ to such
 subbases can be replaced with computations by Lemmas
 \ref{mnogo} and \ref{open}.

 In the remaining cases, there are no rules in $h$
 locking the sectors $P_i^{(j)}Q_{i}^{(j)}$ or there are no rules locking the sectors $Q_{i-1}^{(j)}P_i^{(j)}$.
 If nothing is locked, then we refer to Lemma \ref{sw}
 again. Otherwise  one can apply Lemma \ref{open}. To complete the proof, we recall that by Lemmas \ref{mnogo} (c) and \ref{open}, the restrictions of the quasi-computation $\cal C$ to the subbases 
 $(Q_{i-1}^{(j)}P_i^{(j)}Q_{i}^{(j)})^{\pm 1}$ with different subscripts $i$
%of the form (1)- (4) 
can be modified one-by-one.
 \endproof

 \section{Revolving quasi-computations}

We say that a reduced computation (a quasi-computation) $W_0\to\dots\to W_s$ with a history $\eta_1\dots\eta_s$ is \index[g]{idling computation/quasi-computation} {\it idling} if $W_{i-1}\cdot\eta_i\equiv W_{i-1}$ for every $i=1,\dots,s$. An admissible word $W$ of $M$ is \index[g]{passive word} {\it passive} if there exists
a rule of $M$ such that the computation $W\to W\cdot\eta$ is idling.

\begin{lm}\label{pas} Let ${\cal C}: W_0\to\dots\to W_t$ be a quasi-computation of $M$ with history $h$ nontrivial modulo the Burnside relations, with all rules from $M_3$ or all rules from $M_4$,  and with revolving base $B({\cal C})$. Assume that $W_0$ and $W_t$ are passive words with equal vectors of state letters.
Then either

(1) $W_0\cong W_t$ or

(2) all the state letters of all words $W_l$ ($0\le l\le t$) are equal to the same $p^{\pm 1}$, where $p\in P_i^{(j)}$ for some $i,j$. Furthermore, for some command $\theta$ of $M_0$, all the transitions $W_j\to W_{j+1}$ are obtained by application of the rules $\xi_1(\theta,i)(*)^{\pm 1}$ of a machine $\overrightarrow Z(\theta,i)$ or all of them are obtained by application of the rules $\xi_3(\theta,i)(*)^{\pm 1}$ of a
 machine $\overleftarrow Z(\theta,i)$, or all of them are obtained by application of the rules $\xi'_1(\theta,i)(*)^{\pm 1}$ of a machine $\overleftarrow Z(\phi_{ij})$, or all of them are obtained by application of the rules $\xi'_3(\theta,i)(*)^{\pm 1}$ of a machine $\overrightarrow Z(\phi_{ij})$.
\end{lm}

\proof If the base $B({\cal C})$ has $P$-letters only, then we  have Property (2) since
a $p$-letter can be changed only by a rule locking a $P_iQ_i$-sector or locking
a $Q_{i-1}P_i$-sector. (We often omit superscripts.) So one may assume that there is $Q$-letter in the base. Since the words
$W_0$ and $W_t$ are passive, every $p$-letter from $P_i^{\pm 1}$ stays in these words either near a letter $Q_i^{\pm 1}$ or near a letter from 
$Q_{i-1}^{\pm 1}$.

Assume first that the history $h$ is simple. Let we have a subword $(Q_{i-1}P_iP_i^{-1}Q_{i-1}^{-1})^{-1}$ (or $Q_{i}^{-1}P_i^{-1}P_iQ_{i}$) in the revolving base. Then the $P_i$-letters cannot start running since they cannot stop running returning to the $Q$-head
by Lemmas \ref{pR}, because the maximal stable subhistory of a simple history cannot
be trivial in the free Burnside group. So the three sectors $Q_{i}^{-1}P_i^{-1}P_iQ_{i}$ of $W_0$ do not change, up to congruence, in $\cal C$. 

By the projection argument (see Remark \ref{proj}), we have for $W_0$ and $W_t$, the same $a$-projections
of the subwords in subsectors $Q_{i-1}P_iQ_i$ unless the machine $\overleftarrow Z(\phi_{ij})$ works.
Note also, that these subwords cannot be of the forms $q_{i-1}u'p_iq_i$ and $q_{i-1}p_iu''q_i$ with non-empty $u'$, $u''$ since $W_0$ and $W_t$ are passive. So we obtain Condition (1) if no $\overleftarrow Z(\phi_{ij})$ works, because the sector $t_0t_1$ is locked by $M_3$ and sectors $t_0t_1$ and $t_1^{-1}t_1$ are impossible
for a computation of $M_4$ having passive words.

It remains to consider the work of $\overleftarrow Z(\phi_{ij})$under the assumption that Property (1) is false.
The $P_i$-head of one of such submachine  has to run through the whole sector from $Q_i$ to $Q_{i-1}$ since $W_0$ and $W_t$ are passive.
(It cannot start running from a $Q$-head and return to the same head by Lemma \ref{Ephi}.)
This submachine starts working of a machine $M(\theta)$ for a command $\theta$ of type (\ref{eq57}).
After the $P_i$-head reaches $Q_{i-1}$, it locks the sectors $Q_{i-1}P_i$ and $P_{i-1}Q_{i-1}$
and gives start to $\overleftarrow Z(\theta,i-1)$, which also completes its canonical work
since $W_0$ and $W_t$ are passive. Similarly, we obtain the work of each $Z$-machine from the definition of $M(\theta)$, and every sector of the base will be locked by at least one rule. Hence the revolving base
has to contain the standard base of $M_1$, which contradicts to the assumption that the words $W_0$ and $W_t$ are passive.

If $h$ is not simple, then we can replace it with a simple history $h'$ by Lemma \ref{si}. If we have
got Condition (2) for the quasi-computation with history $h'$, the same property holds for $\cal C$ since
only a stable quasi-computation can have a base consisting of $P$-letters.
\endproof

\begin{df} We say that a regular computation of $M$ or a quasi-computation ${\cal C}: W_1\to\dots\to W_t$ is \index[g]{revolving (quasi-)computation} {\it revolving} if

(a) the base $B({\cal C})$ is revolving,

(b) $W_0\cong W_t$,

(c) the history of $h$ is simple and non-trivial modulo the Burnside relations, and

(d) no subword of the (cyclic) word $W_0^{\pm 1}$ is congruent to
an $M$-accepted configuration.
\end{df}

\begin{lm}\label{conidl} Let ${\cal C}: W_0\to\dots\to W_s$ be an idling  revolving quasi-computation with cyclically reduced
history $h$, and a quasi-computation ${\cal C'}: W'_0\to\dots\to W'_{s'}$  have history $h'$ freely equal to $ghg^{-1}$ for some word $g$. If
$W'_0\equiv W_0$, then there is a computation $W_0\to\dots\to W_t$ with history $g$ satisfying
either Condition (1) or Condition (2) of Lemma \ref{pas}. In the latter case both ${\cal C}$ and
${\cal C'}$ satisfy Condition (2) as well.
\end{lm}

\proof Note that either all the rules of $h$ belong to $M_3$ or all of them belong to $M_4$
since the connecting command $\xi(43)^{\pm 1}$ changes all state letters an so it cannot occur in idling
computations.

We have $h'\equiv fh_0f^{-1}$, where the right-hand side is reduced and $h_0\equiv h_1h_2$ is a cyclic permutation of $h\equiv h_2h_1$ with non-empty $h_1$. So in the free group, we have $ghg^{-1}=fh_1hh_1^{-1}f^{-1}$. Therefore
$g^{-1}fh_1$ belongs to the centralizer of $h$, and if $h$ is a power of $r$, where $r$ is not a proper power,
then we obtain $g=fh_1r^l$ for some integer $l$.

Since $h'\equiv fh_0f^{-1}$, we have a quasi-computation $W_0'\to \dots \to V$ with history
$f$ extended to a quasi-computation $W_0'\to \dots \to V\to\dots V'$ with history $fh_0$,
which starts with a quasi-computation $W_0'\to \dots \to V\to\dots V''$ with history $fh_1$.

Recall that by Lemma \ref{unll}, every sector of the base is either locked by every rule of $h$ or is
unlocked by every rule of $h$, because $h$ is idling. Since the rules of non-empty $h_1$ and  $r$
are contained in $h$, one can extend the above quasi-computation and obtain a quasi-computation
${\cal C'}: W'_0\to\dots\to U$ with history $fh_1r^l=g$.

Note that the first word $W'_0\equiv W_0$ is passive since the nontrivial idling computation $\cal C$ starts with $W_0$. The word $U$ is also passive since the last rule of $r^{-1}$ (i.e. the last rule of $h$) or the last rule of $h_1$ (if $l=0$)
is applicable to $U$, and it changes nothing since the word $U$ has the same base as the words of $\cal C$.
Now the reference to Lemma \ref{pas} completes the proof of the lemma, because Condition (2) for a non-trivial $g$ implies the same condition for $\cal C$ and $\cal C'$ since all these quasi-computations have the same base.
\endproof

\begin{lm} \label{prun} Let ${\cal C}: W_0\to W_1\to \dots\to W_t$ be a non-idling revolving quasi-computation of $M$ with a base $B({\cal C})$ and  cyclically minimal in $B(\Xi^+)$ history $h$ corresponding to the rules of the machine $M_3$.
Then ${\cal C}$ satisfies Condition (2) of Lemma \ref{pas}.
\end{lm}

\proof Assume first that the history $h$ is stable. Then every rule of $h$ belongs to some $\overleftarrow Z(\theta,i)$ (or to $\overrightarrow Z(\theta,i)$, or to
$\overleftarrow Z(\phi_{ij})$, or to $\overrightarrow Z(\phi_{ij})$). Then the letters
staying from the right and from the left of $P_i=P_i^{(j)}$ in the revolving base  $B({\cal C})$ cannot be either $Q_{i-1}$- or $Q_i$-letters
by the assumption $W_t\cong W_0$ and Lemmas \ref{pR}
and \ref{Ephi}, because the stable history $h$ is not congruent to the empty word.
Hence by the definition of admissible word, we should
have the subword $P_i^{-1}P_iP_i^{-1}$ in $B({\cal C})$. (The base $B({\cal C})$, being revolving, can be regarded as cyclic word).
Using the same argument, we obtain $P_iP_i^{-1}P_iP_i^{-1}P_i$, and the entire base $B({\cal C})$ is $(P_iP_i^{-1}P_iP_i^{-1}\dots P_i)^{\pm 1}$. Therefore $p$ will never be replaced with other state letter
during the quasi-computation since only a rule locking the sector $Q_{i-1}P_i$ or the sector $P_iQ_i$  can change the states of the head $P_i$. Thus,
under the above assumption, the lemma is proved.

We may now assume that $h$ is not stable.
Since $\cal C$ is revolving, one can construct a quasi-computation ${\cal C}^m$ with history
$h^m$ for arbitrary $m\ge 1$.  Every stable subhistory of $h^m$ is a subhistory of a cyclic permutation of $h$,
because $h$ is not stable. Therefore $h^m$ is simple history since $h$ is a cyclically minimal word.
Moreover, we can construct an infinite in both direction ''quasi-computation'' $\cal D$ whose history
is periodic with period $h$.

Since the history $h$ is not stable,  the $p$-letter  of a machine $ \overleftarrow Z(\theta,i)$
(or $\overrightarrow Z(\theta,i)$, or$\overleftarrow Z(\theta,i)$, or $\overleftarrow Z(\phi_{ij})$), 
starts working with some $W_l$. 
The $P_i$-head
eventually stops working since the histories of $\cal C$ and $\cal D$ are not stable.
It stops running when
the next rule locks the sector $Q_{i-1}P_i$ of the revolving base (or the sector $P_iQ_i$) taken up to cyclic permutations.

Then by the definition of $M_1$ and $M$, the sector $P_{i-1}Q_{i-1}$ is also locked by the last rule, we have $P_{i-1}$ in the revolving base, and the machine $\overleftarrow Z(\theta,i-1)$ has to start working. (There is no difference if we obtain $P_{i+1}$ and $\overleftarrow Z(\theta,i+1)$ starts working.) 
In turn, the machine $\overleftarrow Z(\theta,i-1)$
will soon or later switch on the next one, and so on. (Here we keep in mind that the $P$-head cannot start
and return to the same $Q$-head by Lemmas \ref{pR}, because the maximal stable subhistory of a simple history cannot
be trivial in the free Burnside group.) 
Reversing the history, we see that the machine
$\overleftarrow Z(\theta,i)$ was switched on right after the machine $\overleftarrow Z(\theta,i+1)$
stopped working. It follows that the quasi-computation with history $h^m $ includes for some $m$, the canonical work
of every $\overleftarrow Z(\theta,i)$, $i=1,\dots,k$.

Therefore every sector of the base $B({\cal C})$
is locked by at least one rule from $h$. We obtain similar conclusion when starting with $\overrightarrow Z(\theta,i)$, or $\overleftarrow Z(\theta,i)$, or $\overleftarrow Z(\phi_{ij})$.
Thus, by Lemma \ref{qqiv}, the revolving base $B({\cal C})^{\pm 1}$ must be (a cyclic permutation) of a power of the standard base. Besides,
it follows that the same machine $M(\theta)$
cannot canonically work permanently, and there are rules in  $h$
erasing the $\theta$-index of the state letters.
So the quasi-computation
$\cal C$ has a word $W_l$ without $\theta$-indices.

Consider the quasi-computation ${\cal E}:W_l\to\dots\to W_0\to\dots\to W_l$ of positive length. By Lemma \ref{beza}, there is a computation $\cal F$ of positive length starting
and ending with
words $W$ and $W'$ congruent to $W_l$
and having  history
$h({\cal F})$ congruent to $h({\cal E})$. Moreover, by Lemma
\ref{beza}, we may assume that all $Q_{i-1}P_i$-sectors of $W$ and $W'$ are either empty or non-trivial in the free Burnside group.
Now by Lemma
\ref{coeq} for $\cal F$, we have $W\equiv W'$, and Lemma \ref{pos} for ${\cal F}$ gives a contradiction.
So $h$ is stable, and the lemma is proved.
\endproof

\begin{lm} \label{trans} Let ${\cal C}: W_0\to\dots\to W_t$ be revolving quasi-computation
whose history  $h$ is cyclically minimal modulo the Burnside relations and contains a rule of the machine $M_4$.
 Then $h$ has the rules of $M_4$ only and either

 (1) ${\cal C}$ is an idling quasi-computation, or

 (2) all the state letters of all words $W_l$ ($0\le l\le t$) are equal to the same $p^{\pm 1}$, where $p\in P_1^{(j)}$ for some $j$, all rules of $h$ are
rules $\xi'_1(a_k)^{\pm 1}$ of a machine $\overleftarrow Z(\phi_{k1})$, or the rules $\xi'_3(a_k)^{\pm 1}$ of a machine $\overrightarrow Z(\phi_{k1})$.
\end{lm}

\proof 
If the maximal $M_4$-quasi-computation $\cal C'$ of  $\cal C$ is idling, then all the rules of $h$ have to belong to $M_4$ since the states of the heads are not changed by $\cal C'$. Then we have option (1). So we may assume further that $\cal C'$ is not idling.

Assume first that every rule of $h$ belongs to the copy of a machine $T(a_l^{\pm 1})$. Since $\cal C$ is not idling,
some $P$-letter is working in this quasi-computation as it follows from the definition of $T$ and $M_4$. (When the rule $\zeta_2(a_j)$ inserts/deletes a letter from the left of $L$
it locks the sector $LP$, and so there is a working $P$-head if $L$ is working head.)

Then (as in Lemma \ref{prun}) the letters
staying from the right and from the left of $P=P_1^{(j)}$ in the revolving base  $B({\cal C})$ cannot be $Q_0^{(j)}$ or $Q_1^{(j)}$
by the assumption $W_t\cong W_0$ and Lemmas \ref{pR}
and \ref{Ephi}, because the history is simple.
Hence we obtain
the word $P^{-1}PP^{-1}$ in $B({\cal C})$,...,
and the entire base $B({\cal C})$ is $(PP^{-1}PP^{-1}\dots P)^{\pm 1}$. Therefore the letter $p$ will never be replaced with other state letter
during the quasi-computation since a rule changing the states of the head $P_1^{(j)}$
locks the sector $Q_0^{(j)}P_1^{(j)}$ or the sector $P_1^{(j)}Q_1^{(j)}$. Thus,
the statement of the lemma follows.

Second, proving by contradiction, we assume now that all the rules of $h$ belong to $M_4$, but do not belong to a copy of the same $T(a_l^{\pm 1})$.
It follows that the $P_1^{(j)}$-head
eventually stops running in $\cal C$ and meets a neighbor $Q_1^{(j)}$-head (or $Q_0^{(j)}$-head). Since  the $P_1^{(j)}$-head
 changes the state when it ends running, and $W_t\cong W_0$,
 $P_1^{(j)}$-head has to meet another neighbor $Q_0^{(j)}$ (resp. $Q_1^{(j)}$).
(The $P$-head cannot start
and return to the same $Q$-head by Lemmas \ref{Ephi}.)

Now we consider the restriction ${\cal D}: V_0\to\dots\to V_t$ of $\cal C$ to the subbase
$Q_0^{(j)}P_1^{(j)}Q_1^{(j)}$,
We have a configuration $V_m$ without $a$-indices
in the computation $\cal D$ 
(Otherwise $\cal D$ has to correspond to the work of a single machine $T(a_k^{\pm 1})$.)
%, and $V_0$ and $V_t$ should have different vectors of state letters.)

There exists a quasi-computation ${\cal E}:V_m\to\dots\to V_0\to\dots\to V_m$. By Lemma \ref{beza}, there is a computation $\cal E'$ of positive length starting
and ending with
words $V$ and $V'$ congruent to $V_m$
and having history
$h({\cal E'})$ congruent to $h({\cal E})$. Moreover, by Lemma
\ref{beza}, we may assume that all $Q_0^{(j)}P_1^{(j)}$-sectors of $W$ and $W'$ are either empty or non-trivial in the free Burnside group.
By Lemma
\ref{coeq}, for $\cal E'$, we have $V\equiv V'$.

The history $h({\cal E'})$
is a product $h_1\dots h_l$, where
every $h_r$ is the history of the canonical work of a submachine
$T(a_{i_r}^{\pm 1})$ for some letter $a_{i_r}$.
Note that $h_{s-1}$ and $h_s$ cannot be associated with $a$ and $a^{-1}$ for some $a\in Y_0^{\pm 1}$ since the computation $\cal C$ is reduced. Hence the word $a_{i_1}^{\pm 1}\dots a_{i_l}^{\pm 1}$ is reduced to.
Recall that the canonical work of every machine $T(a_{i_r})$ corresponds to the application
of the mapping $\phi_{1,i_r}$ to the sector $Q_0^{(j)}P_1^{(j)}$. Hence we have a contradiction
with Lemma \ref{produ}.

Assume finally that there are rules of the machine $M_3$ in $h$, and so there is the rule $\xi^{\pm 1}(43)$ connecting
the machines $M_4$ and $M_3$ This rule locks all sectors except for the input sectors of the machine $M_3$, and therefore there is a subbase $(P_1^{(i)}Q_1^{(i)})^{\pm 1}$ in the revolving base.

A rule of $M_3$ next to $\xi(43)^{\pm 1}$  starts moving the $P_1^{(i)}$-head, and this head has to meet $Q_0^{(i)}$ and so on, which again means that the base of $\cal C$ (without the first or the last letter) is a cyclic permutation of a power of the standard subbase.

Since the rule $\xi(43)$ changes the state letters, we
have both $\xi(43)$ and $\xi(43)^{-1}$ in the history.
So there is a subhistory  $\xi(43)h'\xi(43)^{-1}$ in
 $h$ or in a cyclic permutation of $h$, where $h'$
is a history of $M_3$-quasi-computation. 

Since the sector $t_0t_1$ is locked by $M_3$, by Lemma \ref{beza}, we have the $M$-computation $V\to\dots\to V'$ of $M_3$ with reduced history between the applications of  $\xi(43)$ and $\xi(43)^{-1}$.
So
only input sectors of these words can be non-empty. Thus the
restrictions $U$ and $U'$ to the standard base of $M_2$ are input words. By Lemma \ref{inin}, both of them are accepted
by $M_2$, and so $V$ and $V'$ are (without the first or the last state letter) accepted by $M_3$, since
the sectors with different superscripts have to be copies
of each other being determined by the history.
Hence the word $W_0$ (without the first or the last state letter) or a cyclic permutation of $W_0^{\pm 1}$ is congruent to a word accepted by $M$
 contrary
to the assumption of the lemma saying that the quasi-computation $\cal C$ is revolving.
\endproof

\begin{lm} \label{tworev} Let ${\cal C}_1: W_0\to W_1\to \dots\to W_s$ be a revolving  quasi-computation  with a history $h$ cyclically minimal in the free Burnside group $B(\Xi^+)$.
Let ${\cal C}_2: W_0\equiv W'_0\to W'_1\to \dots\to W'_{s'}$ be another revolving quasi-computation whose history $h'$ is freely
equal to $ghg^{-1}$ for a reduced word $g$. Then there is computation with a history $g'$: $W_0\to\dots\to W_0\cdot g' \cong W_0$, where $g'$ is equal to $g$ modulo the Burnside relations.
\end{lm}

\proof If ${\cal C}_1$ is an idling
quasi-computations, then by Lemma \ref{conidl}, either we obtain the required property
of $g$ or ${\cal C}_1$ satisfies Condition (2) of Lemma \ref{pas}. If ${\cal C}_1$ has a rule of $M_4$, then we come to the same conclusion
%but Condition (2) should be replaced with the similar condition from 
by Lemma \ref{trans}. 
If ${\cal C}_1$ is not idling computation of $M_3$, then again we obtain Condition (2)
by Lemma \ref{prun}. It remains to assume that ${\cal C}_1$ satisfies Condition (2).
Since ${\cal C}_2$ has the same base, it satisfies Condition (2) as well.

Let ${\cal D}: p^{-1}u_0p\to \dots\to p^{-1}u_tp$ (the sector index is omitted) be the restriction of ${\cal C}_1$ to a $P^{-1}P$-sector. (The case $PP^{-1}$ is similar.) Note that this computation satisfies the assumption of Lemma \ref{gen2}, and
using notation of that lemma we obtain that
$\lambda(h)$ belongs to the centralizer $\cal Z$ of $u_0$ in $B$.
Similarly, $\cal Z$ contains $\lambda(ghg^{-1})$. The centralizer of nontrivial element $u_0$ is the unique
cyclic subgroup of order $n$ containing $u_0$ (see \cite{book}, Theorem 19.5), and so $\lambda(g)\in \cal Z$ by \cite{book}, Theorem 19.6.
%i.e., $g$ commutes with $h$ in $\cal B$ since $\lambda$ is a monomorphism. 
Observe that $\cal Z$ is also
the centralizer of the non-trivial $\lambda(h)$. It follows that the images $\lambda(g)$ and $\lambda(h)$
commute, and so do $g$ and $h$ since $\lambda$ is a monomorphism. 

The projection of the equality $ghg^{-1}=h$ to the free Burnside subgroup ${\cal B}_0$ generated
by the rules of the stable computation with history $h$ gives us
$g'h(g')^{-1}=h$, where the rules of $g'$ are among the generators
of ${\cal B}_0$, and therefore  there is a stable computation $W_0\to\dots\to W_0\cdot g'$ with history
$g'$. Since $\lambda(g')$ commutes with $\lambda(h)$, 
we have $\lambda(g') \in \cal Z$ too, and therefore
$W_0\cdot g'\cong W_0$ by Lemma \ref{gen2}. 
Besides, 
$\lambda(g')$ belongs to the centralizer $\cal Z'\le \cal Z$ of $\lambda(h)$ in $\lambda({\cal B}_0)$, and since both subgroups
$\cal Z'$ and $\cal Z$ have order $n$, they are equal, and therefore $\lambda(g)\in \lambda({\cal B}_0)$. Since $\lambda$ is a monomorphism, we have $g\in {\cal B}_0$, and so $g'$ is equal to $g$ modulo the Burnside relations.
\endproof

\section{Groups and diagrams related to the machine $M$}\label{gd}

\subsection{Construction of the embedding}\label{conemb}

We are going to introduce a group $M$  associated with the machine $M.$ Since we will use the main properties of $M$ (given by Lemmas \ref{exten}, \ref{beza} and \ref{tworev}) only and shall not use many details, we re-denote the machine $M$ accepting the language $\cal L$ of positive relators of the group $G$ by \index[g]{$\bf M$} $\bf M$ and simplify notation as follows.

(1) For the set of (state) $q$-letters $\{t^{(0)}\}\sqcup \{t^{(1)}\}\sqcup
Q_0^{(1)}\sqcup  P_1^{1}\sqcup Q_1^{(1)}\sqcup\dots\sqcup \{t^{(L)}\}   $ introduced in Subsection \ref{M}, we
also will use the uniform notation ${\bf Q}=\sqcup_{j=0}^N {\bf Q}_j$, where ${\bf Q}_0=\{t^{(0)}\}$, ${\bf Q}_1=\{t^{(1)}\}$, ${\bf Q}_2= Q_0^{(1)}$,...,${\bf Q}_N=\{t^{(L)}\}=\{t^{(0)}\}$. The number of $t$-letters $L$ is chosen large enough, so that $N\ge n$.

(2) The set of  $a$-letters of $\bf M$ is${\bf Y}=\sqcup_{j=1}^N {\bf Y}_j$
including the input alphabet $A\subset {\bf Y}_1$.

(3) The set of rules of $\bf M$ is now denoted by \index[g]{$\bf\Theta$} $\bf\Theta$.

(4) If $u$ is an input word in the alphabet $A^{\pm 1}$, then the corresponding input configuration of $\bf M$ is denoted by \index[g]{$\Sigma(u)$} $\Sigma(u)$.

The finite set of generators of the group $M$ consists of {\em $q$-letters} from $\bf Q$, {\em $a$-letters} from $\bf Y$, and $\theta$-letters from $N$ copies
 \index[g]{$\bf\Theta^+_j$} $\bf\Theta^+_j$ of
$\bf\Theta^+$, i.e., for every $\theta\in \bf\Theta^+$, we have $N$
generators $\theta_0,\theta_1,\dots,\theta_N$,
where $\theta_N=\theta_0$.

The relations of the \index[g]{group $M$} group $M$ correspond to the rules of the machine $M$;
for every $\theta=[U_0\to V_0,\dots U_{N}\to V_{N}]\in \bf\Theta^+$, we have
\begin{equation}\label{rel1}
U_i\theta_{i+1}=\theta_i V_i,\,\,\,\, \qquad \theta_i a=a\theta_i, \,\,\,\, i=1,...,N
\end{equation}
for every word $a$ from the set of free generators of the subgroup ${\bf Y}_i(\theta)\le F({\bf Y}_i)$.
The first type of relations will be
called \index[g]{$(\theta,q)$-relations} $(\theta,q)$-{\em relations}, the second type of relations are
\index[g]{$(\theta,a)$-relations}$(\theta,a)$-{\em relations}.

Finally, the \index[g]{group $\tilde G$} group $\tilde G$ is given by the generators and
relations of the group $M$  and two more
relation, namely the \index[g]{hub-relations} {\it hub}-relations
\begin{equation}\label{rel3}
\Sigma_0=1, \;\; and \;\; \Sigma_1=1,
\end{equation}
where
\index[g]{$\Sigma_0$} $\Sigma_0$ is the accept word
of the machine $\bf M$ and \index[g]{$\Sigma_1$} $\Sigma_1$ is the  start
word $\Sigma(1)$ (where the input subword $u$ is empty).

Suppose a configuration $W'$ of $\bf M$ is
obtained from a configuration $W$ by an application of a rule
$\theta: \; [U_0\to V_0,\dots,U_{N}\to V_{N}]$. This implies that
$W = U_0w_1U_1\dots w_N U_{N},$ $W' = V_0w_1V_1\dots w_N V_{N},$
where  every $w_i$ is a word in the alphabet ${\bf Y_i}(\theta)$, and therefore
$W'=\theta_0^{-1}W\theta_0$ in $M$ by Relations (\ref{rel1}), since
$\theta_0\equiv\theta_N.$

Now suppose $u\in \cal L$.
Then $u$ is recognized by the machine $\bf M$ by Lemma \ref{inH},
and so the word $\Sigma(u)$
is accepted by the machine $\bf M$ and therefore it is conjugate to $\Sigma_0$
in the group $M.$ Consequently, we have  that
$\Sigma(u)=1$ in $\tilde G$  by (\ref{rel1}, \ref{rel3}).

Note that the word $\Sigma(u)$ is obtained from $\Sigma_1$
by inserting  the subword $u$. It follows from (\ref{rel3})
that $u=1$ in $\tilde G$ too.
Since we identify the alphabet $A$ of the input sector $t_0t_1$ with
the alphabet of the generators of $\tilde G$,  Dyck's lemma implies the following.

\begin{lm} \label{homo} The mapping  $a \mapsto a$ ($a\in A$) extends to a
homomorphism of the group $G$ to $\tilde G$.
\end{lm}

\subsection{Bands and trapezia}
To study (van Kampen) diagrams
over the group $M$ we shall use their simpler subdiagrams such as bands and trapezia, as in \cite{O97},  \cite{SBR}, \cite{BORS}, etc.

We will consider
\index[g]{$q$-band}
$q$-{\it bands}, where the set $\mathcal Z$ from the definition in Subsection \ref{bds} is one of the sets ${\bf Q}_i$ of state letters
for the $S$-machine $\mmm$, and
\index[g]{$\theta$-band}$\theta$-{\it bands} for every $\theta\in\bf\Theta$, where ${\mathcal Z}=\{\theta_0.\theta_1,\dots,\theta_N\}$.

The following is a consequence of the definition of the relations of $M$.

\begin{lm}\label{tba} (\cite{OS19}, Lemma 5.4). Let ${\mathbf e}\iv {\mathbf q}_1{\mathbf f} {\mathbf q}_2\iv$ be the boundary of a $\theta$-band $\mathcal T$ with bottom ${\mathbf q}_1$
and top ${\mathbf q}_2$ in a reduced diagram.

(1) If the start and the end edges ${\mathbf e}$
and ${\mathbf f}$ have different labels, then $\mathcal T$ has $(\theta,q)$-cells.

(2) For every $(\theta,q)$-cell $\pi_i$ of $\mathcal T$, one of its
boundary $q$-edges belongs to ${\mathbf q}_1$ and another one
belongs to ${\mathbf q}_2$.   $\Box$
\end{lm}

To construct the top (or bottom) path of a $\mathcal Z$-band $\mathcal T$, at the beginning
one can just form a product ${\mathbf x}_1\dots {\mathbf x}_n$ of the top paths ${\mathbf x}_i$-s of the cells $\pi_1,\dots,\pi_n$ (where each $\pi_i$ is a $\mathcal Z$-bands of length $1$).
If this product is not freely reduces, then one obtains the
freely reduced form after a number of folding of the subsequent edges with mutually inverse labels.

\begin{remark} \label{tbd}
(1) No  $\theta$-letter is being canceled in the word
$W\equiv Lab({\mathbf x}_1)\dots Lab({\mathbf x}_n)$ if $\mathcal T$ is  a $q$-band since
otherwise two neighbor cells of the band would make the diagram non-reduced. So for every subword $\theta(1)u\theta(2)$ of the
reduced top/bottom label $V$ of a $q$-band, the subword $u$ in $a$-letters is uniquely defined by the pair $(\theta(1),\theta(2))$ of $\theta$-letters.

We will call an arbitrary reduced word $V$ in $\theta$- and $a$-letters
$\theta$-{\it defined} if for every subwords  $\theta(1)u\theta(2)$ and $\theta(3)u'\theta(4)$, the equalities
$\theta(1)=\theta(3)$ and $\theta(2)=\theta(4)$ imply the equality
of the middle $a$-words: $u'\equiv u$.

(2) By Lemma \ref{tba} (2), there are no cancellations of
$q$-letters of
$W$ if $\mathcal T$ is  a $\theta$-band.

If $\mathcal T$ is a $\theta$-band then a few cancellations of $a$-letters (but not $q$-letters) are possible in $W.$
We will always assume
that the top/bottom label of a $\theta$-band is a reduced form of the word $W$.
\end{remark}

In the next lemma, by \index[g]{trimmed word} trimmed word $w$ we mean the maximal subword of $w$ starting and ending with $q$-letters.

\begin{lm} \label{perm} (1) The trimmed bottom and top labels $W_1$ and $W_2$ of any reduced $\theta$-band $\mathcal T$ containing at least one $(\theta,q)-cell$ are $\bf M$-admissible
and $W_2\equiv W_1\cdot\theta$.

(2) If $W_1$ is a $\theta$-admissible word and $W_2= W_1\cdot\theta$,
then there is a reduced $\theta$-band with  trimmed bottom
label $W_1$ and trimmed top label $W_2$.
\end{lm}

\proof (1) By Lemma \ref{tba} (2), we have $W_1\equiv q_1^{\pm 1}u_1q_2^{\pm 1}\dots u_k q_{k+1}^{\pm 1}$, where $q_j^{\pm 1}$ and $q_{j+1}^{\pm 1}$ are
the labels of $q$-edges of some cells $\pi(j)$ and $\pi(j+1)$
such that the subband connecting these cells has no $(\theta, q)$-cells. Therefore by Lemma \ref{tba} (1), all the $\theta$-edges of $\mathcal T$
between $\pi(j)$ and $\pi(j+1)$ have the same label. It follows
from the list of $(\theta,a)$-relations that the subwords $q_iu_iq_{i+1}$ of the word $W_1$ satisfy the definition of admissible word with $u_i\in F({\bf Y}_j(i))$ if $q_i\in {\bf Q}_j(i)$ or $q_i^{-1}\in {\bf Q}_{j(i)+1}$. Hence $W_1$ is admissible, and so is $W_2$.
The word $u_i$ is a reduced form of a product  $bc_1\dots c_s a$, where the words $c_1,\dots, c_s$ are the bottom labels of $(\theta, a)$-cells of the band $\mathcal T$, and so
$c_1\dots c_s\in {\bf Y}_j(i)(\theta)$, $b$ is the subword in the part $*q_ib\to *$ (or in $b^{-1}q_i^{-1}*\to *$) of the rule $\theta$ and $a$ is the subword of the part
$aq_{i+1}*\to *$ (or in the part $*q_{i+1}^{-1}a^{-1}\to*$) of $\theta$. In all these cases,
the word $u_i$ satisfies the requirement for applying of the rule $\theta$.
(For example, if $q_i\in {\bf Q}_{j(i)}$ and $q_{i+1}\in {\bf Q}_{j(i)+1}$, then
$b^{-1}u_ia^{-1}\in {\bf Y}_j(i)(\theta)$, as required; the other cases are similar.)
Thus, the word $W_1$ belongs to the domain of the rule $\theta$. Likewise $W_2$ belongs
to the domain of $\theta^{-1}$ and finally, $W_2=W_1\cdot\theta$
since the bottom and top labels of the $(\theta,q)$-cells of $\mathcal T$ are, resp., the left-hand sides and the right-hand sides
of the parts of $\theta$.

(2) If we have $q_i\in {\bf Q}_j(i)$ and $q_{i+1}\in {\bf Q}_{j(i)+1}$, then one can construct a part of $\mathcal T$
taking two cells $\Pi_i$ and $\Pi_{i+1}$ corresponding to
the $(\theta,q)$-relations involving the letters $q_i$ and
$q_{i+1}$, respectively, and connecting $\Pi_1$ and $\Pi_2$
by a band of $(\theta,a)$-cells, such that the bottom
labels $c_1,\dots,c_s$ of these $(\theta,a)$-cells form a product freely equal to the word $b^{-1}u_ic^{-1}$ (we use the notation of item (1) here). This is possible since by the definition of $\theta$-application, we have $b^{-1}u_ic^{-1}\in {\bf Y}_j(i)(\theta)$. The trimmed bottom label of the obtained $\theta$-band
(including $\Pi_1$ and $\Pi_1$) is equal to $q_iu_iq_{i+1}$.

The reader can easily verify that the same conclusion is
true if one or both letters $q_i$ and $q_{i+1}$ are negative.
Since such a construction works for every $i=1,\dots, k$, one
obtain a longer $\theta$-band with trimmed bottom and top labels $W_1$ and $W_2$, resp.
\endproof

We say that a reduced $\theta$-band is \index[g]{regular $\theta$-band} {\it regular} if the trimmed top and bottom labels
of it are regular admissible words.

We shall consider the
projections of words in the generators of $M$ onto
$\bf\Theta$ (all $\theta$-letters map to the
corresponding element of $\bf\Theta$,
all other letters map to $1$), and the projection onto the
alphabet $\{{\bf Q}_0\sqcup \dots \sqcup {\bf Q_{N-1}}\}$ (every
$q$-letter maps to the corresponding ${\bf Q}_i$, all other
letters map to $1$).

\begin{df}\label{dfsides}
{\rm  The projection of the label
of a side of a $q$-band onto the alphabet $\bf\Theta$ is
called the \index[g]{history of $q$-band}{\em history} of the band.
The projection of the label
of a side of a $\theta$-band onto the alphabet $\{{\bf Q}_0,...,{\bf Q}_{N-1}\}$
is called the \index[g]{base of $\theta$-band} {\em base} of the band, i.e., the base of a $\theta$-band
is equal to the base of the label of its top or bottom.}
\end{df}

We  call a $\mathcal Z$-band
\index[g]{maximal band} {\it maximal} if it is not contained in
any other $\mathcal Z$-band.
Counting the number of maximal $\mathcal Z$-bands
 in a diagram we will not distinguish the bands with boundaries
 ${\mathbf e}\iv {\mathbf q}_1{\mathbf f} {\mathbf q}_2\iv$ and ${\mathbf f} {\mathbf q}_2\iv {\mathbf e}\iv {\mathbf q}_1,$ and
 so every $\mathcal Z$-edge belongs to a unique maximal $\mathcal Z$-band.

We say that a ${\mathcal Z}_1$-band and a ${\mathcal Z}_2$-band \index[g]{crossing bands}
 {\em cross} if
they have a common cell and ${\mathcal Z}_1\cap {\mathcal Z}_2=\emptyset.$

Sometimes we specify the types of bands as follows.
A $q$-band corresponding to one
of the letter ${\bf Q}_i$ of the base is called a \index[g]{$Q_i$-band}
${\bf Q}_i$-band. For example, we will consider \index[g]{$t$-band}
$t_i$-{\it band} (or just $t$-band) corresponding to the part $\{t_i\}$.

By definition, the group \index[g]{$H(\bf Y)$} $H(\bf Y)$ is the factor group of $M$
obtaining from the presentation of $M$ by adding of \index[g]{$a$-relations} $a$-relations, i.e. all relations $w=1$, where the word $w$ in the alphabet ${\bf Y}^{\pm 1}$ is trivial in the free Burnside group $B({\bf Y})$, and all relations $w=1$, where $w$ is a word in  the subalphabet $A={\bf Y}_1$ trivial in the group $G$.
So there is a homomorphism $G\to H(\bf Y)$. The $H(\bf Y)$-diagrams are considered over this presentation.
Let us call such a disk diagram $\Delta$ \index[g]{minimal $H(\bf Y)$-diagram}  {\it minimal} $H(\bf Y)$-diagram if it is reduced, every $\theta$-band in it is regular and $\Delta$
contains no $\theta$-annuli. A reduced diagram
over the presentation of $M$ (i.e. an $H(\bf Y)$-diagram without $a$-cells) is said to be minimal if it has no $\theta$-annuli.

\begin{lm}\label{NoAnnul}
 For every disk diagram $\Delta$ over $H({\bf Y})$ (over $M$), there
is a minimal diagram $\Delta'$ over $H({\bf Y})$ (resp., over $M$) with the same
boundary label and the following properties.
Every $\theta$-band of $\Delta$ shares at most one cell with any
$q$-band; $\Delta$ has no $q$-annuli.
\end{lm}

\proof Proving by contradiction, we consider a counter-example
$\Delta$ with minimal number of cells.

Assume that  a $\theta$-band
$\mathcal T$ and a $q$-band $\mathcal Q$  cross each other two
times. By the minimality of the counter-example,
these bands have exactly two common cells $\pi$ and $\pi'$, and $\Delta$ has no cells outside the region bounded by $\mathcal T$ and $\mathcal Q$. Then $\mathcal Q$ has exactly two cells since otherwise a
maximal $\theta$-band starting with a cell $\pi''$ of $\mathcal Q$,
where $\pi''\notin \{\pi, \pi'\}$, has to end on $\mathcal Q$,
bounding with a part of $\mathcal Q$ a smaller counter-example.
(We use that a $\theta$-band cannot end on an $a$-cell.)
For the similar reason, $\mathcal T$ has no $(\theta,q)$-cells except for $\pi$ and $\pi'$.
Therefore by Lemma \ref{tba} (2), these two
cells have the same labels of $\theta$-edges, so these two neighbor in $\cal Q$ cells are just mirror copies
of each other, since the label of a $q$-edge together with
the label of a $\theta$-edge completely determine the boundary
label of a $(\theta,q)$-cell. Thus, the diagram is not reduced, a contradiction.

If $\Delta$ has a $q$-annulus $\mathcal Q$, then the boundary $\partial\Delta$ is the outer boundary component of $\mathcal Q$, and there must be a $\theta$-band starting and ending on
$\partial\Delta$. It crosses $\mathcal Q$ two times, a contradiction.

Assume that $\Delta$ has a $\theta$-annulus $\cal T$.
Then $\cal T$ has no $(\theta,q)$-cells since otherwise
some $q$-band crosses the annulus $\cal T$ twice.
Notice that every $(\theta,a)$-cell corresponds to a
commutator relation (\ref{rel1}), and so the inner label
and the outer label of $\cal T$ are equal. Hence one
can identify the two boundary components of $\cal T$ and
remove all the $(\theta,a)$-cells of $\cal T$. Such a surgery
preserves the boundary label of $\Delta$ but decreases the
number of $(\theta,a)$-cells. So after several
surgeries of this type, we obtain the required minimal diagram $\Delta'$.

Assume now that a diagram $\Delta$ over $H(\bf Y)$ has a non-regular $\theta$-band $\cal T$. Then by Lemma \ref{tba}, there is a subband
$\cal T'$, whose the first and the last cells $\pi$ and $\pi'$ have mutual inverse boundary
labels, they are connecting with a subband
$\cal T''$ having no $(\theta,q)$-cells, and
the top/bottom labels of $\cal T''$ are congruent to $1$. Then $\cal T''$ can be replaced with a subdiagram with $a$-cells only,
so that $\pi$ and $\pi'$ get adjacent $\theta$-edges. Hence this pair of cells becomes cancellable. Decreasing the number
of $(\theta,q)$-cells, we obtain the desired
diagram $\Delta'$.
\endproof

\begin{cy} \label{GY} The homomorphism $G\to H(\bf Y)$ given by the mapping $a\mapsto a$ is injective.
\end{cy}

\proof  One should prove that an equality $w=1$ in $H(\bf Y)$
for a word $w$ in the alphabet $A^{\pm 1}$ implies the same
equality in $G$.

By Lemma \ref{NoAnnul}, there exists a minimal $H(\bf Y)$-diagram $\Delta$ with
boundary label $w$. Then neither $\theta$-band nor $q$-band can
start on $\partial\Delta$ since $w$ has neither $\theta$- nor $q$-edges. It follows from Lemma \ref{NoAnnul} that $\Delta$
has neither $(\theta,q)$- nor $(\theta,a)$-cells, and so the
equality $w=1$ follows from $a$-relations only.

Note that the replacement of letters from ${\bf Y}\backslash A$
with $1$ in an $a$-relation provides us with a $a$-relation again. So the equality $w=1$ is a consequence of
$a$-relations depending on letters from $A$ only. These
$a$-relations hold in $G$ since $G$ is a group of exponent $n$. Hence we get $w=1$ in $G$.
\endproof

\begin{df}\label{dftrap}\index[g]{trapezium} \index[g]{quasi-trapezium} 
 Let $\Delta$ be a reduced  diagram over $M$ (resp., a minimal diagram over $H(\bf Y)$),
which has  boundary path of the form ${\bf p}_1\iv {\bf q}_1{\bf p}_2{\bf q}_2\iv,$ where
${\bf p}_1$ and ${\bf p}_2$ are sides of $q$-bands, and
${\bf q}_1$, ${\bf q}_2$ are maximal parts of the sides of
$\theta$-bands such that $Lab({\bf q}_1)$, $Lab({\bf q}_2)$ start and end
with $q$-letters. 

\begin{figure}[h!]
%TeXCAD (http://texcad.sf.net/) Picture. File: [pic9.pic]. Options on following lines.
%\grade{\on}
%\emlines{\off}
%\epic{\off}
%\beziermacro{\on}
%\reduce{\on}
%\snapping{\off}
%\pvinsert{% Your \input, \def, etc. here}
%\quality{8.000}
%\graddiff{0.005}
%\snapasp{1}
%\zoom{4.0000}
\unitlength 1mm % = 2.845pt
\linethickness{0.4pt}
\ifx\plotpoint\undefined\newsavebox{\plotpoint}\fi % GNUPLOT compatibility
\begin{picture}(148.25,40)(0,115)
\put(76.5,148){\line(1,0){64.25}}
\put(74.75,143.75){\line(1,0){66.25}}
\put(72.75,140){\line(1,0){68.25}}
\put(71,136.5){\line(1,0){69.75}}
\put(69.5,133.25){\line(1,0){71.75}}
\put(68,130.25){\line(1,0){73.25}}
\put(66.75,127.25){\line(1,0){74.5}}
\put(65.5,124){\line(1,0){75.75}}
\put(76.75,147.75){\line(0,-1){23.5}}
\put(74.25,143.75){\line(0,-1){19.75}}
\put(72.25,140){\line(0,-1){16}}
\put(70,133.5){\line(0,-1){9.5}}
\put(67.5,130.5){\line(0,-1){6.5}}
\put(73.5,148){\line(1,0){3.25}}
%\emline(73.5,148.25)(71.75,147.25)
\multiput(73.5,148.25)(-.0583333,-.0333333){30}{\line(-1,0){.0583333}}
%\end
\put(72,147.25){\line(0,-1){2.75}}
\put(72,144.5){\line(0,1){.25}}
\put(72,144.75){\line(0,-1){.75}}
\put(72,144){\line(1,0){2.25}}
%\emline(72,144.25)(68.5,139.75)
\multiput(72,144.25)(-.033653846,-.043269231){104}{\line(0,-1){.043269231}}
%\end
\put(68.5,139.75){\line(1,0){1.75}}
%\emline(70,140)(67.75,136.5)
\multiput(70,140)(-.03358209,-.05223881){67}{\line(0,-1){.05223881}}
%\end
\put(67.75,136.5){\line(1,0){3.25}}
%\emline(69.5,139.75)(72,140)
\multiput(69.5,139.75)(.3125,.03125){8}{\line(1,0){.3125}}
%\end
\put(67.75,136.5){\line(0,-1){3.25}}
\put(67.75,133.25){\line(1,0){2}}
\put(65.75,133.25){\line(1,0){2}}
%\emline(65.5,133.25)(64,130.25)
\multiput(65.5,133.25)(-.03333333,-.06666667){45}{\line(0,-1){.06666667}}
%\end
\put(64,130.25){\line(1,0){3.25}}
\put(64,130.5){\line(0,-1){6.5}}
\put(64,124){\line(0,1){0}}
\put(64,124){\line(1,0){1.75}}
\put(79.25,147.75){\line(0,-1){23.75}}
%\dottedline(73.75,147.75)(72.25,145.75)
\multiput(73.68,147.68)(-.375,-.5){5}{{\rule{.4pt}{.4pt}}}
%\end
%\dottedline(75,148)(70.5,140.5)
\multiput(74.93,147.93)(-.45,-.75){11}{{\rule{.4pt}{.4pt}}}
%\end
%\dottedline(76.25,147.75)(67.75,134.5)
\multiput(76.18,147.68)(-.5,-.77941){18}{{\rule{.4pt}{.4pt}}}
%\end
%\dottedline(76.5,146)(75,144.25)
\multiput(76.43,145.93)(-.5,-.5833){4}{{\rule{.4pt}{.4pt}}}
%\end
%\dottedline(73.75,141.5)(64,128)
\multiput(73.68,141.43)(-.54167,-.75){19}{{\rule{.4pt}{.4pt}}}
%\end
%\dottedline(66.75,133.25)(64.5,130.75)
\multiput(66.68,133.18)(-.45,-.5){6}{{\rule{.4pt}{.4pt}}}
%\end
%\dottedline(72,137.25)(64.25,126.75)
\multiput(71.93,137.18)(-.55357,-.75){15}{{\rule{.4pt}{.4pt}}}
%\end
%\dottedline(72,135.75)(70.5,133.5)
\multiput(71.93,135.68)(-.375,-.5625){5}{{\rule{.4pt}{.4pt}}}
%\end
%\dottedline(69.75,132.25)(68.5,130.5)
\multiput(69.68,132.18)(-.4167,-.5833){4}{{\rule{.4pt}{.4pt}}}
%\end
%\dottedline(67.5,129)(64.25,125.5)
\multiput(67.43,128.93)(-.54167,-.58333){7}{{\rule{.4pt}{.4pt}}}
%\end
%\dottedline(67,127.25)(64.5,124.5)
\multiput(66.93,127.18)(-.5,-.55){6}{{\rule{.4pt}{.4pt}}}
%\end
%\dottedline(67.25,125.75)(65.5,124)
\multiput(67.18,125.68)(-.5833,-.5833){4}{{\rule{.4pt}{.4pt}}}
%\end
\put(95,148){\line(0,-1){24.25}}
\put(98.25,148){\line(0,-1){24.25}}
%\dottedline(97,148)(95.25,145.5)
\multiput(96.93,147.93)(-.4375,-.625){5}{{\rule{.4pt}{.4pt}}}
%\end
%\dottedline(98,147)(95.25,143)
\multiput(97.93,146.93)(-.45833,-.66667){7}{{\rule{.4pt}{.4pt}}}
%\end
%\dottedline(97.75,144.25)(95.25,140.75)
\multiput(97.68,144.18)(-.5,-.7){6}{{\rule{.4pt}{.4pt}}}
%\end
%\dottedline(98,142.5)(95,138.5)
\multiput(97.93,142.43)(-.42857,-.57143){8}{{\rule{.4pt}{.4pt}}}
%\end
%\dottedline(97.75,140.5)(95.25,137)
\multiput(97.68,140.43)(-.5,-.7){6}{{\rule{.4pt}{.4pt}}}
%\end
%\dottedline(98,138.5)(95.25,134.75)
\multiput(97.93,138.43)(-.55,-.75){6}{{\rule{.4pt}{.4pt}}}
%\end
%\dottedline(97.75,136.25)(95.25,132.75)
\multiput(97.68,136.18)(-.5,-.7){6}{{\rule{.4pt}{.4pt}}}
%\end
%\dottedline(97.75,134)(95.25,130.25)
\multiput(97.68,133.93)(-.5,-.75){6}{{\rule{.4pt}{.4pt}}}
%\end
%\dottedline(97.75,131.25)(95.25,128)
\multiput(97.68,131.18)(-.5,-.65){6}{{\rule{.4pt}{.4pt}}}
%\end
%\dottedline(95.25,128)(95.25,128)
\multiput(95.18,127.93)(0,0){3}{{\rule{.4pt}{.4pt}}}
%\end
%\dottedline(95.25,128)(94.75,128)
\multiput(95.18,127.93)(-.25,0){3}{{\rule{.4pt}{.4pt}}}
%\end
%\dottedline(98,129.25)(95,126.25)
\multiput(97.93,129.18)(-.5,-.5){7}{{\rule{.4pt}{.4pt}}}
%\end
%\dottedline(98.25,127)(95.25,124.5)
\multiput(98.18,126.93)(-.6,-.5){6}{{\rule{.4pt}{.4pt}}}
%\end
%\dottedline(98,125.25)(97,124.25)
\multiput(97.93,125.18)(-.3333,-.3333){4}{{\rule{.4pt}{.4pt}}}
%\end
\put(140.5,148.25){\line(0,-1){8.5}}
\put(140.5,139.75){\line(0,1){0}}
%\emline(140.5,140)(139.75,136.5)
\multiput(140.5,140)(-.0326087,-.1521739){23}{\line(0,-1){.1521739}}
%\end
\put(139.75,136.5){\line(0,1){.25}}
%\emline(139.75,136.75)(141.5,133.25)
\multiput(139.75,136.75)(.03365385,-.06730769){52}{\line(0,-1){.06730769}}
%\end
\put(141.5,133.25){\line(0,1){0}}
\put(141.5,133.25){\line(0,-1){9.25}}
\put(141.5,124){\line(-1,0){.25}}
\put(140.25,148){\line(1,0){4.75}}
%\emline(145.25,148)(144,143.5)
\multiput(145.25,148)(-.03289474,-.11842105){38}{\line(0,-1){.11842105}}
%\end
\put(144.25,144){\line(1,0){1.5}}
\put(145.75,144){\line(0,-1){8}}
\put(145.75,136.5){\line(-1,0){1.5}}
%\emline(144.5,136.25)(146.25,133)
\multiput(144.5,136.25)(.03365385,-.0625){52}{\line(0,-1){.0625}}
%\end
\put(146,133.75){\line(0,-1){10}}
\put(146,124.25){\line(1,0){.25}}
%\emline(140.75,124.25)(141.25,124)
\multiput(140.75,124.25)(.0625,-.03125){8}{\line(1,0){.0625}}
%\end
\put(141.25,124.25){\line(1,0){4.75}}
%\dottedline(143,148)(140.75,145.5)
\multiput(142.93,147.93)(-.45,-.5){6}{{\rule{.4pt}{.4pt}}}
%\end
%\dottedline(144.5,147.75)(141,143.5)
\multiput(144.43,147.68)(-.5,-.60714){8}{{\rule{.4pt}{.4pt}}}
%\end
%\dottedline(143.75,144.5)(140.75,140.75)
\multiput(143.68,144.43)(-.5,-.625){7}{{\rule{.4pt}{.4pt}}}
%\end
%\dottedline(145.25,143)(145.75,144)
\multiput(145.18,142.93)(.25,.5){3}{{\rule{.4pt}{.4pt}}}
%\end
%\dottedline(145,144)(145.25,144.5)
\multiput(144.93,143.93)(.125,.25){3}{{\rule{.4pt}{.4pt}}}
%\end
%\dottedline(145.25,144.5)(140.25,138)
\multiput(145.18,144.43)(-.5,-.65){11}{{\rule{.4pt}{.4pt}}}
%\end
%\dottedline(145.5,142.25)(140,135.75)
\multiput(145.43,142.18)(-.55,-.65){11}{{\rule{.4pt}{.4pt}}}
%\end
%\dottedline(145.25,140)(140.75,135)
\multiput(145.18,139.93)(-.5625,-.625){9}{{\rule{.4pt}{.4pt}}}
%\end
%\dottedline(140.75,135)(140.5,134.75)
\multiput(140.68,134.93)(-.125,-.125){3}{{\rule{.4pt}{.4pt}}}
%\end
%\dottedline(145.5,137.5)(141.5,133.25)
\multiput(145.43,137.43)(-.5,-.53125){9}{{\rule{.4pt}{.4pt}}}
%\end
%\dottedline(145,134.75)(141.5,131.25)
\multiput(144.93,134.68)(-.58333,-.58333){7}{{\rule{.4pt}{.4pt}}}
%\end
%\dottedline(145.75,133.25)(141.5,129.75)
\multiput(145.68,133.18)(-.60714,-.5){8}{{\rule{.4pt}{.4pt}}}
%\end
%\dottedline(145.75,131.75)(141.75,128.5)
\multiput(145.68,131.68)(-.57143,-.46429){8}{{\rule{.4pt}{.4pt}}}
%\end
%\dottedline(145.75,130)(141.5,126.25)
\multiput(145.68,129.93)(-.60714,-.53571){8}{{\rule{.4pt}{.4pt}}}
%\end
%\dottedline(145.75,127.75)(142,125.25)
\multiput(145.68,127.68)(-.75,-.5){6}{{\rule{.4pt}{.4pt}}}
%\end
%\dottedline(145.75,126.5)(143.5,124.75)
\multiput(145.68,126.43)(-.5625,-.4375){5}{{\rule{.4pt}{.4pt}}}
%\end
\put(140.25,143.75){\line(1,0){4.5}}
\put(140.25,140){\line(1,0){5.5}}
\put(140,136.5){\line(1,0){4.75}}
\put(141.5,133.25){\line(1,0){4.5}}
\put(141,130.5){\line(1,0){5.25}}
\put(141,127.25){\line(1,0){4.75}}
\put(137.75,148.25){\line(0,-1){24}}
\put(134.75,148){\line(0,-1){23.75}}
\put(131.75,147.75){\line(0,-1){23.5}}
\put(63,138.75){$p_1$}
\put(125.75,121.5){$q_1$}
\put(148.25,135.5){$p_2$}
\put(110.75,151){$q_2$}
\put(87.5,120){Trapezium}
%\emline(11,129.75)(49.25,150.75)
\multiput(11,129.75)(.0613964687,.0337078652){623}{\line(1,0){.0613964687}}
%\end
%\emline(49.25,150.75)(52.75,144.25)
\multiput(49.25,150.75)(.033653846,-.0625){104}{\line(0,-1){.0625}}
%\end
%\emline(52.75,144.25)(14,123.25)
\multiput(52.75,144.25)(-.0621990369,-.0337078652){623}{\line(-1,0){.0621990369}}
%\end
\put(11.25,129.5){\line(1,-2){3}}
%\emline(16,132.5)(18.75,126.5)
\multiput(16,132.5)(.03353659,-.07317073){82}{\line(0,-1){.07317073}}
%\end
%\emline(44.25,147.75)(47.25,142.25)
\multiput(44.25,147.75)(.03370787,-.06179775){89}{\line(0,-1){.06179775}}
%\end
\put(28.5,139){\line(1,-2){3.25}}
\put(31.75,132.5){\line(0,1){0}}
%\emline(32.75,142.25)(35.75,135.75)
\multiput(32.75,142.25)(.03370787,-.07303371){89}{\line(0,-1){.07303371}}
%\end
\put(15,128){$\pi_1$}
\put(47.5,145.5){$\pi_n$}
\put(31.25,137.75){$\pi_i$}
\put(10.25,125.75){$e_0$}
\put(27.5,135.5){$e_{i-1}$}
\put(35.5,140){$e_i$}
\put(41,135.25){$q_1$}
\put(17.25,137){$q_2$}
\put(28,127.5){Band}
\put(52.5,148.5){$e_n$}
\end{picture}

\begin{center}
\nopagebreak[4] Figure \theppp
\end{center}
\addtocounter{ppp}{1}	

\end{figure}

\begin{figure}
% This is a LaTeX picture output by TeXCAD.
% File name: [Gtr.pic].
% Version of TeXCAD: 4.3
% Reference / build: 30-Jun-2012 (rev. 105)
% For new versions, check: http://texcad.sf.net/
% Options on the following lines.
%\grade{\on}
%\emlines{\off}
%\epic{\off}
%\beziermacro{\on}
%\reduce{\on}
%\snapping{\off}
%\pvinsert{% Your \input, \def, etc. here}
%\quality{8.000}
%\graddiff{0.005}
%\snapasp{1}
%\zoom{4.0000}
\unitlength 1mm % = 2.845pt
\linethickness{0.4pt}
\ifx\plotpoint\undefined\newsavebox{\plotpoint}\fi % GNUPLOT compatibility
\begin{picture}(128,27.75)(0,0)
%\emline(25,2.25)(128,2.75)
\multiput(25,2.25)(6.8666667,.0333333){15}{\line(1,0){6.8666667}}
%\end
\put(25.5,7){\line(1,0){13.75}}
\put(53,7){\line(1,0){34.5}}
\put(113.25,7){\line(1,0){14.25}}
\put(25.5,11.75){\line(1,0){28}}
\put(63.75,11.5){\line(1,0){14.5}}
\put(90.75,11.75){\line(1,0){36.5}}
\put(30,16.5){\line(1,0){11.5}}
\put(54,16.5){\line(1,0){44.25}}
\put(114.5,16.25){\line(1,0){7.5}}
\put(33.25,20.75){\line(1,0){4.75}}
\put(68.5,21){\line(1,0){55}}
\put(33.5,24.5){\line(1,0){52}}
\put(110.25,25){\line(1,0){14.5}}
\put(53.5,7.5){\circle{.707}}
\put(46.625,7.625){\oval(14.25,1.75)[]}
\put(100.375,7.5){\oval(26.25,2)[]}
\put(58.75,11.875){\oval(11,2.25)[]}
\put(84.625,12){\oval(12.25,1.5)[]}
\put(54,16.375){\oval(.5,.75)[]}
\put(47.5,16.75){\oval(13,2.5)[]}
\put(53.25,21.25){\oval(30.5,2)[]}
\put(106.25,16.625){\oval(17,2.75)[]}
\put(98.125,24.625){\oval(25.25,2.25)[]}
\put(68.5,24.5){\line(0,-1){22}}
\put(71.25,24.5){\line(0,-1){22}}
\put(34,24.5){\line(0,-1){3.5}}
\put(34,21.25){\line(-1,-1){4.5}}
%\emline(29.75,16.75)(25.75,12)
\multiput(29.75,16.75)(-.033613445,-.039915966){119}{\line(0,-1){.039915966}}
%\end
\put(25.75,12){\line(0,-1){10}}
\put(37.75,24.5){\line(0,-1){7.75}}
\put(35,16.25){\line(0,-1){5}}
\put(30.5,11.75){\line(-1,0){.5}}
\put(30,11.5){\line(0,-1){9.25}}
%\emline(125,25)(123.25,20.5)
\multiput(125,25)(-.03365385,-.08653846){52}{\line(0,-1){.08653846}}
%\end
\put(121,25){\line(0,-1){3.75}}
%\emline(121,21)(118.5,16.5)
\multiput(121,21)(-.03333333,-.06){75}{\line(0,-1){.06}}
%\end
\put(123.75,21){\line(-1,-2){2.25}}
%\emline(121.75,15.5)(126.5,12)
\multiput(121.75,15.5)(.045673077,-.033653846){104}{\line(1,0){.045673077}}
%\end
\put(119,16.25){\line(0,-1){4.5}}
\put(126.75,11.75){\line(0,-1){8.25}}
\put(123.25,11.75){\line(0,-1){8.75}}
%\dottedline(34.5,23)(36.25,24.5)
\multiput(34.43,22.93)(.5833,.5){4}{{\rule{.4pt}{.4pt}}}
%\end
%\dottedline(34.5,21)(37.75,23.5)
\multiput(34.43,20.93)(.65,.5){6}{{\rule{.4pt}{.4pt}}}
%\end
%\dottedline(30.75,16.5)(37.75,22.5)
\multiput(30.68,16.43)(.63636,.54545){12}{{\rule{.4pt}{.4pt}}}
%\end
%\dottedline(31.5,16.5)(26.25,11.75)
\multiput(31.43,16.43)(-.65625,-.59375){9}{{\rule{.4pt}{.4pt}}}
%\end
%\dottedline(26.25,11.75)(26.25,11.75)
\multiput(26.18,11.68)(0,0){3}{{\rule{.4pt}{.4pt}}}
%\end
%\dottedline(26.25,10.75)(26.25,11)
\multiput(26.18,10.68)(0,.125){3}{{\rule{.4pt}{.4pt}}}
%\end
%\dottedline(25.75,9.75)(37.25,20.75)
\multiput(25.68,9.68)(.63889,.61111){19}{{\rule{.4pt}{.4pt}}}
%\end
%\dottedline(26.25,8.5)(37,19.25)
\multiput(26.18,8.43)(.67188,.67188){17}{{\rule{.4pt}{.4pt}}}
%\end
%\dottedline(25.75,6)(25.75,7)
\multiput(25.68,5.93)(0,.5){3}{{\rule{.4pt}{.4pt}}}
%\end
%\dottedline(26.75,6.75)(26.25,6.75)
\multiput(26.68,6.68)(-.25,0){3}{{\rule{.4pt}{.4pt}}}
%\end
%\dottedline(26.25,6.75)(29.75,10)
\multiput(26.18,6.68)(.58333,.54167){7}{{\rule{.4pt}{.4pt}}}
%\end
%\dottedline(31.5,12.5)(34.5,14.5)
\multiput(31.43,12.43)(.6,.4){6}{{\rule{.4pt}{.4pt}}}
%\end
%\dottedline(36,16.75)(37.75,17.5)
\multiput(35.93,16.68)(.875,.375){3}{{\rule{.4pt}{.4pt}}}
%\end
%\dottedline(25.5,5)(30.25,9.5)
\multiput(25.43,4.93)(.59375,.5625){9}{{\rule{.4pt}{.4pt}}}
%\end
%\dottedline(32.75,12)(34.75,13.5)
\multiput(32.68,11.93)(.5,.375){5}{{\rule{.4pt}{.4pt}}}
%\end
%\dottedline(26.25,3.75)(29.75,7.75)
\multiput(26.18,3.68)(.5,.57143){8}{{\rule{.4pt}{.4pt}}}
%\end
%\dottedline(34.5,12)(34.75,12.5)
\multiput(34.43,11.93)(.125,.25){3}{{\rule{.4pt}{.4pt}}}
%\end
%\dottedline(26.5,2.25)(29.5,6.25)
\multiput(26.43,2.18)(.42857,.57143){8}{{\rule{.4pt}{.4pt}}}
%\end
%\dottedline(28,2)(29.75,4.25)
\multiput(27.93,1.93)(.4375,.5625){5}{{\rule{.4pt}{.4pt}}}
%\end
%\dottedline(68.5,23)(71.5,24.75)
\multiput(68.43,22.93)(.6,.35){6}{{\rule{.4pt}{.4pt}}}
%\end
%\dottedline(69.25,20.75)(71.5,23.25)
\multiput(69.18,20.68)(.45,.5){6}{{\rule{.4pt}{.4pt}}}
%\end
%\dottedline(68.75,19)(70.75,21.5)
\multiput(68.68,18.93)(.4,.5){6}{{\rule{.4pt}{.4pt}}}
%\end
%\dottedline(68.25,17)(71,20.5)
\multiput(68.18,16.93)(.55,.7){6}{{\rule{.4pt}{.4pt}}}
%\end
%\dottedline(68.5,15.75)(71.25,19)
\multiput(68.43,15.68)(.55,.65){6}{{\rule{.4pt}{.4pt}}}
%\end
%\dottedline(68.75,15)(71.25,17.75)
\multiput(68.68,14.93)(.5,.55){6}{{\rule{.4pt}{.4pt}}}
%\end
%\dottedline(69,14)(71.25,15.75)
\multiput(68.93,13.93)(.5625,.4375){5}{{\rule{.4pt}{.4pt}}}
%\end
%\dottedline(69,12.75)(71,14.5)
\multiput(68.93,12.68)(.5,.4375){5}{{\rule{.4pt}{.4pt}}}
%\end
%\dottedline(68.75,11.5)(71.25,13.5)
\multiput(68.68,11.43)(.5,.4){6}{{\rule{.4pt}{.4pt}}}
%\end
%\dottedline(68.5,9.75)(71,11.75)
\multiput(68.43,9.68)(.5,.4){6}{{\rule{.4pt}{.4pt}}}
%\end
%\dottedline(68.5,8.75)(71.5,11)
\multiput(68.43,8.68)(.6,.45){6}{{\rule{.4pt}{.4pt}}}
%\end
%\dottedline(68.5,6.5)(71.25,8.75)
\multiput(68.43,6.43)(.55,.45){6}{{\rule{.4pt}{.4pt}}}
%\end
%\dottedline(69,7.75)(71.25,9.5)
\multiput(68.93,7.68)(.5625,.4375){5}{{\rule{.4pt}{.4pt}}}
%\end
%\dottedline(68.5,5.25)(71.5,7.75)
\multiput(68.43,5.18)(.6,.5){6}{{\rule{.4pt}{.4pt}}}
%\end
%\dottedline(68.75,4.5)(71,6.25)
\multiput(68.68,4.43)(.5625,.4375){5}{{\rule{.4pt}{.4pt}}}
%\end
%\dottedline(68.5,2.25)(71,4.75)
\multiput(68.43,2.18)(.5,.5){6}{{\rule{.4pt}{.4pt}}}
%\end
%\dottedline(70.25,2.75)(71.5,3.75)
\multiput(70.18,2.68)(.4167,.3333){4}{{\rule{.4pt}{.4pt}}}
%\end
%\dottedline(121.25,23.5)(123.5,25.25)
\multiput(121.18,23.43)(.5625,.4375){5}{{\rule{.4pt}{.4pt}}}
%\end
%\dottedline(121.25,22)(125,25)
\multiput(121.18,21.93)(.625,.5){7}{{\rule{.4pt}{.4pt}}}
%\end
%\dottedline(120.5,20)(124.5,22.5)
\multiput(120.43,19.93)(.66667,.41667){7}{{\rule{.4pt}{.4pt}}}
%\end
%\dottedline(119.75,18)(119.5,18.25)
\multiput(119.68,17.93)(-.125,.125){3}{{\rule{.4pt}{.4pt}}}
%\end
%\dottedline(120.25,19.5)(123.5,21.5)
\multiput(120.18,19.43)(.65,.4){6}{{\rule{.4pt}{.4pt}}}
%\end
%\dottedline(119.75,18.25)(123,19.75)
\multiput(119.68,18.18)(.65,.3){6}{{\rule{.4pt}{.4pt}}}
%\end
%\dottedline(119.25,16.25)(120,16.5)
\multiput(119.18,16.18)(.375,.125){3}{{\rule{.4pt}{.4pt}}}
%\end
%\dottedline(119,16.5)(122.25,17.75)
\multiput(118.93,16.43)(.65,.25){6}{{\rule{.4pt}{.4pt}}}
%\end
%\dottedline(119,14.5)(121.5,16.5)
\multiput(118.93,14.43)(.5,.4){6}{{\rule{.4pt}{.4pt}}}
%\end
%\dottedline(119.25,13)(122.25,15.5)
\multiput(119.18,12.93)(.6,.5){6}{{\rule{.4pt}{.4pt}}}
%\end
%\dottedline(119.75,11.75)(123.25,14.5)
\multiput(119.68,11.68)(.7,.55){6}{{\rule{.4pt}{.4pt}}}
%\end
%\dottedline(122,11.75)(124.5,13.75)
\multiput(121.93,11.68)(.5,.4){6}{{\rule{.4pt}{.4pt}}}
%\end
%\dottedline(124.5,13.75)(124.5,13.75)
\multiput(124.43,13.68)(0,0){3}{{\rule{.4pt}{.4pt}}}
%\end
%\dottedline(124.5,13.75)(125.75,12)
\multiput(124.43,13.68)(.4167,-.5833){4}{{\rule{.4pt}{.4pt}}}
%\end
%\dottedline(123.5,11)(125.5,12.5)
\multiput(123.43,10.93)(.5,.375){5}{{\rule{.4pt}{.4pt}}}
%\end
%\dottedline(123.5,9.25)(127.25,12.25)
\multiput(123.43,9.18)(.625,.5){7}{{\rule{.4pt}{.4pt}}}
%\end
%\dottedline(123.5,7.5)(126.75,10.25)
\multiput(123.43,7.43)(.65,.55){6}{{\rule{.4pt}{.4pt}}}
%\end
%\dottedline(123.25,5.75)(126.75,8.5)
\multiput(123.18,5.68)(.7,.55){6}{{\rule{.4pt}{.4pt}}}
%\end
%\dottedline(123.25,4)(126.75,7)
\multiput(123.18,3.93)(.58333,.5){7}{{\rule{.4pt}{.4pt}}}
%\end
%\dottedline(124,3)(126.75,5.5)
\multiput(123.93,2.93)(.55,.5){6}{{\rule{.4pt}{.4pt}}}
%\end
%\dottedline(125.25,2.5)(127,4.5)
\multiput(125.18,2.43)(.4375,.5){5}{{\rule{.4pt}{.4pt}}}
%\end
%\emline(39.75,22.25)(39.5,22)
\multiput(39.75,22.25)(-.03125,-.03125){8}{\line(0,-1){.03125}}
%\end
\put(46.75,-2.5){Quasi-trapezium}
\put(24.5,16.75){${\bf p}_1$}
\put(58,27.75){${q}_2$}
\put(126.5,17.25){${\bf p}_2$}
\put(95.25,-1){${q}_1$}
\end{picture}

\begin{center}
\nopagebreak[4] Figure \theppp
\end{center}
\addtocounter{ppp}{1}	

\end{figure}

Then $\Delta$ is called a \label{trapez}{\em trapezium} (respectively, quasi-trapezium). The path ${\bf q}_1$ is
called the \index[g]{bottom of trapezium}{\em bottom}, the path ${\bf q}_2$ is called the \index[g]{top of trapezium}{\em top} of
the trapezium, the paths ${\bf p}_1$ and ${\bf p}_2$ are called the \index[g]{side of trapezium}{\em left
and right sides} of the (quasi-)trapezium. The history of the $q$-band
whose side is ${\bf p}_2$ is called the \index[g]{history of trapezium}{\em history} of the (quasi-)trapezium;
%the length of the history is called the \index[g]{height of trapezium}{\em height}  of the
%(quasi-)trapezium. 
The base of $Lab ({\bf q}_1)$ is called the \index[g]{base of trapezium}{\em base} of the
(quasi-)trapezium.
\end{df}

\begin{remark} Notice that the top (bottom) side of a
$\theta$-band $\ttt$ does not necessarily coincide with the top
(bottom) side ${\bf q}_2$ (side ${\bf q}_1$) of the corresponding trapezium of height $1$, and ${\bf q}_2$
(${\bf q}_1$) is
obtained from $\topp(\ttt)$ (resp. $\bott(\ttt)$) by trimming the
first and the last $a$-edges
if these paths start and/or end with $a$-edges.
\end{remark}

By Lemma \ref{NoAnnul}, any maximal  $\theta$-band of a
(quasi-)trapezium $\Delta$ connects the left side and the right side of $\Delta$. So one can enumerate them from the bottom to the top of $\Delta$: $\ttt_1,...,\ttt_h$. The following lemma claims that every trapezium (quasi-trapezium)
simulates a computation of $\bf M$ (resp., quasi-computation).
(Similar statements can be found in \cite{OS03}. For the formulations  (1) and (3) below, it is important
that $\bf M$  is an $S$-machine. The analog of
this statement is false for Turing machines. - See \cite{OS01} for a discussion.)

\begin{lm}\label{simul} (1) Let $\Delta$ be a trapezium
with history $\theta_1\dots\theta_d$ ($d\ge 1$).
Assume that $\Delta$
has consecutive maximal $\theta$-bands  ${\cal T}_1,\dots
{\cal T}_d$, and the words
$U_j$
and $V_j$
are the  trimmed bottom and the
trimmed top labels of ${\cal T}_j,$ ($j=1,\dots,d$).
Then the history of $\Delta$ is a reduced word, $U_j$, $V_j$ are admissible
words for $M,$ and
$$V_1= U_1\cdot \theta_1,\; U_2\equiv V_1,\dots, U_d \equiv V_{d-1},\; V_d= U_d\cdot \theta_d.$$

(2) For every reduced computation $U\to\dots\to U\cdot h \equiv V$ of $\bf M$
with $|h|\ge 1$
there exists
a trapezium $\Delta$ with bottom label $U$, top label $V$, and with history $h.$

(3) Let $\Delta$ be a quasi-trapezium
with history $\theta_1\dots\theta_d$ ($d\ge 1$)
and consecutive maximal $\theta$-bands  ${\cal T}_1,\dots
{\cal T}_d$. Let the words
$U_j$
and $V_j$
are the  trimmed bottom and the
trimmed top labels of regular $\theta$-bands ${\cal T}_j,$ ($j=1,\dots,d$).
Then the history of $\Delta$ is a reduced word, $U_j$, $V_j$ are regular admissible
words for $M,$ and
$$V_1= U_1\cdot \theta_1, \;U_2\cong V_1,\dots, U_d \cong V_{d-1},\; V_d = U_d\cdot \theta_d.$$

(4) For every reduced quasi-computation $U_0\to U_1\to\dots\to U_d$ with history $h=\theta(1)\dots \theta(d)$
where $d\ge 1$
there exists
a quasi-trapezium $\Delta$ with bottom label $U_0$, top label $U_d$, and with history $h.$
\end{lm}

\proof (1) The equalities $V_i= U_i\cdot \theta_i$ follow from Lemma \ref{perm} (1). The equalities $V_{i-1}\equiv U_i$ are true
because every cell of $\Delta$ is a $\theta$-cell, and so there
are no cells between the top of the band ${\mathcal T}_{i-1}$
and the bottom of ${\mathcal T}_{i-1}$.

(2) Given a reduced computation ${\cal C}: U\equiv W_0\to\dots\to W_s\equiv V$, then by Lemma \ref{perm} (2) one can construct a reduced band ${\mathcal T}_{i}$ for every subcomputation $W_{i-1}\to W_i$.
Since the trimmed top label $W_i$ of ${\mathcal T}_{i}$ coincides with the trimmed bottom label of ${\mathcal T}_{i+1}$, the obtained bands can be glue up. The obtained diagram $\Delta$ is reduced since the history of $\cal C$ is reduced, and so
the cells from different $\theta$-band cannot form reducible pairs. Clearly, the history of $\cal C$ is equal to the history
of the constructed trapezium.

(3) We have regular $\theta$-bands ${\mathcal T}_{i}$-s corresponding to the letters $\theta_i$ of the history and the equalities $V_i= U_i\cdot \theta_i$ as in the proof of Statement (1), but now there
are $a$-{\it cells} corresponding to $a$-relations
between the neighbor bands ${\mathcal T}_{i}$ and ${\mathcal T}_{i+1}$. Since the $q$-edges of the $(\theta,q)$-cells from
the neighbor $\theta$-bands belong to the same $q$-bands,
we have the equalities of the form $v_{ij}=u_{i+1,j}$ modulo the
$a$-relations
for the corresponding sector words of the words $V_i$ and $U_{i+1}$. However $v_{ij}$ and $u_{i+1,j}$ are the words in the same subalphabet ${\bf Y}_k\subset \bf Y$ since the words
$V_i$ and $U_{i+1}$ are admissible by Lemma \ref{perm} (1). They are regular since every $\theta$-band of a quasi-trapesium is regular.
Hence the homomorphism mapping all the letters from $\bf Y\backslash {\bf Y}_k$ to $1$ shows that the words $v_{ij}$ and $u_{i+1,j}$ are equal modulo the Burnside relations of the
group $B({\bf Y}_k)$ (and modulo the relations of $G$ if the sector is the input one for $\bf M$). Hence $U_{i+1}\cong V_i$,
as required.

(4) As in the proof of Statement (2), one can  start constructing  $\Delta$ with the subsequent regular $\theta$-bands ${\mathcal T}_1, \dots$ identifying the $q$-edges of the neighbor bands. The congruence $W_{i}\cong
W_{i-1}\cdot\theta_i$ implies that for the corresponding
sector subwords $v_{i-1,j}$ and $u_{ij}$ of these two words, we have equalities modulo
the Burnside relations of $B({\bf Y}_k)$, $k=k(j)$ (or modulo $G$-relations). It remains to glue up each hole having boundary label $v_{i-1,j}u_{ij}^{-1}$ with the corresponding $a$-cell.
\endproof

\section{More Burnside relations}\label{more}

\subsection{Burnside relations without $q$-letters}

In this subsection, we study presentations of some auxiliary groups. Let us start with the group \index[g]{$H(a)$} $H(a)$ generated by the set $\cup_{j=1}^N\bf\Theta^+_j\cup\bf Y$ and defined by all $(\theta,a)$-relations and all $a$-relations.
 It follows from these
relations that
there is a retraction of $H(a)$ onto the free
group $F(\cup_j\bf\Theta^+_j)$.
%and so every element $g$ has a projection $\bar g$ in $F(\cup_j\bf\Theta^+_j)$. 
Furthermore,
$H(a)$ is the multiple HNN-extension of the
subgroup $K(a)$ generated by $\bf Y$ with stable letters from $\cup_j\bf\Theta^+_j$. If $\theta_i$ belongs
to some $\bf\Theta^+_i$, then the subgroup  $\tilde{\bf Y}_i(\theta)$ of $K(a)$ associated with $\theta_i$ is the image of the subgroup ${\bf Y}_i(\theta)$ in the free Burnside group $B({\bf Y}_i)$ (in $G$ if $i=1$) and the conjugation by
$\theta_i$ is identical on $\tilde{\bf Y}_i(\theta)$.

\begin{lm} \label{man} Every subgroup $\tilde{\bf Y}_i(\theta)$ is malnormal in the group $K(a)$.
\end{lm}
\proof The group $K(a)$ is generated by the set ${\bf Y}$ of all $a$-letters. The $a$-relations can be imposed 
as follows. At first we impose the relations on the words in the subalphabet ${\bf Y}_1=A$
and impose the Burnside relations on the words in the subalphabet ${\bf Y}_j$ fro every
$j=2,\dots, N$. We obtain the free product $P$ of $G$ and free Burnside groups $B({\bf Y}_j)$, $j=2,\dots, N$.
Then imposing all the Burnside relations on $P$ we obtain the group $K(a)\cong P/P^n$, which is called
the free product of the groups $G$ and $B({\bf Y}_j)$ ($j=2,\dots, N$) in the Burnside variety ${\cal B}_n$. A few properties of free multiplication in ${\cal B}_n$ can be found in Subsection 36.1 of \cite{book}. 
The operation $o^{n}_{i\in J} G_j$ studied there coincides with the free product operation in ${\cal B}_n$
if the factors $G_j$-s belong to ${\cal B}_n$. In particular, we have $G_j\cap xG_jx^{-1}=\{1\}$ if
$x\notin G_j$, i.e. the factors $G_j$ of a free product in ${\cal B}_n$ are malnormal subgroups. In particular
we have that the subgroups $G$ and $B({\bf Y}_j)$ ($j=2,\dots, N$) are malnormal in $K(a)$.

It follows from the definition of the machine $\bf M$, that every subgroup ${\bf Y}_j(\theta)$ is either
equal to the free group $F({\bf Y}_j)$ or trivial, or is the copy of some $H(ij)$ (we have no latter variant
if $j=1$, i.e. if ${\bf Y}_j(\theta)=G$). Therefore its image
$\tilde{\bf Y}_i(\theta)$ is either obviously malnormal in $G$ or in $B({\bf Y}_j)$ or it is malnormal
in one of these groups by Lemma \ref{maln}. Now the statement of the lemma follows from the previous paragraph since the malnormality is a transitive property.
\endproof

We will verify the conditions (Z1), (Z2), and (Z3)
from Section \ref{axioms} for the group $H(a)$ under the assumption
that $a$-letters are zero-letters and $\theta$-letters are non-zero
letters. So one can talk of reduced, cyclically reduced and essential words as they were defined in Section \ref{axioms}. Now ${\cal R}_0={\cal S}_0$ is the union 
of the set $(\theta,a)$-relations and the set of $a$-relations; the set ${\cal S}_{1/2}$ is empty.

\begin{lm} \label{HaZ} The presentation of $H(a)$ satisfies the conditions (Z1),
(Z2) and (Z3).
\end{lm}
\proof Since the set ${\cal S}_{1/2}$ is empty, Condition (Z2) holds. Condition (Z1) follows from the definition of  $H(a)$.
It remains to verify Condition (Z3).

Assume that $g$ is an essential  cyclically reduced word and for some $x\in K(a)\backslash\{1\}$,
$g^{-4}xg^4$ is also a $0$-elements, that is $g^{-4}xg^4\in K(a)$. Then there are no
pinches in $g^4$ regarded as a product in the HNN-extension $H(a)$, and every time one obtains 
only one pinch
in the middle when reducing the product $g^{-4}xg^4$ to the normal form
in the HNN-extension. Hence the product $g^{-1}xg=x'$ is a $0$-element too.

Let $g=y_0\theta(1)y_1\dots\theta(k)y_k$ be the normal form of  $g$ in the HNN extension $H(a)$, where $\theta(1),\dots, \theta(k)$
are stable letters or their inverses ($k\ge 1$). Then the
element $y_0^{-1}xy_0$ must belong to the subgroup $\tilde{\bf Y}_i(\theta)$ for some $i$ and $\theta$, and so $\theta(1)$ commutes
with $y_0xy_0^{-1}$ in $H(a)$. Similarly we obtain that $\theta(2)$
commutes with $y_1^{-1}y_0^{-1}xy_0y_1$, and so on,
and $x'=y^{-1}xy$, where $y=y_0y_1\dots y_k$.

We also have the next pinch $\theta(1)^{-1}y_0^{-1}x'y_0\theta(1)$,
and therefore the product $$y_0^{-1}x'y_0=
y_0^{-1}y^{-1}xyy_0=
(y_0^{-1}y^{-1}y_0)(y_0^{-1}xy_0)(y_0^{-1}yy_0)$$
is also in $\tilde{\bf Y}_i( \theta)$.
By Lemma \ref{man}, we have $y_0^{-1}yy_0=y_1\dots y_ky_0\in \tilde{\bf Y}_i(\theta)$ as well. Hence
$\theta(1)$ commutes with $y(1)^{-s}(y_0^{-1}xy_0)y(1)^s$
for every $s$, where $y(1) = y_1\dots y_ky_0$.

Changing $g$ by a cyclic permutations, we similarly
obtain that arbitrary $\theta(m)$ commutes with
$y(m)^{-s}(y_0y_1\dots y_{m-1})^{-1}x(y_0y_1\dots y_{m-1})y(m)^s$
for every $s$, where $y(m) = y_m\dots y_ky_0\dots y_{m-1}$. Using these commutativity
relations $k$ times and taking into account that every $y(m)$ is a cyclic permutation of the product $y$,
we have $$g^{-1}y^{1-s}xy^{s-1}g=(y_0\theta(1)y_1\dots\theta(k)y_k)^{-1}y^{1-s}xy^{s-1}(y_0\theta(1)y_1\dots\theta(k)y_k)=
$$ $$ (y_1\theta(2)\dots\theta(k)y_k)^{-1}y(1)^{1-s}y_0^{-1}xy_0y(1)^{s-1} (y_1\theta(2)\dots\theta(k)y_k)=\dots =y^{-s}xy^s,$$
because $y_0\dots y_{k-1}y(k)^{s-1}y_k= y^s$.
Hence by induction on $s$, we obtain equalities
$g^{-s}xg^s=
y^{-s}xy^s$ for every $s\ge 0$.
Since $y\in K(0)$, Conditions $(Z3.1)$ holds. Condition $(Z3.2)$ follows
as well since the definition of $y$ does not depend on the choice of $x$ in ${\bf 0}(g)$. \endproof

Lemma \ref{HaZ} and Lemma \ref{mainprop}  provide us with

\begin{cy} \label{Hainf} There exists a graded presentation of the factor group $H(a,\infty)$
of $H(a)$ over the subgroup generated by all $n$-th powers of elements of $H(a)$,
such that every g-reduced diagram over this presentation satisfies Property A
from Section \ref{bcs}.
\end{cy} $\Box$

We denote by \index[g]{$H(\theta,a)$} $H(\theta,a)$ the group
$H(a,\infty)$ from Corollary \ref{Hainf}, that is the group
generated by $\cup_j\bf\Theta_j^+\cup\bf Y$, which is subject to all Burnside relations $w^n=1$ in this alphabet (and so $H(\theta,a)\in {\cal B}_n$), all relations of $G$ in the subalphabet $A$, and all $(\theta,a)$-relations.

\begin{lm} \label{retr} The mapping identical on the set of generators $A$
of the group $G$ and trivial on other generators from $\cup\bf\Theta_j^+\cup\bf Y$ is a retraction of the group $H(\theta,a)$
onto the subgroup $G$.
\end{lm}
\proof Indeed, this mapping preserves all defining relations
of $H(\theta,a)$ since $G^n=\{1\}$.\endproof

\begin{lm} \label{thea} Let $\theta(1),\dots,\theta(m)$ be different rules from $\Theta^+$ and
$\theta_j(1),\dots,\theta_j(m)$ be their copies from some set $\bf\Theta_j^+$. Assume that
$x(i)= c(i)\theta_j(i)d(i)$ for some group words $c_i$, $d_i$ over $\bf Y$ ($i=1,\dots m$). Then

(a) the words $x(i)$ freely generate
a free Burnside subgroup, which is the retract of $H(\theta,a)$,

(b) The mapping  $x(i)\mapsto \theta(i)$ ($i=1,\dots m$) defines an isomorphism between the free Burnside groups $B(x(1),\dots,x(m))$ and $B(\theta(1),\dots,\theta(m))$.
\end{lm}

\proof We preserve the set of defining relations of $H(\theta,a)$
if we replace all the generators, except for $\theta_j(1),\dots,\theta_j(m)$, with $1$. Hence there is a homomorphism of $H(\theta,a)$ onto the free Burnside group
$B(\theta_j(1),\dots,\theta_j(m))$ such that $x(i)\mapsto \theta_j(i)$.
The mapping $\theta_j(i)\mapsto x(i)$ ($i=1,\dots,m$) also defines
a homomorphism since all the relations $w^n=1$ hold in $H(\theta,a)$, which proves te lemma.
\endproof

\begin{cy} \label{tb} 
(a) A product of several letters $x(i)^{\pm 1}$ ($i=1,\dots,m$)
is minimal (cyclically minimal) word $w$ in $H(\theta,a)$
if and only if the image of $w$ under the mapping $x(i)\to \theta(i)$ ($i=1,\dots m$) is minimal (cyclically minimal) in $B(\bf\Theta)$.

(b) The top label of a $q$-band is trivial in the group
$H(\theta,a)$ if and only if the bottom label is trivial in $H(\theta,a)$.
\end{cy}

\proof (a) This follows from Lemma \ref{thea}, because $a$-letters have length $0$.

(b) The top label is freely equal to a product of the words having the form
$x(i)^{\pm 1}=(c(i)\theta(i)d(i))^{\pm 1}$ for $\theta\in \bf\Theta^+_j$, where $c(i), d(i)$ are words over $\bf Y$. The bottom label is a product of the words
$y(i)$ of the similar form with middle terms from $\bf\Theta^+_{j+1}$, or from $\bf\Theta^+_{j-1}$ and the corresponding factors of these
products have the same $\theta$ in the middle. It follows
from Lemma \ref{thea} (a) that the top label is trivial in $H(\theta,a)$ if and only if the bottom one is trivial.
\endproof

We say that a reduced (disk or annular) diagram over $H(a)$ is \index[g]{minimal diagram over $H(a)$} {\it minimal} if it contains no $\theta$-annuli. Similarly, a reduced diagram $\Delta$ over $H(\theta,a)$ is {\it minimal} if it is $g$-reduced and
has no $\theta$-annuli.

\begin{lm} \label{notha} (1) For every (disk or annular) diagram $\Delta$ over $H(a)$ (over $H(\theta,a)$), there is a minimal diagram $\Delta'$ over $H(a)$ (over $H(\theta,a)$) with the same
boundary label(s).

(2) If $\Delta$ is an annular diagram with
boundary components $\bf p$ and $\bf q$,  a path $\bf s$ connects the vertices ${\bf p}_-$
and ${\bf q}_-$ in $\Delta$, and the word $Lab({\bf p})$ is not a conjugate of a word
of $\theta$-length $0$ in $H(a)$ (resp., in $H(\theta,a)$), then there is a path $\bf s'$ in the diagram $\Delta'$,
such that $Lab({\bf s})= Lab({\bf s'})$ in $H(\theta,a)$.
\end{lm}

\proof In a disk diagram $\Delta$, $\theta$-bands can be eliminated
exactly as in the proof of Lemma \ref{NoAnnul}
since their top and bottom have the same label.
The same trick works for annular diagrams,
since if a $\theta$-annulus surrounds the hole
of $\Delta$, then we have an annular subdiagram
with boundary components $\bf p$ and the top (or the bottom) of this annulus. It follows that $Lab(\bf p)$ is conjugate to a word of
$\theta$-length $0$, a contradiction.

Note that the modifications of $\Delta$
(i.e. eliminating of $\theta$-band and $g$-reductions) are 'local', i.e. the path $\bf s$ can be replaced with a homotopic path which
does not cross the modified subdiagram. Thus,
Statement (2) follows.
\endproof

Let now \index[g]{$H(\bf\Theta)$} $H(\bf\Theta)$ be the factor group of $H(\bf Y)$
over the all relations $w^n=1$, where $w$ is a
word in generates $\bf\Theta\cup\bf Y$.
In other words, $H(\bf\Theta)$
is the factor group of the free product $H(\theta,a)\star F({\bf Q})$ over all $(\theta,q)$-relations.

We will call a reduced disk diagram $\Delta$ over $H(\bf\Theta)$
\index[g]{minimal diagram over $H(\bf\Theta)$} {\it minimal} if
it has no $q$-annuli and every subdiagram
of $\Delta$ over $H(\theta,a)$ or over $H(\bf Y)$ is minimal.

\begin{lm} \label{noq} For every disk diagram $\Delta$ over $H(\bf\Theta)$, there exists a minimal diagram $\Delta'$ with the same boundary label.
\end{lm}

\proof Assume $\Delta$ contains a $q$-annulus $\cal Q$. Proving by contradiction we may assume that the inner subdiagram $\Delta_1$  of $\cal Q$ has no $q$-annuli. Since the top/bottom of $\cal Q$ has no $q$-edges, $\Delta_1$ has no $(\theta,q)$-cells, i.e., it is a diagram over $H(\theta,a)$. So the label of the inner boundary component of $\cal Q$ has trivial in $H(\theta,a)$ label. By Corollary \ref{tb} (b), the label of the outer
boundary component of the annulus ${\cal Q}$ is also trivial in $H(\theta,a)$. So the annulus ${\cal Q}$ together with $\Delta_1$ can be replaced with a diagram without $(\theta,q)$-cell. One can continue decreasing
the number of $(\theta,q)$-cells and obtain the required
minimal diagram $\Delta'$ using Lemmas \ref{noq} and \ref{NoAnnul}.
\endproof

\begin{lm} \label{toHT} The canonical mapping $a\mapsto a$ defines an
embedding of the group $G$ in $H(\bf\Theta)$.
\end{lm}

\proof This mapping defines a homomorphism $G\to H(\bf\Theta)$ by Corollary \ref{GY}. Assume that a word $v$ in the generators $A$ of $G$ is trivial in $H(\bf\Theta)$ and consider a minimal diagram $\Delta$
for this equality.
By Lemma \ref{noq}, it has no $q$-annuli, and therefore $\Delta$ has no $(\theta,q)$-cells at all since the boundary $\partial\Delta$ has no $q$-edges.
Hence $v=1$ in the group $H(\theta,a)$. By lemma \ref{retr},
$v=1$ in $G$, which completes the proof. \endproof

\begin{lm}\label{tohist} The mapping sending every $\theta$-generator $\theta_j\in \bf\Theta^+_j$ of $H(\bf\Theta)$ (of $H(\bf Y)$) to the corresponding history letter $\theta\in \bf\Theta^+$
and every $q$- and $a$-generator of $H(\bf\Theta)$ (of $H(\bf Y)$)  to $1$, extends to a homomorphism of $H(\bf\Theta)$ (resp., of $H(\bf Y)$) to the
free Burnside group $B(\bf\Theta^+)$ (to the free group $F(\bf\Theta^+)$, resp.).

For every $j$, the restriction of this mapping to the subgroup generated by $\bf \Theta^+_j$ is an isomorphism.
\end{lm}

\proof The image of every relator of $H(\bf\Theta)$ under
this mapping is either $1$ or a word of the form $u^n$. Therefore the first statement follows from Dyck's lemma. For the restriction
to $H(\bf\Theta^+_j)$, there is an obvious
inverse homomorphism given by $\theta\mapsto \theta^{(j)}$, which proves the second statement
of the lemma for $H(\bf\Theta)$. The proof
for the group $H(\bf Y)$ is similar.
\endproof

\subsection{Diagrams with one maximal $q$-band}

\begin{lm} \label{through} (1) Let $\Delta$ be a $g$-reduced annular diagram
over $H(\bf\Theta)$ with boundary components ${\bf s}_1$ and ${\bf t}_1{\bf t}_2{\bf s}_2$, containing a unique maximal $q$-band $\cal Q$ with history $h$,
start/end edges ${\bf t}_1$ and ${\bf t}_2$ and boundary ${\bf t}_1{\bf y}_1{\bf t}_2{\bf y}_2$, where ${\bf s}_1$ is contained in the region bounded by the closed path
${\bf y}_1$ (\setcounter{pdeleven}{\value{ppp}} left part of Figure
\thepdeleven). Assume that
$h\ne 1$ in the free Burnside group $B(\Theta^+)$ and
the annular subdiagram $\Delta_1$ with boundary components  ${\bf s}_1$ and ${\bf y}_1$ is minimal over
$H(\theta,a)$ and contains a positive
cell $\pi$ with contiguity degree $\ge\varepsilon$ to ${\bf y}_1$ via a contiguity
subdiagram $\Gamma$ of rank $0$, where
${\bf y}_1$ is regarded as the path starting
and ending at the vertex $({\bf t}_1)_+$.

Then there exists a minimal annular diagram $\Delta'$
with boundary components ${\bf s}'_1$ and ${\bf t}'_1{\bf t}'_2{\bf s}'_2$, where $Lab ({\bf s}_1')\equiv Lab ({\bf s}_1)$,
$Lab ({\bf s}_2')\equiv Lab ({\bf s}_2)$, $Lab ({\bf t}_1')\equiv Lab ({\bf t}_1)$, $Lab ({\bf t}_2')\equiv Lab ({\bf t}_2)$ (\setcounter{pdeleven}{\value{ppp}} right part of Figure
\thepdeleven), such that

(a) $\Delta'$ has a unique maximal $q$-band $\cal Q'$ with boundary
${\bf t}'_1{\bf y}'_1{\bf t}'_2{\bf y}'_2$,

(b) the history $h'$ of $\cal Q'$ is
equal to $h$ modulo the Burnside relations,

(c) the type of the annular subdiagram $\Delta'_1$ over $H(\theta,a)$ with boundary components  ${\bf s}'_1$ and ${\bf y}'_1$
is less than $\tau(\Delta_1)$.

(d) If a simple path ${\bf z}$ connects the
vertices $({\bf s}_1)_-$ and $({\bf s}_2)_-$ in $\Delta_1$, then there exists a simple path ${\bf z}'$ connecting the vertices $({\bf s}'_1)_-$ and $({\bf s}'_2)_-$ in $\Delta'_1$
such that the words $Lab({\bf z}')$ and $Lab({\bf z})$ are equal in the group $H(\theta,a)$.

%(2) Similar statement with disk diagrams $\Delta$, $\Delta_1$ and $\Delta'$, where $\Delta$ has the boundary ${\bf t}_1{\bf s}_1{\bf t}_2{\bf s}_2$ %(\setcounter{pdeleven}{\value{ppp}}  Figure
%\thepdeleven).

%\end{lm}

%This is a LaTeX picture output by TeXCAD.
% File name: [thr.pic].
% Version of TeXCAD: 4.3
% Reference / build: 30-Jun-2012 (rev. 105)
% For new versions, check: http://texcad.sf.net/
% Options on the following lines.
%\grade{\on}
%\emlines{\off}
%\epic{\off}
%\beziermacro{\on}
%\reduce{\on}
%\snapping{\off}
%\pvinsert{% Your \input, \def, etc. here}
%\quality{8.000}
%\graddiff{0.005}
%\snapasp{1}
%\zoom{4.0000}
\unitlength 1mm % = 2.845pt
\linethickness{0.4pt}
\ifx\plotpoint\undefined\newsavebox{\plotpoint}\fi % GNUPLOT compatibility


\begin{center}
\nopagebreak[4] Figure \theppp
\end{center}
\addtocounter{ppp}{1}	

%\end{figure}

(2) Similar statement holds for disk diagrams $\Delta$, $\Delta_1$ and $\Delta'$ and a single maximal $q$-band $\cal Q$, where $\Delta$ has the boundary ${\bf t}_1{\bf s}_1{\bf t}_2{\bf s}_2$ and $\cal Q$ has boundary ${\bf t}_1{\bf y}_1{\bf t}_2{\bf y}_2$ (\setcounter{pdeleven}{\value{ppp}}  Figure
\thepdeleven).
\end{lm}

% This is a LaTeX picture output by TeXCAD.
% File name: [dthr.pic].
% Version of TeXCAD: 4.3
% Reference / build: 30-Jun-2012 (rev. 105)
% For new versions, check: http://texcad.sf.net/
% Options on the following lines.
%\grade{\on}
%\emlines{\off}
%\epic{\off}
%\beziermacro{\on}
%\reduce{\on}
%\snapping{\off}
%\pvinsert{% Your \input, \def, etc. here}
%\quality{8.000}
%\graddiff{0.005}
%\snapasp{1}
%\zoom{4.0000}
\unitlength 1mm % = 2.845pt
\linethickness{0.4pt}
\ifx\plotpoint\undefined\newsavebox{\plotpoint}\fi % GNUPLOT compatibility
\begin{picture}(161.25,55.5)(10,0)
\put(69.125,10.875){\oval(.25,.25)[]}
\put(157.875,12.125){\oval(.25,.25)[]}
\put(40.25,37.25){\oval(58.5,32.5)[]}
\put(129,38.5){\oval(58.5,32.5)[]}
\put(11.25,33.75){\line(1,0){58.25}}
\put(11,30){\line(1,0){58.75}}
\put(50.25,40.375){\oval(0,.25)[]}
\put(40.75,42.625){\oval(18.5,7.25)[]}
\put(34.5,39){\line(0,-1){8.75}}
\put(47.25,39){\line(0,-1){8.75}}
\put(22.75,33.75){\line(0,-1){3.5}}
\put(15.25,33.5){\line(0,-1){3.25}}
\put(18.75,33.25){\line(0,-1){2.75}}
\put(26.5,34){\line(0,-1){3.5}}
\put(30.5,33.75){\line(0,-1){3.25}}
\put(40.75,33.75){\line(0,-1){3.5}}
\put(44,33.5){\line(0,-1){3}}
\put(37.75,31.5){\line(0,-1){.75}}
\put(58.25,33.75){\line(0,-1){4}}
\put(55,33.75){\line(0,-1){3.5}}
\put(62,33.5){\line(0,-1){3}}
\put(66,33.75){\line(0,-1){4}}
%\emline(37.75,33.25)(38.5,33.5)
\multiput(37.75,33.25)(.09375,.03125){8}{\line(1,0){.09375}}
%\end
\put(38,34){\line(0,-1){4.25}}
\put(51,33.5){\line(0,-1){3.5}}
\put(99.75,33){\line(1,0){17}}
%\emline(116.75,33)(119,34.75)
\multiput(116.75,33)(.04326923,.03365385){52}{\line(1,0){.04326923}}
%\end
%\emline(119,34.75)(118.75,41.5)
\multiput(119,34.75)(-.03125,.84375){8}{\line(0,1){.84375}}
%\end
%\emline(119.25,41)(121,42.75)
\multiput(119.25,41)(.03365385,.03365385){52}{\line(0,1){.03365385}}
%\end
\put(121,42.75){\line(1,0){16.25}}
%\emline(137.25,42.75)(139,41)
\multiput(137.25,42.75)(.03365385,-.03365385){52}{\line(0,-1){.03365385}}
%\end
\put(139,41){\line(0,-1){5.25}}
\put(139,35.75){\line(0,1){0}}
%\emline(139,35.75)(140.25,33.25)
\multiput(139,35.75)(.03289474,-.06578947){38}{\line(0,-1){.06578947}}
%\end
\put(140.25,33.25){\line(1,0){17.75}}
\put(100,29.75){\line(1,0){17.25}}
%\emline(117.25,29.75)(121.25,33)
\multiput(117.25,29.75)(.041237113,.033505155){97}{\line(1,0){.041237113}}
%\end
\put(121.25,33){\line(0,1){6}}
%\emline(121.25,39)(122.75,39.75)
\multiput(121.25,39)(.0652174,.0326087){23}{\line(1,0){.0652174}}
%\end
\put(122.75,39.75){\line(1,0){12.25}}
%\emline(135,39.75)(136.25,38)
\multiput(135,39.75)(.03289474,-.04605263){38}{\line(0,-1){.04605263}}
%\end
%\emline(136.25,38)(136,34)
\multiput(136.25,38)(-.03125,-.5){8}{\line(0,-1){.5}}
%\end
%\emline(136,34)(138.5,30.5)
\multiput(136,34)(.03333333,-.04666667){75}{\line(0,-1){.04666667}}
%\end
%\emline(138.5,30.5)(137.75,30.25)
\multiput(138.5,30.5)(-.09375,-.03125){8}{\line(-1,0){.09375}}
%\end
\put(137.75,30.25){\line(1,0){1.75}}
\put(139.5,30.25){\line(1,0){18.25}}
\put(104,33){\line(0,-1){3.75}}
\put(108.25,33){\line(0,-1){3.25}}
\put(112.5,33.25){\line(0,-1){3.5}}
\put(116.75,32.75){\line(0,-1){2.75}}
%\emline(119.25,34.25)(121.25,32.75)
\multiput(119.25,34.25)(.04444444,-.03333333){45}{\line(1,0){.04444444}}
%\end
\put(119.25,37){\line(1,0){2}}
%\emline(119.25,41)(121.75,39.5)
\multiput(119.25,41)(.05555556,-.03333333){45}{\line(1,0){.05555556}}
%\end
\put(124.5,42){\line(0,-1){2.5}}
\put(128.75,42.25){\line(0,-1){2.75}}
%\emline(133,42.75)(132.75,39.75)
\multiput(133,42.75)(-.03125,-.375){8}{\line(0,-1){.375}}
%\end
%\emline(138,42)(135.25,39)
\multiput(138,42)(-.03353659,-.03658537){82}{\line(0,-1){.03658537}}
%\end
\put(135.75,37){\line(1,0){3.25}}
%\emline(139.5,35.25)(137,32.75)
\multiput(139.5,35.25)(-.03333333,-.03333333){75}{\line(0,-1){.03333333}}
%\end
%\emline(141.5,33.25)(141.25,30.25)
\multiput(141.5,33.25)(-.03125,-.375){8}{\line(0,-1){.375}}
%\end
%\emline(145.75,33.25)(145.5,30.5)
\multiput(145.75,33.25)(-.03125,-.34375){8}{\line(0,-1){.34375}}
%\end
%\emline(149.75,33.25)(150,30.5)
\multiput(149.75,33.25)(.03125,-.34375){8}{\line(0,-1){.34375}}
%\end
%\emline(153.75,33.25)(154,29.75)
\multiput(153.75,33.25)(.03125,-.4375){8}{\line(0,-1){.4375}}
%\end
\put(78,40){\line(1,0){9.75}}
\put(78.25,37){\line(1,0){9}}
%\emline(86,41.75)(89,38.25)
\multiput(86,41.75)(.03370787,-.03932584){89}{\line(0,-1){.03932584}}
%\end
%\emline(89,38.5)(86,35.75)
\multiput(89,38.5)(-.03658537,-.03353659){82}{\line(-1,0){.03658537}}
%\end
\put(11.75,18.25){$\Delta$}
\put(99,21.25){$\Delta'$}
\put(41,43){$\pi$}
\put(38.5,35){$\Gamma$}
\put(58,36){$\cal Q$}
\put(140.25,44.5){$\cal Q'$}
\put(20.25,42){$\Delta_1$}
\put(108.75,40.25){$\Delta'_1$}
\put(55.75,50.5){${\bf s}_1$}
\put(121.75,51.75){${\bf s'}_1$}
\put(143,25){${\bf s'}_2$}
\put(56,23.75){${\bf s}_2$}
\put(19.25,36.5){${\bf y}_1$}
\put(30.25,27.5){${\bf y}_2$}
\put(148.75,36.25){${\bf y'}_1$}
\put(114,27.25){${\bf y'}_2$}
\put(6.25,32.25){${\bf t}_1$}
\put(71.25,32){${\bf t}_2$}
\put(96,31.5){${\bf t'}_1$}
\put(161.25,32){${\bf t'}_2$}
%\emline(46,54.75)(48.25,53.75)
\multiput(46,54.75)(.075,-.0333333){30}{\line(1,0){.075}}
%\end
%\emline(46.25,52.25)(48.5,53.5)
\multiput(46.25,52.25)(.05921053,.03289474){38}{\line(1,0){.05921053}}
%\end
%\emline(130.5,55.5)(133.5,54.75)
\multiput(130.5,55.5)(.1304348,-.0326087){23}{\line(1,0){.1304348}}
%\end
%\emline(131,53.75)(133,55)
\multiput(131,53.75)(.05263158,.03289474){38}{\line(1,0){.05263158}}
%\end
%\emline(43.75,21)(47,22.25)
\multiput(43.75,21)(.08552632,.03289474){38}{\line(1,0){.08552632}}
%\end
%\emline(44.5,20.75)(47,19.75)
\multiput(44.5,20.75)(.0833333,-.0333333){30}{\line(1,0){.0833333}}
%\end
%\emline(134.75,23.75)(132.5,22.25)
\multiput(134.75,23.75)(-.05,-.03333333){45}{\line(-1,0){.05}}
%\end
%\emline(132.5,22.25)(135,21)
\multiput(132.5,22.25)(.06578947,-.03289474){38}{\line(1,0){.06578947}}
%\end
%\emline(24.5,35.25)(27.25,33.5)
\multiput(24.5,35.25)(.05288462,-.03365385){52}{\line(1,0){.05288462}}
%\end
%\emline(24.5,33)(24.25,32.75)
\multiput(24.5,33)(-.03125,-.03125){8}{\line(0,-1){.03125}}
%\end
%\emline(24.25,32.75)(27.25,33.75)
\multiput(24.25,32.75)(.1,.0333333){30}{\line(1,0){.1}}
%\end
%\emline(38,30.25)(40.5,31.75)
\multiput(38,30.25)(.05555556,.03333333){45}{\line(1,0){.05555556}}
%\end
%\emline(38.25,30.25)(40,29)
\multiput(38.25,30.25)(.04605263,-.03289474){38}{\line(1,0){.04605263}}
%\end
%\emline(38,30.25)(40,28.5)
\multiput(38,30.25)(.03846154,-.03365385){52}{\line(1,0){.03846154}}
%\end
%\emline(120.25,35.75)(121.25,33)
\multiput(120.25,35.75)(.0333333,-.0916667){30}{\line(0,-1){.0916667}}
%\end
%\emline(122.5,35.5)(121.25,33.5)
\multiput(122.5,35.5)(-.03289474,-.05263158){38}{\line(0,-1){.05263158}}
%\end
%\emline(125.75,44)(128.25,43)
\multiput(125.75,44)(.0833333,-.0333333){30}{\line(1,0){.0833333}}
%\end
%\emline(126,41.5)(128.5,42.75)
\multiput(126,41.5)(.06578947,.03289474){38}{\line(1,0){.06578947}}
%\end
\end{picture}

\begin{center}
\nopagebreak[4] Figure \theppp
\end{center}
\addtocounter{ppp}{1}	

%\end{figure} 
	\proof (1) Let $\partial(\pi, \Gamma, {\bf y}_1)={\bf p}_1{\bf q}_1{\bf p}_2{\bf q}_2$; then $Lab ({\bf q}_1)$ is a $u^{\pm 1}$-periodic word, where $u^n$ is the boundary label of $\partial\pi$. We will assume that it is $u$-periodic. So $|{\bf q}_1|_{\theta}>\varepsilon n |u|_{\theta}$ and $|{\bf p}_1|_{\theta}=|{\bf p}_2|_{\theta}=0$
since $r(\Gamma)=0$ (when regarding
$a$-edges as $0$-edges).

Since $u$ is a simple period of some rank $k>0$, no $\theta$-band can start
and also end on ${\bf q}_1$. The same property holds for the subpath ${\bf q}_2$ of ${\bf y}_1$ by Lemma \ref{NoAnnul} (1). Therefore every maximal
$\theta$-band of $\Gamma$ starts on ${\bf q}_1$ and ends on ${\bf q}_2$. It follows
that the words $Lab({\bf q}_1)$ and $Lab({\bf q}_2)$ have equal $\theta$-projections; in particular,
$|{\bf q}_2|_{\theta}=|{\bf q}_1|_{\theta}$ and $Lab({\bf q}_2^{-1})$ is a $v$-periodic word with $|v|_{\theta}=|u|_{\theta}$, because the side
label of a reduced $q$-band is $\theta$-defined
(see Remark \ref{tbd} (1)).
Thus, the word $v$ can be chosen so that the $\theta$-projection of $v$ is the $\theta$-projection of $u$.

Replacing now $\Gamma$ with a smaller
contiguity
diagram $\Gamma'$ with $\partial\Gamma'={\bf p}'_1{\bf q}'_1{\bf p}'_2{\bf q}'_2$,
we may assume that $Lab({\bf q}'_1)\equiv u^l$,
$Lab({\bf q}'_2)^{-1}\equiv v^l$, $|{\bf q}'_1|_{\theta}=|{\bf q}'_2|_{\theta}$  and $l >\frac{\varepsilon n}{2}$. Therefore
\begin{equation}\label{ul}
x_1u^l= v^lx_2
\end{equation}
in the group $H(a)$, where $|x_1|_{\theta}=|x_2|_{\theta}=0$.
The word $v$ is cyclically reduced since the history $h$ is reduced.

\begin{remark}\label{same} If ${\cal Q}_0$ is a subband of $\cal Q$ with history $h_0^2$ corresponding to $v^2$ and ${\cal Q}_0={\cal Q}_1{\cal Q}_2$,
where ${\cal Q}_1$ has history $h_0$,
then the top and the bottom of ${\cal Q}_1$ have equal labels since both occurrences of $h_0$ start with the same rule. This makes possible the following surgery: one can construct
a $q$-band by inserting, between ${\cal Q}_1$ and ${\cal Q}_2$, a $q$-band with
history of the form $h_0^d$.
\end{remark}

Assume that the word $v$ is not simple in rank $0$.  Then
By Lemmas \ref{3.25}, there is a simple in
rank $0$ word $w$ such that for some word $x$, we have
\begin{equation}\label{vw}
v=xw^bgx^{-1}
\end{equation}
in rank $0$, where $g\in {\bf 0}(w)$ and $|w|_{\theta}\le |v|_{\theta}$. Since the word $v$ is cyclically
reduced in rank $0$, one can choose $x$ with $\theta$-length $0$. Indeed, otherwise we have that either the product $xw$ or the product $wgx^{-1}$ (and therefore $wx^{-1}$ since $g\in {\bf 0}(w)$)
has a 2-letter subword equal to a 0-word in rank $0$. It follows that $xwx^{-1}$ is equal
to $yw'y^{-1}$, where $|w'|_{\theta}=|w|_{\theta}$ and $|y|_{\theta}<|x|_{\theta}$. By Remark \ref{0simple}, the word $w'$ is also simple.
We obtain $v=yw'^bg'y^{-1}$, with $g'= (y^{-1}x)g(y^{-1}x)^{-1}$ and $g'\in {\bf 0}(w')$ (see Condition (Z3)). So the element $x$ in (\ref{vw}) can be shorten.

Now it follows from (\ref{vw}) that  $v^l=xw^{bl}x'x^{-1}$ for some $x'\in  {\bf 0}(w)$ since $w$ normalizes the subgroup ${\bf 0}(w)$.
Together with (\ref{ul}), this implies
the equality
\begin{equation}\label{uw}
v^lx_2=xw^{bl}x'x^{-1}x_2=x_1u^l
\end{equation}
in rank $0$, where both $u$ and $w$ are simple in rank $0$ and $x_2, x, x', x_1$ have  $\theta$-length $0$.
Since $l>\frac{\varepsilon n}{2}>\delta^{-1}$, we can
apply Lemma \ref{3.32} to equality (\ref{uw}) and conclude that $u=x''w^{\pm 1}x_3(x'')^{-1}$ in rank $0$ for some word $x''$ and some $x_3\in{\bf 0}(w)$. We obtain a contradiction because
$u$ is simple in rank $0$ but $|w|_{\theta}<|v|_{\theta}=|u|_{\theta}$.

Thus, $v$ is a simple word of rank $0$, and we may apply Lemma \ref{3.32} to the diagram $\Gamma'$ with periodic
labels $u^l$ and $v^{-l}$ of the contiguity arcs. It follows that
$v=Lab({\bf p}'_1)u^{\pm 1}zLab({\bf p}'_1)^{-1}$ in rank $0$ for some $z\in {\bf 0}(u)$. By Lemma \ref{star1}  for $H(a)$, we have $(u^{\pm 1}z)^n=u^{\pm n}$.
Hence
\begin{equation}\label{vn}
v^n= Lab({\bf p}'_1)u^{\pm n}Lab({\bf p}'_1)^{-1}
\end{equation}
in $H(a)$.

The right-hand side of (\ref{vn}) is the label of the closed path starting at the vertex $({\bf p}'_1)_-$ and
bounding $\partial\pi$. So one can insert it in ${\bf y}_1$ to bypath $\pi$ and
obtain an annular diagram $E$ of smaller type than
$\Delta_1$. To obtain $\Delta_1'$ we
glue up the diagram of rank $0$ corresponding to the equality (\ref{vn}) to $E$ and obtain the desired $\Delta'_1$.

However passing from $\Delta_1$ to $\Delta'_1$, we insert the subpath labeled by $v^n$ in ${\bf y}_1$ . So we have to change the q-band $\cal Q$ accordingly, i.e. using Remark \ref{same}, we insert in $\cal Q$ a
$q$-band with history $h_0^{\pm n}$
with possible reduction of the history
of the obtained $q$-band $\cal Q'$.
Now we see that the path ${\bf y}_2$ is modified too. But the label of ${\bf y}'_2$ is equal to the label of ${\bf y}_2$ modulo
the relations of $H(\theta,a)$, and so
one can glue up a number of cells corresponding to these relations to
${\bf y}'_2$ so that these additional cells
and the diagram $\Delta_2$ with boundary ${\bf s}_2{\bf y}_2$ form a diagram $\Delta'_2$.

The diagram $\Delta'$ is built from
$\Delta'_1$, $\cal Q'$ and $\Delta'_2$
(with subsequent reductions if needed). Property (d) follows from Lemma \ref{notha} (2). Indeed, if $Lab({\bf s}_1)$ is conjugate
of a word of length $0$ in $H(\theta,a)$, then
so is $Lab({\bf q}_1)$. Then Lemmas \ref{thea} (b) and \ref{tohist} imply that $h$ is a conjugate
of $1$ in $B(\Theta^+)$, despite  of the assumption
of the lemma.

(2) The proof is similar.
\endproof

\begin{lm} \label{min} Let $\cal D$ be a reduced $q$-band
with boundary ${\bf x}_1{\bf y}_1{\bf x}_2{\bf y}_2$,
where ${\bf y}_1$ and ${\bf y}_2$ are the
top and the bottom of $\cal D$.
Then

(a) there exists a $q$-band
$\cal D'$
with boundary ${\bf x}'_1{\bf y}'_1{\bf x}'_2{\bf y}'_2$, where $Lab ({\bf x}'_i)=Lab ({\bf x}_i)$ ($i=1,2$),
the top ${\bf y}'_1$ of $\cal D'$
has label equal to $Lab({\bf y}_1)$ in $H(\theta,a)$ and
$Lab({\bf y}'_1)$ is minimal in $H(\theta,a)$ word.

(b)
there exists a  reduced $q$-band ${\cal D}^0$ with boundary ${\bf x}_1^0{\bf y}_1^0{\bf x}_2^0{\bf y}_2^0$, where $Lab({\bf x}_j^{0})\equiv Lab ({\bf x}_j)$ ($j=1,2$),  such that the top ${\bf y}_1^0$ having label equal to $Lab({\bf y}_1)$ in $H(\theta,a)$ and $Lab ({\bf y}^0_1)\equiv uvu^{-1}$, where $v$ is cyclically
minimal in $H(\theta,a)$ word.
\end{lm}

\proof
(a) Consider the first version. One may assume that the history $h$ of $\cal D$ is nontrivial in $B(\Theta)$ since otherwise by Lemmas \ref{thea} (b) and \ref{tohist},  $Lab({\bf y}_1)=1$ in $H(\theta,a)$, and the statement
of the lemma is obviously true.

Let $S_1$ be a minimal word equal to $Lab({\bf y}_1)$ in $H(\theta,a)$ and $\Delta_1$ a minimal diagram over $H(\theta,a)$ with boundary ${\bf s}_1{\bf y}_1$, where $Lab({\bf s}_1)\equiv S_1^{-1}$.  We attach $\Delta_1$ to $\cal D$ along the path ${\bf y}_1$,
and denote the obtained diagram by $\Delta$.

If $\Delta_1$ has a positive cell, then by Corollary \ref{Hainf} and Lemma \ref{four}, $\Delta_1$ has a cell $\Pi$ with the sum of contiguity degrees
to ${\bf s}_1$ and to ${\bf y}_1$ at least $\bar\gamma$. However the degree
of contiguity to ${\bf s}_1$ is less than $\bar\alpha$ by Lemma \ref{bara} and Remark \ref{geo},
since ${\bf s}_1$ is a geodesic path in $\Delta_1$. Hence the
contiguity degree to ${\bf y}_1$ is at least $\bar\gamma -\bar\alpha>1/3$.  Therefore by Lemma \ref{0cont}, $\Delta_1$ has a positive  cell $\pi$ and a diagram of contiguity $\Gamma$ of
rank $0$ of $\pi$ to ${\bf y}_1$ with $(\pi,\Gamma,{\bf y}_1)>\varepsilon$.
Using now Lemma \ref{through} (2), we can replace $\Delta_1$ with
a diagram $\Delta'_1$ of smaller type with boundary ${\bf s}_1{\bf y}'_1$
and replace $\cal D$ with a reduced $q$-band $\cal D'$  having the same labels of the start and end edges and having top label $Lab ({\bf y}'_1)$ equal to $Lab({\bf y}_1)$ in $H(\theta,a)$.
We can continue decreasing the types
of $\Delta'_1,\Delta''_1,\dots $ until we get a diagram $\Delta^{(k)}$ of rank $0$ and a reduced $q$-band ${\cal D}_k$ whose
top label $Lab({\bf y}_1^{(k)})$ is is equal to $Lab ({\bf y}_1)$ in $H(\theta,a)$ and equal to $S_1$ in rank $0$. Hence
$|{\bf y}_1^{(k)}|_{\theta}=|{\bf s}_1|_{\theta}$ and the top label of ${\cal D'}={\cal D}_k$ is a minimal in $H(\theta,a)$ word, as required.

(b) The word $Lab ({\bf y}_1)$ is conjugate in $H(\theta,a)$ to a cyclically minimal word $S_1$,
and so there exists a minimal  annular
diagram $\Delta_1$ over $H(\theta,a)$ whose inner boundary component is
labeled by $S_1^{-1}$ and the outer component of the boundary
has label $Lab({\bf y}_1)$. Therefore one can obtain an annular
diagram $\Delta$ identifying the boundary segments of $\Delta_1$ and $\cal D$ labeled by $Lab({\bf y}_1)$. Then we can reduce the type of $\Delta_1$ as in item (a) (but using the statement (1) of Lemma \ref{through}), and the modifications of ${\bf y}_1$
do not change its label modulo the relations of $H(\theta,a)$.
Finally we will replace $\Delta_1$ with an annular diagram
$E$ without positive cells, where one of the boundary label is the cyclically
minimal word $S_1^{-1}$ and therefore every $\theta$-edges of these boundary is adjacent to an edge
of ${\bf y}^0_1$. (See \setcounter{pdeleven}{\value{ppp}} Figure
\thepdeleven .)

\begin{figure}
% This is a LaTeX picture output by TeXCAD.
% File name: [uv.pic].
% Version of TeXCAD: 4.3
% Reference / build: 30-Jun-2012 (rev. 105)
% For new versions, check: http://texcad.sf.net/
% Options on the following lines.
%\grade{\on}
%\emlines{\off}
%\epic{\off}
%\beziermacro{\on}
%\reduce{\on}
%\snapping{\off}
%\pvinsert{% Your \input, \def, etc. here}
%\quality{8.000}
%\graddiff{0.005}
%\snapasp{1}
%\zoom{4.0000}
\unitlength 1mm % = 2.845pt
\linethickness{0.4pt}
\ifx\plotpoint\undefined\newsavebox{\plotpoint}\fi % GNUPLOT compatibility
\begin{picture}(82,41.75)(-10,10)
\put(43.5,30.875){\oval(20,10.75)[]}
%\emline(58.75,32.75)(54.75,38.25)
\multiput(58.75,32.75)(-.033613445,.046218487){119}{\line(0,1){.046218487}}
%\end
\put(54.75,38.25){\line(-1,0){21}}
%\emline(33.75,38.25)(30,31.25)
\multiput(33.75,38.25)(-.033482143,-.0625){112}{\line(0,-1){.0625}}
%\end
\put(30,31.25){\line(1,-2){3.5}}
\put(33.5,24.25){\line(1,0){21}}
%\emline(54.5,24.25)(58.5,28.5)
\multiput(54.5,24.25)(.033613445,.035714286){119}{\line(0,1){.035714286}}
%\end
%\emline(58.5,28.5)(80.25,30.75)
\multiput(58.5,28.5)(.32462687,.03358209){67}{\line(1,0){.32462687}}
%\end
%\emline(80.25,30.75)(58.75,32.5)
\multiput(80.25,30.75)(-.41346154,.03365385){52}{\line(-1,0){.41346154}}
%\end
%\emline(81,35.75)(79.75,30.5)
\multiput(81,35.75)(-.03289474,-.13815789){38}{\line(0,-1){.13815789}}
%\end
%\emline(79.75,30.5)(81.25,27)
\multiput(79.75,30.5)(.03333333,-.07777778){45}{\line(0,-1){.07777778}}
%\end
%\emline(81.25,27)(60.25,25.75)
\multiput(81.25,27)(-.55263158,-.03289474){38}{\line(-1,0){.55263158}}
%\end
\put(60.25,25.75){\line(-5,-6){3.75}}
%\emline(56.5,21.25)(32.5,21.5)
\multiput(56.5,21.25)(-3,.03125){8}{\line(-1,0){3}}
%\end
%\emline(27.5,30.5)(32.25,21.75)
\multiput(27.5,30.5)(.033687943,-.062056738){141}{\line(0,-1){.062056738}}
%\end
%\emline(28,30.5)(32.75,41.5)
\multiput(28,30.5)(.033687943,.078014184){141}{\line(0,1){.078014184}}
%\end
\put(32.75,41.5){\line(1,0){24}}
%\emline(56.75,41.5)(60.75,36)
\multiput(56.75,41.5)(.033613445,-.046218487){119}{\line(0,-1){.046218487}}
%\end
%\emline(60.75,36)(80.5,34.5)
\multiput(60.75,36)(.43888889,-.03333333){45}{\line(1,0){.43888889}}
%\end
\put(71.5,33){$u$}
%\emline(71,32)(73.5,32.5)
\multiput(71,32)(.1666667,.0333333){15}{\line(1,0){.1666667}}
%\end
\put(71,31.75){\line(1,0){.5}}
%\emline(71.5,31.75)(72,30.75)
\multiput(71.5,31.75)(.0333333,-.0666667){15}{\line(0,-1){.0666667}}
%\end
%\emline(68.25,30.25)(70.25,29.75)
\multiput(68.25,30.25)(.1333333,-.0333333){15}{\line(1,0){.1333333}}
%\end
%\emline(70.25,29.75)(68.75,28.75)
\multiput(70.25,29.75)(-.05,-.0333333){30}{\line(-1,0){.05}}
%\end
\put(71.25,26.75){$u'$}
\put(55.5,29.5){$E$}
\put(35.5,28.25){$S_1$}
%\emline(36.75,27.25)(37.5,26.5)
\multiput(36.75,27.25)(.0326087,-.0326087){23}{\line(0,-1){.0326087}}
%\end
%\emline(36.25,25.75)(36.5,26.5)
\multiput(36.25,25.75)(.03125,.09375){8}{\line(0,1){.09375}}
%\end
%\emline(35.75,25.25)(37.5,26.25)
\multiput(35.75,25.25)(.0583333,.0333333){30}{\line(1,0){.0583333}}
%\end
\put(55.25,37.75){\line(0,-1){2}}
%\emline(55.5,37.75)(56.25,37.25)
\multiput(55.5,37.75)(.05,-.0333333){15}{\line(1,0){.05}}
%\end
\put(57.25,36){$v$}
\put(82,33.25){${\bf x}_1^0$}
\put(82,28.75){${\bf x}_2^0$}
%\emline(45.25,21.75)(47,22.5)
\multiput(45.25,21.75)(.076087,.0326087){23}{\line(1,0){.076087}}
%\end
%\emline(45,21.25)(47,20.5)
\multiput(45,21.25)(.0869565,-.0326087){23}{\line(1,0){.0869565}}
%\end
\put(48.5,16.75){${\bf y}_2^0$}
%\emline(76.25,34.75)(76,31.5)
\multiput(76.25,34.75)(-.03125,-.40625){8}{\line(0,-1){.40625}}
%\end
\put(69,35.75){\line(0,-1){4}}
\put(61.25,36.25){\line(0,-1){3}}
\put(54.25,41.5){\line(0,-1){3.25}}
%\emline(48.5,41.5)(48.75,38)
\multiput(48.5,41.5)(.03125,-.4375){8}{\line(0,-1){.4375}}
%\end
\put(43.25,41.75){\line(0,-1){3.75}}
\put(37.75,41.5){\line(0,-1){3.25}}
%\emline(32.5,41.5)(33.75,38.25)
\multiput(32.5,41.5)(.03289474,-.08552632){38}{\line(0,-1){.08552632}}
%\end
\put(64.75,35.75){\line(0,-1){3.75}}
%\emline(30.75,37.25)(32.5,36.25)
\multiput(30.75,37.25)(.0583333,-.0333333){30}{\line(1,0){.0583333}}
%\end
%\emline(29.25,35)(31.25,34)
\multiput(29.25,35)(.0666667,-.0333333){30}{\line(1,0){.0666667}}
%\end
\put(28,31.25){\line(1,0){1.5}}
%\emline(30.75,28.5)(29.25,27.75)
\multiput(30.75,28.5)(-.0652174,-.0326087){23}{\line(-1,0){.0652174}}
%\end
%\emline(32.5,25.75)(31,24.75)
\multiput(32.5,25.75)(-.05,-.0333333){30}{\line(-1,0){.05}}
%\end
\put(34.75,24){\line(0,-1){2}}
\put(39,24.25){\line(0,-1){3}}
%\emline(43.25,24.5)(43,21.25)
\multiput(43.25,24.5)(-.03125,-.40625){8}{\line(0,-1){.40625}}
%\end
\put(48.25,24){\line(0,-1){2.75}}
%\emline(53,24.25)(52.75,21.75)
\multiput(53,24.25)(-.03125,-.3125){8}{\line(0,-1){.3125}}
%\end
%\emline(55.75,25.5)(57.75,23.75)
\multiput(55.75,25.5)(.03846154,-.03365385){52}{\line(1,0){.03846154}}
%\end
%\emline(59.75,28.5)(60.75,26.5)
\multiput(59.75,28.5)(.0333333,-.0666667){30}{\line(0,-1){.0666667}}
%\end
%\emline(64.25,29.25)(64.75,27.25)
\multiput(64.25,29.25)(.0333333,-.1333333){15}{\line(0,-1){.1333333}}
%\end
\put(68.75,29.5){\line(0,-1){3}}
%\emline(75.25,30)(75.5,27.25)
\multiput(75.25,30)(.03125,-.34375){8}{\line(0,-1){.34375}}
%\end
%\vector(24.25,37.25)(28.5,35)
\put(28.5,35){\vector(2,-1){.07}}\multiput(24.25,37.25)(.06343284,-.03358209){67}{\line(1,0){.06343284}}
%\end
\put(19.5,38.5){${\cal D}^0$}
%\emline(33.75,31.75)(38.75,36.25)
\multiput(33.75,31.75)(.037313433,.03358209){134}{\line(1,0){.037313433}}
%\end
%\emline(35.5,31.5)(41,36)
\multiput(35.5,31.5)(.041044776,.03358209){134}{\line(1,0){.041044776}}
%\end
%\emline(38,31.5)(43.25,36)
\multiput(38,31.5)(.039179104,.03358209){134}{\line(1,0){.039179104}}
%\end
\put(38.75,29.75){\line(1,0){.75}}
%\emline(38.75,30.75)(45.5,36)
\multiput(38.75,30.75)(.043269231,.033653846){156}{\line(1,0){.043269231}}
%\end
%\emline(38.5,29)(47.5,36)
\multiput(38.5,29)(.043269231,.033653846){208}{\line(1,0){.043269231}}
%\end
%\emline(38,27)(49,35.5)
\multiput(38,27)(.0436507937,.0337301587){252}{\line(1,0){.0436507937}}
%\end
%\emline(38.75,25.5)(39.25,26.25)
\multiput(38.75,25.5)(.0333333,.05){15}{\line(0,1){.05}}
%\end
%\emline(39.25,26.25)(51,35.25)
\multiput(39.25,26.25)(.0440074906,.0337078652){267}{\line(1,0){.0440074906}}
%\end
%\emline(41,25.5)(52.25,34.25)
\multiput(41,25.5)(.0432692308,.0336538462){260}{\line(1,0){.0432692308}}
%\end
%\emline(43.5,25.5)(53,32.75)
\multiput(43.5,25.5)(.044186047,.03372093){215}{\line(1,0){.044186047}}
%\end
%\emline(46,25.5)(53.5,30.75)
\multiput(46,25.5)(.048076923,.033653846){156}{\line(1,0){.048076923}}
%\end
%\emline(49.25,25.75)(49.5,25.5)
\multiput(49.25,25.75)(.03125,-.03125){8}{\line(0,-1){.03125}}
%\end
%\emline(49.5,25.5)(49.25,25.25)
\multiput(49.5,25.5)(-.03125,-.03125){8}{\line(0,-1){.03125}}
%\end
\put(48.5,25.75){\line(5,3){5}}
%\emline(49.75,36.25)(53.5,30.75)
\multiput(49.75,36.25)(.033482143,-.049107143){112}{\line(0,-1){.049107143}}
%\end
%\emline(48.25,36.25)(53,29.25)
\multiput(48.25,36.25)(.033687943,-.04964539){141}{\line(0,-1){.04964539}}
%\end
%\emline(47,36.5)(53,27.25)
\multiput(47,36.5)(.033707865,-.051966292){178}{\line(0,-1){.051966292}}
%\end
%\emline(45.75,36.25)(51.5,27)
\multiput(45.75,36.25)(.033625731,-.054093567){171}{\line(0,-1){.054093567}}
%\end
%\emline(44,36.5)(50.25,26.5)
\multiput(44,36.5)(.033602151,-.053763441){186}{\line(0,-1){.053763441}}
%\end
%\emline(42.5,36.25)(48.5,26.5)
\multiput(42.5,36.25)(.033707865,-.054775281){178}{\line(0,-1){.054775281}}
%\end
%\emline(41.25,35.75)(47,26.25)
\multiput(41.25,35.75)(.033625731,-.055555556){171}{\line(0,-1){.055555556}}
%\end
%\emline(39.25,36.25)(45.25,25.5)
\multiput(39.25,36.25)(.033707865,-.060393258){178}{\line(0,-1){.060393258}}
%\end
%\emline(37.75,36.25)(43.75,25.25)
\multiput(37.75,36.25)(.033707865,-.061797753){178}{\line(0,-1){.061797753}}
%\end
%\emline(36.5,35.75)(42.25,26)
\multiput(36.5,35.75)(.033625731,-.057017544){171}{\line(0,-1){.057017544}}
%\end
%\emline(35.25,35)(37.25,31.75)
\multiput(35.25,35)(.03333333,-.05416667){60}{\line(0,-1){.05416667}}
%\end
%\emline(39.25,28.75)(40.75,26)
\multiput(39.25,28.75)(.03333333,-.06111111){45}{\line(0,-1){.06111111}}
%\end
%\emline(34.75,33.75)(35.25,32.25)
\multiput(34.75,33.75)(.0333333,-.1){15}{\line(0,-1){.1}}
%\end
%\emline(38.5,27.75)(39.75,25.25)
\multiput(38.5,27.75)(.03289474,-.06578947){38}{\line(0,-1){.06578947}}
%\end
\end{picture}
\begin{center}
\nopagebreak[4] Figure \theppp
\end{center}
\addtocounter{ppp}{1}	

\end{figure} 

Hence $Lab({\bf y}^0_1)$ is of the form $uvu'$, where $v$ is equal to a cyclic permutation of $S_1$
% it follows that $uu'=1$. 
 $u$ and $u'$
are the labels of side arcs of some subbands of the band ${\cal D}^0$ obtained after the modifications of $\cal D$, and the $\theta$-edges of these arcs have to be mutually  adjacent in the diagram $E$ of rank $0$. It follows that $u'\equiv u^{-1}$ since side labels of $q$-bands are $\theta$-defined. Statement (b) is proved.
\endproof

\subsection{Axioms for $H({\bf\Theta})$}\label{ax}

In this section we start with the group ${\bf G}(0)=H(\bf\Theta)$. The $q$-letters are now non-zero letters,
while $\theta$- and $a$-letters are regarded as $0$-letters. The group $H(\bf\Theta)$ satisfies all Burnside relation $u^n=1$ of $q$-length zero,
so we have Property (Z1.2) defined in
Section \ref{axioms}. Property (Z1.1) follows
from the definition of $(\theta,q)$-relations.

By definition, the group \index[g]{${\bf G}(1/2)$} ${\bf G}(1/2)$ is the factor group of ${\bf G}(0)$ modulo the hub relations.

\begin{lm} \label{Z2} Property (Z2) holds in the group ${\bf G}(1/2)$.
\end{lm}

\proof Since we have chosen $N\ge n$ for the number of $q$-letters in the hubs, Property (Z2.1) holds. Since all $q$-letters of the hubs are pairwise distinct, we also obtain Property (Z2.2).

Then we will use the notation and the diagram $\Delta$ from the definition of (Z2.3). Let us denote by $\Gamma$
the subdiagram of $\Delta$ with boundary ${\bf p}_1{\bf q}_1{\bf p}_2{\bf q}_2$, where $Lab({\bf p}_1)\equiv u_1$, $Lab({\bf q_1})\equiv v_1$, $Lab({\bf p}_2)\equiv u_2^{-1}$,
$Lab({\bf q}_2)\equiv v_2^{-1}$ (see Fig. 5 in Subsection \ref{axioms}). We may assume that $\Gamma$ is a minimal diagram over the group ${\bf G}(0)=H(\bf\Theta)$.

Since we have $|v_1|_q\ge \varepsilon |v_1w_1|_q$ and one can choose $L>3\varepsilon n$, at least three maximal $t$-bands start
in $\Gamma$ on ${\bf q}_1$. By Property (Z2.2) they cannot end on ${\bf q}_1$. They cannot end on ${\bf p}_1$ or ${\bf p}_2$ either since $|{\bf p}_1|_q=|{\bf p}_2|_q=0$. So they
end on ${\bf q}_2$, and therefore there are two
maximal $t$-bands ${\cal Q}_i$ and ${\cal Q}_{i+1}$ connecting ${\bf q}_1$ and ${\bf q}_2$ and corresponding to $t$-letters
$t^{(i)}$ and $t^{(i+1)}$ with $\{i, i+1\} \ne \{0, 1\}$.

To prove that the boundary label of $\Delta$ is trivial in ${\bf G}(0)=H(\bf\Theta)$, one may cut off subdiagrams over $H(\bf\Theta)$ from $\Delta$ or glue up such diagrams to $\Delta$.
We will define such surgeries below.
At the end, we will obtain a spherical
diagram, whose boundary label is the empty word obviously equal to $1$  in
$H(\bf\Theta)$, as desired.

First of all we replace the subdiagram $\Gamma$ with
the smaller subdiagram $\Gamma_i$ bounded by ${\cal Q}_i$, ${\cal Q}_{i+1}$, and the subpaths ${\bf \bar q}_1$ and
${\bf \bar q}_2$ of ${\bf q}_1$ and ${\bf q}_2$,resp. By Lemma \ref{noq}, $\Gamma_i$ contains no $q$-annuli. So every maximal $q$-band of $\Gamma_i$ connects
${\bf \bar q}_1$ and ${\bf \bar q}_2$. Let ${\cal D}_0={\cal Q}_i,\dots, {\cal D}_m={\cal Q}_{i+1}$ be
all these consecutive $q$-bands starting on ${\bf \bar q}_1$.

By Lemma \ref{min}, it is possible to replace the band ${\cal D}_0$
with a band ${\cal D'}_0$ having minimal in $H(\theta,a)$ label
of the top. (This operation can also change the boundary label of the whole diagram modulo the relations of $H(\theta,a)$.) By Lemma
\ref{thea} (b), the bottom label of ${\cal D'}_0$ is also a minimal word. We may keep notation ${\cal D}_0$ for the modified $q$-band.

Consider now the minimal subdiagram $E$ over $H(\theta,a)$  between ${\cal D}_0$ and ${\cal D}_1$ bounded by these two $q$-bands and by subpaths ${\bf x}_1$, ${\bf x}_2$ of ${\bf \bar q}_1$ and ${\bf \bar q}_2$. If $E$ has a positive cell, then
by Lemma \ref{four}, it has a cell $\Pi$ with the sum of contiguity degrees to ${\bf x}_1$, ${\bf x}_2$ and to the sides (i.e. top/bottoms) of
 ${\cal D}_0$ and ${\cal D}_1$ at least $\bar\gamma$. However the first two contiguity degrees are $0$ since ${\bf x}_1$ and ${\bf x}_2$ have
 no $\theta$-edges, and the degree of contiguity to the side of ${\cal D}_0$ is $<\bar \alpha$ by Lemma \ref{bara} and Remark \ref{geo} since the label of this side is minimal in $H(\theta,a)$. Hence the degree of contiguity of $\Pi$ to the side of ${\cal D}_1$ is $>\bar\gamma-\bar\alpha>1/3$, and therefore by Lemma \ref{0cont}, there exists a positive cell $\pi$ in
 $E$ with contiguity diagram  of rank $0$ to the side of ${\cal D}_1$ and wit degree of contiguity at least $\varepsilon$. Now Lemma \ref{through} helps us to reduce the type of $E$ (but the type
 of diagram over $H(\theta,a)$ between ${\cal D}_1$ and ${\cal D}_2$ may increase). Repeating such operation we may obtain a modified subdiagram $E$ without cells of positive ranks. Thus,
 we may assume that $E$ is a diagram over the group $H(a)$.

 It follows that the projections of side labels of the bands
 ${\cal D}_0$ and modified ${\cal D}_1$ (for which we keep the same notation) on the alphabet $\bf\Theta$ are equal, and therefore
 the side labels of ${\cal D}_1$ are also minimal over $H(\theta,a)$ by Lemma \ref{thea} (b). Therefore one can eliminate now all positive cell from the subdiagram settled between ${\cal D}_1$ and ${\cal D}_2$ at the expense of the modifications between ${\cal D}_2$ and ${\cal D}_3$, and so on.
 Thus one may assume now that the whole $\Gamma_i$ is a diagram over the group $H({\bf Y})$.

 By Lemma \ref{NoAnnul}, every maximal $\theta$-band of $\Gamma_i$
 crosses every $q$-band of $\Gamma_i$ exactly once, and so we obtain a quasi-trapezium $\Delta_i$ bounded by ${\cal Q}_i$, ${\cal Q}_{i+1}$ and containing all $\theta$-bands of $\Gamma_i$. By Lemma \ref{simul}, $\Delta_0$ corresponds to some quasi-computation
 $\cal C$, and $\cal C$ should start/end with the start/end rules of $\bf M$ since the labels of the paths ${\bf \bar q}_1$ and ${\bf \bar q}_2$ are the subwords of hub relation. Since these rules
 lock all sectors between $t^{(i)}$ and $t^{(i+1)}$, there are
 no $a$-cells between ${\bf \bar q}_1$ (${\bf \bar q}_2$) and the nearest
 $\theta$-band, i.e. $\Delta_i=\Gamma_i$.

 Since the side label  of ${\cal D}_0={\cal Q}_i$ is minimal
 in $H(\theta,a)$, the history of $\Gamma_i$ is simple by Corollary \ref{tb} (b), and
 by Lemmas \ref{beza} and \ref{NoAnnul}, there is a trapezium
 $\Gamma'_i$ with the same top and bottom labels  whose history is 
congruent to the history of $\Gamma_i$, and so the side
 labels are equal in $H(\theta,a)$ to the side labels of $\Gamma_i$ by Lemma \ref{thea} (a).
 Hence the replacement $\Gamma_i$ by $\Gamma'_i$ does not
 change the boundary label of the whole diagram modulo the relations of $H(\bf\Theta)$. So we further assume that $\Gamma_i$
 is a trapezium and $\cal C$ is a computation.
(See \setcounter{pdeleven}{\value{ppp}} Figure
\thepdeleven.)

\begin{figure}
% This is a LaTeX picture output by TeXCAD.
% File name: [hubs.pic].
% Version of TeXCAD: 4.3
% Reference / build: 30-Jun-2012 (rev. 105)
% For new versions, check: http://texcad.sf.net/
% Options on the following lines.
%\grade{\on}
%\emlines{\off}
%\epic{\off}
%\beziermacro{\on}
%\reduce{\on}
%\snapping{\off}
%\pvinsert{% Your \input, \def, etc. here}
%\quality{8.000}
%\graddiff{0.005}
%\snapasp{1}
%\zoom{4.0000}
\unitlength 1mm % = 2.845pt
\linethickness{0.4pt}
\ifx\plotpoint\undefined\newsavebox{\plotpoint}\fi % GNUPLOT compatibility
\begin{picture}(88,72.5)(-35,0)
\put(47.375,63.375){\oval(65.25,18.25)[]}
\put(78.25,16.5){\oval(0,.5)[]}
\put(47.625,19.875){\oval(65.75,18.75)[]}
\put(35.25,53.5){\line(0,-1){24.25}}
\put(37.25,54){\line(0,-1){25}}
\put(47.25,54.25){\line(0,-1){26}}
\put(49.5,53.75){\line(0,-1){24.75}}
%\dashline{1}(59,54)(59,29.25)
\put(58.93,53.93){\line(0,-1){.99}}
\put(58.93,51.95){\line(0,-1){.99}}
\put(58.93,49.97){\line(0,-1){.99}}
\put(58.93,47.99){\line(0,-1){.99}}
\put(58.93,46.01){\line(0,-1){.99}}
\put(58.93,44.03){\line(0,-1){.99}}
\put(58.93,42.05){\line(0,-1){.99}}
\put(58.93,40.07){\line(0,-1){.99}}
\put(58.93,38.09){\line(0,-1){.99}}
\put(58.93,36.11){\line(0,-1){.99}}
\put(58.93,34.13){\line(0,-1){.99}}
\put(58.93,32.15){\line(0,-1){.99}}
\put(58.93,30.17){\line(0,-1){.99}}
%\end
%\dashline{1}(61.25,54)(61.25,28.75)
\put(61.18,53.93){\line(0,-1){.9712}}
\put(61.18,51.987){\line(0,-1){.9712}}
\put(61.18,50.045){\line(0,-1){.9712}}
\put(61.18,48.103){\line(0,-1){.9712}}
\put(61.18,46.16){\line(0,-1){.9712}}
\put(61.18,44.218){\line(0,-1){.9712}}
\put(61.18,42.276){\line(0,-1){.9712}}
\put(61.18,40.334){\line(0,-1){.9712}}
\put(61.18,38.391){\line(0,-1){.9712}}
\put(61.18,36.449){\line(0,-1){.9712}}
\put(61.18,34.507){\line(0,-1){.9712}}
\put(61.18,32.564){\line(0,-1){.9712}}
\put(61.18,30.622){\line(0,-1){.9712}}
%\end
%\dashline{1}(25,54.25)(25,28.5)
\put(24.93,54.18){\line(0,-1){.9904}}
\put(24.93,52.199){\line(0,-1){.9904}}
\put(24.93,50.218){\line(0,-1){.9904}}
\put(24.93,48.237){\line(0,-1){.9904}}
\put(24.93,46.257){\line(0,-1){.9904}}
\put(24.93,44.276){\line(0,-1){.9904}}
\put(24.93,42.295){\line(0,-1){.9904}}
\put(24.93,40.314){\line(0,-1){.9904}}
\put(24.93,38.334){\line(0,-1){.9904}}
\put(24.93,36.353){\line(0,-1){.9904}}
\put(24.93,34.372){\line(0,-1){.9904}}
\put(24.93,32.391){\line(0,-1){.9904}}
\put(24.93,30.41){\line(0,-1){.9904}}
%\end
%\dashline{1}(22.75,54)(22.75,29.5)
\put(22.68,53.93){\line(0,-1){.98}}
\put(22.68,51.97){\line(0,-1){.98}}
\put(22.68,50.01){\line(0,-1){.98}}
\put(22.68,48.05){\line(0,-1){.98}}
\put(22.68,46.09){\line(0,-1){.98}}
\put(22.68,44.13){\line(0,-1){.98}}
\put(22.68,42.17){\line(0,-1){.98}}
\put(22.68,40.21){\line(0,-1){.98}}
\put(22.68,38.25){\line(0,-1){.98}}
\put(22.68,36.29){\line(0,-1){.98}}
\put(22.68,34.33){\line(0,-1){.98}}
\put(22.68,32.37){\line(0,-1){.98}}
\put(22.68,30.41){\line(0,-1){.98}}
%\end
%\dashline{1}(70,54)(70.5,48.75)
\put(69.93,53.93){\line(0,-1){.875}}
\put(70.096,52.18){\line(0,-1){.875}}
\put(70.263,50.43){\line(0,-1){.875}}
%\end
%\dashline{1}(71.75,53.75)(72,50.25)
\put(71.68,53.68){\line(0,-1){.875}}
\put(71.805,51.93){\line(0,-1){.875}}
%\end
%\dashline{1}(70.25,35.25)(70.25,28.75)
\put(70.18,35.18){\line(0,-1){.9286}}
\put(70.18,33.323){\line(0,-1){.9286}}
\put(70.18,31.465){\line(0,-1){.9286}}
\put(70.18,29.608){\line(0,-1){.9286}}
%\end
%\dashline{1}(72.25,35.25)(72,28.25)
\put(72.18,35.18){\line(0,-1){.875}}
\put(72.117,33.43){\line(0,-1){.875}}
\put(72.055,31.68){\line(0,-1){.875}}
\put(71.992,29.93){\line(0,-1){.875}}
%\end
%\dashline{1}(78.75,58)(81,52.25)
\multiput(78.68,57.93)(.032143,-.082143){10}{\line(0,-1){.082143}}
\multiput(79.323,56.287)(.032143,-.082143){10}{\line(0,-1){.082143}}
\multiput(79.965,54.644)(.032143,-.082143){10}{\line(0,-1){.082143}}
\multiput(80.608,53.001)(.032143,-.082143){10}{\line(0,-1){.082143}}
%\end
%\dashline{1}(79.5,59.25)(82.25,53.5)
\multiput(79.43,59.18)(.032738,-.068452){12}{\line(0,-1){.068452}}
\multiput(80.215,57.537)(.032738,-.068452){12}{\line(0,-1){.068452}}
\multiput(81.001,55.894)(.032738,-.068452){12}{\line(0,-1){.068452}}
\multiput(81.787,54.251)(.032738,-.068452){12}{\line(0,-1){.068452}}
%\end
%\dashline{1}(80.25,35)(78.25,27.25)
\multiput(80.18,34.93)(-.031746,-.123016){7}{\line(0,-1){.123016}}
\multiput(79.735,33.207)(-.031746,-.123016){7}{\line(0,-1){.123016}}
\multiput(79.291,31.485)(-.031746,-.123016){7}{\line(0,-1){.123016}}
\multiput(78.846,29.763)(-.031746,-.123016){7}{\line(0,-1){.123016}}
\multiput(78.402,28.041)(-.031746,-.123016){7}{\line(0,-1){.123016}}
%\end
%\dashline{1}(82.5,34.25)(79.25,26)
\multiput(82.43,34.18)(-.0325,-.0825){10}{\line(0,-1){.0825}}
\multiput(81.78,32.53)(-.0325,-.0825){10}{\line(0,-1){.0825}}
\multiput(81.13,30.88)(-.0325,-.0825){10}{\line(0,-1){.0825}}
\multiput(80.48,29.23)(-.0325,-.0825){10}{\line(0,-1){.0825}}
\multiput(79.83,27.58)(-.0325,-.0825){10}{\line(0,-1){.0825}}
%\end
%\dashline{1}(79.75,67.5)(85.5,66.25)
\multiput(79.68,67.43)(.136905,-.029762){6}{\line(1,0){.136905}}
\multiput(81.323,67.073)(.136905,-.029762){6}{\line(1,0){.136905}}
\multiput(82.965,66.715)(.136905,-.029762){6}{\line(1,0){.136905}}
\multiput(84.608,66.358)(.136905,-.029762){6}{\line(1,0){.136905}}
%\end
%\dashline{1}(78.75,69.25)(86.5,67.75)
\multiput(78.68,69.18)(.161458,-.03125){6}{\line(1,0){.161458}}
\multiput(80.617,68.805)(.161458,-.03125){6}{\line(1,0){.161458}}
\multiput(82.555,68.43)(.161458,-.03125){6}{\line(1,0){.161458}}
\multiput(84.492,68.055)(.161458,-.03125){6}{\line(1,0){.161458}}
%\end
%\dashline{1}(81,19.5)(88,23.75)
\multiput(80.93,19.43)(.0538462,.0326923){13}{\line(1,0){.0538462}}
\multiput(82.33,20.28)(.0538462,.0326923){13}{\line(1,0){.0538462}}
\multiput(83.73,21.13)(.0538462,.0326923){13}{\line(1,0){.0538462}}
\multiput(85.13,21.98)(.0538462,.0326923){13}{\line(1,0){.0538462}}
\multiput(86.53,22.83)(.0538462,.0326923){13}{\line(1,0){.0538462}}
%\end
%\dashline{1}(80.5,16)(79.75,16.5)
\multiput(80.43,15.93)(-.046875,.03125){8}{\line(-1,0){.046875}}
%\end
%\dashline{1}(80.25,16.25)(80.25,16.5)
\put(80.18,16.18){\line(0,1){.125}}
%\end
%\dashline{1}(80,17.25)(87.75,21.5)
\multiput(79.93,17.18)(.0596154,.0326923){13}{\line(1,0){.0596154}}
\multiput(81.48,18.03)(.0596154,.0326923){13}{\line(1,0){.0596154}}
\multiput(83.03,18.88)(.0596154,.0326923){13}{\line(1,0){.0596154}}
\multiput(84.58,19.73)(.0596154,.0326923){13}{\line(1,0){.0596154}}
\multiput(86.13,20.58)(.0596154,.0326923){13}{\line(1,0){.0596154}}
%\end
%\dashline{1}(15.75,58.5)(11,51.5)
\multiput(15.68,58.43)(-.0316667,-.0466667){15}{\line(0,-1){.0466667}}
\multiput(14.73,57.03)(-.0316667,-.0466667){15}{\line(0,-1){.0466667}}
\multiput(13.78,55.63)(-.0316667,-.0466667){15}{\line(0,-1){.0466667}}
\multiput(12.83,54.23)(-.0316667,-.0466667){15}{\line(0,-1){.0466667}}
\multiput(11.88,52.83)(-.0316667,-.0466667){15}{\line(0,-1){.0466667}}
%\end
%\dashline{1}(15,60.25)(10.25,53)
\multiput(14.93,60.18)(-.0316667,-.0483333){15}{\line(0,-1){.0483333}}
\multiput(13.98,58.73)(-.0316667,-.0483333){15}{\line(0,-1){.0483333}}
\multiput(13.03,57.28)(-.0316667,-.0483333){15}{\line(0,-1){.0483333}}
\multiput(12.08,55.83)(-.0316667,-.0483333){15}{\line(0,-1){.0483333}}
\multiput(11.13,54.38)(-.0316667,-.0483333){15}{\line(0,-1){.0483333}}
%\end
%\dashline{1}(12,34.25)(15.75,27.75)
\multiput(11.93,34.18)(.0334821,-.0580357){14}{\line(0,-1){.0580357}}
\multiput(12.867,32.555)(.0334821,-.0580357){14}{\line(0,-1){.0580357}}
\multiput(13.805,30.93)(.0334821,-.0580357){14}{\line(0,-1){.0580357}}
\multiput(14.742,29.305)(.0334821,-.0580357){14}{\line(0,-1){.0580357}}
%\end
%\dashline{1}(10.25,33)(15.75,25.25)
\multiput(10.18,32.93)(.0323529,-.0455882){17}{\line(0,-1){.0455882}}
\multiput(11.28,31.38)(.0323529,-.0455882){17}{\line(0,-1){.0455882}}
\multiput(12.38,29.83)(.0323529,-.0455882){17}{\line(0,-1){.0455882}}
\multiput(13.48,28.28)(.0323529,-.0455882){17}{\line(0,-1){.0455882}}
\multiput(14.58,26.73)(.0323529,-.0455882){17}{\line(0,-1){.0455882}}
%\end
\put(31.5,49){${\cal Q}_i$}
\put(51.25,49){${\cal Q}_{i+1}$}
\put(62.5,49.75){${\cal Q}_{i+2}$}
\put(40,39){$\Gamma_i$}
\put(51.75,38.5){$\Gamma_{i+1}$}
\put(35.5,51.5){\line(1,0){1.5}}
\put(35.5,48.25){\line(1,0){1.5}}
\put(35.5,45.25){\line(1,0){1.75}}
\put(35.25,42){\line(1,0){2}}
\put(35.5,38){\line(1,0){1.75}}
%\emline(35.5,35.25)(37,35)
\multiput(35.5,35.25)(.1875,-.03125){8}{\line(1,0){.1875}}
%\end
\put(35.25,33){\line(1,0){2}}
\put(47,51.75){\line(1,0){3}}
\put(47.5,48.5){\line(1,0){2}}
%\emline(47.5,46.25)(50,46)
\multiput(47.5,46.25)(.3125,-.03125){8}{\line(1,0){.3125}}
%\end
\put(47.5,43.25){\line(1,0){2}}
\put(47.25,40.75){\line(1,0){2.5}}
\put(47.5,37.75){\line(1,0){1.75}}
%\emline(47.75,35)(49.5,35.25)
\multiput(47.75,35)(.21875,.03125){8}{\line(1,0){.21875}}
%\end
\put(47.25,32.25){\line(1,0){2.5}}
\end{picture}
\begin{center}
\nopagebreak[4] Figure \theppp
\end{center}
\addtocounter{ppp}{1}	

\end{figure}

	We can now construct a copy $\Gamma_{i+1}$ of $\Gamma_i$ by adding $1$ to all superscripts of edges if $i+2\ne 0$ ($mod\; L$).
 This diagram over $H({\bf Y})$ can be inserted between two the hubs (with identification the $t$-band ${\cal Q}_{i+1}$ in
 $\Gamma_i$ and $\Gamma_{i+1}$).  Similarly we can connect
 the hubs by $t$-bands ${\cal Q}_0$,..., ${\cal Q}_L$, and the
 holes between them will be filled in with the diagrams $\Gamma_l$-s, but
 with one exception: to obtain a spherical diagram one should
 fill in the hole between the $t$-bands ${\cal Q}_0$ and ${\cal Q}_1$ with bases $t^{(0)}$ and $t^{(1)}$. Indeed this is possible
 by Lemma \ref{simul} since the computation $\cal C$ can be extended by Lemma \ref{exten} so that the whole between
${\cal Q}_0$ and ${\cal Q}_1$ is glued up with a quasi-trapezium
too. This complete the construction of the sherical diagram and
the proof of  Lemma \ref{Z2}.
\endproof

\begin{lm} \label{Z3} Let $w$ be an essential element represented by a cyclically minimal in rank $1/2$ word $W$ and $x$ be a $0$-element
such that $w^{-4}xw^4$ is a $0$-element $y$ in $H(\bf\Theta)$. Then $x\in {\bf 0}(w)$ and there is a $0$-element $c$ in $H({\bf \Theta})$
such that the product $wc^{-1}$ commutes with every element from ${\bf 0}(w)$.

\end{lm}
\proof We subdivide the proof in several steps.

\medskip

{\bf 0.} Note that if $w'=ywy^{-1}$ and $x'=yxy^{-1}$ for some
$0$-element $y$, then  the statement of the lemma holds for
the pair $(w,x)$ iff it holds for $(w',x')$ since ${\bf 0}(ywy^{-1})=y{\bf 0}(w)y^{-1}$.

\medskip

{\bf 1.} It follows that we may assume that the word $W$ starts
with some $q$-letter $q_0$. Let $\Delta$ be a minimal diagram over $H(\bf\Theta)$ with boundary ${\bf p}_1{\bf q}_1{\bf p}_2{\bf q}_2$, where
the section ${\bf p}_1$  has labels representing
$x$, $Lab({\bf q}_1)\equiv Lab({\bf q}_2)^{-1}\equiv W^4$, and $|{\bf p}_1|_q=|{\bf p}_2|_q=0$.

Since the word $W$ is cyclically minimal in $H(\bf\Theta)$, there is no maximal $q$-band in $\Delta$, which starts and ends on ${\bf q}_1$ (on ${\bf q}_2$). (Otherwise $W$ were conjugate of a
word of length $|W|_q-2$.)
Hence each of them starts on ${\bf q}_1$, ends on ${\bf q}_2$, and one
can enumerate them counting from ${\bf p}_1$:  ${\cal D}_1,\dots,{\cal D}_m$ ($m\ge 4$). So we have a diagram over $H(\theta,a)$ between
${\bf p}_1$ and ${\cal D}_1$, it can be assumed empty, and so $x$ is represented by the side label of ${\cal D}_1$. We may assume that the
history of ${\cal D}_1$ is non-trivial in $B(\bf\Theta)$, since otherwise $x=1$ in $H(\bf \Theta)$ by Lemmas \ref{thea} (a) and  \ref{tohist}, and so $x\in {\bf 0}(w)$.

\medskip

{\bf 2.} Let ${\bf y}_1{\bf x}_1{\bf y}_2{\bf x}_2$ be the boundary of ${\cal D}_1$ with
sides ${\bf y}_1={\bf p}_1$ and ${\bf y}_2$. Since the labels of ${\bf q}_1$ and ${\bf q}_2^{-1}$ start with the same letter, we have $Lab({\bf x}_2)\equiv Lab({\bf x}_1)^{-1}$. So by Lemma \ref{min} (b), one can construct a
$q$-band ${\cal D}^0_1$ with boundary ${\bf y}^0_1{\bf x}^0_1{\bf y}^0_2{\bf x}^0_2$, where $Lab({\bf x}^0_2)\equiv Lab({\bf x}^0_1)^{-1}$, $Lab({\bf y}_2^0)$ is equal
to $Lab({\bf y}_2)$ in $H(\theta,a)$,
$Lab({\bf y}_2^0)$ is
$vsv^{-1}$ for some word $v$, with cyclically
minimal $s$
in the group $H(\theta,a)$.
Similarly, we have $Lab({\bf y}_1^0)\equiv utu^{-1}$ with cyclically minimal $t$ by Lemma \ref{thea} (b).
Therefore we preserves the boundary label of $\Delta$, when
replacing the band ${\cal D}_1$ with the band ${\cal D}_1^0$
(and inserting an auxiliary diagrams over $H(\theta,a)$ to the top/bottom of ${\cal D}_1^0$). Let $E$ be the modified $\Delta$. We may assume that all $q$-bands of $E$ are
reduced, all the subdiagrams over $H(\theta,a)$ are minimal, and $E$ has no $q$-annuli by Lemma \ref{noq}.

Note that $|u|_q=0$, and  now we can replace the partition ${\bf p}_1{\bf q}_1{\bf p}_2{\bf q}_2$
of $\partial\Delta$ with the partition ${\bf y}({\bf z}_1{\bf q}_1){\bf p}_2({\bf q}_2{\bf z}_2)$,
where $Lab({\bf y})\equiv t$ and $Lab({\bf z}_1^{-1})\equiv Lab({\bf z}_2)\equiv u$. Furthermore,
adding edges to the boundary, we have a diagram $E'$ with boundary ${\bf y}({\bf z}_1{\bf q}_1{\bf z}_3)({\bf z}_3^{-1}{\bf p}_2{\bf z}_4)({\bf z}_4^{-1}{\bf q}_2{\bf z}_2)$,
where $Lab({\bf z}_3)\equiv Lab({\bf z}_4)\equiv u$. This partition corresponds
to simultaneous conjugation of the element $x$ and $w$ by the
word $u$ of $q$-length $0$. Therefore we may replace $(x,w)$
with the conjugate pair $(x',w')$, where $x'=t$ (see item $\bf 0$). Finaly,
the beginning ${\bf z}_1{\bf x}_1$ of the path ${\bf z}_1{\bf q}_1{\bf z}_3$ is homotopic in $E'$ to the path ${\bf l r}$, where ${\bf l}$ is a $q$-edge cutting the band ${\cal D}_1^0$ and ${\bf r}$ is a subpath of ${\bf y}_2^{-1}$.
Let $E''$ be the diagram, where the boundary subpath ${\bf z}_1{\bf x}_2$
is replaced by ${\bf lr}$, and similarly we replace the subpath
${\bf x}_2{\bf z}_2$. (See \setcounter{pdeleven}{\value{ppp}} Figure
\thepdeleven .)
The replacement of $E'$ with $E''$ does not change
the pair $(x', w')$, but in $E''$ the element $x'$ is represented by a cyclically minimal word and the word $W'$ representing $w'$ starts with
the $q$-letter $Lab({\bf l})$.

%\bigin{figure}
% This is a LaTeX picture output by TeXCAD.
% File name: [to2.pic].
% Version of TeXCAD: 4.3
% Reference / build: 30-Jun-2012 (rev. 105)
% For new versions, check: http://texcad.sf.net/
% Options on the following lines.
%\grade{\on}
%\emlines{\off}
%\epic{\off}
%\beziermacro{\on}
%\reduce{\on}
%\snapping{\off}
%\pvinsert{% Your \input, \def, etc. here}
%\quality{8.000}
%\graddiff{0.005}
%\snapasp{1}
%\zoom{4.0000}
\unitlength 1mm % = 2.845pt
\linethickness{0.4pt}
\ifx\plotpoint\undefined\newsavebox{\plotpoint}\fi % GNUPLOT compatibility
\begin{picture}(86.25,67.75)(-5,5)
\put(35,64.75){\line(0,-1){57.25}}
%\emline(35.25,64.75)(34.75,64.5)
\multiput(35.25,64.75)(-.0625,-.03125){8}{\line(-1,0){.0625}}
%\end
\put(35.25,64.75){\vector(1,0){47.25}}
%\vector(81.75,8.25)(35.5,7.75)
\put(35.5,7.75){\vector(-1,0){.07}}\multiput(81.75,8.25)(-3.0833333,-.0333333){15}{\line(-1,0){3.0833333}}
%\end
\thicklines
\put(35.25,48.25){\line(1,0){11.5}}
\put(35,24.5){\line(1,0){11.25}}
\put(46.25,64.5){\line(0,-1){16.25}}
\put(46.25,24.25){\line(0,-1){16.25}}
\thinlines
\put(46.25,47.75){\line(0,-1){23.25}}
\put(30.5,55.75){${\bf z}_1$}
\put(29,36.5){$\bf y$}
\put(30,16.25){${\bf z}_2$}
\put(39.25,67.75){${\bf x}_1$}
\put(40,5.25){${\bf x}_2$}
\put(40,44.25){$\bf l$}
\put(67.75,36.5){$E''$}
\put(86.25,64.75){${\bf q}_1$}
\put(85,8){${\bf q}_2$}
\put(47.75,55.5){$\bf r$}
%\emline(71.75,9.25)(70,8)
\multiput(71.75,9.25)(-.04605263,-.03289474){38}{\line(-1,0){.04605263}}
%\end
%\emline(70,8)(71.75,7.25)
\multiput(70,8)(.076087,-.0326087){23}{\line(1,0){.076087}}
%\end
%\dashline{1}(39.25,65)(35.25,61.75)
\multiput(39.18,64.93)(-.0408163,-.0331633){14}{\line(-1,0){.0408163}}
\multiput(38.037,64.001)(-.0408163,-.0331633){14}{\line(-1,0){.0408163}}
\multiput(36.894,63.073)(-.0408163,-.0331633){14}{\line(-1,0){.0408163}}
\multiput(35.751,62.144)(-.0408163,-.0331633){14}{\line(-1,0){.0408163}}
%\end
%\dashline{1}(41.75,65)(34.75,59.25)
\multiput(41.68,64.93)(-.0388889,-.0319444){18}{\line(-1,0){.0388889}}
\multiput(40.28,63.78)(-.0388889,-.0319444){18}{\line(-1,0){.0388889}}
\multiput(38.88,62.63)(-.0388889,-.0319444){18}{\line(-1,0){.0388889}}
\multiput(37.48,61.48)(-.0388889,-.0319444){18}{\line(-1,0){.0388889}}
\multiput(36.08,60.33)(-.0388889,-.0319444){18}{\line(-1,0){.0388889}}
%\end
%\dashline{1}(45.75,63.5)(35,53.75)
\multiput(45.68,63.43)(-.0358333,-.0325){20}{\line(-1,0){.0358333}}
\multiput(44.246,62.13)(-.0358333,-.0325){20}{\line(-1,0){.0358333}}
\multiput(42.813,60.83)(-.0358333,-.0325){20}{\line(-1,0){.0358333}}
\multiput(41.38,59.53)(-.0358333,-.0325){20}{\line(-1,0){.0358333}}
\multiput(39.946,58.23)(-.0358333,-.0325){20}{\line(-1,0){.0358333}}
\multiput(38.513,56.93)(-.0358333,-.0325){20}{\line(-1,0){.0358333}}
\multiput(37.08,55.63)(-.0358333,-.0325){20}{\line(-1,0){.0358333}}
\multiput(35.646,54.33)(-.0358333,-.0325){20}{\line(-1,0){.0358333}}
%\end
%\dashline{1}(46.75,60.5)(47,61)
\multiput(46.68,60.43)(.03125,.0625){4}{\line(0,1){.0625}}
%\end
%\dashline{1}(46,62)(35.25,52)
\multiput(45.93,61.93)(-.0353618,-.0328947){19}{\line(-1,0){.0353618}}
\multiput(44.586,60.68)(-.0353618,-.0328947){19}{\line(-1,0){.0353618}}
\multiput(43.242,59.43)(-.0353618,-.0328947){19}{\line(-1,0){.0353618}}
\multiput(41.898,58.18)(-.0353618,-.0328947){19}{\line(-1,0){.0353618}}
\multiput(40.555,56.93)(-.0353618,-.0328947){19}{\line(-1,0){.0353618}}
\multiput(39.211,55.68)(-.0353618,-.0328947){19}{\line(-1,0){.0353618}}
\multiput(37.867,54.43)(-.0353618,-.0328947){19}{\line(-1,0){.0353618}}
\multiput(36.523,53.18)(-.0353618,-.0328947){19}{\line(-1,0){.0353618}}
%\end
%\dashline{1}(46.25,59.25)(46,59.75)
\multiput(46.18,59.18)(-.03125,.0625){4}{\line(0,1){.0625}}
%\end
%\dashline{1}(45.75,60)(34.75,49.75)
\multiput(45.68,59.93)(-.0359477,-.0334967){18}{\line(-1,0){.0359477}}
\multiput(44.386,58.724)(-.0359477,-.0334967){18}{\line(-1,0){.0359477}}
\multiput(43.091,57.518)(-.0359477,-.0334967){18}{\line(-1,0){.0359477}}
\multiput(41.797,56.312)(-.0359477,-.0334967){18}{\line(-1,0){.0359477}}
\multiput(40.503,55.106)(-.0359477,-.0334967){18}{\line(-1,0){.0359477}}
\multiput(39.209,53.9)(-.0359477,-.0334967){18}{\line(-1,0){.0359477}}
\multiput(37.915,52.694)(-.0359477,-.0334967){18}{\line(-1,0){.0359477}}
\multiput(36.621,51.489)(-.0359477,-.0334967){18}{\line(-1,0){.0359477}}
\multiput(35.327,50.283)(-.0359477,-.0334967){18}{\line(-1,0){.0359477}}
%\end
%\dashline{1}(46.5,58.75)(35.5,49)
\multiput(46.43,58.68)(-.0361842,-.0320724){19}{\line(-1,0){.0361842}}
\multiput(45.055,57.461)(-.0361842,-.0320724){19}{\line(-1,0){.0361842}}
\multiput(43.68,56.242)(-.0361842,-.0320724){19}{\line(-1,0){.0361842}}
\multiput(42.305,55.023)(-.0361842,-.0320724){19}{\line(-1,0){.0361842}}
\multiput(40.93,53.805)(-.0361842,-.0320724){19}{\line(-1,0){.0361842}}
\multiput(39.555,52.586)(-.0361842,-.0320724){19}{\line(-1,0){.0361842}}
\multiput(38.18,51.367)(-.0361842,-.0320724){19}{\line(-1,0){.0361842}}
\multiput(36.805,50.148)(-.0361842,-.0320724){19}{\line(-1,0){.0361842}}
%\end
%\dashline{1}(46,56.75)(37.25,48.75)
\multiput(45.93,56.68)(-.0354251,-.0323887){19}{\line(-1,0){.0354251}}
\multiput(44.584,55.449)(-.0354251,-.0323887){19}{\line(-1,0){.0354251}}
\multiput(43.237,54.218)(-.0354251,-.0323887){19}{\line(-1,0){.0354251}}
\multiput(41.891,52.987)(-.0354251,-.0323887){19}{\line(-1,0){.0354251}}
\multiput(40.545,51.757)(-.0354251,-.0323887){19}{\line(-1,0){.0354251}}
\multiput(39.199,50.526)(-.0354251,-.0323887){19}{\line(-1,0){.0354251}}
\multiput(37.853,49.295)(-.0354251,-.0323887){19}{\line(-1,0){.0354251}}
%\end
%\dashline{1}(46,55)(38.75,48)
\multiput(45.93,54.93)(-.0335648,-.0324074){18}{\line(-1,0){.0335648}}
\multiput(44.721,53.763)(-.0335648,-.0324074){18}{\line(-1,0){.0335648}}
\multiput(43.513,52.596)(-.0335648,-.0324074){18}{\line(-1,0){.0335648}}
\multiput(42.305,51.43)(-.0335648,-.0324074){18}{\line(-1,0){.0335648}}
\multiput(41.096,50.263)(-.0335648,-.0324074){18}{\line(-1,0){.0335648}}
\multiput(39.888,49.096)(-.0335648,-.0324074){18}{\line(-1,0){.0335648}}
%\end
%\dashline{1}(46.25,53.25)(41,48.25)
\multiput(46.18,53.18)(-.0343137,-.0326797){17}{\line(-1,0){.0343137}}
\multiput(45.013,52.069)(-.0343137,-.0326797){17}{\line(-1,0){.0343137}}
\multiput(43.846,50.957)(-.0343137,-.0326797){17}{\line(-1,0){.0343137}}
\multiput(42.68,49.846)(-.0343137,-.0326797){17}{\line(-1,0){.0343137}}
\multiput(41.513,48.735)(-.0343137,-.0326797){17}{\line(-1,0){.0343137}}
%\end
%\dashline{1}(46.25,51.25)(42.75,48)
\multiput(46.18,51.18)(-.0343137,-.0318627){17}{\line(-1,0){.0343137}}
\multiput(45.013,50.096)(-.0343137,-.0318627){17}{\line(-1,0){.0343137}}
\multiput(43.846,49.013)(-.0343137,-.0318627){17}{\line(-1,0){.0343137}}
%\end
%\dashline{1}(46.25,49.5)(44.5,48.25)
\multiput(46.18,49.43)(-.0448718,-.0320513){13}{\line(-1,0){.0448718}}
\multiput(45.013,48.596)(-.0448718,-.0320513){13}{\line(-1,0){.0448718}}
%\end
%\dashline{1}(42.75,65.25)(43,65.25)
\put(42.68,65.18){\line(1,0){.125}}
%\end
%\dashline{1}(43,65.25)(43.25,65.25)
\put(42.93,65.18){\line(1,0){.125}}
%\end
%\dashline{1}(43.25,65.25)(43.75,65.25)
\put(43.18,65.18){\line(1,0){.25}}
%\end
%\dashline{1}(43.75,65)(42.75,64.75)
\multiput(43.68,64.93)(-.125,-.03125){4}{\line(-1,0){.125}}
%\end
%\dashline{1}(43.5,65.25)(35.25,58)
\multiput(43.43,65.18)(-.0381944,-.0335648){18}{\line(-1,0){.0381944}}
\multiput(42.055,63.971)(-.0381944,-.0335648){18}{\line(-1,0){.0381944}}
\multiput(40.68,62.763)(-.0381944,-.0335648){18}{\line(-1,0){.0381944}}
\multiput(39.305,61.555)(-.0381944,-.0335648){18}{\line(-1,0){.0381944}}
\multiput(37.93,60.346)(-.0381944,-.0335648){18}{\line(-1,0){.0381944}}
\multiput(36.555,59.138)(-.0381944,-.0335648){18}{\line(-1,0){.0381944}}
%\end
%\dashline{1}(44.75,64.5)(35,56.25)
\multiput(44.68,64.43)(-.0386905,-.0327381){18}{\line(-1,0){.0386905}}
\multiput(43.287,63.251)(-.0386905,-.0327381){18}{\line(-1,0){.0386905}}
\multiput(41.894,62.073)(-.0386905,-.0327381){18}{\line(-1,0){.0386905}}
\multiput(40.501,60.894)(-.0386905,-.0327381){18}{\line(-1,0){.0386905}}
\multiput(39.108,59.715)(-.0386905,-.0327381){18}{\line(-1,0){.0386905}}
\multiput(37.715,58.537)(-.0386905,-.0327381){18}{\line(-1,0){.0386905}}
\multiput(36.323,57.358)(-.0386905,-.0327381){18}{\line(-1,0){.0386905}}
%\end
%\dashline{1}(37.5,24.5)(35.25,22.75)
\multiput(37.43,24.43)(-.0432692,-.0336538){13}{\line(-1,0){.0432692}}
\multiput(36.305,23.555)(-.0432692,-.0336538){13}{\line(-1,0){.0432692}}
%\end
%\dashline{1}(39,24)(35.25,20.75)
\multiput(38.93,23.93)(-.0367647,-.0318627){17}{\line(-1,0){.0367647}}
\multiput(37.68,22.846)(-.0367647,-.0318627){17}{\line(-1,0){.0367647}}
\multiput(36.43,21.763)(-.0367647,-.0318627){17}{\line(-1,0){.0367647}}
%\end
%\dashline{1}(41.25,24.25)(34.75,18.25)
\multiput(41.18,24.18)(-.0361111,-.0333333){18}{\line(-1,0){.0361111}}
\multiput(39.88,22.98)(-.0361111,-.0333333){18}{\line(-1,0){.0361111}}
\multiput(38.58,21.78)(-.0361111,-.0333333){18}{\line(-1,0){.0361111}}
\multiput(37.28,20.58)(-.0361111,-.0333333){18}{\line(-1,0){.0361111}}
\multiput(35.98,19.38)(-.0361111,-.0333333){18}{\line(-1,0){.0361111}}
%\end
%\dashline{1}(43.25,24.75)(35,16.25)
\multiput(43.18,24.68)(-.0334008,-.034413){19}{\line(0,-1){.034413}}
\multiput(41.91,23.372)(-.0334008,-.034413){19}{\line(0,-1){.034413}}
\multiput(40.641,22.064)(-.0334008,-.034413){19}{\line(0,-1){.034413}}
\multiput(39.372,20.757)(-.0334008,-.034413){19}{\line(0,-1){.034413}}
\multiput(38.103,19.449)(-.0334008,-.034413){19}{\line(0,-1){.034413}}
\multiput(36.834,18.141)(-.0334008,-.034413){19}{\line(0,-1){.034413}}
\multiput(35.564,16.834)(-.0334008,-.034413){19}{\line(0,-1){.034413}}
%\end
%\dashline{1}(45,24.5)(35,14)
\multiput(44.93,24.43)(-.0328947,-.0345395){19}{\line(0,-1){.0345395}}
\multiput(43.68,23.117)(-.0328947,-.0345395){19}{\line(0,-1){.0345395}}
\multiput(42.43,21.805)(-.0328947,-.0345395){19}{\line(0,-1){.0345395}}
\multiput(41.18,20.492)(-.0328947,-.0345395){19}{\line(0,-1){.0345395}}
\multiput(39.93,19.18)(-.0328947,-.0345395){19}{\line(0,-1){.0345395}}
\multiput(38.68,17.867)(-.0328947,-.0345395){19}{\line(0,-1){.0345395}}
\multiput(37.43,16.555)(-.0328947,-.0345395){19}{\line(0,-1){.0345395}}
\multiput(36.18,15.242)(-.0328947,-.0345395){19}{\line(0,-1){.0345395}}
%\end
%\dashline{1}(45.5,23.5)(34.75,11.75)
\multiput(45.43,23.43)(-.0332817,-.0363777){19}{\line(0,-1){.0363777}}
\multiput(44.165,22.047)(-.0332817,-.0363777){19}{\line(0,-1){.0363777}}
\multiput(42.9,20.665)(-.0332817,-.0363777){19}{\line(0,-1){.0363777}}
\multiput(41.636,19.283)(-.0332817,-.0363777){19}{\line(0,-1){.0363777}}
\multiput(40.371,17.9)(-.0332817,-.0363777){19}{\line(0,-1){.0363777}}
\multiput(39.106,16.518)(-.0332817,-.0363777){19}{\line(0,-1){.0363777}}
\multiput(37.841,15.136)(-.0332817,-.0363777){19}{\line(0,-1){.0363777}}
\multiput(36.577,13.753)(-.0332817,-.0363777){19}{\line(0,-1){.0363777}}
\multiput(35.312,12.371)(-.0332817,-.0363777){19}{\line(0,-1){.0363777}}
%\end
%\dashline{1}(45.75,21.75)(34.25,9.5)
\multiput(45.68,21.68)(-.0336257,-.0358187){19}{\line(0,-1){.0358187}}
\multiput(44.402,20.319)(-.0336257,-.0358187){19}{\line(0,-1){.0358187}}
\multiput(43.124,18.957)(-.0336257,-.0358187){19}{\line(0,-1){.0358187}}
\multiput(41.846,17.596)(-.0336257,-.0358187){19}{\line(0,-1){.0358187}}
\multiput(40.569,16.235)(-.0336257,-.0358187){19}{\line(0,-1){.0358187}}
\multiput(39.291,14.874)(-.0336257,-.0358187){19}{\line(0,-1){.0358187}}
\multiput(38.013,13.513)(-.0336257,-.0358187){19}{\line(0,-1){.0358187}}
\multiput(36.735,12.152)(-.0336257,-.0358187){19}{\line(0,-1){.0358187}}
\multiput(35.457,10.791)(-.0336257,-.0358187){19}{\line(0,-1){.0358187}}
%\end
%\dashline{1}(46,20.25)(35,8.75)
\multiput(45.93,20.18)(-.0321637,-.0336257){19}{\line(0,-1){.0336257}}
\multiput(44.707,18.902)(-.0321637,-.0336257){19}{\line(0,-1){.0336257}}
\multiput(43.485,17.624)(-.0321637,-.0336257){19}{\line(0,-1){.0336257}}
\multiput(42.263,16.346)(-.0321637,-.0336257){19}{\line(0,-1){.0336257}}
\multiput(41.041,15.069)(-.0321637,-.0336257){19}{\line(0,-1){.0336257}}
\multiput(39.819,13.791)(-.0321637,-.0336257){19}{\line(0,-1){.0336257}}
\multiput(38.596,12.513)(-.0321637,-.0336257){19}{\line(0,-1){.0336257}}
\multiput(37.374,11.235)(-.0321637,-.0336257){19}{\line(0,-1){.0336257}}
\multiput(36.152,9.957)(-.0321637,-.0336257){19}{\line(0,-1){.0336257}}
%\end
%\dashline{1}(46,18)(36.75,8.25)
\multiput(45.93,17.93)(-.0324561,-.0342105){19}{\line(0,-1){.0342105}}
\multiput(44.696,16.63)(-.0324561,-.0342105){19}{\line(0,-1){.0342105}}
\multiput(43.463,15.33)(-.0324561,-.0342105){19}{\line(0,-1){.0342105}}
\multiput(42.23,14.03)(-.0324561,-.0342105){19}{\line(0,-1){.0342105}}
\multiput(40.996,12.73)(-.0324561,-.0342105){19}{\line(0,-1){.0342105}}
\multiput(39.763,11.43)(-.0324561,-.0342105){19}{\line(0,-1){.0342105}}
\multiput(38.53,10.13)(-.0324561,-.0342105){19}{\line(0,-1){.0342105}}
\multiput(37.296,8.83)(-.0324561,-.0342105){19}{\line(0,-1){.0342105}}
%\end
%\dashline{1}(46,16)(38.5,8.25)
\multiput(45.93,15.93)(-.0328947,-.0339912){19}{\line(0,-1){.0339912}}
\multiput(44.68,14.638)(-.0328947,-.0339912){19}{\line(0,-1){.0339912}}
\multiput(43.43,13.346)(-.0328947,-.0339912){19}{\line(0,-1){.0339912}}
\multiput(42.18,12.055)(-.0328947,-.0339912){19}{\line(0,-1){.0339912}}
\multiput(40.93,10.763)(-.0328947,-.0339912){19}{\line(0,-1){.0339912}}
\multiput(39.68,9.471)(-.0328947,-.0339912){19}{\line(0,-1){.0339912}}
%\end
%\dashline{1}(46,13.75)(40.25,8)
\multiput(45.93,13.68)(-.0336257,-.0336257){19}{\line(0,-1){.0336257}}
\multiput(44.652,12.402)(-.0336257,-.0336257){19}{\line(0,-1){.0336257}}
\multiput(43.374,11.124)(-.0336257,-.0336257){19}{\line(0,-1){.0336257}}
\multiput(42.096,9.846)(-.0336257,-.0336257){19}{\line(0,-1){.0336257}}
\multiput(40.819,8.569)(-.0336257,-.0336257){19}{\line(0,-1){.0336257}}
%\end
%\dashline{1}(46,11.75)(42.25,8.25)
\multiput(45.93,11.68)(-.0347222,-.0324074){18}{\line(-1,0){.0347222}}
\multiput(44.68,10.513)(-.0347222,-.0324074){18}{\line(-1,0){.0347222}}
\multiput(43.43,9.346)(-.0347222,-.0324074){18}{\line(-1,0){.0347222}}
%\end
%\dashline{1}(46.25,11)(44,7.75)
\multiput(46.18,10.93)(-.0321429,-.0464286){14}{\line(0,-1){.0464286}}
\multiput(45.28,9.63)(-.0321429,-.0464286){14}{\line(0,-1){.0464286}}
\multiput(44.38,8.33)(-.0321429,-.0464286){14}{\line(0,-1){.0464286}}
%\end
\end{picture}
\begin{center}
\nopagebreak[4] Figure \theppp
\end{center}
\addtocounter{ppp}{1}	
%\end{figure} 

\begin{remark} \label{redu} The replacement of the pair
$(x,w)$ with the pair $(x',w')$ (and the diagram $\Delta$ with $E''$) is
possible if the history of ${\cal D}_1$ is reduced but not necessarily  cyclically reduced.
Let us call such a modification of the
pair $(x,w)$ and the diagram $\Delta$
the \index[g]{cyclic reduction of $q$-band} cyclic reduction of the band ${\cal D}_1$ in $\Delta$.
\end{remark}

The advantage of the obtained diagram is that the band ${\cal D}_1$ is replaced with the band having cyclically minimal in $H(\theta, a)$ labels of the top and the bottom and the modified paths ${\bf q}_1$ and ${\bf q}_2^{-1}$ start with $q$-edges. Now we may re-establish the original notation $x, w, {\bf p}_1,\dots$ but assume
that the first maximal $q$-band ${\cal D}_1$ of $\Delta$ has
labels of the top and bottom {\bf cyclically minimal} in the group $H(\theta,a)$. 

If the $\theta$-letters of $Lab({\bf p}_1)$ belong to a subalphabet
${\bf Q}_j$, then the label of the opposite side of the band ${\cal D}_1$ has $\theta$-letters from ${\bf Q}_{j'}$, where $j'=j\pm 1$ $ (mod N)$ by the definition of $(\theta, q)$-relations.

\medskip

{\bf 3.} Let ${\bf \bar q}_1$ (resp., ${\bf \bar q}_2^{-1}$) be the subpath of ${\bf q}_1$ (of ${\bf q}_2^{-1}$) starting with the vertex $({\bf q}_1)_-$ (with $({\bf q}_2)_+$) and labeled by $Wq_0$, where $q_0$ is the first letter of $W$. All $q$-edges of $\bar q_1$ are connected with $q$-edges of the subpath ${\bf \bar q}_2$ by
the $q$-bands ${\cal D}_1,\dots, {\cal D}_s$, where $s-1=|W|_q$. Let $\Gamma$
be the subdiagram of $\Delta$ bounded by ${\bf \bar q}_1$, ${\bf \bar q}_2$, ${\cal D}_1$ and ${\cal D}_s$
(${\cal D}_1$ and ${\cal D}_s$ are included in $\Gamma$).

 Since $Lab({\bf \bar q}_1)\equiv Lab({\bf \bar q}_2)^{-1}$, identifying ${\bf \bar q}_1$ and ${\bf \bar q}_2^{-1}$ in $\partial\Gamma$,
 except for the last edges of these paths labeled by $q_0$,
 one obtains an annular diagram $\Gamma(0)$. The $q$-bands
${\cal D}_1,\dots, {\cal D}_{s-1}$ becomes $q$-annuli in $\Gamma_0$.
We denote by ${\cal Q}_1,\dots, {\cal Q}_{s-1}$ the reduced forms of
these annuli, where ${\cal Q}_1={\cal D}_1$ by item ${\bf 3}$. We do not
change ${\cal D}_s$, and so there is a cutting path ${\bf p}$ connecting two vertices on ${\cal Q}_1$ and ${\cal D}_s$, whose label is equal to $Wq_0$ in $H(\bf\Theta)$.

The  {\it chamber} $\Gamma_1$ bounded by ${\cal Q}_1$
and ${\cal Q}_2$ can be made  minimal annular diagram
over $H(\theta,a)$ according to Lemma \ref{notha} (2). Since the side labels of ${\cal Q}_1$ are cyclically  minimal in $H(\theta,a)$, there are no positive (over the presentation of $H(\theta,a)$) cells in this chamber with contiguity degree $\ge \bar\alpha$ to the side of ${\cal Q}_1$ by Lemma \ref{bara}. Using Lemma \ref{through} (as in Lemma \ref{min}), one can reduce
the type of $\Gamma_1$ and  finally transform it in an annular diagram over
$H(a)$ (at the expense of the modification of the next chamber $\Gamma_2$). This transformation does
not change the band ${\cal D}_1$, and
according to Lemma \ref{through} (d), it preserves the label of a simple path homotopic to ${\bf p}$
modulo the relations of $H(\Theta)$.

Recall that by item {\bf 1}, the history of ${\cal Q}_1$ is nontrivial in $B({\Theta}$). Therefore
the retraction from Lemma \ref{thea} (a) shows that the side labels of ${\cal Q}_1$ are not conjugate
in $H(\theta,a)$ to a word of $\theta$-length $0$. Hence applying Lemma \ref{notha} (2), we may assume
that $\Gamma_1$ is a minimal diagram over $H(a)$.

Since both ${\cal Q}_1$ and ${\cal Q}_2$ are reduced annuli, every maximal $\theta$-band of $\Gamma_1$ must start on the boundary of ${\cal Q}_1$ and end on the boundary of ${\cal Q}_2$. Since the boundary labels
of ${\cal Q}_1$ are cyclically minimal in $H(\theta,a)$ the same is true for the boundary labels of ${\cal Q}_2$ by Lemma \ref{thea}(b). Therefore the next chamber $\Gamma_2$ between ${\cal Q}_2$ and ${\cal Q}_3$ can be transformed in a minimal annular diagram over $H(a)$ too. The iteration of this procedure
makes the boundary labels of all $q$-annuli of $\Gamma(0)$ cyclically minimal over $H(\theta,a)$ and all the chambers $\Gamma_i$ minimal annular diagrams
over the group $H(a)$.
According to Lemma \ref{through} (d), the procedure preserves, modulo the relations of $H(\bf\Theta)$, the label of a cutting simple path ${\bf p}'$ connecting ${\bf p}_-$ and ${\bf p}_+$ in the modified annular diagram $\Gamma(0)$,
i.e. the word $Lab({\bf p}')$ represents the
element $wq_0$, and $|{\bf p}'|_q=|{\bf p}|_q$.

\medskip

{\bf 4.}  By Lemma \ref{NoAnnul}, no $\theta$-band crosses a $q$-annulus ${\cal Q}_i$ twice in $\Gamma(0)$
 ($i=1,\dots, s-1$), and so 
the maximal $\theta$-bands starting on ${\cal Q}_1$ end on  ${\cal Q}_s$. Therefore one can connect the vertex ${\bf p}'_-=({\bf p}_1)_-$ with the vertex on the band
${\cal Q}_s$ along a $\theta$-band of $\Gamma(0)$. So
${\bf p}'$ is homotopic in $\Gamma(0)$ to a path ${\bf p}''$, which is a product ${\bf p}(1){\bf p}(2)$ of a side ${\bf p}(1)$ of a maximal $\theta$-band and a path ${\bf p}(2)$ whose
label $S$ is freely equal to a subword of a power of the (outer) side label of the $q$-band ${\cal Q}_s$.
(We use that the two boundary $q$-edges of ${\cal Q}_s$ have the same label $q_0$.)

Cutting the diagram $\Gamma(0)$ along the path ${\bf p}''$, we obtain a modification $\Gamma'$ of $\Gamma$ with the additional properties:

(1) $\Gamma'$ has no
positive cells over $H(\theta,a)$, i.e.
$\Gamma'$ is a minimal diagram over $H(\bf Y)$,

(2) Every maximal $\theta$-band of $\Gamma'$ starting on ${\bf p}_1$ crosses
each of the $q$-bands ${\cal D}_1,\dots, {\cal D}_{s}$ and

(3) $Wq_0=W(1)W(2)$, where, $|W(1)|_q=s$, $W(1)$ has
no $\theta$-letters, and $W(2)$ is a product of the words $x(i)$-s defined in Lemma \ref{thea} for the set $\bf\Theta^+_{j'}$ from the end of item {\bf 2}.

Note that the factorization of $Wq_0$ with these properties is unique, that is the factor $W(2)$ is uniquely defined modulo the relations of $H(\bf Y)$ by Lemma \ref{tohist},
since the words $Wq_0$ and $W(2)$ have
the same image under the homomorpism
to $F(\bf\Theta)$ defined in Lemma \ref{tohist}.

Simplifying the notation, we will assume that the diagram $\Gamma$ enjoys the properties (1) - (3) itself.

\medskip

{\bf 5.} Since all the $\theta$-bands
starting on ${\bf p}_1$ cross ${\cal D}_s$ in $\Gamma$ and the band ${\cal D}_s$ is reduced we conclude that ${\cal D}_s$ is a reduced form of a product ${\cal F}={\cal F}_1{\cal F}_2{\cal F}_3$, where the band ${\cal F}_2$ is a copy of ${\cal D}_1$ and the history of the band ${\cal F}_3$ is inverse to the history of  ${\cal F}_1$, because the side labels of these two bands are $W(2)^{\pm 1}$.
Futhermore, one can construct a new diagram $E$ as follows.

Take the subdiagram $E_1$ formed by all maximal $\theta$-bands of $\Gamma$ starting with ${\bf p}_1$ and by $a$-cells between them (if any).
The  top/bottom
label of $E_1$ are equal to $W(1)$, i.e. they
are minimal words in rank $1/2$. Hence they
cannot contain a subword of $q$-length $2$ equal in $H({\bf Y})$ to a word of $q$-length $0$. It follows
by induction that all the $\theta$-bands of $E_1$ are regular, and therefore $E_1$ is a quasi-trapezium.
(See \setcounter{pdeleven}{\value{ppp}} Figure
\thepdeleven .)

% This is a LaTeX picture output by TeXCAD.
% File name: [to5.pic].
% Version of TeXCAD: 4.3
% Reference / build: 30-Jun-2012 (rev. 105)
% For new versions, check: http://texcad.sf.net/
% Options on the following lines.
%\grade{\on}
%\emlines{\off}
%\epic{\off}
%\beziermacro{\on}
%\reduce{\on}
%\snapping{\off}
%\pvinsert{% Your \input, \def, etc. here}
%\quality{8.000}
%\graddiff{0.005}
%\snapasp{1}
%\zoom{4.0000}
\unitlength 1mm % = 2.845pt
\linethickness{0.4pt}
\ifx\plotpoint\undefined\newsavebox{\plotpoint}\fi % GNUPLOT compatibility
\begin{picture}(135,71.5)(10,5)
\put(43.25,58.25){\line(0,-1){27.5}}
\put(43.75,57.5){\line(1,0){19.75}}
\put(43.5,31){\line(1,0){19}}
\put(50.75,57.5){\line(0,-1){26.25}}
%\dashline{1}(66.5,57.25)(89.5,57.5)
\put(66.43,57.18){\line(1,0){.9583}}
\put(68.346,57.201){\line(1,0){.9583}}
\put(70.263,57.221){\line(1,0){.9583}}
\put(72.18,57.242){\line(1,0){.9583}}
\put(74.096,57.263){\line(1,0){.9583}}
\put(76.013,57.284){\line(1,0){.9583}}
\put(77.93,57.305){\line(1,0){.9583}}
\put(79.846,57.326){\line(1,0){.9583}}
\put(81.763,57.346){\line(1,0){.9583}}
\put(83.68,57.367){\line(1,0){.9583}}
\put(85.596,57.388){\line(1,0){.9583}}
\put(87.513,57.409){\line(1,0){.9583}}
%\end
%\dashline{1}(65.25,30.75)(88.75,30.75)
\put(65.18,30.68){\line(1,0){.9792}}
\put(67.138,30.68){\line(1,0){.9792}}
\put(69.096,30.68){\line(1,0){.9792}}
\put(71.055,30.68){\line(1,0){.9792}}
\put(73.013,30.68){\line(1,0){.9792}}
\put(74.971,30.68){\line(1,0){.9792}}
\put(76.93,30.68){\line(1,0){.9792}}
\put(78.888,30.68){\line(1,0){.9792}}
\put(80.846,30.68){\line(1,0){.9792}}
\put(82.805,30.68){\line(1,0){.9792}}
\put(84.763,30.68){\line(1,0){.9792}}
\put(86.721,30.68){\line(1,0){.9792}}
%\end
\put(90.25,57){\line(1,0){41.5}}
\put(89.25,30.5){\line(1,0){43.25}}
%\emline(130.75,70.25)(131.25,18.25)
\multiput(130.75,70.25)(.0333333,-3.4666667){15}{\line(0,-1){3.4666667}}
%\end
\put(125,70){\line(0,-1){51.75}}
\put(124.75,69.75){\line(1,0){6.5}}
\put(125,19.25){\line(0,1){.75}}
\put(125.25,18.75){\line(1,0){5.75}}
\put(98.75,57.25){\line(0,-1){28.25}}
\put(106,57.25){\line(0,-1){26.75}}
\put(43,52.25){\line(1,0){19.75}}
\put(43.5,42.75){\line(1,0){19}}
\put(43.25,47.25){\line(1,0){18.75}}
\put(43.5,37){\line(1,0){18.5}}
\put(90.25,51){\line(1,0){40.75}}
\put(90.5,41.5){\line(1,0){40.75}}
\put(90.5,35.75){\line(1,0){40.5}}
\put(111,60.25){\line(1,0){5.75}}
%\emline(116.75,60.25)(120.75,63.5)
\multiput(116.75,60.25)(.041237113,.033505155){97}{\line(1,0){.041237113}}
%\end
%\vector(120.25,63.75)(121,68.75)
\put(121,68.75){\vector(1,4){.07}}\multiput(120.25,63.75)(.0326087,.2173913){23}{\line(0,1){.2173913}}
%\end
\put(121,20.25){\line(0,1){3.25}}
%\emline(121,23.5)(118.25,26)
\multiput(121,23.5)(-.03666667,.03333333){75}{\line(-1,0){.03666667}}
%\end
%\vector(119.25,25.25)(114.25,26)
\put(114.25,26){\vector(-4,1){.07}}\multiput(119.25,25.25)(-.2173913,.0326087){23}{\line(-1,0){.2173913}}
%\end
\put(106,59.75){${\bf q'}_1$}
\put(108.75,26.25){${\bf q'}_2$}
\put(38,45.5){${\bf p}_1$}
\put(46.5,26.25){${\cal D}_1$}
\put(99.25,26){${\cal D}_{s-1}$}
\put(125.25,22){${\cal F}_1$}
\put(125,65.25){${\cal F}_3$}
\put(126.5,45.75){${\cal F}_2$}
\put(133.5,40){${\bf p'}_2$}
\put(127.75,71.5){$q_0$}
\put(127.5,16.75){$q_0$}
%\dashline{1}(43.75,52.75)(50,57.5)
\multiput(43.68,52.68)(.0434028,.0329861){16}{\line(1,0){.0434028}}
\multiput(45.069,53.735)(.0434028,.0329861){16}{\line(1,0){.0434028}}
\multiput(46.457,54.791)(.0434028,.0329861){16}{\line(1,0){.0434028}}
\multiput(47.846,55.846)(.0434028,.0329861){16}{\line(1,0){.0434028}}
\multiput(49.235,56.902)(.0434028,.0329861){16}{\line(1,0){.0434028}}
%\end
%\dashline{1}(43.75,49)(50.75,54)
\multiput(43.68,48.93)(.0466667,.0333333){15}{\line(1,0){.0466667}}
\multiput(45.08,49.93)(.0466667,.0333333){15}{\line(1,0){.0466667}}
\multiput(46.48,50.93)(.0466667,.0333333){15}{\line(1,0){.0466667}}
\multiput(47.88,51.93)(.0466667,.0333333){15}{\line(1,0){.0466667}}
\multiput(49.28,52.93)(.0466667,.0333333){15}{\line(1,0){.0466667}}
%\end
%\dashline{1}(43.5,44.5)(51,50.75)
\multiput(43.43,44.43)(.040107,.0334225){17}{\line(1,0){.040107}}
\multiput(44.793,45.566)(.040107,.0334225){17}{\line(1,0){.040107}}
\multiput(46.157,46.702)(.040107,.0334225){17}{\line(1,0){.040107}}
\multiput(47.521,47.839)(.040107,.0334225){17}{\line(1,0){.040107}}
\multiput(48.884,48.975)(.040107,.0334225){17}{\line(1,0){.040107}}
\multiput(50.248,50.112)(.040107,.0334225){17}{\line(1,0){.040107}}
%\end
%\dashline{1}(43.75,40)(43.5,40)
\put(43.68,39.93){\line(-1,0){.125}}
%\end
%\dashline{1}(43.5,41.25)(50.5,46.75)
\multiput(43.43,41.18)(.0411765,.0323529){17}{\line(1,0){.0411765}}
\multiput(44.83,42.28)(.0411765,.0323529){17}{\line(1,0){.0411765}}
\multiput(46.23,43.38)(.0411765,.0323529){17}{\line(1,0){.0411765}}
\multiput(47.63,44.48)(.0411765,.0323529){17}{\line(1,0){.0411765}}
\multiput(49.03,45.58)(.0411765,.0323529){17}{\line(1,0){.0411765}}
%\end
%\dashline{1}(43.25,38.25)(50.5,43.75)
\multiput(43.18,38.18)(.0426471,.0323529){17}{\line(1,0){.0426471}}
\multiput(44.63,39.28)(.0426471,.0323529){17}{\line(1,0){.0426471}}
\multiput(46.08,40.38)(.0426471,.0323529){17}{\line(1,0){.0426471}}
\multiput(47.53,41.48)(.0426471,.0323529){17}{\line(1,0){.0426471}}
\multiput(48.98,42.58)(.0426471,.0323529){17}{\line(1,0){.0426471}}
%\end
%\dashline{1}(43.25,34.5)(50.75,39.75)
\multiput(43.18,34.43)(.046875,.0328125){16}{\line(1,0){.046875}}
\multiput(44.68,35.48)(.046875,.0328125){16}{\line(1,0){.046875}}
\multiput(46.18,36.53)(.046875,.0328125){16}{\line(1,0){.046875}}
\multiput(47.68,37.58)(.046875,.0328125){16}{\line(1,0){.046875}}
\multiput(49.18,38.63)(.046875,.0328125){16}{\line(1,0){.046875}}
%\end
%\dashline{1}(43.75,31.5)(50.5,36.75)
\multiput(43.68,31.43)(.0421875,.0328125){16}{\line(1,0){.0421875}}
\multiput(45.03,32.48)(.0421875,.0328125){16}{\line(1,0){.0421875}}
\multiput(46.38,33.53)(.0421875,.0328125){16}{\line(1,0){.0421875}}
\multiput(47.73,34.58)(.0421875,.0328125){16}{\line(1,0){.0421875}}
\multiput(49.08,35.63)(.0421875,.0328125){16}{\line(1,0){.0421875}}
%\end
%\dashline{1}(48.25,31)(47.5,30.5)
\multiput(48.18,30.93)(-.046875,-.03125){8}{\line(-1,0){.046875}}
%\end
%\dashline{1}(47.5,30.5)(48,30)
\multiput(47.43,30.43)(.03125,-.03125){8}{\line(0,-1){.03125}}
%\end
%\dashline{1}(47.5,30.25)(47,30.25)
\put(47.43,30.18){\line(-1,0){.25}}
%\end
%\dashline{1}(46.75,31)(50.5,33.25)
\multiput(46.68,30.93)(.0535714,.0321429){14}{\line(1,0){.0535714}}
\multiput(48.18,31.83)(.0535714,.0321429){14}{\line(1,0){.0535714}}
\multiput(49.68,32.73)(.0535714,.0321429){14}{\line(1,0){.0535714}}
%\end
%\dashline{1}(99.25,53.75)(103.5,57)
\multiput(99.18,53.68)(.0433673,.0331633){14}{\line(1,0){.0433673}}
\multiput(100.394,54.608)(.0433673,.0331633){14}{\line(1,0){.0433673}}
\multiput(101.608,55.537)(.0433673,.0331633){14}{\line(1,0){.0433673}}
\multiput(102.823,56.465)(.0433673,.0331633){14}{\line(1,0){.0433673}}
%\end
%\dashline{1}(98.75,50.5)(106.25,56)
\multiput(98.68,50.43)(.0441176,.0323529){17}{\line(1,0){.0441176}}
\multiput(100.18,51.53)(.0441176,.0323529){17}{\line(1,0){.0441176}}
\multiput(101.68,52.63)(.0441176,.0323529){17}{\line(1,0){.0441176}}
\multiput(103.18,53.73)(.0441176,.0323529){17}{\line(1,0){.0441176}}
\multiput(104.68,54.83)(.0441176,.0323529){17}{\line(1,0){.0441176}}
%\end
%\dashline{1}(99,47.25)(105.5,51.75)
\multiput(98.93,47.18)(.0481481,.0333333){15}{\line(1,0){.0481481}}
\multiput(100.374,48.18)(.0481481,.0333333){15}{\line(1,0){.0481481}}
\multiput(101.819,49.18)(.0481481,.0333333){15}{\line(1,0){.0481481}}
\multiput(103.263,50.18)(.0481481,.0333333){15}{\line(1,0){.0481481}}
\multiput(104.707,51.18)(.0481481,.0333333){15}{\line(1,0){.0481481}}
%\end
%\dashline{1}(99,43.5)(99.5,43.5)
\put(98.93,43.43){\line(1,0){.25}}
%\end
%\dashline{1}(99.25,44.25)(106,49)
\multiput(99.18,44.18)(.046875,.0329861){16}{\line(1,0){.046875}}
\multiput(100.68,45.235)(.046875,.0329861){16}{\line(1,0){.046875}}
\multiput(102.18,46.291)(.046875,.0329861){16}{\line(1,0){.046875}}
\multiput(103.68,47.346)(.046875,.0329861){16}{\line(1,0){.046875}}
\multiput(105.18,48.402)(.046875,.0329861){16}{\line(1,0){.046875}}
%\end
%\dashline{1}(98.75,40)(98.5,40.5)
\multiput(98.68,39.93)(-.03125,.0625){4}{\line(0,1){.0625}}
%\end
%\dashline{1}(98.5,40.5)(105.5,46)
\multiput(98.43,40.43)(.0411765,.0323529){17}{\line(1,0){.0411765}}
\multiput(99.83,41.53)(.0411765,.0323529){17}{\line(1,0){.0411765}}
\multiput(101.23,42.63)(.0411765,.0323529){17}{\line(1,0){.0411765}}
\multiput(102.63,43.73)(.0411765,.0323529){17}{\line(1,0){.0411765}}
\multiput(104.03,44.83)(.0411765,.0323529){17}{\line(1,0){.0411765}}
%\end
%\dashline{1}(99.25,36.75)(106,42)
\multiput(99.18,36.68)(.0421875,.0328125){16}{\line(1,0){.0421875}}
\multiput(100.53,37.73)(.0421875,.0328125){16}{\line(1,0){.0421875}}
\multiput(101.88,38.78)(.0421875,.0328125){16}{\line(1,0){.0421875}}
\multiput(103.23,39.83)(.0421875,.0328125){16}{\line(1,0){.0421875}}
\multiput(104.58,40.88)(.0421875,.0328125){16}{\line(1,0){.0421875}}
%\end
%\dashline{1}(99.5,32.75)(105.75,38.5)
\multiput(99.43,32.68)(.0347222,.0319444){18}{\line(1,0){.0347222}}
\multiput(100.68,33.83)(.0347222,.0319444){18}{\line(1,0){.0347222}}
\multiput(101.93,34.98)(.0347222,.0319444){18}{\line(1,0){.0347222}}
\multiput(103.18,36.13)(.0347222,.0319444){18}{\line(1,0){.0347222}}
\multiput(104.43,37.28)(.0347222,.0319444){18}{\line(1,0){.0347222}}
%\end
%\dashline{1}(102,30.5)(102.25,30.75)
\multiput(101.93,30.43)(.03125,.03125){4}{\line(0,1){.03125}}
%\end
%\dashline{1}(101.25,31)(106,35)
\multiput(101.18,30.93)(.0395833,.0333333){15}{\line(1,0){.0395833}}
\multiput(102.367,31.93)(.0395833,.0333333){15}{\line(1,0){.0395833}}
\multiput(103.555,32.93)(.0395833,.0333333){15}{\line(1,0){.0395833}}
\multiput(104.742,33.93)(.0395833,.0333333){15}{\line(1,0){.0395833}}
%\end
%\dashline{1}(105.25,30.5)(105,31.25)
\multiput(105.18,30.43)(-.03125,.09375){4}{\line(0,1){.09375}}
%\end
%\dashline{1}(105,31.25)(105.25,32.5)
\multiput(104.93,31.18)(.03125,.15625){4}{\line(0,1){.15625}}
%\end
%\dashline{1}(125,62.25)(127.5,64.75)
\multiput(124.93,62.18)(.0333333,.0333333){15}{\line(0,1){.0333333}}
\multiput(125.93,63.18)(.0333333,.0333333){15}{\line(0,1){.0333333}}
\multiput(126.93,64.18)(.0333333,.0333333){15}{\line(0,1){.0333333}}
%\end
%\dashline{1}(125,59.5)(130.5,64.75)
\multiput(124.93,59.43)(.0339506,.0324074){18}{\line(1,0){.0339506}}
\multiput(126.152,60.596)(.0339506,.0324074){18}{\line(1,0){.0339506}}
\multiput(127.374,61.763)(.0339506,.0324074){18}{\line(1,0){.0339506}}
\multiput(128.596,62.93)(.0339506,.0324074){18}{\line(1,0){.0339506}}
\multiput(129.819,64.096)(.0339506,.0324074){18}{\line(1,0){.0339506}}
%\end
%\dashline{1}(125,56.5)(130.5,61.75)
\multiput(124.93,56.43)(.0339506,.0324074){18}{\line(1,0){.0339506}}
\multiput(126.152,57.596)(.0339506,.0324074){18}{\line(1,0){.0339506}}
\multiput(127.374,58.763)(.0339506,.0324074){18}{\line(1,0){.0339506}}
\multiput(128.596,59.93)(.0339506,.0324074){18}{\line(1,0){.0339506}}
\multiput(129.819,61.096)(.0339506,.0324074){18}{\line(1,0){.0339506}}
%\end
%\dashline{1}(124.75,53.5)(131.25,59.25)
\multiput(124.68,53.43)(.0361111,.0319444){18}{\line(1,0){.0361111}}
\multiput(125.98,54.58)(.0361111,.0319444){18}{\line(1,0){.0361111}}
\multiput(127.28,55.73)(.0361111,.0319444){18}{\line(1,0){.0361111}}
\multiput(128.58,56.88)(.0361111,.0319444){18}{\line(1,0){.0361111}}
\multiput(129.88,58.03)(.0361111,.0319444){18}{\line(1,0){.0361111}}
%\end
%\dashline{1}(125,51)(131.25,56)
\multiput(124.93,50.93)(.0416667,.0333333){15}{\line(1,0){.0416667}}
\multiput(126.18,51.93)(.0416667,.0333333){15}{\line(1,0){.0416667}}
\multiput(127.43,52.93)(.0416667,.0333333){15}{\line(1,0){.0416667}}
\multiput(128.68,53.93)(.0416667,.0333333){15}{\line(1,0){.0416667}}
\multiput(129.93,54.93)(.0416667,.0333333){15}{\line(1,0){.0416667}}
%\end
%\dashline{1}(125.5,25.5)(131.75,32.25)
\multiput(125.43,25.43)(.0328947,.0355263){19}{\line(0,1){.0355263}}
\multiput(126.68,26.78)(.0328947,.0355263){19}{\line(0,1){.0355263}}
\multiput(127.93,28.13)(.0328947,.0355263){19}{\line(0,1){.0355263}}
\multiput(129.18,29.48)(.0328947,.0355263){19}{\line(0,1){.0355263}}
\multiput(130.43,30.83)(.0328947,.0355263){19}{\line(0,1){.0355263}}
%\end
%\dashline{1}(125.25,29.75)(131,35.75)
\multiput(125.18,29.68)(.0319444,.0333333){18}{\line(0,1){.0333333}}
\multiput(126.33,30.88)(.0319444,.0333333){18}{\line(0,1){.0333333}}
\multiput(127.48,32.08)(.0319444,.0333333){18}{\line(0,1){.0333333}}
\multiput(128.63,33.28)(.0319444,.0333333){18}{\line(0,1){.0333333}}
\multiput(129.78,34.48)(.0319444,.0333333){18}{\line(0,1){.0333333}}
%\end
%\dashline{1}(125.25,34.75)(125.75,34.75)
\put(125.18,34.68){\line(1,0){.25}}
%\end
%\dashline{1}(125.75,34.75)(125.25,34.5)
\multiput(125.68,34.68)(-.0625,-.03125){4}{\line(-1,0){.0625}}
%\end
%\dashline{1}(125.25,34)(125,34.75)
\multiput(125.18,33.93)(-.03125,.09375){4}{\line(0,1){.09375}}
%\end
%\dashline{1}(125.5,33.5)(131,39)
\multiput(125.43,33.43)(.0321637,.0321637){19}{\line(0,1){.0321637}}
\multiput(126.652,34.652)(.0321637,.0321637){19}{\line(0,1){.0321637}}
\multiput(127.874,35.874)(.0321637,.0321637){19}{\line(0,1){.0321637}}
\multiput(129.096,37.096)(.0321637,.0321637){19}{\line(0,1){.0321637}}
\multiput(130.319,38.319)(.0321637,.0321637){19}{\line(0,1){.0321637}}
%\end
%\dashline{1}(125.25,37.5)(131,42.75)
\multiput(125.18,37.43)(.0354938,.0324074){18}{\line(1,0){.0354938}}
\multiput(126.457,38.596)(.0354938,.0324074){18}{\line(1,0){.0354938}}
\multiput(127.735,39.763)(.0354938,.0324074){18}{\line(1,0){.0354938}}
\multiput(129.013,40.93)(.0354938,.0324074){18}{\line(1,0){.0354938}}
\multiput(130.291,42.096)(.0354938,.0324074){18}{\line(1,0){.0354938}}
%\end
%\dashline{1}(125,41.75)(132,46.25)
\multiput(124.93,41.68)(.05,.0321429){14}{\line(1,0){.05}}
\multiput(126.33,42.58)(.05,.0321429){14}{\line(1,0){.05}}
\multiput(127.73,43.48)(.05,.0321429){14}{\line(1,0){.05}}
\multiput(129.13,44.38)(.05,.0321429){14}{\line(1,0){.05}}
\multiput(130.53,45.28)(.05,.0321429){14}{\line(1,0){.05}}
%\end
%\dashline{1}(126.75,19)(131.5,23.5)
\multiput(126.68,18.93)(.0349265,.0330882){17}{\line(1,0){.0349265}}
\multiput(127.867,20.055)(.0349265,.0330882){17}{\line(1,0){.0349265}}
\multiput(129.055,21.18)(.0349265,.0330882){17}{\line(1,0){.0349265}}
\multiput(130.242,22.305)(.0349265,.0330882){17}{\line(1,0){.0349265}}
%\end
%\dashline{1}(125,47.75)(130.75,52.5)
\multiput(124.93,47.68)(.0399306,.0329861){18}{\line(1,0){.0399306}}
\multiput(126.367,48.867)(.0399306,.0329861){18}{\line(1,0){.0399306}}
\multiput(127.805,50.055)(.0399306,.0329861){18}{\line(1,0){.0399306}}
\multiput(129.242,51.242)(.0399306,.0329861){18}{\line(1,0){.0399306}}
%\end
%\dashline{1}(129,25.5)(131.25,27.75)
\multiput(128.93,25.43)(.0321429,.0321429){14}{\line(0,1){.0321429}}
\multiput(129.83,26.33)(.0321429,.0321429){14}{\line(0,1){.0321429}}
\multiput(130.73,27.23)(.0321429,.0321429){14}{\line(0,1){.0321429}}
%\end
%\dashline{1}(128.25,67.75)(130.25,70.5)
\multiput(128.18,67.68)(.033333,.045833){12}{\line(0,1){.045833}}
\multiput(128.98,68.78)(.033333,.045833){12}{\line(0,1){.045833}}
\multiput(129.78,69.88)(.033333,.045833){12}{\line(0,1){.045833}}
%\end
\put(135,38.25){\vector(0,-1){5.75}}
%\vector(40.5,50.25)(40.75,55.5)
\put(40.75,55.5){\vector(0,1){.07}}\multiput(40.5,50.25)(.03125,.65625){8}{\line(0,1){.65625}}
%\end
\put(75,21.5){$E$}
\put(69,43){$E_1$}
\end{picture}
\begin{center}
\nopagebreak[4] Figure \theppp
\end{center}
\addtocounter{ppp}{1}

The maximal $q$-bands of $E_1$ are ${\cal D}_1,\dots,{\cal D}_{s-1}, {\cal F}_2$. 
Then we attach to ${\cal F}_2$ two $q$-bands
${\cal F}_1$ and ${\cal F}_3$. The boundary path
of $E$ is ${\bf p}_1{\bf q}_1'{\bf p}'_2{\bf q}'_2$, where
${\bf p}'_2$ is the side label of $\cal F$ (equal  to the outer side label of ${\cal D}_s$). The labels of ${\bf q}'_1$ and $({\bf q}'_2)^{-1}$ represent the same element $wq_0$ in $H(\bf Y)$. So (changing notation) we may assume
that these paths are labeled by $Wq_0$ and $\Gamma=E$. 

The clock-wise boundary label of ${\cal F}_1$ is $q_0^{-1}Sq_0S'$, where $S$ and $S'$ are side labels
having $\theta$-letters from the alphabets ${\bf Q}_j$ and ${\bf Q}_{j'}$, respectively.

By Lemma \ref{simul}, there is a revolving quasi-computation $\cal C$ corresponding to the quasi-trapezium $E_1$ with top/bottom labels $W(1)$
and history $h$ equals the history of ${\cal D}_1$. Removing the $q$-band ${\cal F}_2$ from $E_1$, we see that the word $W(1)q_0^{-1}$
commutes with $x$ in $H(\bf Y)$.

\medskip

{\bf 6.} Consider now the subpath $ {\bf\tilde q}_1$ of ${\bf q}_1$ starting with the last edge of ${\bf \bar q}_1$ and having
label $W^3$. The maximal $q$-bands ${\cal Q}_s,\dots, {\cal Q}_m$ connect it with the
subpath ${\bf \tilde q}_2^{-1}$ of ${\bf q}_2^{-1}$. The diagram $\tilde\Gamma$
is bounded by ${\cal Q}_s, {\cal Q}_m$,
${\bf \tilde q}_1$ and ${\bf \tilde q}_2$.

It corresponds to the pair $(y,w)$, where $y$ is the label of the outer
(in $\tilde\Gamma$) side of the band
${\cal D}_s$ with $\theta$-edges labeled in ${\bf Q}_j$. We first replace ${\cal D}_s$ with $\cal F$, and then make the
reduction of $\cal F$, removing ${\cal F}_1$ and ${\cal F}_3$ with mutually inverse histories $h_1$ and $h_3$ according Remark
\ref{redu}. So we replace the pair
$(y,w)$ with $(x, u)$, where $u=S^{-1}wS$ in $H({\bf Y})$ 
%$S$ is the label of inner (in $\tilde\Gamma$) side of ${\cal F}_1$,
and $x$ is represented by a cyclically minimal label of the side of ${\cal F}_2$ (which is a copy of ${\cal D}_1$).

As in item {\bf 4}, for the word $Uq_0$ representing $uq_0$ in $H(\bf\Theta)$, we have a decomposition $U= U(1)U(2)$, where $U$ and $U(2)$  have equal images in $F(\bf\Theta^+)$. Since $U=S^{-1}WS$, the image of $U(2)$
is equal to the image of $S^{-1}W(2)S$
in $F(\bf\Theta^+)$. Therefore  $U(2)=(S')^{-1}W(2)S'$ in $H(\bf Y)$ by Lemmas \ref{thea} and \ref{tohist}. Taking into account that $Sq_0(S')^{-1}=q_0$ in $H({\bf Y})$ (since the boundary of ${\cal F}_1$ is trivial in $H({\bf Y})$), 
%and ${\cal F}_2$ is a copy of ${\cal D}_1$), 
we have

\begin{equation}\label{UW}
U(1)=Uq_0U(2)^{-1}=S^{-1}WSq_0(S')^{-1}W(2)^{-1}S'=S^{-1}(Wq_0 W(2)^{-1})S'=S^{-1}W(1)S'
\end{equation}
in $H(\bf Y)$.

As in {\bf 5}, there exists a quasi-trapezium with top/bottom label equal to $U(1)$ and history $h$. Removing the last $q$-band (which is a copy of the first
$q$-band of this quasi-trapezium), we obtain
that $U(1)q_0^{-1}$ represents an element commuting with $x$. It follows from (\ref{UW}) that $W(1)q_0^{-1}=SU(1)(S')^{-1}q_0^{-1}=S(U(1)q_0^{-1})S^{-1}$ commutes with $SxS^{-1}$ in $H(\bf Y)$.

\medskip

{\bf 7.} Thus, there is a minimal diagram over $H(\bf Y)$ with a boundary ${\bf x}_1{\bf y}_1{\bf x}_2{\bf y}_2$, where the labels of ${\bf x}_1$  and ${\bf x}_2^{-1}$ are the reduced forms of  $SxS^{-1}$ and
the labels of ${\bf y}_1$ and ${\bf y}_2^{-1}$ are equal to
$W(1)q_0^{-1}$. The $q$ edges of ${\bf y}_1$ are connected with corresponding edges of
${\bf y}_2$ by $q$-bands, and one can add
one more $q$-band (which is a copy of
the first $q$-band) obtaining, by Lemma \ref{NoAnnul}, a quasi-trapezium with top/bottom label $W(1)$. Since a side label of this quasi-trapezium is equal to $SxS^{-1}$,
the history $g$ of it is the reduced form of the product $h_1h_2h_1^{-1}$.

So, by Lemma \ref{simul} we have two revolving quasi-computations starting with the same  $W(1)$ and histories
$h_2$ (the history of the quasi-trapezium $E_1$ from item {\bf 5})
and $g$. Note that that the word $W(1)$ contains no cyclic subwords $V$ congruent to $\bf M$-accepting configurations. Indeed, otherwise Lemma \ref{simul} allows us  to construct a quasi-trapezium with top label $V$ and bottom label equal to $\Sigma(0)$. Gluing up the hub to the bottom, we see that $V$ is equal to $1$
in rank $1/2$. Hence the word $W$ is not cyclically minimal in rank $1/2$, a contradiction.

It follows from Lemma \ref{tworev} that there is a revolving
computation $\cal C$ with base $W(1)$
and history $h'$ congruent to $h_1$.

The side label of the corresponding trapezium is equal in $H(\bf Y)$ to $S$, and so $W'=W(1)q_0^{-1}$ commutes with $S$ in $H(\bf Y)$. 
Therefore $W^n=(W'S^{-1})^n=(W')^nS^n=(W')^n$ in $H(\Theta)$,
because $S$ has $q$-length $0$. Since $W'$ commutes with $x$, this implies
$W^nx=xW^n$ in $H(\bf \Theta)$.

To prove that $w^kxw^{-k}$ is a $0$-element, it suffices now to consider
the maximal $q$-bands in a diagram
of the equality $w^{ns}x=xw^{sn}$ for
$sn>k$ (as in item {\bf 1}.)

{\bf 8}. Recall that the value of $W(1)$, and therefore the value of $W'$ in $H({\bf\Theta})$ does not
depend on an element $x$ from ${\bf 0}(w)$  (see the end of item {\bf 4}). We have that the $0$-word $S$ equal to $(W')^{-1}W$
in $H({\bf \Theta})$,  does not depend on $x$ too. Hence the word $WS^{-1}=W'$ commutes with {\it every} element
from ${\bf 0}(w)$, and choosing an element $c\in H(\Theta)$ represented by the word $S$, we complete the proof.
\endproof

\subsection{End of the proof}

{\bf Proof of Theorem \ref{HE}}. Let a finitely generated group $G$ be a subgroup of a group $H$  finitely
presented in ${\cal B}_n$. Note that the set $R$ of all words in the generators $a_1,...,a_k$ of $H$, which are equal to $1$ in $H$, recursively enumerable since in class of all groups the set of defining relations of $H$
consists of a finite set and the recursively enumerable set of powers $w^n=1$ for group words $w$
in the alphabet   $\{a_1,...,a_k\}$. Each generator $b_i$ of $G$ from the finite set of generators $\{b_1,...,b_l\}$ is equal in $H$ to a word $w_i=w_i(a_1,\dots a_k)$, and to obtain the list of words in $b_i$-s
equal to $1$, it suffices to obtain the list of the corresponding products of $w_i$-s equal to $1$ in $G$.
Since the procedure of extracting such products from the list $R$ is effective, the group $G$ is recursively presented, i.e. it is given by a recursively enumerable set of relations in the generators $b_1,...,b_l$.

For the nontrivial part of Theorem \ref{HE}, we assume that there is a recursively enumerable set of defining relations for a finitely
 generated group $G$ of exponent $n$, and we have a presentation $G=\langle A\mid {\cal L}\rangle$ as in the beginning of Section \ref{interf}. The axioms (Z1), (Z2), (Z3) were checked for the group ${\bf G}(0)=H(\bf\Theta)$ and for the factor group ${\bf G}(1/2)=\bf G$ of ${\bf G}(0)$ in Subsection \ref{ax}. Therefore the group ${\bf G}(\infty)$
defined in Subsection \ref{construction} satisfies the identity $x^n=1$ by Lemma \ref{mainprop}. It is finitely presented in the variety ${\cal B}_n$ since the initial group $\tilde G$ given by relations (\ref{rel1}), (\ref{rel3}) is finitely presented, and the additional relation used in the constructions of groups
$H(a)$, $H(\theta,a)$, $H(\bf\Theta)$ and $H={\bf G}(\infty)$ are either $a$-relations (which
follow from the relations of $\tilde G$ by Lemma \ref{homo}) or relations of the form $u^n=1$.

It remains to prove that the canonical homomorphism $G\to H$ is injective. For this goal, we assume that a word $w$ in generators from $A$ is equal to $1$ in $H$ and prove that $w=1$ in $G$. So there is a $g$-reduced diagram $\Delta$ over $H$ with boundary $\bf p$ labeled by $w$. If $\Delta$ has a positive
cell, then by Lemmas \ref{four} and \ref{0cont}
(and Lemma \ref{mainprop}), there is a positive cell $\pi$ and a contiguity subdiagram $\Gamma$ of rank $0$ such that $(\pi, \Gamma,{\bf p})\ge \varepsilon$. Let ${\bf q}_1=\Gamma\bigwedge\pi$ and ${\bf q}_2=\Gamma\bigwedge{\bf p}$. By Condition (A2) (see the definition of A-map), the path ${\bf q}_1$ is reduced in 
$\Gamma$, and so every $q$-edge ${\bf e}$ of it is adjacent with a $q$-edge of $q_2^{-1}$. Such a $q$-edge
does exists in ${\bf q}_1$ since $\varepsilon n>1$, but does not exist in ${\bf q}_2$, a contradiction. 
%Then by Lemma \ref{2.2}, $$|{\bf q}_2|_q> (1+2\beta)^{-1} |{\bf q}_1|_q\ge (1+2\beta)^{-1}\varepsilon %|\partial\pi|_q> 0,$$ an so there is a $q$-edge in the path $\bf p$, a contradiction. 
It follows that $\Delta$ has no positive cells, and so it is a diagram over $H(\bf\Theta)$. So $w$ is trivial in the group $G$ by Lemma \ref{toHT}, as required. $\Box$

\medskip

{\bf Proof of Proposition \ref{prop12}}
(1) Let $G'$ be a homomorphic image of the group $G$. We can modify the definition of
congruent admissible words, replacing
equality of the words in the input sector $\{t_0\}\{t_1\}$ (and also in $\{t_0\}\{t_0\}^{-1}$ and $\{t_1\}^{-1}\{t_1\}$) modulo the relations of $G$ with equality modulo the bigger set  relations of $G'$.  The wording of other  definitions do not change (although the content of other concepts depending on the definition of congruence changes, we obtain different set
of regular admissible word, and so on).

Respectively, we enlarge the set of $a$-relations now. Let $M'$ be the factor group of $M$ by
all relations of $G'$.
%, which are now form the set of  $a$-relations. 
(So $M'$ is not necessarily finitely presented.) Respectively we obtain
the modified groups $\hat M'$, $H'(a)$, $H'({\bf Y})$, and so on.

In the modified Statement of Lemma \ref{homo} we have a homomorphism of $G'$ to $\hat M'$.
This just follows from the new definition of $a$-relations. Respectively, we obtain injective homomorphism $G'\to H'({\bf Y})$,
$G'\to H'({\bf \Theta})$, and $G'\to H'$ in Corollary \ref{GY},
Lemma \ref{toHT}, and the above proof of Theorem \ref{HE}. There are no changes in the formulations and proofs of other statements.

Thus, an arbitrary homomorphism $G\to G'$
extends to a homomorphism $H\to H'$. This means
that the embedding $G\hookrightarrow H$ is a CEP-embedding, and Proposition \ref{prop12} is proved for the embedding given by Theorem \ref{HE}. The proofs of Corollaries \ref{twog} and \ref{all} given in Introduction
also provide us with the CEP for the embeddings in $G\hookrightarrow E$ and $G_i\hookrightarrow E$, because the CEP is transitive and the retract of a group is a CEP subgroup.

\medskip

(2)Assume that two words $v$ and $w$ in the alphabet $A$ are conjugate in the group $H={\bf G}(\infty)$, and $v\ne 1$. So we have a $g$-reduced annular diagram $\Delta$ over $H$. Then the assumption that $\Delta$ has a positive cell
can be disproved exactly as in the proof of
Theorem \ref{HE}. Thus, $\Delta$ is a diagram over $H(\bf \Theta)$. By Lemma \ref{noq}, a $q$-annulus
(if any exists in $\Delta$) surrounds the hole of the annulus. So we have an annular diagram over $H(\theta,a)$
for the conjugacy of $v$ and the boundary label $V$ of such an annulus. However the retraction from Lemma \ref{thea} maps the word $v$ of $\theta$-length zero to $1$. It follows that $V=1$ in $H(\theta,a)$, and
therefore $v=1$ too, contrary the above assumption.

%one may assume that $\Delta$ is a minimal diagram,
%and therefore it contains no $q$-annuli. 
Hence
$\Delta$ has no $q$-annuli and  $(\theta,q)$-cells since $q$-bands cannot start/end on the boundary on $\Delta$.
Therefore $\Delta$ is a diagram over $H(\theta,a)$. By Lemma \ref{notha}, one may assume that $\Delta$ is minimal, and so it has
no $\theta$-annuli. If $\Delta$ has a positive
cell $\Pi$, then we obtain a contradiction
using Lemmas \ref{mainprop}, \ref{four}. \ref{0cont}, exactly as we used these lemmas in the proof of Theorem \ref{HE} (but
applying them to the graded diagram over $H(\theta,a)$ now), because there are no $\theta$-edges on the boundary of $\Delta$.
Hence $\Delta$ is a diagram over $H(a)$ without
$(\theta,a)$-cells, since the boundary has no
$\theta$-edges. Thus, every nontrivial cell of $\Delta$ is an $a$-cell. Therefore the words
$v$ and $w$ are conjugate in the group $G$.
The statement (2) of Proposition \ref{prop12} follows for the embedding $G\hookrightarrow H$.
We also have the Frattini property for the embeddings in  Corollaries \ref{twog} and \ref{all}
since following their proofs in Introduction, one should just keep in mind that the Frattini property of
subgroups is transitive and every retract has Frattini property.
$\Box$

\addcontentsline{toc}{section}{Index}

\printindex[g]

\vskip .5 in

\begin{minipage}[t]{3.5 in}
\noindent Alexander Yu. Ol'shanskii\\ Department of Mathematics\\
Vanderbilt University
and\\ Department of
Higher Algebra,\\ MEHMAT,
 Moscow State University
 \\ alexander.olshanskiy@vanderbilt.edu\\
\end{minipage}
\end{document}